\documentclass[12pt,twoside,reqno]{amsart}

\usepackage{amsmath, amssymb, amsthm, enumitem, esint}
\usepackage[normalem]{ulem}
\usepackage{soul}
\usepackage{xcolor}
\usepackage{mathrsfs}
\usepackage[hyperindex]{hyperref}

\usepackage{hyperref}
\hypersetup{bookmarksdepth=3}

\newcommand{\lto}{\longrightarrow}

\newcommand{\cyl}{\textnormal{cyl}}

\newcommand{\Spec}{\textnormal{Spec}}

\DeclareMathOperator{\graft}{graft}
\DeclareMathOperator{\sph}{sph}
\DeclareMathOperator{\aff}{aff}
\DeclareMathOperator{\con}{conf}
\newcommand{\sV}{\mathscr{V}}
\newcommand{\IR}{\mathbb{R}}
\newcommand{\lb}{\linebreak[1]}

\newcommand{\IN}{\mathbb{N}}
\newcommand{\IB}{\mathbb{B}}
\newcommand{\IS}{\mathbb{S}}
\newcommand{\IZ}{\mathbb{Z}}

\newcommand{\LL}{\mathcal{L}}
\newcommand{\IF}{\mathbb{F}}

\newcommand{\Qu}{\mathsf{Q}}
\newcommand{\DD}{\mathcal{D}}

\newcommand{\PP}{\mathcal{P}}
\newcommand{\UU}{\mathcal{U}}

\newcommand{\MM}{\mathcal{M}}
\newcommand{\RR}{\mathcal{R}}
\newcommand{\bY}{\mathbf{Y}}
\newcommand{\bO}{\mathbf{0}}

\newcommand{\bb}{\mathbf{b}}

\newcommand{\by}{\mathbf{y}}
\newcommand{\bx}{\mathbf{x}}

\newcommand{\XX}{\mathcal{X}}
\newcommand{\WW}{\mathcal{W}}

\newcommand{\eps}{\varepsilon}
\newcommand{\la}{\lambda}

\newcommand{\ov}[1]{\overline{#1}}
\newcommand{\td}[1]{\widetilde{#1}}

\DeclareMathOperator*{\osc}{osc}
\DeclareMathOperator{\Bry}{Bry}

\DeclareMathOperator{\Image}{Image}

\DeclareMathOperator{\rot}{rot}

\DeclareMathOperator{\Var}{Var}
\DeclareMathOperator{\Int}{Int}
\DeclareMathOperator{\Ric}{Ric}
\DeclareMathOperator{\Rm}{Rm}
\DeclareMathOperator{\tr}{tr}
\DeclareMathOperator{\id}{id}

\DeclareMathOperator{\diam}{diam}

\DeclareMathOperator{\conf}{conf}

\DeclareMathOperator{\inj}{inj}

\DeclareMathOperator{\vol}{vol}

\newcommand{\dotcup}{\ensuremath{\mathbin{\mathaccent\cdot\cup}}}
\newcommand{\EMPTY}[1]{}

\newtheorem{Theorem}[equation]{Theorem}
\newtheorem{Lemma}[equation]{Lemma}
\newtheorem{Corollary}[equation]{Corollary}
\newtheorem{Proposition}[equation]{Proposition}

\newtheorem{Claim}[equation]{Claim}

\theoremstyle{definition}
\newtheorem{Definition}[equation]{Definition}

\theoremstyle{remark}
\newtheorem{Remark}[equation]{Remark}

\numberwithin{equation}{section}

\title{Cylindrical Ricci flows and the PDE--ODI principle}
\author{Richard H. Bamler}
\address{Department of Mathematics, UC Berkeley, CA 94720, USA}
\email{rbamler@berkeley.edu}
\author{Yi Lai}
\address{Department of Mathematics, UC Irvine, CA 92697, USA}
\email{ylai25@uci.edu}
\thanks{R.~H. Bamler was supported by NSF grants DMS-2204364 and DMS-2604326. Y. Lai is supported by NSF grants DMS-2506832 and DMS-2543311 (CAREER)}
\date{\today}

\begin{document}

\begin{abstract}
We derive precise asymptotics for ancient asymptotically cylindrical Ricci flows. 
We establish a dichotomy in which the leading deviation from the cylinder is either linear or quadratic in the axial variables. 
In the linear case, we show that the flow is exponentially asymptotic to a Bryant soliton times a Euclidean factor. 
For $\kappa$-solutions, we further show that the flow is even isometric to such a product, up to parabolic rescaling.
In the quadratic case, we obtain higher-order asymptotics governed by a finite-dimensional invariant. We also establish propagation of cylindrical regions and uniqueness of cylindrical tangent flows.

These results closely parallel those in our earlier work on cylindrical mean curvature flows.
Both analyses follow the same PDE--ODI principle and analytic strategies, yielding a unified treatment of cylindrical Ricci flow and mean curvature flow.
The essential new ingredient specific to the Ricci flow setting is a coherent harmonic map heat flow gauge, which allows this principle to be applied even in the absence of assumptions on the sign of the curvature.
The local analytic tools developed here also form part of a program aimed at constructing Ricci flow with surgery in dimension~4.
\end{abstract}

\maketitle
\tableofcontents

\section{Introduction}

\subsection{Overview}

Cylindrical models play a central role in the analysis of singularities in Ricci flow. A fundamental example is Perelman's construction of Ricci flow with surgery on $3$-manifolds~\cite{Perelman1,Perelman2}, which led to the resolution of the Poincar\'e and Geometrization Conjectures. In this construction, regions that are nearly cylindrical are identified, cut and capped in order to continue the flow through singularities. In dimension~$3$, this analysis benefitted from strong curvature control and a detailed understanding of the possible models near singularities (via the Hamilton--Ivey pinching estimate and Perelman's theory of $\kappa$-solutions). No comparably strong tools are available for general Ricci flows in higher dimensions. This is one of the main obstacles to constructing Ricci flows with surgery in dimensions~$4$ and higher. New techniques are therefore needed to analyze flows that are only locally close to a cylinder, without assuming these special $3$-dimensional inputs. The purpose of this paper is to develop a comprehensive analytic theory for this problem, applicable to Ricci flows that are locally close to round cylinders $\IR^k\times\IS^{n-k}$ on finite or infinite time-intervals.

From an analytic viewpoint, the main difficulty is that the flow is only \emph{locally} close to a cylindrical model. The analogous localization problem arises in cylindrical mean curvature flow~\cite{Bamler_Lai_MCF1,Bamler_Lai_MCF2}. In that setting, the geometry allows a graphical representation over the cylinder only on part of the hypersurface, and the relevant cylindrical region may change with time. Any effective approach must therefore begin with a careful localization step. This step is delicate, because the linearized operator at the cylinder has neutral and unstable modes. The crucial ingredient in overcoming these obstacles was the \emph{PDE--ODI principle}, developed in~\cite{Bamler_Lai_MCF1}. The principle converts the original parabolic PDE into a system of ordinary differential inequalities for a finite number of leading modes, while controlling the remaining modes and the localization error. High-order information for the original PDE can then be extracted using standard ODE arguments, even though the flow is only known to be cylindrical on a large but finite region. In the present paper, we show that the PDE--ODI principle overcomes the same localization and mode-control problems for Ricci flow. The two settings thus fit into a common analytic framework. Much of our discussion can therefore be organized in parallel with the mean curvature flow theory, and we review the PDE--ODI principle in more detail later in the introduction.

The passage from mean curvature flow to Ricci flow does, however, require some new ingredients.
In mean curvature flow, the fixed ambient geometry provides a local graphical representation over the cylinder. However, Ricci flow has no fixed ambient geometry, and its equation is parabolic only modulo diffeomorphisms. To apply the PDE--ODI principle, we therefore first construct a harmonic map heat flow (HMHF) from an open spacetime domain of the solution to the model cylinder, defined throughout the prescribed time interval. In the induced coordinates, the metric satisfies the strictly parabolic Ricci--DeTurck equation. Establishing this global gauge is the central analytic prerequisite for our analysis; its construction and the additional difficulties specific to Ricci flow are discussed later in the introduction.

Using this global HMHF gauge, we obtain results that closely parallel those for cylindrical mean curvature flow in~\cite{Bamler_Lai_MCF1}, subject to several qualifications specific to the Ricci flow setting.
We establish a dichotomy according to whether the linear or the quadratic mode is dominant and control all remaining modes to arbitrarily high order.
In the dominant linear mode case, we prove exponential convergence to the Bryant soliton times a Euclidean factor at cylindrical infinity and classify asymptotically cylindrical $\kappa$-solutions with nonnegative curvature operator.
In the dominant quadratic mode case, we construct a well-defined, continuous invariant $\Qu(\MM)$ that determines the asymptotics of the full rescaled metric, in suitable HMHF coordinates, to arbitrary polynomial order in inverse rescaled time.
We also construct a family of symmetric $\kappa$-solutions realizing every admissible value of $\Qu$, consisting of ancient ovals, their products with Euclidean factors, and the round shrinking cylinder.
The same approach also yields propagation of cylindrical regions on finite time intervals and uniqueness of cylindrical tangent flows.

We keep the discussion in this introduction focused on the main conclusions and on the new ideas that are needed to pass from mean curvature flow to Ricci flow. The statements below are simplified. Precise formulations appear in the body of the paper.

\subsection{Main results I: Ancient flows}

We begin with ancient Ricci flows $\MM=(M,(g_t)_{t<T})$, which are defined for all sufficiently negative times. We call $\MM$ \emph{asymptotically $(n,k,\Gamma)$-cylindrical} if, for some fixed point $x\in M$, the rescaled pointed manifolds $(M,|t|^{-1}g_t,x)$ converge as $t\to-\infty$ in the smooth Cheeger--Gromov topology to a round cylinder $\IR^k\times(\IS^{n-k}/\Gamma)$, where $1\leq k\leq n-2$ and $\Gamma\subset O(n-k+1)$ is a finite subgroup acting freely on the sphere.
{For noncollapsed $\IF$-limits, this is equivalent to having the corresponding round shrinking cylinder as a tangent flow at infinity in the sense of $\IF$-convergence~\cite{Bamler_RF_compactness,Bamler_HK_RF_partial_regularity} (see Lemma~\ref{lem:asymptotically_cylindrical_tangent_flows}).}

Theorems~\ref{thm:intro_neck_stability} and~\ref{thm:intro_cylindrical_limits} also apply to flows that may have singularities, including noncollapsed $\IF$-limits. In this introduction we take $\Gamma=\{1\}$.

Our first result describes the leading deviation from the cylinder. Apart from the cylinder itself, there are two possible behaviors: the leading change in the spherical metric is either linear or quadratic in the Euclidean directions of the cylinder. This is the Ricci flow counterpart of the linear--quadratic dichotomy discussed in our earlier work on mean curvature flow~\cite{Bamler_Lai_MCF1} and provides the starting point for our analysis here.

\begin{Theorem}[Asymptotic dichotomy; Propositions~\ref{Prop_ancient_summary}, \ref{p:coefficient}, and~\ref{Prop_RF_dom_qu_asymp}]\label{thm:intro_dichotomy}
Let $\MM=(M,(g_t)_{t<T})$ be an ancient asymptotically $(n,k)$-cylindrical Ricci flow  with complete time-slices and bounded curvature on compact time-intervals. 
Then, after translating time if necessary, there are time-dependent coordinates, expressed by maps $\phi_t : \Omega_t \to \IR^k \times \IS^{n-k}$ defined on larger and larger neighborhoods $\Omega_t \subset M$ and evolving by rescaled HMHF, such that we can locally express the rescaled flow as a perturbation of the cylinder
\[
  e^{\tau} (\phi_{-e^{-\tau}})_* g_{-e^{-\tau}}=g_{\cyl}+u_\tau, \qquad \tau = - \log(-t),
\]
and the perturbation $u_\tau$ shows one of the following behaviors:
\begin{enumerate}[label=(\alph*)]
\item We have $u_\tau = 0$ and $\MM$ is isometric to a round shrinking cylinder.
\item There is a nonzero vector $\bb=(b_1,\ldots,b_k)\in\IR^k$ such that
\[
 u_\tau
 =e^{\tau/2}\sum_{i=1}^k b_i\mathfrak p_i^{(1)}(\bx)g_{\IS^{n-k}}
  +o(e^{\tau/2}).
\]
In this case we say that the \emph{linear mode is dominant}.
\item There is a nonzero symmetric matrix-valued function $(c_{ij}(\tau))$ with $C^{-1}|\tau|^{-1}\leq\|(c_{ij}(\tau))\|\leq C|\tau|^{-1}$ for some $C\geq1$ such that, in such coordinates,
\begin{equation}\label{eq:intro_quadratic_expansion}
 u_\tau
 =\sum_{i,j=1}^k c_{ij}(\tau)\mathfrak p_{ij}^{(2)}(\bx)g_{\IS^{n-k}}
 -\frac{n-k}{2\sqrt2}\sum_{i,j=1}^k c_{ij}(\tau)\,d\bx_i\otimes d\bx_j
  +o(|\tau|^{-1}).
\end{equation}
In this case we say that the \emph{quadratic mode is dominant}.
\end{enumerate}
Here $\mathfrak p_i^{(1)}$ and $\mathfrak p_{ij}^{(2)}$ denote the linear and quadratic Hermite polynomials in the axial variable $\bx\in\IR^k$. The displayed little-$o$ estimates hold in $C^m$ on every fixed compact subset of the model cylinder, for every $m$. 
\end{Theorem}

We remark that using the PDE--ODI principle~\cite{Bamler_Lai_MCF1} we essentially automatically obtain higher order expansions, up to arbitrary polynomial order in $e^{\tau/2}$ (linear case) or in $|\tau|^{-1}$ (quadratic case).
In the linear case, these expansions provide the basis for the sharper asymptotic description in Theorem~\ref{thm:intro_bryant_infinity}.
In the quadratic case, they are expressed by the asymptotic characterization in Theorem~\ref{thm:intro_quadratic}.

\medskip
We first discuss the dominant linear mode case.
In our mean curvature flow paper~\cite{Bamler_Lai_MCF1}, we proved that the flow must be a rescaling of the bowl soliton times a Euclidean factor.
The corresponding model for Ricci flow is the Bryant soliton times a Euclidean factor $\MM'=\IR^{k-1}\times\MM^{n-k+1}_{\Bry}$.
The lack of a comparison principle prevents us from obtaining the same global classification in this generality.
We can, however, prove that the flow is uniformly asymptotic to the Bryant model with an exponentially decaying error.
This improves substantially on the polynomial error bound obtained in~\cite{Bamler_Lai_MCF1} before applying the comparison principle.

\begin{Theorem}[Exponential Bryant asymptotics at infinity; Theorem~\ref{thm:linear_bryant_asymptotics}]\label{thm:intro_bryant_infinity}
Let $\MM=(M,(g_t)_{t<T})$ be an ancient asymptotically $(n,k)$-cylindrical Ricci flow with dominant linear mode and  complete time-slices and bounded curvature on compact time-intervals. 
Write $\MM' = (M', (g')_{t \in \IR})$ and denote by $y(\cdot,t) : M' \to \IR_{\geq 0}$ the time-$t$ distance to $\IR^{k-1} \times \{ x_{\Bry} \}$, where $x_{\Bry}$ denotes the tip of $\MM_{\Bry}^{n-k+1}$.
Then for a dimensional choice of $Y_0$ the following is true.

After parabolically rescaling $\MM$, there is a smooth family of diffeomorphisms onto their images
\[ \psi_t : \{ y(\cdot,t) > Y_0 \} \to M \]
such that $\psi_{t_1} |_{\{ y(\cdot,t_2) > Y_0 \}} = \psi_{t_2}$ {for $t_1\leq t_2<T$}. For every integer $m\geq0$ and all $t<T$, we have
\[
 \big|\nabla^{m}_{g'_t}(\psi_t^*g_t-g'_t)\big|_{g'_t}
 \leq C_m e^{-c_m y(\cdot,t)}
 \qquad\textnormal{on }\{y(\cdot,t)>Y_0\}.
\]
Here $C_m, c_m > 0$ are constants that only depend on $m$ and the dimension of $M$.
\end{Theorem}

If we assume in addition that the flow is a $\kappa$-solution, then we obtain a full classification in the dominant linear case. Recall that a $\kappa$-solution is a nonflat ancient Ricci flow with complete time-slices, bounded curvature on compact time-intervals, nonnegative curvature operator, and $\kappa$-noncollapsing on all scales for a fixed $\kappa>0$.
The non-collapsing condition could be dropped, as it is implied by asymptotic cylindricality.

\begin{Theorem}[Corollary~\ref{thm:unique_bryant}]\label{thm:intro_linear}
Let $\MM=(M,(g_t)_{t<T})$ be an ancient asymptotically $(n,k)$-cylindrical $\kappa$-solution with dominant linear mode.
Then, up to parabolic rescaling and restriction of the model to $(-\infty,T)$,
\[
 \MM\cong\IR^{k-1}\times\MM^{n-k+1}_{\Bry}.
\]
\end{Theorem}

The uniqueness of the Bryant soliton has been studied extensively. 
Earlier uniqueness results for steady gradient solitons were established by Brendle in dimension $3$ under noncollapsing assumptions and, in higher dimensions, assuming positive sectional curvature and cylindrical asymptotics~\cite{Brendle_steady_3d,Brendle_steady_higher_dim}.
For general ancient flows, Brendle~\cite{Brendle_3dBryant_unique} showed that any noncompact, 3-dimensional $\kappa$-solution must be a Bryant soliton or a quotient of a cylinder (see also \cite{MR4319063} for an alternate approach).
This was generalized to higher dimensions by Brendle--Naff~\cite[Corollary~1.6]{Brendle_Naff_ancient} under the additional assumption that the flow has uniformly positive isotropic curvature.
All these past results address the case $k=1$, i.e., when the flow is asymptotic to a cylinder of the form $\IR \times \IS^{n-1}$.

\medskip

Let us now discuss the case in which the quadratic mode dominates. 
Here we obtain quantitative asymptotic control analogous to that established in our first mean curvature flow paper~\cite{Bamler_Lai_MCF1}. 
As in that paper, these results do not give a complete classification. Such a classification in mean curvature flow is obtained using the leading mode condition in our second paper~\cite{Bamler_Lai_MCF2}; transferring those ideas to Ricci flow lies beyond the scope of the present work. 
We also obtain the same type of asymptotic invariant $\Qu(\MM)$ as in \cite{Bamler_Lai_MCF1}. It records finer asymptotic information beyond the leading $|\tau|^{-1}$ term of the quadratic coefficient matrix. {We have $\Qu(\MM)=0$ if and only if $\MM$ is a round shrinking cylinder or has dominant linear mode.}
In the mean curvature flow setting, these considerations formed the basis of our subsequent work \cite{Bamler_Lai_MCF2}; in the Ricci flow setting they are likely to lead to a full classification of asymptotically cylindrical $\kappa$-solutions.

Theorem~\ref{thm:intro_quadratic} establishes the existence of the asymptotic invariant $\Qu(\MM)$ and discusses basic properties.

\begin{Theorem}[Propositions~\ref{Prop_RF_Q_basic} and~\ref{Prop_RF_Q_continuity}]\label{thm:intro_quadratic}
Let $\MM=(M,(g_t)_{t<T})$ be an ancient asymptotically $(n,k)$-cylindrical Ricci flow with complete time-slices and bounded curvature on compact time-intervals. 
There is a well-defined invariant $\Qu(\MM)\in\IR^k_{\uparrow,\geq0}$,
where $\IR^k_{\uparrow,\geq0}:=\{(q_1,\ldots,q_k)\in\IR^k:0\leq q_1\leq\cdots\leq q_k\}$.
This vector is determined by the asymptotic behavior of $\MM$.
In the dominant quadratic case, it is determined by the asymptotics of the time-dependent matrix $(c_{ij}(\tau))$ from the representation \eqref{eq:intro_quadratic_expansion}.
This invariant has the following properties:
\begin{enumerate}[label=(\alph*)]
\item If $\Qu(\MM) = 0$, then $\MM$ is a round shrinking cylinder or has dominant linear mode. {In the latter case, Theorem~\ref{thm:intro_bryant_infinity} applies.}
\item {$\Qu(\MM)$ is invariant under time-shifts. If} $\lambda \MM = (M, (\lambda^2 g_{\lambda^{-2}t})_{t < \lambda^2T})$ denotes the parabolically rescaled flow, then $\Qu(\lambda \MM) = \lambda \Qu(\MM)$.
\item The invariant $\Qu$ is continuous under pointed smooth Hamilton convergence within the class of regular parts of noncollapsed $\IF$-limits in Proposition~\ref{Prop_RF_Q_continuity}, as long as the limit is again asymptotically $(n,k)$-cylindrical.
\end{enumerate}
\end{Theorem}

The next result shows that $\Qu(\MM)$ determines the asymptotics of the rescaled flow to any prescribed polynomial order in $|\tau|^{-1}$. 
It is the analogue of~\cite[Theorem~1.6]{Bamler_Lai_MCF1} from the mean curvature flow setting.

\begin{Theorem}[Proposition~\ref{Prop_RF_same_Q_close}]\label{thm:intro_higher_order_asymptotics}
Suppose that $$\Qu(\MM^0) = \Qu(\MM^1)$$ for two ancient asymptotically $(n,k)$-cylindrical Ricci flows $\MM^0 = (M^0,(g^0_t)_{t<T^0})$ and $\MM^1= (M^1,(g^1_t)_{t<T^1})$ with complete time-slices and bounded curvature on compact time-intervals.
For every $A\geq0$ and integer $m\geq0$, the time-dependent coordinate maps $\phi^i_t : \Omega_t \to \IR^k \times \IS^{n-k}$ from Theorem~\ref{thm:intro_dichotomy} can be chosen so that, after possibly composing $\phi^0_t$ with a fixed rotation in the $\IR^k$-factor, we have for all sufficiently negative $\tau$
\[ \Bigl\Vert e^\tau \bigl( (\phi_{-e^{-\tau}}^1)_* g^1_{-e^{-\tau}} - (\phi_{-e^{-\tau}}^0)_* g^0_{-e^{-\tau}} \bigr) \Bigr\Vert_{C^m(\IB^k_{A\sqrt{\log |\tau|}} \times \IS^{n-k})} \leq O( |\tau|^{-A}). \]
\end{Theorem}

On the other hand, we show that each admissible value of $\Qu(\MM)$ can be realized by an ancient oval, which is a $\kappa$-solution.
These solutions can therefore serve as comparison models in Theorem~\ref{thm:intro_higher_order_asymptotics}.
The following theorem is the analogue of our existence result for mean curvature flow~\cite[Theorem~1.7]{Bamler_Lai_MCF1}.

\begin{Theorem}[Theorem~\ref{Thm_RF_Q_realization}]\label{thm:intro_realization}
For any $\Qu' \in \IR^k_{\uparrow, \geq 0}$ there is a $\kappa$-solution $\MM_{\Qu'} = (M_{\Qu'}, (g_{\Qu',t})_{t< 0})$ with $\Qu(\MM_{\Qu'}) = \Qu'$.

Moreover, we have $\MM_{\Qu'} = \IR^l \times \MM_{\Qu''}$, where $l$ is the number of vanishing entries in $\Qu'$ and $\Qu'' \in \IR^{k-l}_{\uparrow,>0}$ consists of the remaining entries. 
We can take $M_{\Qu''}=\IS^{n-l}$ such that the flow
$(g_{\Qu'',t})$ is invariant under the standard
$\IZ_2^{k-l}\times SO(n-k+1)$-action on
$\IS^{n-l}\subset\IR^{k-l}\times\IR^{n-k+1}$.
If $l = k$, then $\MM_{\Qu'}$ is the round shrinking cylinder.
\end{Theorem}

The solutions in Theorem \ref{thm:intro_realization} are called ancient ovals.  Perelman~\cite[Section~1.4]{Perelman2} constructed the first ancient oval in dimension 3. Its asymptotics and uniqueness were established in~\cite{Angenent_Brendle_Daskalopoulos_Sesum_asymptotics,BDS_2021}. For compact noncollapsed ancient flows on $\IS^n$ with bounded curvature, Brendle--Daskalopoulos--Naff--\v{S}e\v{s}um obtained the higher-dimensional classification assuming that the curvature is uniformly $\textnormal{PIC}$ and weakly $\textnormal{PIC2}$~\cite{Brendle_Daskalopoulos_Naff_Sesum_ancient}. 
In dimension $4$, Haslhofer constructed a one-parameter family of ancient ovals with asymptotic cylinder $\IR^2\times\IS^2$~\cite{Haslhofer_4d_kappa}. 
Asymptotics for higher-dimensional $SO(k)\times SO(n-k+1)$-symmetric ancient ovals are obtained in~\cite{Daskalopoulos_Du_Sesum_Zhao_RF_ovals}.

\subsection{Main results II: Rigidity, neck propagation, and uniqueness of tangent flows}

The methods developed in this paper also apply beyond the ancient asymptotically cylindrical setting of the previous subsection.
They recover the cylindrical shrinker rigidity results of Li--Wang~\cite{Li_Wang_cylinder_rigidity} and Colding--Minicozzi~\cite{Colding_Minicozzi_RF_2025} and yield smooth convergence to cylindrical tangent flows, for which Fang--Li~\cite{Fang_Li_strong_uniqueness} established strong uniqueness results in the smooth setting.
We also obtain propagation of cylindrical regions on finite time-intervals (Theorem~\ref{thm:intro_neck_stability}).
Our approach places these results in a common framework based on the PDE--ODI principle.
The local nature of the analysis also allows us to treat singular flows, including noncollapsed $\IF$-limits arising in blow-up analysis~\cite{Bamler_RF_compactness,Bamler_HK_RF_partial_regularity}.

Let us now present these results in a simplified form; the precise statements are again indicated in parentheses.
For simplicity, the rigidity and neck propagation results below are stated only for conventional Ricci flows and cylinders without quotients.

The first result below recovers the cylindrical shrinker rigidity results of Li--Wang~\cite{Li_Wang_cylinder_rigidity} for $\IR\times\IS^{n-1}$ and of Colding--Minicozzi~\cite{Colding_Minicozzi_RF_2025} for $\IR^k\times\IS^{n-k}$. The same result also gives rigidity for Ricci flows defined on $(-\infty,0)$ that remain close to the cylinder at scale $\sqrt{-t}$ around a fixed point, without assuming self-similarity.

\begin{Theorem}[Rigidity; Corollary~\ref{Cor_eternal}]\label{thm:intro_cylindrical_rigidity}
There is $\delta=\delta(n,k)>0$ with the following properties.
\begin{enumerate}[label=(\alph*)]
\item Let $(M,g,f)$ be a complete connected gradient shrinking Ricci soliton with $\Ric_g+\nabla_g^2f=\frac12g$. If $(M,g)$ contains an $(n,k,\delta)$-neck at scale $1$ centered at a critical point of $f$, then $(M,g)$ is isometric to the normalized round cylinder.
\item Let $\MM=(M,(g_t)_{t<0})$ be a connected Ricci flow with complete time-slices and bounded curvature on compact time-intervals. If, for some fixed $p\in M$, the point $(p,t)$ is a center of an $(n,k,\delta)$-neck at scale $\sqrt{-t}$ for every $t<0$, then $\MM$ is the round shrinking cylindrical flow.
\end{enumerate}
\end{Theorem}

Next, Theorem~\ref{thm:intro_neck_stability} propagates cylindricality from the initial and final times to all intermediate times under an $H_n$-center condition. The required $\delta=\delta(n,k,\eps)$ is independent of the length of the time-interval.

\begin{Theorem}[Neck propagation; Theorem~\ref{Thm_propagation_necks}]\label{thm:intro_neck_stability}
For every $\eps>0$ there is $\delta=\delta(n,k,\eps)>0$ with the following properties.
Let $\MM=(M,(g_t)_{t\in[-T_-,-T_+]})$ be a complete Ricci flow with uniformly bounded curvature, where $0<T_+<T_-$. Suppose that $p_\pm\in M$ and that $(p_\pm,-T_\pm)$ are centers of $(n,k,\delta)$-necks at scales $\sqrt{T_\pm}$. Assume that $(p_-,-T_-)$ is an $H_n$-center of $(p_+,-T_+)$.
Then $(p_+,t)$ is a center of an $(n,k,\eps)$-neck at scale $\sqrt{-t}$ for every $t\in[-T_-,-T_+]$, and $d_{-T_-}(p_-,p_+)\leq C(n)\sqrt{T_-}$.
All constants are independent of the location and length of the time-interval.
\end{Theorem}

Here $H_n$-centers are defined as in~\cite[Definition~3.10]{Bamler_HK_entropy_estimates}.
Intuitively, $p_-$ is an approximate center of mass at time $-T_-$ of the conjugate heat flow evolving backward from $(p_+,-T_+)$.

We also establish the uniqueness and full convergence whenever a subsequential cylindrical limit occurs, both at a finite singularity and at $t=-\infty$.
Colding--Minicozzi~\cite[Theorem~10.8]{Colding_Minicozzi_RF_2025} established uniqueness of cylindrical tangent flows at Type~I singularities of closed Ricci flows.
Fang--Li~\cite[Theorems~2.23, 7.7, and~7.11]{Fang_Li_strong_uniqueness} proved strong uniqueness using a gauge generated by a conjugate heat-kernel potential.
Our Corollaries~\ref{Cor_unique_tangent} and~\ref{Cor_unique_tangent_infinity} give pointed smooth convergence at every scale with basepoints along a single worldline of the original Ricci flow.

\begin{Theorem}[Uniqueness of tangent flows; Corollaries~\ref{Cor_unique_tangent} and~\ref{Cor_unique_tangent_infinity}]\label{thm:intro_cylindrical_limits}
Let $\XX$ be a metric flow arising as a noncollapsed $\IF$-limit of
$n$-dimensional Ricci flows in the sense of~\cite{Bamler_RF_compactness,Bamler_HK_RF_partial_regularity}.
\begin{enumerate}[label=(\alph*)]
\item If a tangent flow at a point $p_0\in\XX_{t_0}$, with $t_0$ strictly later than the initial time, is the round shrinking cylinder $\IR^k\times\IS^{n-k}$, then there is a regular worldline $p(t)$ for $t<t_0$ sufficiently close to $t_0$, with $p(t)\to p_0$, along which all parabolic blow-ups converge smoothly to the cylindrical flow as the scale tends to zero.
\item Suppose that $\XX$ is ancient. If a tangent flow at infinity based at the conjugate heat flow from $p_0\in\XX_{t_0}$ is the round shrinking cylinder $\IR^k\times\IS^{n-k}$, then there is a regular worldline $p(t)$ defined for all sufficiently negative times along which all parabolic blow-downs converge smoothly to the cylindrical flow as the scale tends to infinity. Hence the regular part of $\XX$ is an ancient asymptotically $(n,k)$-cylindrical Ricci flow spacetime in the sense of Definition~\ref{defn:asymptotically_cylindrical}; see also Lemma~\ref{lem:asymptotically_cylindrical_tangent_flows}.
\end{enumerate}
In both assertions, smooth convergence means pointed smooth Cheeger--Gromov convergence on compact subsets of $\IR^k\times\IS^{n-k}\times(-\infty,0)$. 
\end{Theorem}

For rigidity results on asymptotically conical shrinkers, shrinkers with cylindrical ends, and generalized cylinders, see~\cite{Kotschwar_Wang_AC_rigidity,Kotschwar_Wang_cylindrical_uniqueness,Li_Zhang_shrinker_rigidity,Li_Yan_Zhang_AC_I,Kim_Ozuch_blowdowns}.
Earlier stability results include backward stability of necks in three-dimensional noncompact $\kappa$-solutions~\cite[Theorem~6.1]{Kleiner_Lott_singular} and weak stability of generalized cylinders under admissible warped-product perturbations~\cite[Theorem~29]{Carson_Isenberg_Knopf_Sesum_cylinders}.
For related uniqueness results for compact or asymptotically conical tangent flows, see~\cite{Sesum_conv_to_soliton,Ache_Ricci_asymptotic_limits,Sun_Wang_KRF,Chan_Ma_Zhang_closed_tangent_flows,Chan_Ma_Zhang_tangent_rate,Choi_Lai_HMHF,Li_compact_ancient_parametrization,Li_Yan_Zhang_AC_II}.

\subsection{The harmonic map heat flow gauge and the PDE--ODI principle}\label{subsec:intro_HMHF_methods}
In our work \cite{Bamler_Lai_MCF1,Bamler_Lai_MCF2}, we developed a PDE--ODI principle for a broad class of parabolic equations.
It converts a solution of a parabolic PDE into a solution of a finite-dimensional system of ordinary differential inequalities, which can then be analyzed using standard ODE methods.
In this way, heuristic perturbative pictures can often be made rigorous.
Crucially, the principle only requires the flow to be \emph{locally} close to the model.
The key to this localization is a new mechanism for absorbing the errors introduced by cutoff functions into the evolution inequalities.
In mean curvature flow, a cylindrical region can be represented as a local graph over the cylinder, and the graph function satisfies a strictly parabolic equation to which the principle applies.
Ricci flow, however, is only weakly parabolic because of its diffeomorphism invariance, so we must first construct a suitable coordinate gauge throughout the time interval under consideration.

A standard choice for such a gauge is the Ricci--DeTurck gauge, which is obtained by pushing the flow forward by diffeomorphisms evolving by harmonic map heat flow (HMHF).
Let us describe this construction in more detail.
Consider a Ricci flow $(M,(g_t))$ close, modulo rescalings and diffeomorphisms, to the shrinking cylindrical flow $(M_{\cyl},(g_{\cyl,t}))$, where $M_{\cyl}=\IR^k\times\IS^{n-k}$ and $g_{\cyl,t}=g_{\IR^k}+(-t)g_{\IS^{n-k}}$ (with $g_{\IS^{n-k}}$ of sectional curvature $\frac{1}{2(n-k-1)}$).
To obtain the Ricci--DeTurck gauge, we first solve the harmonic map heat flow equation for maps $\chi_t:M\to M_{\cyl}$:
\begin{equation}\label{eq:intro_HMHF}
 \partial_t\chi_t=\Delta_{g_t,g_{\cyl,t}}\chi_t.
\end{equation}
As long as the maps are diffeomorphisms onto their images, the pushed-forward metrics $g'_t:=(\chi_t)_*g_t$ satisfy the Ricci--DeTurck flow equation
\begin{equation}\label{eq:intro_RDTF}
 \partial_t g'_t=-2\Ric_{g'_t}-\LL_{X_{g_{\cyl,t}}[g'_t]}g'_t,
\end{equation}
where $X_h[g]:=\Delta_{g,h}\id$ is the DeTurck vector field.
This equation is strictly parabolic in $g'_t$.

To compare with the cylindrical model, we next rescale so that the cylindrical flow becomes stationary.
To do this, set $\tau=-\log(-t)$ and let $\xi_\lambda$ be the dilation $\xi_\lambda(\bx,\by)=(\lambda\bx,\by)$ of the Euclidean factor in $M_{\cyl}=\IR^k\times\IS^{n-k}$.
The rescaled metrics
\[
 \td g'_\tau:=e^\tau(\xi_{e^{\tau/2}})_*g'_{-e^{-\tau}}
\]
are identically equal to $g_{\cyl}:=g_{\cyl,-1}$ when $g'_t=g_{\cyl,t}$, and satisfy the rescaled Ricci--DeTurck flow equation
\begin{equation}\label{eq:intro_full_rescaled_RDTF}
 \partial_\tau\td g'_\tau=-2\Ric_{\td g'_\tau}+\td g'_\tau
 -\LL_{X_{g_{\cyl}}[\td g'_\tau]+\nabla f}\td g'_\tau,
\end{equation}
where $f=\frac14|\bx|^2$ and the gradient is taken with respect to $g_{\cyl}$.
Thus the perturbation $u_\tau:=\td g'_\tau-g_{\cyl}$ satisfies a parabolic equation of the form
\begin{equation}\label{eq:intro_rescaled_RDTF}
 \partial_\tau u_\tau=Lu_\tau+Q[u_\tau],
 \qquad
 Lu=\Delta_fu+2\Rm_{g_{\cyl}}(u),
\end{equation}
where $\Delta_f=\Delta_{g_{\cyl}}-\nabla_{\nabla f}$ and $Q$ vanishes to second order at the origin.

We have thus transformed the weakly parabolic Ricci flow equation into the strictly parabolic equation~\eqref{eq:intro_rescaled_RDTF}, in a form suitable for applying the PDE--ODI principle.
However, the HMHF is itself nonlinear, and we must ensure that its maps remain diffeomorphisms throughout the prescribed interval.
Moreover, since the original flow is only locally close to a cylinder, the HMHF must be constructed \emph{locally.}
We therefore must construct the harmonic map heat flow solution in tandem with the analysis of the associated Ricci--DeTurck flow perturbation $u_\tau$ via the PDE--ODI principle.
The following theorem illustrates this simultaneous construction and control.

\begin{Theorem}[Theorem~\ref{thm:HMHF}, simplified]\label{thm:intro_HMHF}
For any $m\geq4$, $J\geq2$, and $\lambda<0$, the following holds if $0<\eta\leq\ov\eta(\lambda,J)$ and $0<\delta\leq\ov\delta(m,J,\lambda,\eta)$.
Let $(M,(g_t)_{t\in I})$ be a Ricci flow with complete time-slices and bounded curvature on a compact time-interval $I\subset(-\infty,0)$, and let $p\in M$.
Assume $(p,t)$ is a center of an $(n,k,\delta)$-neck at scale $\sqrt{-t}$ for all $t\in I$.
Then there is a harmonic map heat flow $\chi_t:\Omega_t\longrightarrow\IR^k\times\IS^{n-k}$, defined on a family of open subsets $\Omega_t\subset M$ containing $p$, for $t\in I$.
The maps $\chi_t$ are diffeomorphisms onto their images and keep the rescaled image of the basepoint close to the center of the cylinder, in the sense that
\[
 \xi_{e^{\tau/2}}\circ\chi_t(p)\in\IB^k_\eta\times\IS^{n-k},\qquad \tau=-\log(-t).
\]
There is also a continuous radius function $R(\tau)\in[100,\infty]$, which is $C^1$ wherever finite.
The perturbation $u_\tau:=e^\tau(\xi_{e^{\tau/2}})_*(\chi_t)_*g_t-g_{\cyl}$ satisfies \eqref{eq:intro_rescaled_RDTF} and
\begin{align*}
 \IB^k_{R(\tau)}\times\IS^{n-k}
 &\subset \xi_{e^{\tau/2}}\circ\chi_t(\Omega_t),\qquad
 \|u_\tau\|_{C^m(\IB^k_{R(\tau)}\times\IS^{n-k})}\leq\eta.
\end{align*}

Let $\sV_{>\lambda}$ be the finite-dimensional sum of eigenspaces of $L$ with eigenvalue greater than $\lambda$. For a suitable small constant $\eps=\eps(m,J,\lambda,\eta)\in(0,1)$ and a corresponding cutoff function $\omega_{R(\tau)}$ supported in the cylindrical region, set
\[
 U^+(\tau):=\PP_{\sV_{>\lambda}}(u_\tau\omega_{R(\tau)}),\qquad
 \mathcal U^-(\tau):=\big\|\PP_{\sV_{\leq\lambda}}(u_\tau\omega_{R(\tau)})\big\|_{L^2_f}
 +\eta e^{-\frac{(1-\eps)^2R^2(\tau)}{8}}.
\]
Here we set $\omega_{R(\tau)} \equiv 1$ if $R(\tau)=\infty$. Then
\begin{align}
 \big\|\partial_\tau U^+-LU^+-Q_J^+(U^+)\big\|_{L^2_f}
 &\leq C(\lambda,J)\|U^+\|_{L^2_f}^{J+1}+C(\lambda)\sqrt\eta\,\mathcal U^-,
 \label{eq:intro_ODI_plus}\\
 \partial_\tau\mathcal U^-
 &\leq(\lambda+C(\lambda)\sqrt\eta)\mathcal U^-
  +\|Q_J^-(U^+)\|_{L^2_f}
  +C(\lambda,J)\|U^+\|_{L^2_f}^{J+1}.
\label{eq:intro_ODI_minus}
\end{align}
\end{Theorem}

\medskip

The inequalities \eqref{eq:intro_ODI_plus}--\eqref{eq:intro_ODI_minus} have the same coarse structure as those in our mean curvature flow analysis~\cite{Bamler_Lai_MCF1} and can be studied using similar standard ODE methods.
At a finer level, however, the equations reflect the different geometric setting of Ricci flow.
For example, we now work with symmetric $2$-tensors, and the spectrum of $L$ and the nonlinear terms $Q_J^\pm$ differ from those in mean curvature flow.
These differences already appear in the decomposition
\[
U^+=U^+_{<0}+U^{++}+U^+_{\aff}+U^+_{\conf}.
\]
Here $U^+_{<0}$ consists of the stable modes, with negative eigenvalues.
The component $U^{++}$, which we view as the leading component of $U^+$, consists of the constant, linear, and quadratic multiples of $g_{\IS^{n-k}}$ in the axial coordinates and will play a central role in the subsequent analysis.
The affine modes $U^+_{\aff}$ are constant symmetric tensors in the Euclidean directions, arising from infinitesimal linear changes of coordinates, and have eigenvalue zero.
The conformal modes $U^+_{\conf}$ arise from infinitesimal conformal changes of the spherical coordinates and are neutral when $n-k=2$ and unstable when $n-k\geq3$.
These modes are not geometrically significant, but reflect coordinate freedom that remains after passing to the Ricci--DeTurck gauge.
Controlling them is crucial for the construction of the harmonic map heat flow.

As in the mean curvature flow setting, we study the resulting ODI by dividing the time interval into subintervals according to which modes dominate.
After an initial interval in which the stable modes and localization error decay exponentially, all other modes, including the affine and conformal components, and the localization error are controlled by $U^{++}$.
The latter interval further divides into an interval dominated by the quadratic modes and one dominated by the constant and linear modes.
This partition is already recorded in Theorem~\ref{thm:HMHF} because it enters the proof, where the initial gauge is chosen so that the affine and conformal components vanish at suitable transition times.
In the ancient case, bootstrapping these ODIs gives higher-order control of the remaining modes and a finite-dimensional approximation with metric and evolution errors bounded by any prescribed power of $\|U^{++}\|$.
The affine and conformal components require additional estimates, but these too come from an ODE analysis.

\medskip

Let us briefly outline the proof of Theorem~\ref{thm:HMHF} (Theorem~\ref{thm:intro_HMHF}).
We first construct the local HMHF using a grafting argument, following Bamler--Kleiner~\cite{Bamler_Kleiner_uniqueness} and Lai~\cite[Section~6]{Lai_O2_steady}.
This allows us to solve and continue the HMHF as long as the perturbation $u_\tau$ remains small on a fixed large cylindrical region and the rescaled image of the basepoint stays in a fixed central region.
However, the affine and conformal modes, although not geometrically significant, may grow and make $u_\tau$ too large for further continuation.
The rescaled image of the basepoint may also escape the controlled region.
Either possibility would prevent us from continuing the HMHF and hence from applying the PDE--ODI principle throughout the prescribed interval.

To address these issues, we consider a finite-dimensional family of initial conditions for the harmonic map heat flow, parametrized by affine and conformal changes of coordinates, together with translations to control the basepoint.
We then apply the PDE--ODI principle to each induced Ricci--DeTurck flow, obtaining \eqref{eq:intro_ODI_plus}--\eqref{eq:intro_ODI_minus} and the resulting mode estimates and transition times.
To select a suitable initial condition, we use a Wa\.zewski-type box argument, as in~\cite{Stolarski_AC_shrinkers}.
More precisely, the mode estimates give a continuous map from the parameter box whose normalized boundary map has degree one.
This map must therefore have a zero, giving the desired normalization of the affine and conformal components at the transition times and of the basepoint.
For this choice, $U^{++}$ controls the other modes after the initial decay interval.
Since $U^{++}$ is geometrically significant, the $\delta$-cylindricality assumption bounds its size.
This bound and the dominance of $U^{++}$, together with the initial decay estimates, keep $u_\tau$ small.
With the basepoint also controlled, this rules out an early exit and ensures that the HMHF exists throughout the prescribed interval.

\medskip

The method of constructing a global HMHF to convert an ancient Ricci flow into Ricci--DeTurck flow is also used by Choi--Lai~\cite[Section~6.1]{Choi_Lai_HMHF} to classify exponentially converging ancient flows asymptotic to compact shrinkers. 
In their gauge construction, they control the perturbation by successively modifying the HMHF to remove the growing geometrically insignificant modes. The box argument in the present paper allows us to select the initial map once and for all, ensuring that the affine and conformal modes remain controlled throughout the interval.
Kim-Ozuch~\cite{Kim_Ozuch_blowdowns} used HMHF, a related box argument, and the PDE--ODI principle to prove dynamical stability of the four-dimensional blowdown shrinker.
For related box arguments in Ricci flow and mean curvature flow, see~\cite{Angenent_Velazquez_MCF_neckpinches,Angenent_Isenberg_Knopf_RF_neckpinches,Stolarski_AC_shrinkers,Lee_Zhao_AC_singularities,Chen_Lee_Sun_Zhao_tangent_flows}.

\subsection{Organization of the paper}

In Section~\ref{sec:preliminaries} we introduce Ricci flow spacetimes, the cylindrical model and the pseudolocality assumptions used throughout the paper.

In Section~\ref{sec:calculations} we derive the rescaled Ricci--DeTurck equation around the round cylinder, compute its linearization and spectrum on symmetric $2$-tensors and calculate the quadratic and cubic terms needed for the mode analysis.

In Section~\ref{sec:gauging} we develop the local HMHF construction. We establish the required local Ricci--DeTurck estimates and the continuous dependence of the HMHF on the initial data.

Section~\ref{sec:pde_odi} contains the main gauging argument. We combine HMHF with the PDE--ODI principle and carry out the finite-dimensional box argument. This produces the global gauge from Theorem~\ref{thm:intro_HMHF}.

In Section~\ref{sec:higher_mode_analysis} we analyze the mode inequalities, obtain higher-order estimates and establish the linear--quadratic alternatives underlying Theorem~\ref{thm:intro_dichotomy}. This section also establishes cylindrical rigidity, neck propagation, and uniqueness of cylindrical blow-ups and blow-downs, summarized in Theorems~\ref{thm:intro_cylindrical_rigidity}, \ref{thm:intro_neck_stability} and~\ref{thm:intro_cylindrical_limits}.

In Section~\ref{sec:dominant_linear} we treat the dominant linear mode. We determine the asymptotic constants, construct the gauge at infinity and prove the exponential Bryant asymptotics in Theorem~\ref{thm:intro_bryant_infinity} and the classification in Theorem~\ref{thm:intro_linear}.

In Section~\ref{sec:dominant_quadratic} we treat the dominant quadratic mode. We define and analyze the invariant $\Qu$ in Theorem~\ref{thm:intro_quadratic}, show that it determines the asymptotics to arbitrary polynomial order in Theorem~\ref{thm:intro_higher_order_asymptotics} and construct the $\kappa$-solutions in Theorem~\ref{thm:intro_realization}.

Throughout the paper we regard $n$ and $k$ as fixed and usually suppress the dependence of constants on them.

\medskip

\noindent\textbf{Statement on the use of AI.}
The main results of this paper were established in preliminary form in 2024, although the proofs were not yet organized in terms of the PDE--ODI principle.
We spent the following months digesting and reorganizing these arguments, which led to the formulation of the general PDE--ODI principle.
We then set this project aside to work out the resolution of the mean convex neighborhood conjecture in~\cite{Bamler_Lai_MCF1,Bamler_Lai_MCF2}, as mean curvature flow offered a more direct setting in which to develop and present this approach.
After completing that work, we returned to the present paper and recast its proofs using the PDE--ODI principle and the language and terminology developed in the mean curvature flow papers.
A draft of the manuscript was nearly complete by April 2026.

In completing the manuscript, AI assistance was used to polish the exposition and restructure a few lemmas.
All revisions were independently verified by the authors, who take full responsibility for the manuscript.

\section{Preliminaries}\label{sec:preliminaries}
\subsection{Ricci flow spacetimes}
In this paper we formulate our results in terms of Ricci flow spacetimes, which include conventional Ricci flows as a special case.
This formulation also allows us to consider singular Ricci flows, i.e., metric flows in the sense of~\cite{Bamler_RF_compactness,Bamler_HK_RF_partial_regularity}, which are of particular interest as they often arise as $\IF$-limits of blow-up sequences.
For our purposes, we do not need the full structure of these limits.
It suffices to work with their regular parts, which are Ricci flow spacetimes.
We now introduce this notion and formulate the additional properties needed for our analysis.

We use the definition of a Ricci flow spacetime from~\cite[Definition~5.1]{Bamler_Kleiner_uniqueness}, allowing arbitrary dimension and time interval.

\begin{Definition}[Ricci flow spacetime]
An $(n+1)$-dimensional \textbf{Ricci flow spacetime} over a time-interval $I \subset \IR$ is a tuple $(\MM, \mathfrak t, \partial_{\mathfrak t}, g)$ with the following properties:
\begin{enumerate}
\item $\MM$ is a smooth $(n+1)$-dimensional manifold with (smooth) boundary $\partial \MM$.
\item $\mathfrak t : \MM \to I$ is a smooth function without critical points (called \textbf{time function}). For any $t \in I$ we denote by $\MM_t := \mathfrak t^{-1} (t) \subset \MM$ the \textbf{time-$t$-slice} of $\MM$.
\item $\mathfrak t^{-1}(\partial I) = \partial \MM$, i.e., the boundary of $\MM$ equals the union of initial and final time-slices (if they exist).
\item $\partial_{\mathfrak t}$ is a smooth vector field (the \textbf{time vector field}) on $\MM$ that satisfies $\partial_{\mathfrak t} \mathfrak t \equiv 1$.
\item $g$ is a smooth inner product on the spatial subbundle $\ker (d\mathfrak t) \subset T \MM$.
For any $t \in I$ we denote by $g_t$ the restriction of $g$ to the time-$t$-slice $\MM_t$; so $(\MM_t, g_t)$ is a Riemannian manifold.
\item $g$ satisfies the Ricci flow equation: $\mathcal L_{\partial_{\mathfrak t}} g = - 2 \Ric_g$.
Here $\Ric_g$ denotes the symmetric $(0,2)$-tensor on $\ker (d\mathfrak t)$ that restricts to the Ricci tensor of $(\MM_t, g_t)$ for all $t \in I$.
\end{enumerate}
If $p \in \MM_{t_0}$ is a point and the trajectory of $\partial_{\mathfrak t}$ is defined at time $t \in I$, then we say that $p$ \textbf{survives} until time $t$ and we denote the corresponding point on this trajectory by $p(t) \in \MM_t$.
The trajectory $t \mapsto p(t)$ is also called a \textbf{worldline.}
We call $\MM$ \textbf{ancient} if $I=(-\infty,T)$ or $I=(-\infty,T]$, \textbf{immortal} if $I=(T,\infty)$ or $I=[T,\infty)$, and \textbf{eternal} if $I=\IR$. 
We refer to \cite{Bamler_Kleiner_uniqueness} for further elementary terminology relating to Ricci flow spacetimes.
\end{Definition}

The coordinate maps used in this paper evolve by harmonic map heat flow.

\begin{Definition}\label{def:HMHF}
Let $\Omega \subset \MM$ be an open subset of a Ricci flow spacetime over $I$, and let $(h_t)_{t\in I}$ be a smooth family of Riemannian metrics on $N$.
A smooth map $\chi:\Omega\to N$ is a \textbf{harmonic map heat flow (HMHF)} if, for every $t$, its restriction $\chi_t := \chi |_{\Omega_t}$  satisfies
\[
\partial_{\mathfrak t} \chi := d\chi(\partial_{\mathfrak t})\big|_{\Omega_t}=\Delta_{g_t,h_t}\chi_t.
\]
\end{Definition}

We recall that a conventional Ricci flow $(M,(g_t)_{t \in I})$ can be viewed as a Ricci flow spacetime by setting $\MM := M \times I$, letting $\mathfrak t$ and $\partial_{\mathfrak t}$ be the projection and unit vector field of the last coordinate and endowing the time-slices $\MM_t = M \times \{ t \} \cong M$ with the metric $g_t$.
Throughout the paper, the flows used to form $\IF$-limits are assumed to be complete Ricci flows with bounded curvature on compact time-intervals; see~\cite[Appendix~A]{Bamler_2021_fundamental_group} for the noncompact case.
A metric flow $\mathcal X$ arising as a non-collapsed $\IF$-limit of a sequence of conventional Ricci flows has a canonical decomposition $\mathcal X = \mathcal R \dotcup \mathcal S$ into its regular and singular parts.
The regular part $\mathcal R$ is a Ricci flow spacetime and determines
$\mathcal X$ as a metric flow on the interior of its time
interval~\cite[Theorems~2.4 and~15.29]{Bamler_HK_RF_partial_regularity}.

We emphasize that working with a \emph{general} Ricci flow spacetime $\MM$ typically requires additional assumptions or geometric context, since $\MM$ may exhibit various forms of incompleteness: its time-slices may be incomplete as Riemannian manifolds, and its worldlines may cease to exist forward or backward in time.
In~\cite{Bamler_Kleiner_uniqueness}, these issues were addressed by imposing additional completeness and canonical neighborhood conditions.
For our purposes, a surprisingly weak pseudolocality assumption suffices.
This assumption holds in many settings of interest, for example for the spacetimes considered in~\cite{Bamler_Kleiner_uniqueness} and for regular parts of non-collapsed $\IF$-limits.

To introduce this property, we first recall the following pseudolocality result for conventional Ricci flows due to Perelman~\cite{Perelman1}.

\begin{Theorem}[Perelman's pseudolocality theorem]\label{thm:Perelman_pseudolocality}
There is a dimensional constant $\eps_{\mathrm{pseudo},n}>0$ with the following property.
Consider an $n$-dimensional Ricci flow $(M,(g_t)_{t\in I})$ on a compact manifold and let $(x_0,t_0)\in M\times I$ and $r>0$ be such that
\[
 \sup_{B_{g_{t_0}}(x_0,r)}|\Rm_{g_{t_0}}|<r^{-2}
 \qquad\text{and}\qquad
 \big|B_{g_{t_0}}(x_0,r)\big|_{g_{t_0}}>\frac12\omega_n r^n,
\]
where $\omega_n=|\IB_1^n|$.
Then
\[
 |\Rm_{g_t}|\leq(\eps_{\mathrm{pseudo},n}r)^{-2}
 \quad\text{on}\quad
 B_{g_t}(x_0,\eps_{\mathrm{pseudo},n}r)
\]
for all $t\in[t_0,t_0+(\eps_{\mathrm{pseudo},n}r)^2]\cap I$.
\end{Theorem}

Let us fix the constant $\eps_{\mathrm{pseudo},n}$ for the remainder of this paper.
We define the analogous property for Ricci flow spacetimes.

\begin{Definition}[Pseudolocality property]\label{def:pseudolocality_property}
Let $\MM$ be an $(n+1)$-dimensional Ricci flow spacetime over the time-interval $I$.
We say that $\MM$ satisfies the \textbf{pseudolocality property} if the following holds.
Suppose that, for some $x_0\in\MM$ with $t_0=\mathfrak t(x_0)$ and some $r>0$,
\[
 \sup_{B_{g_{t_0}}(x_0,r)}|\Rm_{g_{t_0}}|<10r^{-2}
 \qquad\text{and}\qquad
 \big|B_{g_{t_0}}(x_0,r)\big|_{g_{t_0}}>\frac1{20}\omega_n r^n.
\]
Then $B_{g_{t_0}}(x_0,r)$ survives until every time in
\[
 [t_0,t_0+(\tfrac1{10}\eps_{\mathrm{pseudo},n}r)^2]\cap I,
\]
and
\[
 |\Rm_{g_t}|\leq(\tfrac1{10}\eps_{\mathrm{pseudo},n}r)^{-2}
\]
on the time-$t$ image of
$B_{g_{t_0}}(x_0,\tfrac1{10}\eps_{\mathrm{pseudo},n}r)$
for every such time $t$.
\end{Definition}

Throughout the paper, unless otherwise stated, Ricci flow spacetimes
are assumed to satisfy the pseudolocality property in
Definition~\ref{def:pseudolocality_property}.

We note that the pseudolocality property is closed under smooth convergence.
Ricci flow spacetimes arising in the standard settings commonly possess this property.
For example, it follows from Perelman's work that the unique $(3+1)$-dimensional Ricci flow spacetime from~\cite{Bamler_Kleiner_uniqueness}, starting from compact---and also most noncompact---initial data, satisfies the pseudolocality property.

The regular part of a noncollapsed $\IF$-limit of Ricci flows is a Ricci flow spacetime by~\cite[Theorem~9.12(a)]{Bamler_RF_compactness}. Convergence is smooth on the regular part~\cite[Theorem~2.5]{Bamler_HK_RF_partial_regularity}, so the pseudolocality property follows by passing to the limit in Theorem~\ref{thm:Perelman_pseudolocality}.

\subsection{The round cylinder and necks}
Throughout the paper we denote the round sphere by
\[
\bigl(\IS^n,g_{\IS^n}\bigr):=\bigl(\IS^n(1),2(n-1)g_{\IS^n(1)}\bigr),
\]
where $g_{\IS^n(1)}$ is the unit round metric.
Thus $\sec_{g_{\IS^n}}=\frac{1}{2(n-1)}$, $\Ric_{g_{\IS^n}}=\frac12 g_{\IS^n}$ and $R_{g_{\IS^n}}=\frac n2$.
With this normalization, $\IS^n$ is the spherical factor of the time-$(-1)$ slice of the round shrinking cylinder that becomes extinct at time zero.

\begin{Definition}\label{def:cylindrical_model}
Fix $1\leq k\leq n-2$ and a finite subgroup $\Gamma\subset O(n-k+1)$ acting freely on $\IS^{n-k}$.
The cylindrical model is
\[
M_{\cyl,\Gamma}:=\IR^k\times(\IS^{n-k}/\Gamma),\qquad g_{\cyl,\Gamma}:=g_{\IR^k}+g_{\IS^{n-k}/\Gamma},
\]
where $g_{\IS^{n-k}/\Gamma}$ is the quotient metric.
The corresponding shrinking Ricci flow is $g_{\cyl,\Gamma,t}:=g_{\IR^k}+(-t)g_{\IS^{n-k}/\Gamma}$ for $t<0$, and the standard cylindrical potential is $f(\bx,\by):=\frac14|\bx|^2$.
We omit $\Gamma$ from the notation when it is trivial.

Unless otherwise indicated, pointwise tensor inner products, covariant derivatives and $C^m$-norms on the model are taken with respect to $g_{\cyl,\Gamma}$. For tensor fields $u,v$ of the same type, set
\[
\langle u,v\rangle_{L^2_f}:=|\Gamma|\int_{M_{\cyl,\Gamma}}\langle u,v\rangle_{g_{\cyl,\Gamma}}e^{-f}\,dg_{\cyl,\Gamma},\qquad \|u\|_{L^2_f}^2:=\langle u,u\rangle_{L^2_f}.
\]
We write $L^2_f$ for the space of tensor fields with finite norm. On a subset of the cylinder, the integral is restricted to that subset. The factor $|\Gamma|$ makes this norm equal to the $L^2_f$-norm of the invariant lift to $\IR^k\times\IS^{n-k}$.
\end{Definition}

We use the following notion of a neck to quantify closeness to the cylindrical model.

\begin{Definition}\label{def:neck}
Let $(M,g)$ be an $n$-dimensional Riemannian manifold, $p\in M$, $r>0$ and $0<\delta<1$. We call $p$ a \textbf{center of an $(n,k,\delta,\Gamma)$-neck at scale $r$} if there is a diffeomorphism
\[
\Phi:\IB^k_{\delta^{-1}}\times(\IS^{n-k}/\Gamma)\longrightarrow U\subset M
\]
onto an open neighborhood of $p$ such that
\[
p\in\Phi\big(\{0\}\times(\IS^{n-k}/\Gamma)\big),\qquad \|r^{-2}\Phi^*g-g_{\cyl,\Gamma}\|_{C^{\delta^{-1}}(\IB^k_{\delta^{-1}}\times(\IS^{n-k}/\Gamma))}<\delta.
\]
Here $C^{\delta^{-1}}$ means $C^{\lfloor\delta^{-1}\rfloor}$.
We call $U$ an \textbf{$(n,k,\delta,\Gamma)$-neck} and omit $\Gamma$ when it is trivial.
A point $p$ in a Ricci flow spacetime $\MM$ is a center of an $(n,k,\delta,\Gamma)$-neck at scale $r$ if it is one in its time-slice $(\MM_t,g_t)$, where $t=\mathfrak t(p)$.
\end{Definition}

\subsection{Asymptotically cylindrical Ricci flows}
\label{subsec:asymptotically_cylindrical}
Our main concern will be Ricci flows that, after parabolic rescaling, converge as $t \to -\infty$ to a round cylinder of the form $\IR^k \times (\IS^{n-k}/\Gamma)$, where $\Gamma \subset O(n-k+1)$ is a discrete subgroup acting freely on $\IS^{n-k}$.
We will use the following terminology.

\begin{Definition}[Asymptotically cylindrical Ricci flow]
\label{defn:asymptotically_cylindrical}
Let $\MM$ be an ancient Ricci flow spacetime with $n$-dimensional time-slices, let
$1\leq k\leq n-2$, and let $\Gamma \subset O(n-k+1)$ be a discrete subgroup acting freely on $\IS^{n-k}$. We call $\MM$ \textbf{asymptotically
$(n,k,\Gamma)$-cylindrical} (omitting $\Gamma$ when it is trivial) if there is a worldline $t \mapsto p(t)$ defined for all
sufficiently negative $t$ such that
\begin{equation}\label{eq:asymptotically_cylindrical_convergence}
\bigl(\MM_t,(-t)^{-1}g_t,p(t)\bigr)
\longrightarrow
\bigl(\IR^k\times(\IS^{n-k}/\Gamma),
g_{\IR^k}+g_{\IS^{n-k}/\Gamma},o\bigr)
\qquad\textnormal{as }t\to-\infty
\end{equation}
in the pointed smooth Cheeger--Gromov topology, for some point $o$
of the limiting cylinder. Here $g_{\IS^{n-k}/\Gamma}$ is the round
metric of sectional curvature $\frac{1}{2(n-k-1)}$.
Moreover, we require $\MM$ to be connected and for any two worldlines $t\mapsto p_1(t),p_2(t)$ satisfying~\eqref{eq:asymptotically_cylindrical_convergence}, we require an asymptotic bound of the form $d_{g_t}(p_1(t),p_2(t))\leq C(p_1,p_2)\sqrt{-t}$ for all sufficiently negative $t$.
\end{Definition}

For ancient non-collapsed $\IF$-limits, Lemma~\ref{lem:asymptotically_cylindrical_tangent_flows} gives an equivalent formulation in terms of cylindrical tangent flows at infinity.

The last condition expresses that there is only one cylindrical end at backward infinity.
It is automatic for connected conventional Ricci flows with complete time-slices and bounded curvature on compact time-intervals.
Indeed, the curvature bounds near $p_1(t),p_2(t)$ supplied by~\eqref{eq:asymptotically_cylindrical_convergence} and Perelman's distance distortion estimate~\cite[Lemma~8.3(b)]{Perelman1} imply for sufficiently negative $t$
\[
\frac{d}{dt}d_{g_t}(p_1(t),p_2(t))\geq-\frac{C_n}{\sqrt{-t}}.
\]
Integrating from $t$ to a fixed sufficiently negative time $t_0$ yields
\[
d_{g_t}(p_1(t),p_2(t))
\leq d_{g_{t_0}}(p_1(t_0),p_2(t_0))+2C_n\bigl(\sqrt{-t}-\sqrt{-t_0}\bigr)
\leq C(p_1,p_2)\sqrt{-t}.
\]

Definition~\ref{defn:asymptotically_cylindrical} is the Ricci flow analogue of the notion of asymptotically cylindrical mean curvature flow studied and classified in~\cite{Bamler_Lai_MCF1, Bamler_Lai_MCF2}.

\section{The rescaled Ricci--DeTurck equation near the cylinder}\label{sec:calculations}
In this section, we collect the calculations for the rescaled Ricci--DeTurck equation near the cylinder. First, Lemma~\ref{Lem_Q_precise_formula} gives explicit formulas for the linearized operator $L$ and the nonlinear remainder $Q$. Lemma~\ref{Lem_cyl_spectrum} identifies the nonnegative eigenspaces of $L$ and the negative oscillatory spectral gap. We then compute the quadratic Taylor polynomial $Q_2$ and its neutral projection for neutral inputs in Lemma~\ref{l:Q_2}, and evaluate it explicitly on rotationally symmetric tensors in Lemma~\ref{Lem_Q2_computation}. Finally, Lemma~\ref{Lem_Q3_computation} computes the neutral spherical projection of $Q_3-Q_2$, which is used in the refined asymptotics of Section~\ref{sec:dominant_quadratic} (quadratic dominance case).

We compute on the round cylinder, identifying quotient tensors with their $\Gamma$-invariant lifts.

Consider Ricci flows $(M,g_t)$ and $(N,h_t)$ on a time-interval $I$, and an HMHF $\chi_t$ between them whose time-slice maps are diffeomorphisms onto their images. The pushed-forward metrics $g'_t:=(\chi_t)_*g_t$ satisfy
\[
\partial_t g'_t=-2\Ric_{g'_t}-\LL_{X_{h_t}[g'_t]}g'_t,
\]
where $X_h[g]:=\triangle_{g,h}\id$ is the DeTurck vector field; see~\cite[Section~6.5 and Appendix~A.1]{Bamler_Kleiner_uniqueness}.

For the cylindrical background, take $(N,h_t)=(M_{\cyl},g_{\cyl,t})$ and $I\subset(-\infty,0)$. Set $\tau=-\log(-t)$, let $\xi_\lambda(\bx,\by):=(\lambda\bx,\by)$, and define
\[
\td g_\tau:=e^\tau(\xi_{e^{\tau/2}})_*g'_{-e^{-\tau}}.
\]
The rescaled metrics satisfy
\begin{equation} \label{eq_resc_RdT}
   \partial_\tau \td g_\tau = - 2 \Ric_{\td g_\tau} + \td g_{\tau} - \LL_{X_{g_{\cyl}}[\td g_\tau] + \nabla f} \td g_\tau 
\end{equation}
where $f=\frac14|\bx|^2$ and the gradient is taken with respect to $g_{\cyl}$. The identity $\Ric_{g_{\cyl}}+\nabla^2f=\frac12g_{\cyl}$ shows that $g_{\cyl}$ is a stationary solution of \eqref{eq_resc_RdT}.

The following lemma gives the equation for the perturbation $u_\tau:=\td g_\tau-g_{\cyl}$.

\begin{Lemma} \label{Lem_Q_precise_formula}
If we set $\td g_\tau = g_{\cyl} + u_\tau$, then \eqref{eq_resc_RdT} is equivalent to
\[ \partial_\tau u_\tau = L u_\tau + Q[u_\tau]. \]
The linear part is
\begin{equation} \label{eq_L_def}
   L u = \triangle_f u + 2\Rm_{g_{\cyl}} (u), 
\end{equation}
where $\triangle_f u = \triangle u - \nabla_{\nabla f} u$ and $(\Rm_{g_{\cyl}}(u))_{ij} = g^{kl} R_{kij}^{\quad s} u_{ls}$.
If $\{ e_i \}_{i=1}^n$ is a local orthonormal frame on $\IR^k \times \IS^{n-k}$ with respect to $g_{\cyl}$ with the property that $e_1, \ldots, e_k$ are parallel to the first factor and $e_{k+1}, \ldots, e_n$ are parallel to the second factor, then
\begin{equation} \label{eq_Lu_explicit}
   (Lu)_{ij} = \begin{cases}
\triangle_f u_{ij} &\text{if $i \leq k$ or $j \leq k$} \\
\triangle_f u_{ij} + \frac1{n-k-1} (\sum_{s=k+1}^n u_{ss}) \delta_{ij} - \frac1{n-k-1}  u_{ij} &\text{if $i > k$ and $j > k$}
\end{cases} 
\end{equation}
The non-linear part is
\begin{align}
(Q[u])_{ij} &= (g+u)^{pq} (g+u)^{st} \big(  \tfrac12  \nabla_i u_{ps} \nabla_j u_{qt} - \nabla_p u_{si} \nabla_q u_{tj} + \nabla_p u_{si} \nabla_t u_{qj} \notag \\
  &\qquad\qquad\qquad\qquad\qquad\qquad\qquad\qquad\qquad   -  \nabla_i u_{ps} \, \nabla_q u_{tj}    - \nabla_j u_{ps} \, \nabla_q u_{ti}     \big) \notag \\
  &\qquad +( (g+u)^{pq} -g^{pq} )  ( \nabla^2_{pq} u_{ij} - R_{ipq}^{\quad s} u_{sj} - R_{ipj}^{\quad s} u_{qs} - R_{jpi}^{\quad s} u_{sq} - R_{jpq}^{\quad s} u_{is}). \label{eq_Q_def}
\end{align}
Here the Laplacian, the curvature terms $R_{ijp}^{\quad l}$ and all covariant derivatives are taken with respect to the metric $g_{\cyl}$. 
\end{Lemma}

\begin{proof}
We verify the lemma at some fixed point $z \in \IR^k \times \IS^{n-k}$.
Choose local normal coordinates $(x^1, \ldots, x^n)$ near $z$ such that $x^1, \ldots, x^k$ depend on the $\IR^k$-component and $x^{k+1}, \ldots, x^n$ depend on the $\IS^{n-k}$-component and set $e_i := \frac{\partial}{\partial x^i}$.
Recall that $\{ e_i \}$ is orthonormal at $z$ and $\nabla_{e_i} e_j = 0$ at $z$, where $\nabla$ denotes covariant derivatives with respect to $g_{\cyl}$.
Let us drop the time-index and write $g := g_{\cyl}$ and $\td g =: g + u$.
The covariant derivative $\td\nabla$ with respect to $\td g$ can be expressed as follows:
\begin{align*} 
\td\nabla_{e_i} e_j 
&= \tfrac12 \td g^{pq} \big( e_i \td g(e_j, e_p) + e_j \td g(e_i , e_p) - e_p \td g(e_i , e_j) \big) e_q \\
&= \tfrac12 \td g^{pq} \big( (\nabla_{e_i} \td g)(e_j,e_p) + \td g( \nabla_{e_i} e_j, e_p) + \td g(e_j, \nabla_{e_i} e_p )\\
&\qquad\qquad + (\nabla_{e_j} \td g)(e_i,e_p) + \td g(\nabla_{e_j} e_i, e_p) + \td g(e_i, \nabla_{e_j} e_p) \\
&\qquad\qquad - (\nabla_{e_p} \td g)(e_i, e_j) - \td g(\nabla_{e_p} e_i, e_j) - \td g(e_i, \nabla_{e_p} e_j ) \big) e_q \\
&= \nabla_{e_i} e_j + \tfrac12 (g+u)^{pq} \big( (\nabla_{e_i} u)(e_j,e_p) + (\nabla_{e_j} u)(e_i,e_p) - (\nabla_{e_p} u)(e_i,e_j) \big) e_q \\&
=: \nabla_{e_i} e_j + (g+u)^{pq} \gamma_{ijp} e_q .
\end{align*}
Note that $(g+u)^{pq} \gamma_{ijq}$ are in fact the differences between the Christoffel symbols of $\td g$ and $g$.
We obtain the following identity at $z$:
\[ \td\nabla_{e_a} \td\nabla_{e_i} e_j = \nabla_{e_a} \nabla_{e_i} e_j + (g+u)^{pq} \partial_a \gamma_{ijp} e_q +  (g+u)^{pq} (g+u)^{st} (  -  \nabla_a u_{ps} \,\gamma_{ijt} +  \gamma_{ijs} \gamma_{atp} ) e_q, \]
so the curvature $R_{\td g}$ of $\td g$ satisfies
\begin{multline*}
  R_{\td g}(e_a,e_i)e_j = R(e_a, e_i)e_j + (g+u)^{pq} (\partial_a \gamma_{ijp}-\partial_i \gamma_{ajp}) e_q \\
 +  (g+u)^{pq} (g+u)^{st} (  -  \nabla_a u_{ps} \,\gamma_{ijt} +   \nabla_i u_{ps} \,\gamma_{ajt}  +  \gamma_{ijs} \gamma_{atp} -   \gamma_{ajs} \gamma_{itp} ) e_q, \end{multline*}
which implies that its Ricci tensor is
\begin{multline} \label{eq_Ricij}
 (\Ric_{\td g})_{ij} = \Ric_{ij} +  (g+u)^{pq} (\partial_q \gamma_{ijp}-\partial_i \gamma_{qjp}) \\
 +  (g+u)^{pq} (g+u)^{st} (  -  \nabla_q u_{ps} \,\gamma_{ijt} +   \nabla_i u_{ps} \,\gamma_{qjt}  +  \gamma_{ijs} \gamma_{qtp} -   \gamma_{qjs} \gamma_{itp} ) . 
\end{multline}

Next, the following is true near $z$:
\[ X:=X_{g}[\td g] = \triangle_{\td g, g} \id = \td g^{pq} 
(  \nabla_{e_p} e_q - \td\nabla_{e_p} e_q )
= - (g+u)^{pq} (g+u)^{st} \gamma_{pqs} e_t.  \]
So we obtain that at $z$
\begin{align} 
(\LL_{X_{g}[\td g]} \td g)(e_i,e_j) &= X ( \td g(e_i, e_j)) - \td g ([X,e_i], e_j) - \td g(e_i, [X, e_j]) \notag \\
&= (\nabla_X \td g) (e_i, e_j) + \td g(\nabla_X e_i - [X, e_i], e_j) + \td g(e_i, \nabla_X e_j - [X,e_j]) \notag \\
&= (\nabla_X \td g) (e_i, e_j) + \td g (\nabla_{e_i} X, e_j) + \td g (e_i , \nabla_{e_j} X) \notag \\
&= -(g+u)^{pq} (g+u)^{st} \gamma_{pqs} \nabla_t u_{ij} \notag \\
&\qquad - \partial_i \big( (g+u)^{pq} (g+u)^{st} \gamma_{pqs} \big) (g+u)_{tj} - \partial_j \big( (g+u)^{pq} (g+u)^{st} \gamma_{pqs} \big) (g+u)_{ti} \notag \\
&= (g+u)^{pq} (g+u)^{st}  (- \nabla_t u_{ij} \, \gamma_{pqs}+  \nabla_i u_{ps} \, \gamma_{qtj} + \nabla_i u_{js} \, \gamma_{pqt} + \nabla_j u_{ps} \, \gamma_{qti} + \nabla_j u_{is} \, \gamma_{pqt})  \notag \\
&\qquad  - (g+u)^{pq} (\partial_i \gamma_{pqj} + \partial_j \gamma_{pqi})  \label{eq_LXij}
\end{align}

Combining \eqref{eq_Ricij} and \eqref{eq_LXij} implies that
\begin{align}
 (- 2 \Ric_{\td g} - \LL_{X_g[\td g]} \td g)_{ij}
&= - 2\Ric_{ij} + (g+u)^{pq} (g+u)^{st} (  2  \nabla_q u_{ps} \,\gamma_{ijt} -2   \nabla_i u_{ps} \,\gamma_{qjt}  -2  \gamma_{ijs} \gamma_{qtp} +2   \gamma_{qjs} \gamma_{itp}  \notag \\
&\qquad\qquad\qquad + \nabla_t u_{ij} \, \gamma_{pqs}-  \nabla_i u_{ps} \, \gamma_{qtj} - \nabla_i u_{js} \, \gamma_{pqt} - \nabla_j u_{ps} \, \gamma_{qti} - \nabla_j u_{is} \, \gamma_{pqt} ) \notag \\
&\qquad + (g+u)^{pq} ( -2 \partial_q \gamma_{ijp} + 2 \partial_i \gamma_{qjp} + \partial_i \gamma_{pqj} + \partial_j \gamma_{pqi} ). \label{eq_RicLXij}
\end{align}
The product terms from the first two lines in which the indices $i$ and $j$ occur in the same factor are
\begin{align}
(g+u)^{pq} &(g+u)^{st} ( 2 \nabla_q u_{ps} \,\gamma_{ijt}  -2  \gamma_{ijs} \gamma_{qtp} 
+ \nabla_t u_{ij} \, \gamma_{pqs} - \nabla_i u_{js} \, \gamma_{pqt}  - \nabla_j u_{is} \, \gamma_{pqt}) \notag \\
&= (g+u)^{pq} (g+u)^{st} \big( (2\nabla_q u_{pt} - 2 \gamma_{qtp}) \gamma_{ijs} + (\nabla_s u_{ij} - \nabla_i u_{js} - \nabla_j u_{is} ) \gamma_{pqt} \big) \notag \\
&= (g+u)^{pq} (g+u)^{st} (2\gamma_{pqt} \gamma_{ijs} - 2\gamma_{ijs} \gamma_{pqt} ) = 0. \label{eq_sameij_0}
\end{align}
The remaining terms in the first two lines of \eqref{eq_RicLXij} are
\begin{align}
  (g+u)^{pq} &(g+u)^{st}  (  -2   \nabla_i u_{ps} \,\gamma_{qjt}   +2   \gamma_{qjs} \gamma_{itp}   -  \nabla_i u_{ps} \, \gamma_{qtj} - \nabla_j u_{ps} \, \gamma_{qti}  ) \notag \\
  &= (g+u)^{pq} (g+u)^{st} (  -   \nabla_i u_{ps} \, \nabla_j u_{qt}   +2   \gamma_{qjs} \gamma_{itp}   -  \nabla_i u_{ps} \, \nabla_q u_{tj} \notag \\
  &\qquad\qquad\qquad\qquad\qquad + \tfrac12  \nabla_i u_{ps} \, \nabla_j u_{qt}   - \nabla_j u_{ps} \, \nabla_q u_{ti}  + \tfrac12 \nabla_j u_{ps} \, \nabla_i u_{qt}   ) \notag \\
  &= (g+u)^{pq} (g+u)^{st} (  2   \gamma_{qjs} \gamma_{itp}   -  \nabla_i u_{ps} \, \nabla_q u_{tj}    - \nabla_j u_{ps} \, \nabla_q u_{ti}     ) \notag \\
  &= (g+u)^{pq} (g+u)^{st} \big(  \tfrac12  \nabla_i u_{ps} \nabla_j u_{qt} - \nabla_p u_{si} \nabla_q u_{tj} + \nabla_p u_{si} \nabla_t u_{qj} \notag \\
  &\qquad\qquad\qquad\qquad\qquad   -  \nabla_i u_{ps} \, \nabla_q u_{tj}    - \nabla_j u_{ps} \, \nabla_q u_{ti}     \big). \label{eq_ijdiff}
\end{align}
The third line of \eqref{eq_RicLXij} can be simplified as follows
\begin{align}
 (g+u)^{pq} &( -2 \partial_q \gamma_{ijp} + 2 \partial_i \gamma_{qjp} + \partial_i \gamma_{pqj} + \partial_j \gamma_{pqi} ) \notag \\
 &=  (g+u)^{pq} \big( \nabla^2_{pq} u_{ij} + \tfrac12 (\nabla^2_{ij} u_{pq} - \nabla^2_{ji} u_{pq} ) + (\nabla^2_{ip} u_{qj} - \nabla^2_{pi} u_{qj}) + (\nabla^2_{jp} u_{iq} - \nabla^2_{pj} u_{iq}) \big) \notag \\
 &=  (g+u)^{pq} \big( \nabla^2_{pq} u_{ij} - \tfrac12 ( R_{ijp}^{\quad s} u_{sq} + R_{ijq}^{\quad s} u_{ps} ) \\
 &\qquad\qquad\qquad - (R_{ipq}^{\quad s} u_{sj} + R_{ipj}^{\quad s} u_{qs}) - (R_{jpi}^{\quad s} u_{sq} + R_{jpq}^{\quad s} u_{is}) \big) \notag \\
 &=   ( (g+u)^{pq} -g^{pq} )  ( \nabla^2_{pq} u_{ij} - R_{ipq}^{\quad s} u_{sj} - R_{ipj}^{\quad s} u_{qs} - R_{jpi}^{\quad s} u_{sq} - R_{jpq}^{\quad s} u_{is}) \notag \\
 &\qquad + g^{pq} \nabla^2_{pq} u_{ij}  + 2 g^{pq} R_{pij}^{\quad s} u_{qs} - g^{pq} \Ric_{ip} u_{qj} - g^{pq} \Ric_{jp} u_{iq}. \label{eq_ij_second}
 \end{align}
To see the last equality, one may assume without loss of generality that $u_{pq}$ and $g_{pq}$ are both diagonal.

Lastly, we calculate as in \eqref{eq_LXij} that
\begin{multline}
 (\LL_{\nabla f} \td g)(e_i, e_j)
= (\nabla_{\nabla f}\td g)(e_i, e_j) + \td g(\nabla_{e_i} \nabla f, e_j) + \td g(e_i, \nabla_{e_j} \nabla f) \\
= (\nabla_{\nabla f} u)_{ij} + 2 \nabla_{ij} f + g^{pq} \nabla^2_{ip} f \, u_{qj} + g^{pq} \nabla^2_{jp} f \, u_{iq}. \label{eq_Lnabf_ij}
\end{multline}
Combining \eqref{eq_RicLXij}, \eqref{eq_sameij_0}, \eqref{eq_ijdiff}, \eqref{eq_ij_second} and \eqref{eq_Lnabf_ij} implies that the RHS of \eqref{eq_resc_RdT} is equal to
\begin{multline*}
 (Lu + Q[u])_{ij} - 2 \Ric_{ij} - g^{pq} \Ric_{ip} u_{qj} - g^{pq} \Ric_{jp} u_{iq}  - 2 \nabla_{ij} f - g^{pq} \nabla^2_{ip} f \, u_{qj} - g^{pq} \nabla^2_{jp} f \, u_{iq} + g_{ij} + u_{ij} \\
= (Lu + Q[u])_{ij} - g^{pq} (\Ric_{ip}+ \nabla^2_{ip} f) u_{qj} - g^{pq} (\Ric_{jp} + \nabla^2_{jp} f) u_{iq}  + u_{ij}  =(Lu + Q[u])_{ij}
 \end{multline*}
where $Lu$ and $Q[u]$ are as in \eqref{eq_L_def} and \eqref{eq_Q_def}.
The first case in \eqref{eq_Lu_explicit} is clear since the curvature terms vanish in this case.
Suppose now that $i,j > k$ and only sum over indices $> k$.
Then
\begin{multline}
   (\Rm_{g_{\cyl}}(u))_{ij} 
= \sum_{s,t=k+1}^n R_{sij}^{\quad t} u_{st}
= \frac1{2(n-k-1)} \sum_{s,t=k+1}^n (g_{ij} g_{st} - g_{it} g_{js}) u_{st}  \\
= \frac1{2(n-k-1)} \Big( \sum_{s=k+1}^n u_{ss} g_{ij}  - u_{ij} \Big),
\end{multline}
which shows the second case.
\end{proof}
\bigskip

For the remainder of this paper, we will denote by $\sV_\la$ the eigenspace of $L$ corresponding to the eigenvalue $\la$.
We write $\sV_{> \la}, \sV_{\geq \la}, \sV_{< \la}, \sV_{\leq \la}$ to denote the closure of the sum of eigenspaces $\sV_{\la'}$ for all $\la' > \la$, $\la' \geq \la$, etc.

We next identify the nonnegative modes and the negative oscillatory spectral gap of $L$.

We use the Hermite polynomials
\[
\mathfrak p^{(0)}=1,\qquad \mathfrak p_i^{(1)}(\bx)=\frac{\bx_i}{\sqrt2},\qquad \mathfrak p_{ij}^{(2)}(\bx)=\frac{\bx_i\bx_j-2\delta_{ij}}{2\sqrt2}.
\]
Let $\sV_{\rot}$ be the closed subspace of tensors invariant under the standard $O(n-k+1)$-action on the spherical factor, and let $\sV_{\osc}:=\sV_{\rot}^{\perp}$ in $L^2_f$. On a quotient cylinder, these spaces are understood through their $\Gamma$-invariant lifts.

\begin{Lemma}\label{Lem_cyl_spectrum}
The only non-negative eigenvalues of $L$ are $1,\frac12, 1-\frac{n-k}{2(n-k-1)}, 0$ (if $n-k \geq 3$) and $1, \tfrac12, 0$ (if $n-k=2$).
We have 
\[ \sV_{\geq 0} = \sV_{1} \oplus \sV_{\frac12} \oplus \sV_{\sph,0} \oplus \sV_{\aff} \oplus \sV_{\con},  \]
where 
\begin{align*}   
\sV_{1} &= \{ c \, g_{\IS^{n-k}} \;\; : \;\; c \in \IR \}, \\
\sV_{\frac12} &= \bigg\{ \sum_{i=1}^k c_i \bx_i \, g_{\IS^{n-k}}  \;\; : \;\; c_i \in \IR \bigg\}, \\ 
\sV_{\sph,0} &= \bigg\{ \sum_{i,j=1}^k c_{ij} \mathfrak p^{(2)}_{ij} (\bx_1, \ldots, \bx_k) \, g_{\IS^{n-k}}  \;\; : \;\; c_{ij} \in \IR \bigg\} \\
 \sV_{\aff} &= \bigg\{ \sum_{i,j=1}^k h_{ij} \, d\bx_i \otimes d\bx_j \;\; : \;\; h_{ij} \in \IR \bigg\} \\
 \sV_{\con} &= \bigg\{ \sum_{i=1}^{n-k+1} u_i \, \by_i \, g_{\IS^{n-k}} \;\; : \;\; u_i \in \IR \bigg\}, 
\end{align*}
where $\by_i$ denotes the restriction of the $i$-th coordinate function to $\IS^{n-k} \subset \IR^{n-k+1}$.
Moreover,
\begin{alignat*}{3}
\sV_{0} &= \sV_{\sph,0} \oplus \sV_{\aff} \oplus \sV_{\con}, &\qquad  {\textstyle 1- \frac{n-k}{2(n-k-1)}} &= 0 \qquad &&\text{if} \quad n-k =2 \\
\sV_{0} &= \sV_{\sph,0} \oplus \sV_{\aff},& \qquad \sV_{1- \frac{n-k}{2(n-k-1)}} &= \sV_{\con} \qquad &&\text{if} \quad n-k \geq 3.
\end{alignat*}
The largest negative eigenvalue of $L$ on the mixed subspace $\sV_{as}$ is
$-\frac{1}{2(n-k-1)}$, realized by the tensors $d\bx_i \odot d\by_a$.
Consequently, the spectral gap of $L$ on $\sV_{\osc,<0}$ is
\[ \gamma_{\osc}=\frac{1}{2(n-k-1)}. \]
\end{Lemma}

\begin{proof}
The product splitting of the tangent bundle gives a decomposition $\sV_{aa}\oplus\sV_{as}\oplus\sV_{ss}$ of the smooth symmetric $2$-tensors into their Euclidean, mixed and spherical components.
The coefficients of each component may depend on both factors.
Using \eqref{eq_Lu_explicit}, it is easy to see that $L$ preserves this decomposition and that $L u = \triangle_f u$ for $u \in \sV_{aa} \oplus \sV_{as}$.
For any $u = u_{ij} dx^i dx^j \in \sV_{aa}$ we have $L u = (\triangle_f u_{ij}) dx^i dx^j$, which reduces the study of $L$ restricted to $\sV_{aa}$ to a scalar problem.
Next, write $\sV_{ss}=\sV_{ss,0}\oplus\sV_{ss,1}$, where $\sV_{ss,0}$ consists of the trace-free spherical tensors and $\sV_{ss,1}$ consists of tensors $u=u'g_{\IS^{n-k}}$ with $u'\in C^\infty(\IR^k\times\IS^{n-k})$.
By \eqref{eq_Lu_explicit}, $Lu=\triangle_f u-\frac1{n-k-1}u$ on $\sV_{ss,0}$ and $L(u'g_{\IS^{n-k}})=(\triangle_f u'+u')g_{\IS^{n-k}}$ on $\sV_{ss,1}$.
So it suffices to show the following two claims.

\begin{Claim}
The eigenvalues of $\triangle_f$ that are $\geq -1$ are $0, -\frac12, - \frac{n-k}{2(n-k-1)}, -1$.
If $n-k \geq 3$, then these eigenvalues are pairwise distinct and the eigenspaces corresponding to the first, second and last eigenvalues are given by Hermite polynomials of degree 0, 1 and 2, respectively.
The eigenspace corresponding to the third eigenvalue is given by linear combinations of the functions $\by_i$ on $\IS^{n-k}$.
If $n-k=2$, then the last two eigenvalues agree and the corresponding eigenspace is the direct sum of the previously mentioned spaces.
\end{Claim}

\begin{proof}
We can carry out a separation of variables procedure as in the proof of \cite[Lemma~5.4]{Bamler_Lai_MCF1}.
Fix an eigenfunction $u \in C^\infty(\IR^k \times \IS^{n-k})$ with $\triangle_f u = \la u$.
Let $b \in C^\infty(\IS^{n-k})$, $\triangle b = \mu b$ be an eigenfunction of the Laplacian on $(\IS^{n-k},g_{\IS^{n-k}})$ and consider the function $w(\bx_1, \ldots, \bx_k) := \int_{\IS^{n-k}} u(\bx_1, \ldots, \bx_k, \ldots) b dg_{\IS^{n-k}}$.
It is not hard to see that $\triangle_f w = (\la-\mu) w$, so the only possible values of $\la$ are sums of eigenvalues of $\triangle_f$ on $\IR^k$ and possible values for $\mu$.
It is well known that the eigenvalues of $\triangle_f$ on $\IR^k$ are $0, -\frac12, -1, \ldots$ and the eigenfunctions are given by Hermite polynomials.
Moreover, the first three eigenvalues of the Laplacian on the unit $(n-k)$-sphere are $0, -(n-k)$ and $-2(n-k+1)$, so the first three possible values of $\mu$ are $0,-\frac{n-k}{2(n-k-1)}$ and $-\frac{n-k+1}{n-k-1}$.
Since the last number is $<-1$, we find that the possible values for $\la$ are sums of elements of $\{ 0,-\frac12, -1 \}$ and $\{ 0, -\frac{n-k}{2(n-k-1)} \}$, which proves the claim.
\end{proof}

\begin{Claim}
The operator $\triangle_f$ restricted to $\sV_{as}$ is strictly negative definite, and its largest eigenvalue is $-\frac{1}{2(n-k-1)}$.
The corresponding eigenspace is spanned by $d\bx_i\odot d\by_a$ and, when $n-k=2$, by $d\bx_i\odot K^\flat$, where $K$ ranges over the Killing vector fields on $\IS^2$.
\end{Claim}

\begin{proof}
Suppose that $\triangle_f u = \la u$ for some $u \in \sV_{as}$.
Fix $i \in \{ 1, \ldots, k \}$ such that $\alpha_i := u(e_i, \cdot)$ does not vanish.
Then $\alpha_i$ is a 1-form with $\alpha_i (v) = 0$ for any tangent vector $v$ that is parallel to the $\IR^k$-factor and we still have $\triangle_f \alpha_i = \la \alpha_i$.
We can now argue as before.
Fix an eigenform $\beta \in \Omega^1(\IS^{n-k})$ such that $\triangle \beta = \mu \beta$.
As before, we can define
$w(\bx_1, \ldots, \bx_k) := \int_{\IS^{n-k}} \langle \alpha_i (\bx_1, \ldots, \bx_k, \ldots), \beta \rangle dg_{\IS^{n-k}}$ and conclude that $\triangle_f w = (\la-\mu) w$.
The largest eigenvalue of the rough Laplacian on $1$-forms on $(\IS^{n-k},g_{\IS^{n-k}})$ is $-\frac{1}{2(n-k-1)}$, realized by $d\by_a$, since
\[ \triangle d\by_a=d\triangle \by_a+\tfrac12d\by_a=-\frac{1}{2(n-k-1)}d\by_a. \]
When $n-k=2$, the Killing vector fields also satisfy
$\triangle K^\flat=-\Ric(K,\cdot)=-\tfrac12K^\flat$.
The remaining eigenvalues are at most this number, while the eigenvalues of $\triangle_f$ in the $\IR^k$-directions are $0,-\frac12,-1,\ldots$.
It follows that $\la\leq -\frac{1}{2(n-k-1)}$, with equality precisely on the span of the tensors listed in the claim.
\end{proof}

This finishes the proof of the lemma.
\end{proof}
\bigskip

We write $\sV_{\rot,\sph}$ for tensors of the form $u'(\bx)g_{\IS^{n-k}}$ and $\sV_{\rot,\aff}$ for tensors of the form $\sum h_{ij}(\bx)\,d\bx_i\otimes d\bx_j$. Additional spectral subscripts denote intersections with the corresponding eigenspaces or their sums. For a tensor $U$, subscripts denote orthogonal projections, and we set
\[
\sV^{++}:=\sV_1\oplus\sV_{\frac12}\oplus\sV_{\sph,0},\qquad U^{++}:=\PP_{\sV^{++}}U,\qquad \la_{\conf}:=1-\frac{n-k}{2(n-k-1)}.
\]
If $\PP_{\sV_{\sph,0}}U=\sum c_{ij}\mathfrak p_{ij}^{(2)}g_{\IS^{n-k}}$ with $(c_{ij})$ symmetric, we denote its smallest and largest eigenvalues by $U_{0,\min}$ and $U_{0,\max}$, respectively.

\begin{Lemma}\label{l:quotient_cylinder_spectrum}
On $M_{\cyl,\Gamma}$ the eigenspaces of $L$ are the $\Gamma$-invariant parts of those in Lemma~\ref{Lem_cyl_spectrum}. If $\Gamma\ne\{1\}$, the nonnegative eigenspaces are unchanged except that $\sV_{\conf}=\{0\}$. The negative oscillatory eigenvalues still satisfy $\lambda\leq-1/(2(n-k-1))$.
\end{Lemma}

\begin{proof}
The operator $L$ commutes with pullback to the cover. A nonzero $\Gamma$-invariant first spherical harmonic would give a common fixed point of $\Gamma$ on the sphere, contradicting freeness. The other nonnegative modes in Lemma~\ref{Lem_cyl_spectrum} are all invariant.
\end{proof}

Fix now some $\la < 0$ and $J \geq 2$ and recall the construction of the Taylor polynomial $Q_J = (Q^+_J, Q^-_J) : \sV_{> \la} \to \sV_{> \la} \oplus \sV_{\leq \la}$ of the non-linear part of \eqref{eq_resc_RdT}.
We will be particularly interested in the second Taylor polynomial.

\begin{Lemma}\label{l:Q_2}
For any $U \in \sV_{> \la}$ we have in a local orthonormal frame 
\begin{multline} \label{eq_Q2}
(Q_2 (U))_{ij} =   \tfrac12  \nabla_i U_{ps} \nabla_j U_{ps} - \nabla_p U_{si} \nabla_p U_{sj} + \nabla_p U_{si} \nabla_s U_{pj}   -  \nabla_i U_{ps} \, \nabla_p U_{sj}    - \nabla_j U_{ps} \, \nabla_p U_{si}      \\
  \qquad - U_{pq}  ( \nabla^2_{pq} U_{ij} - R_{ipq}^{\quad s} U_{sj} - R_{ipj}^{\quad s} U_{qs} - R_{jpi}^{\quad s} U_{sq} - R_{jpq}^{\quad s} U_{is}).
\end{multline}
Here we use the Einstein summation convention and write $u_{ij} := u(e_i,e_j)$, $\nabla_s u_{ij} := (\nabla_{e_s} u)(e_i,e_j)$ and $\nabla^2_{st} u_{ij} := (\nabla^2_{e_s,e_t} u)(e_i,e_j)$ for a local orthonormal frame $\{ e_i \}$.
Moreover, if $U \in \sV_{0}$, then
\begin{equation} \label{eq_PV0Q2}
   \PP_{\sV_0} \big( Q_2(U) \big) =\PP_{\sV_0} \big(  Q_2\big(\PP_{ \sV_{\sph,0}} U \big) \big). 
\end{equation}
\end{Lemma}

\begin{proof}
Identity \eqref{eq_Q2} is a direct consequence of \eqref{eq_Q_def}.
To see \eqref{eq_PV0Q2} let $Q(U', U'') := \frac12 (Q_2(U'+U'') - Q_2(U') - Q_2(U''))$ be the polarization of $Q_2$, so that $Q_2(U') = Q(U',U')$ and write $U = U_{\sph,0} + U_{\aff} + U_{\con}$, where $U_{\conf} = 0$ if $n-k \geq 3$.
We obtain that
\begin{equation} \label{eq_Q2_sac}
   Q_2(U) = Q_2 (U_{\sph,0}) + Q_2(U_{\aff} ) + Q_2(U_{\con}) + 2Q(U_{\sph,0}, U_{\aff}) + 2Q(U_{\sph,0}, U_{\con}) + 2Q(U_{\aff}, U_{\con}). 
\end{equation}
Since all terms in \eqref{eq_Q2} involve a factor that is a covariant derivative of $U$ or a curvature term, we find that $Q_2(U_{\aff}) =0$.
By the same reasoning we obtain for $U' := U_{\sph,0} + U_{\con}$
\[ (Q(U_{\aff}, U'))_{ij} = - \tfrac12 (U_{\aff})_{pq} \nabla^2_{pq} U'_{ij}  = - \tfrac12 (U_{\aff})_{pq} \nabla^2_{pq} (U_{\sph,0})_{ij}. \]
This term vanishes unless $i, j > k$, in which case it is of the form $c \delta_{ij}$ for some $c \in \IR$.
So $Q(U_{\aff}, U') \in \sV_1$.

Next, consider the term $Q(U_{\sph,0}, U_{\con})$.
Note that $\nabla_p (U_{\sph,0})_{ij}$ vanishes unless $p \leq k$ and $i,j > k$.
On the other hand, $\nabla_p (U_{\con})_{ij}$ vanishes unless $p,i,j > k$.
Write $U_{\sph,0}=u'g_{\IS^{n-k}}$ and $U_{\con}=u''g_{\IS^{n-k}}$.
Only the mixed components of the derivative products in $2Q(U_{\sph,0},U_{\con})$ can be nonzero.
For $i\leq k<a$, the first and fourth terms in \eqref{eq_Q2} give a total of $\big(\frac{n-k}{2}-1\big)(\nabla_i u')(\nabla_a u'')$, which vanishes since $U_{\con}=0$ unless $n-k=2$.
Thus for $i,j>k$ we obtain, summing over spherical indices,
\begin{align*} 
   2(Q(U_{\sph,0}, U_{\con}))_{ij} &= - u' \delta_{pq} (\nabla^2_{pq} u'' \, \delta_{ij} - R_{ipq}^{\quad s}  u'' \delta_{sj} - R_{ipj}^{\quad s} u'' \delta_{qs} - R_{jpi}^{\quad s} u'' \delta_{sq} - R_{jpq}^{\quad s} u'' \delta_{is}) \\
 &\qquad  - u'' \delta_{pq} (\nabla^2_{pq} u' \, \delta_{ij} - R_{ipq}^{\quad s}  u' \delta_{sj} - R_{ipj}^{\quad s} u' \delta_{qs} - R_{jpi}^{\quad s} u' \delta_{sq} - R_{jpq}^{\quad s} u' \delta_{is}) \\
 &= - u' \delta_{pq} \nabla^2_{pq} u'' \, \delta_{ij} = c_{n,k} u' u'' \delta_{ij}.
\end{align*}
Since $u'$ depends only on the $\IR^k$-component and $u''$ depends only on the $\IS^{n-k}$-component, it is not hard to see that $Q(U_{\sph,0}, U_{\con})$ is an eigenfunction of $L$ corresponding to a negative eigenvalue, so it is perpendicular to $\sV_0$.

Lastly, we need to analyze the term $Q_2(U_{\conf})$ if $n-k=2$.
Recall that $\nabla_p (U_{\con})_{ij}$ vanishes unless $p,i,j > k$.
So if we set $U = U_{\conf}$, then the first line in \eqref{eq_Q2} vanishes unless $i, j > k$.
It is not hard to see that the same is true for the second line.
Thus $\PP_{\sV_{\aff}} Q_2(U_{\con}) = 0$.
Suppose now that $i, j > k$.
These terms are constant in the Euclidean directions, so $\PP_{\sV_{\sph,0}} Q_2(U_{\con}) = 0$.
Next, consider the map
\[ \sV_{\con} \to \sV_{\con}, \qquad U_{\con} \mapsto \PP_{\sV_{\conf}} Q_2(U_{\con}) \]
This map is equivariant under the natural $O(3)$-action on $\sV_{\con} \cong \IR^{3}$, so 
\[ \PP_{\sV_{\conf}} Q_2(U_{\con}) = - \PP_{\sV_{\conf}} Q_2(-U_{\con}) =  - \PP_{\sV_{\conf}} Q_2(U_{\con}) \]
and thus $\PP_{\sV_{\conf}} Q_2(U_{\con}) = 0$.

This shows that the projections to $\sV_0$ of all terms  on the RHS of \eqref{eq_Q2_sac}, except for $Q_2(U_{\sph,0})$, vanish and finishes the proof.
\end{proof}
\bigskip

We compute $Q_2$ on the rotationally symmetric tensors used in the mode analysis.

\begin{Lemma}\label{Lem_Q2_computation}
Suppose that $U \in \sV_{> \la}$ is of the form
\[ U = \sum_{i,j=1}^k U'_{ij} (\bx) dx^i \otimes dx^j + u'(\bx) g_{\IS^{n-k}}, \]
where $\bx$ is the coordinate on the first factor of $\IR^k \times \IS^{n-k}$ and $U'_{ij}=U'_{ji},u'\in C^\infty(\IR^k)$.
Then
\begin{equation} \label{eq_Q2UQpqp}
Q_2(U) = \sum_{i,j=1}^k Q'_{ij} (\bx) d\bx_i \otimes d\bx_j + q'(\bx) g_{\IS^{n-k}}, 
\end{equation}
where (the indices $p$ and $s$ run from $1$ to $k$)
\begin{align}
Q'_{ij} &=   \tfrac{n-k}2  \nabla_i u' \nabla_j u' +
\tfrac12  \nabla_i U'_{ps} \nabla_j U'_{ps}  
- \nabla_p U'_{si} \nabla_p U'_{sj} 
+ \nabla_p U'_{si} \nabla_s U'_{pj} \notag \\ &\qquad -  \nabla_i U'_{ps}  \nabla_p U'_{sj}     - \nabla_j U'_{ps}  \nabla_p U'_{si}  - U'_{ps}  \nabla^2_{ps} U'_{ij}     
    \label{eq_Qpij} \\
q' &= - |\nabla u' |^2      
   - U'_{pq}   \nabla^2_{pq} u' . \label{eq_qp}
\end{align}
If $U'_{ij} = U'_{ji}, a, b_i, c_{ij} = c_{ji} \in \IR$ and
\[ U = \sum_{i,j=1}^k U'_{ij} d\bx_i \otimes d\bx_j + \bigg( a + \sum_{i=1}^k b_i \mathfrak p^{(1)}_i (\bx) + \sum_{i,j=1}^k c_{ij} \mathfrak p^{(2)}_{ij}(\bx) \bigg) g_{\IS^{n-k}}, \]
then
\[ Q_2(U) = \sum_{i,j=1}^k Q'_{ij}(\bx) d\bx_i \otimes d\bx_j + \bigg( \ov a + \sum_{i=1}^k \ov b_i \mathfrak p^{(1)}_i (\bx) + \sum_{i,j=1}^k \ov c_{ij} \mathfrak p^{(2)}_{ij}(\bx) \bigg) g_{\IS^{n-k}}, \]
where
\begin{align} 
Q'_{ij} &= \frac{n-k}{4}b_i b_j+\frac{n-k}{2 \sqrt{2}} \sum_{s=1}^k (b_i c_{js} + b_j c_{is}) \mathfrak p^{(1)}_s(\bx)  + \frac{n-k}{\sqrt 2} \sum_{s,t=1}^k c_{is} c_{jt} \mathfrak p^{(2)}_{st}(\bx) + \frac{n-k}2 \sum_{s=1}^k c_{is} c_{js} \label{eq_Us2QU_1} \\
\ov a &= - \frac12 \sum_{s=1}^k b_s^2 - \sum_{s,t=1}^k c_{st}^2 - \frac1{\sqrt 2}\sum_{s,t=1}^k U'_{st} c_{st} \label{eq_Us2QU_2} \\
 \ov b_i &= - \sqrt 2 \sum_{s=1}^k  b_s c_{is}  \label{eq_Us2QU_3} \\
 \ov c_{ij} &=  - \sqrt 2 \sum_{s=1}^k c_{is}c_{sj}. \label{eq_Us2QU_4} 
\end{align}
If
\[ U = \sum_{i,j=1}^k \big(v_{ij} + w_{ijst} \mathfrak p^{(2)}_{st}(\bx) \big) d\bx_i \otimes d\bx_j + \bigg( a  + \sum_{i,j=1}^k c_{ij} \mathfrak p^{(2)}_{ij}(\bx) \bigg) g_{\IS^{n-k}}, \]
then
\[ \PP_{\sV_{\sph,0}} \big( Q_2(U) \big) = \sum_{i,j=1}^k \ov c_{ij} \mathfrak p_{ij}^{(2)}(\bx)g_{\IS^{n-k}},\]
where
\begin{equation} \label{eq_Uh2QU}
\ov c_{ij} =  - \sqrt 2 \sum_{s=1}^k c_{is}c_{sj} - \frac1{\sqrt 2} \sum_{s,t=1}^k w_{stij} c_{st}.
\end{equation}
\end{Lemma}

\begin{proof}
The fact that $Q_2(U)$ is of the form \eqref{eq_Q2UQpqp} follows by symmetry.
Identity \eqref{eq_Qpij} follows from \eqref{eq_Q2}.
Note that all curvature terms vanish.
Likewise, evaluating \eqref{eq_Q2} for $i \in \{ k +1, \ldots, n \}$,
\begin{align*}
 q' &=  - |\nabla u' |^2        
   - U'_{pq}   \nabla^2_{pq} u'  - u' \delta_{pq} ( - R_{ipq}^{\quad s} \, u' \, \delta_{sj} - R_{ipj}^{\quad s} \, u' \, \delta_{qs} - R_{jpi}^{\quad s} \, u' \, \delta_{sq} - R_{jpq}^{\quad s} \, u' \, \delta_{is}) \\
  &=  - |\nabla u' |^2      
   - U'_{pq}   \nabla^2_{pq} u'   - u'  ( - R_{ipp}^{\quad i} \, u'  - R_{ipi}^{\quad p} \, u'  - R_{ipi}^{\quad p} \, u'  - R_{ipp}^{\quad i} \, u' ) \\
   &=  - |\nabla u' |^2      
   - U'_{pq}   \nabla^2_{pq} u'  ,
\end{align*}
proving \eqref{eq_qp}.
Suppose now that $U'_{ij} \equiv const$ and $u'(\bx) =  a + \sum_{i=1}^k b_i \mathfrak p^{(1)}_i (\bx) + \sum_{i,j=1}^k c_{ij} \mathfrak p^{(2)}_{ij}(\bx)$.
Then, with repeated Euclidean indices summed,
\begin{align*}
Q'_{ij}&=\frac{n-k}{2}\nabla_i u'\nabla_j u'
=\frac{n-k}{4}(b_i+c_{is}\bx_s)(b_j+c_{jt}\bx_t),\\
q'&=-\frac12(b_s+c_{is}\bx_i)(b_s+c_{sj}\bx_j)
-\frac1{\sqrt2}U'_{st}c_{st}\\
&=-\frac12 b_s^2-c_{st}^2-\frac1{\sqrt2}U'_{st}c_{st}
-\sqrt2 b_s c_{is}\mathfrak p_i^{(1)}
-\sqrt2 c_{is}c_{sj}\mathfrak p_{ij}^{(2)}.
\end{align*}
This confirms \eqref{eq_Us2QU_1}--\eqref{eq_Us2QU_4}.
Next, suppose that $U'_{ij} = v_{ij} + w_{ijst} \mathfrak p^{(2)}_{st}(\bx)$ and $u' = a + c_{ij} \mathfrak p^{(2)}_{ij}(\bx)$.
Then
\begin{align*}
q' &= - \frac12  c_{is} \bx_i c_{sj} \bx_j  - \frac1{\sqrt{2}} (v_{pq} + w_{pqst} \mathfrak p^{(2)}_{st}) c_{pq} \\
q' &= - \sqrt 2 c_{is} c_{sj}\mathfrak p^{(2)}_{ij}(\bx) - c_{st}^2  - \frac1{\sqrt 2} v_{pq} c_{pq} - \frac1{\sqrt 2} w_{stij} c_{st} \mathfrak p^{(2)}_{ij}(\bx).
\end{align*}
This shows \eqref{eq_Uh2QU}.
\end{proof}
\medskip

The following calculation gives the cubic contribution needed in Section~\ref{sec:dominant_quadratic}.

\begin{Lemma} \label{Lem_Q3_computation}
Suppose that $U'_{ij}=U'_{ji},c_i\in\IR$ and
\[
U=\sum_{i,j=1}^k U'_{ij}\,d\bx_i\otimes d\bx_j
+\sum_{i=1}^k c_i\mathfrak p_{ii}^{(2)}g_{\IS^{n-k}}
\in\sV_{\aff}\oplus\sV_{\sph,0}.
\]
Then
\[
\PP_{\sV_{\sph,0}}(Q_3(U)-Q_2(U))
=\sum_{i,j=1}^k\ov c_{ij}\mathfrak p_{ij}^{(2)}g_{\IS^{n-k}},
\]
where
\begin{equation} \label{eq_ovci82}
\ov c_{ij}
=\delta_{ij}\bigg(4c_i^3+c_i\sum_{s=1}^k c_s^2\bigg)
+\sqrt2\,c_i c_j U'_{ij}.
\end{equation}
\end{Lemma}

\begin{proof}
Write $u'=\sum_i c_i\mathfrak p_{ii}^{(2)}$ and regard
$U'=(U'_{ij})$ as a constant symmetric matrix.
By \eqref{eq_Q_def}, the spherical block of $Q[U]$ is
\[
\bigg[-\frac{((I+U')^{-1})^{pq}}{1+u'}
\partial_pu'\partial_qu'
+\big(((I+U')^{-1})^{pq}-\delta^{pq}\big)
\partial^2_{pq}u'\bigg]g_{\IS^{n-k}}.
\]
Its homogeneous cubic part is
\[
\big(u'|\nabla u'|^2
+U'_{pq}\partial_pu'\partial_qu'
+(U'^2)_{pq}\partial^2_{pq}u'\big)g_{\IS^{n-k}}.
\]
Since $\partial_pu'=c_p\bx_p/\sqrt2$ and
$\partial^2_{pq}u'=c_p\delta_{pq}/\sqrt2$, the last term is
constant. With respect to the normalized Gaussian inner product,
\[
\langle(\mathfrak p_{ii}^{(2)})^2,1\rangle=1,\qquad
\langle(\mathfrak p_{ii}^{(2)})^2\bx_i^2,1\rangle=10,\qquad
\langle(\mathfrak p_{ii}^{(2)})^2\bx_s^2,1\rangle=2
\quad(s\ne i).
\]
Consequently, the degree-two Hermite projections are
\begin{align*}
\PP_2(u'|\nabla u'|^2)
&=\sum_i\bigg(4c_i^3+c_i\sum_s c_s^2\bigg)
\mathfrak p_{ii}^{(2)},\\
\PP_2(U'_{pq}\partial_pu'\partial_qu')
&=\sqrt2\sum_{p,q}U'_{pq}c_pc_q\mathfrak p_{pq}^{(2)}.
\end{align*}
This proves \eqref{eq_ovci82}.
\end{proof}

\section{A continuation criterion for harmonic map heat flow}\label{sec:gauging}
Proposition~\ref{Prop_gauging} constructs an HMHF from a cylindrical chart by extending the metric to the whole cylinder. If the construction stops before the prescribed final time, either the rescaled perturbation fails the required $L^2$ bound or the coordinate image does not contain the required cylindrical region. The proof of Theorem~\ref{thm:HMHF} uses continuous dependence on the initial chart to choose an initial map for which the HMHF reaches the final time, provided the Ricci flow remains sufficiently cylindrical along a fixed worldline.

From this section onward we suppress $\Gamma$ from the notation for the chosen cylindrical model and its mode spaces. Constants may depend on $\Gamma$. Unless specified otherwise, differences of maps close to the identity are taken on their lifts near the identity in the round cover embedded in $\IR^{n+1}$.

We will frequently consider the dilation diffeomorphism
\[ \xi_{\la} : \IR^{k} \times (\IS^{n-k}/\Gamma) \lto \IR^{k} \times (\IS^{n-k}/\Gamma), \qquad  (\bx,\by) \mapsto (\la \bx, \by) \]

\begin{Proposition} \label{Prop_gauging}
If $m \geq 10$, $\eta \leq \ov\eta$ and $\sigma \leq \ov\sigma(m, \eta)$, then the following is true.

Let $\MM$ be a Ricci flow spacetime satisfying the pseudolocality property in Definition~\ref{def:pseudolocality_property}, defined over a time-interval $[-e^{-\tau_0}, -e^{-\tau_1}]$ for $\tau_0<\tau_1$.
Consider an open subset $\Omega' \subset \MM_{-e^{-\tau_0}}$ and a smooth diffeomorphism
$$\chi' : \Omega' \lto \IB_{1.1e^{-\tau_0/2}R'}^k \times (\IS^{n-k}/\Gamma)$$
with the property that $\Vert e^{\tau_0}(\xi_{e^{\tau_0/2}})_* \chi'_* g_{-e^{-\tau_0}} - g_{\cyl} \Vert_{C^{m+10}(\IB^k_{1.1R'}\times(\IS^{n-k}/\Gamma))} < \sigma$ and $R' > 2\sigma^{-1}$.

Then there are a time $\tau_*\in(\tau_0,\tau_1]$, an open subset $\Omega\subset\MM_{[-e^{-\tau_0},-e^{-\tau_*}]}$, and a smooth map $\chi:\Omega\to\IR^k\times(\IS^{n-k}/\Gamma)$ whose time-slice maps are diffeomorphisms onto their images. The time $\tau_*$ is either $\tau_1$ or $\tau_0+j$ for some positive integer $j$. The following hold:
\begin{enumerate}[label=(\alph*)]
\item \label{Prop_gauging_a} $\Omega_{-e^{-\tau_0}} \subset \Omega'$ and $\chi_{-e^{-\tau_0}} = \chi'|_{\Omega_{-e^{-\tau_0}}}$.
\item \label{Prop_gauging_c} The map $\chi$, viewed as a family of diffeomorphisms, satisfies the following harmonic map heat flow equation
\begin{equation} \label{eq_HMHF_in_chi_prop}
   \partial_t \chi_t = \triangle_{g_{t}, g_{\cyl,t}} \chi_t . 
\end{equation}
Consequently, the push-forwards $g'_t := (\chi_t)_* g_t$
 satisfy the Ricci DeTurck flow equation with background metrics $(g_{\cyl,t})$ and 
\begin{equation} \label{eq_mod_RDTF}
   \partial_t  g'_t = - 2 \Ric_{ g'_t} -  \mathcal{L}_{X_{g_{\cyl,t}} [ g'_t]  }  g'_t .
\end{equation}
\item \label{Prop_gauging_d}
Define the domain $\DD \subset \IR^k \times (\IS^{n-k}/\Gamma) \times [\tau_0, \tau_1]$ and the tensor fields $\td g' = g_{\cyl} + u$ on $\DD$ such that for all $\tau \in [\tau_0, \tau_*]$
\begin{align}
 \DD_\tau &:= \Image (\xi_{e^{\tau/2}} \circ \chi_{-e^{-\tau}}), \\
  \td g'_{\tau} &:= g_{\cyl} +  u_{\tau} := e^\tau (\xi_{e^{\tau/2}})_* g'_{-e^{-\tau}}   . \label{eq:h'_tau}
\end{align}
Then $\td g'$ satisfies the rescaled Ricci DeTurck flow equation
\begin{equation} \label{eq_rescale_mod_RDTF}
 \partial_\tau \td g'_\tau = - 2 \Ric_{\td g'_\tau} + \td g'_\tau -\LL_{X_{g_{\cyl}} [\td g'_{\tau}] +  \nabla f_{\cyl}} \td g'_\tau. 
\end{equation}
\item \label{Prop_gauging_e} $\Vert u_\tau \Vert_{C^{m}(\DD_\tau) } < \eta$ for all $\tau \in [\tau_0, \tau_*]$. 
\item \label{Prop_gauging_f} $\IB^k_{R'/2 + 1} \times (\IS^{n-k}/\Gamma) \subset \DD_\tau$ for all $\tau \in [\tau_0, \min\{\tau_0+1,\tau_1\}]$.
\item\label{Prop_gauging_g}
There is $\sigma'=\sigma'(m,\eta)>0$ such that the tensor field $u$
(and its domain) satisfies the $(\sigma',\eta,m')$-pseudolocality
property in \cite[Definition~4.10]{Bamler_Lai_MCF1} for every
integer $4\leq m'\leq m$.
\item \label{Prop_gauging_h} If $\tau_* < \tau_1$, then  $\Vert u_{\tau_*} \Vert_{L^2 (\IB^k_{\sigma^{ -1}} \times (\IS^{n-k}/\Gamma))} \geq \sigma$ or $\IB^k_{\sigma^{-1}} \times (\IS^{n-k}/\Gamma) \not\subset \DD_{\tau_*}$.
\item \label{Prop_gauging_i} There is a well-defined process with input $\chi'$ and output being the time $\tau_*$, domain $\DD$ and tensor field $u$.
Moreover we have the following continuity property:
Consider a sequence of inputs $\chi^{\prime,(i)}$, for $i \leq \infty$, which satisfy the assumptions of the proposition and such that $\chi^{\prime,(i)} \to  \chi^{\prime,(\infty)}$ locally smoothly.
Let $\tau^{*,(i)}$, $\DD^{(i)}$ and $u^{ (i)}$ be the resulting times, domains, and tensor fields.
Then $\liminf_{i\to\infty}\tau^{*,(i)}\geq\tau^{*,(\infty)}$. Before $\tau^{*,(\infty)}$, each compact subset of a limiting domain lies in the corresponding domains for all large $i$, and $u^{(i)}\to u^{(\infty)}$ and $\chi^{(i)}\to\chi^{(\infty)}$ locally smoothly.
At the terminal time the same convergence holds on compact subsets that remain in the corresponding domains for all large $i$.
\end{enumerate}
\end{Proposition}
\begin{Remark}
The construction in Proposition~\ref{Prop_gauging} also applies to any fixed
complete gradient shrinking Ricci soliton $(N,\bar g,f)$ with bounded
curvature, normalized by $\Ric_{\bar g}+\nabla^2_{\bar g}f=\tfrac12\bar g$.
The pseudolocality assumption on $\MM$ is unchanged. Writing $\Phi_s$ for
the flow of $\nabla_{\bar g}f$, use the background
$(-t)\Phi_{-\log(-t)}^*\bar g$ and replace $\xi_{e^{\tau/2}}$ by $\Phi_\tau$.
Replace cylindrical product balls by sufficiently large sublevel sets
$\{f<R^2/4\}$; the constants and numerical factors in the domain inclusions
may depend on the shrinker. The grafting and HMHF estimates require only
bounded geometry, while pointed compactness replaces Euclidean
translations in the compactness arguments.
\end{Remark}
\medskip

\subsection{Estimates for the modified Ricci DeTurck flow equation}

The following estimates control derivatives and propagate smallness of the Ricci--DeTurck perturbation.

\begin{Lemma} \label{Lem_RdT_estimates}
If $m \geq 2$, $\eps > 0$, $\eta \leq \ov\eta$ and $R \geq \underline R (m, \eps)$, then the following is true.

Let $(g'_t = g_{\cyl,t} + u'_t)_{t \in [-1,-e^{-1}]}$ be a solution to the Ricci DeTurck flow equation on $\IB^k_R \times (\IS^{n-k}/\Gamma)$.
Suppose that 
\[ \Vert u'_t \Vert_{C^2(\IB^k_R \times (\IS^{n-k}/\Gamma))} \leq \eta . \]
Then the following is true:
\begin{enumerate}[label=(\alph*)]
\item \label{Lem_RdT_estimates_a} We have the following bound for all  $t \in [-e^{-1/2}, -e^{-1}]$
\[ \Vert u'_t \Vert_{C^m(\IB^k_{(1-\eps)R} \times (\IS^{n-k}/\Gamma))} \leq C(m, \eps) \eta. \]
\item \label{Lem_RdT_estimates_b} If $\delta > 0$ and $\sigma \leq \ov\sigma(m, \eps, \delta)$ and $R \geq \underline R(m, \eps, \delta)$, then the following is true.
If for some $t \in [-1, -e^{-1}]$ we have $\Vert u'_t \Vert_{C^2(\IB^k_R \times (\IS^{n-k}/\Gamma))} \leq \sigma$, then $\Vert u'_{-e^{-1}} \Vert_{C^m(\IB^k_{(1-\eps)R} \times (\IS^{n-k}/\Gamma))} \leq \delta$.
\end{enumerate}
\end{Lemma}

\begin{proof}
Assertion~\ref{Lem_RdT_estimates_a} follows from standard local derivative estimates (see, for example, \cite[Lemma~A.14]{Bamler_Kleiner_uniqueness}).

Assertion~\ref{Lem_RdT_estimates_b} follows via a limit argument:
Assume that $\eta \leq \ov\eta$ so that we can apply local derivative estimates and fix $m, \eps, \delta$.
If the statement was wrong, then we could find a sequence of Ricci DeTurck flows $(g'_{i,t} = g_{\cyl, t} + u'_{i,t})_{t \in [-1,-e^{-1}]}$ on $\IB^k_{R_i} \times (\IS^{n-k}/\Gamma)$ for $R_i \to \infty$ such that $\Vert u'_{i,t_i} \Vert_{C^2(\IB^k_{R_i} \times (\IS^{n-k}/\Gamma))} \leq \sigma_i \to 0$ for $t_i \in [-1,-e^{-1}]$, yet $\Vert u'_{i,-e^{-1}} \Vert_{C^m(\IB^k_{(1-\eps)R_i} \times (\IS^{n-k}/\Gamma))} > \delta$.
Suppose that $|u'_{i,-e^{-1}}|$ or the norm of one of its  derivatives up to order $m$ exceeds $c(m)\delta$ at some point $(\bx_i,\by_i) \in \IB^k_{(1-\eps)R_i} \times (\IS^{n-k}/\Gamma)$.
Since the background metric $g_{\cyl, t}$ is invariant under translations along the $\IR^k$-factor, we can translate the sequence of flows along $\IR^k$ so that the points $\bx_i$ are mapped to the origin. Let us denote the spatially translated sequence by $\td{u}'_{i,t}$.

The translated domains contain $\IB^k_{\eps R_i}\times(\IS^{n-k}/\Gamma)$
throughout $[-1,-e^{-1}]$. The local derivative estimates used in
Assertion~\ref{Lem_RdT_estimates_a} therefore give uniform higher-order
bounds on every compact spatial set and every compact interval in
$(-1,-e^{-1}]$. Passing to a subsequence,
we obtain $t_i\to t_\infty\in[-1,-e^{-1}]$ and a $C^m$-limit
$\td u'_{\infty,t}$ on $(-1,-e^{-1}]$. The assumed $C^2$-bound and the equation
also give $\|\partial_t\td u'_{i,t}\|_{C^0}\leq C\eta$.
Thus the limit extends continuously to $t=-1$, and the $C^2$ bound
$\sigma_i\to0$ at $t_i$ gives $\td u'_{\infty,t_\infty}=0$, including
when $t_\infty$ is either endpoint.

Uniqueness for the Ricci--DeTurck equation among uniformly
$C^2$-small solutions gives $\td u'_{\infty,t}=0$ for
$t\in[t_\infty,-e^{-1}]$.
Local $C^m$-convergence at $t=-e^{-1}$ contradicts
the lower bound $c(m)\delta$ on $\{0\}\times(\IS^{n-k}/\Gamma)$ for sufficiently large $i$.
\end{proof}

\bigskip

\subsection{Solving the modified harmonic map heat flow equation}
\label{subsec_modified_hmhf}
{The following lemma constructs HMHF over an interval of one unit in rescaled time. In its proof, we first use pseudolocality to control the Ricci flow in the initial coordinates. We then graft its metric to the cylindrical metric outside a large region, solve HMHF on the complete cylinder, and restrict the solution to a smaller region.
To prove Proposition~\ref{Prop_gauging}, we iterate this construction, using the terminal map as the initial condition for the next step whenever the perturbation remains sufficiently small on a sufficiently large region.}

\begin{Lemma} \label{Lem_one_step}
If $m \geq 10$, $\eta \leq \ov\eta$ and $\delta \leq \ov\delta(m,\eta)$, then the following is true.

Set $T_{0} := - e^{-\tau_0}$, $T_{1} := - e^{-\tau_1}$ for $\tau_1 = \tau_0 + 1 \in \IZ$ and let $\MM$ be a Ricci flow spacetime satisfying the pseudolocality property in Definition~\ref{def:pseudolocality_property}, defined over the time-interval $[T_{0}, T_{1}]$.
Let $R_0, R_1 > \delta^{ -1}$ and let $\chi' : \Omega' \to \IB^k_{e^{-\tau_0/2}R_0} \times (\IS^{n-k}/\Gamma)$ be a smooth diffeomorphism, where $\Omega' \subset \MM_{T_0}$ is open.
Suppose that
\begin{equation} \label{eq_R1R0xig}
 1.1 R_{1} \leq R_{0}, \qquad  \Vert e^{\tau_0} (\xi_{e^{\tau_0/2}})_* \chi'_* g_{T_{0}} - g_{\cyl} \Vert_{C^{m+10}(\IB^k_{1.1 R_{1}} \times (\IS^{n-k}/\Gamma) )} \leq \delta.
\end{equation}
Then there is an open subset $\Omega \subset \chi^{\prime -1}(\IB^k_{1.1e^{-\tau_0/2}R_1} \times (\IS^{n-k}/\Gamma)) \subset \Omega' \subset \MM_{T_0}$, which survives until time $T_1$, and a smooth family of diffeomorphisms onto their images
\[ \big( \chi_t : \Omega(t) \lto \IR^k \times (\IS^{n-k}/\Gamma) \big)_{t \in [T_0, T_1]} \]
such that
\begin{enumerate}[label=(\alph*)]
\item \label{Lem_one_step_a} $\chi_{T_0} = \chi' |_{\Omega}$.
\item \label{Lem_one_step_b} $(\chi_t)$ satisfies the  harmonic map heat flow equation \eqref{eq_HMHF_in_chi_prop} and $g'_t :=  (\chi_t)_* g_t$ satisfies the Ricci DeTurck flow equation \eqref{eq_mod_RDTF} and $u_\tau := e^\tau (\xi_{e^{\tau/2}})_* g'_{-e^{-\tau}} - g_{\cyl}$ satisfies the rescaled Ricci DeTurck flow equation.
Moreover, $u_\tau$ is defined on the domain $\DD_\tau := \Image (\xi_{e^{\tau/2}} \circ \chi_{-e^{-\tau}})$.
\item \label{Lem_one_step_c} $\IB^k_{e^{(\tau-\tau_0)/2}R_1} \times (\IS^{n-k}/\Gamma) \subset \DD_{\tau} \subset \IB^k_{1.1e^{(\tau-\tau_0)/2}R_1} \times (\IS^{n-k}/\Gamma)$ for all $\tau \in [\tau_0, \tau_1]$.
\item \label{Lem_one_step_d} $\Vert u_\tau \Vert_{C^{m}(\DD_\tau)} < \eta$ for all $\tau \in [\tau_0, \tau_1]$ and $\Vert u_{\tau_1} \Vert_{C^{m+20}(\DD_{\tau_1})} < \eta$.
\item \label{Lem_one_step_e} The family $(\chi_t)$ depends continuously on $\chi'$ and $R_1$ in the local smooth topology.
\end{enumerate}
\end{Lemma}

\begin{proof}
Since the lemma is scaling invariant, we may assume in the following that $[T_0, T_1] = [-1, -e^{-1}]$ and $[\tau_0, \tau_1] = [0,1]$.
After decreasing $R_0$, we may moreover assume that $R_0 = 1.1 R_1 =:1.1 R$.
Let $\delta' > 0$ be a constant whose value we will determine later.
The next claim allows us to reduce the setting to the case in which $\MM$ is a product spacetime.

\begin{Claim}
If $\delta \leq \ov\delta(m,\eta,\delta')$, then all points in 
\[ \Omega'_0 := \chi^{\prime -1} (\IB^k_{1.09R} \times (\IS^{n-k}/\Gamma)) \subset \Omega' \subset \MM_{-1} \]
survive until time $-e^{-1}$ and if we identify $\Omega'_0 ([-1,-e^{-1}])$ with $\Omega'_0 \times [-1, -e^{-1}]$, then
\begin{equation} \label{eq_chi0gtgcylt}
 \Vert \chi'_* g_t - g_{\cyl,t} \Vert_{C^{m+5}(\IB^k_{1.09R} \times (\IS^{n-k}/\Gamma))} \leq \delta'.
\end{equation}
\end{Claim}

\begin{proof}
Fix $\eta,\delta'$ and choose counterexamples $\MM_i$, $\chi'_i:W'_i\to\IB^k_{1.1R_i}\times(\IS^{n-k}/\Gamma)$ with $\delta_i\to0$ and hence $R_i>\delta_i^{-1}\to\infty$.
Choose $p_i\in(\chi'_i)^{-1}(\IB^k_{1.09R_i}\times(\IS^{n-k}/\Gamma))$ whose worldline does not survive until $-e^{-1}$ or near which \eqref{eq_chi0gtgcylt} fails at some $t_i\in[-1,-e^{-1}]$.
Let $(\bx_i,\by_i) := \chi'_i (p_i)$ and consider the translation 
\[ \varphi_{-\bx_i} : \IR^k \times (\IS^{n-k}/\Gamma) \to \IR^k \times (\IS^{n-k}/\Gamma), \qquad (\bx,\by) \mapsto (\bx-\bx_i, \by). \]
The recentered initial domains contain $\IB^k_{0.01R_i}\times(\IS^{n-k}/\Gamma)$, and the initial metrics converge locally in $C^{m+10}$ to $g_{\cyl,-1}$.
Choose $r>0$, depending only on the dimension and $\Gamma$, such that every radius-$r$ ball of $g_{\cyl,t}$, $t\in[-1,-e^{-1}]$, satisfies the two bounds in Definition~\ref{def:pseudolocality_property} with strict inequalities, and set
\[
\theta:=\bigl(\tfrac1{10}\eps_{\mathrm{pseudo},n}r\bigr)^2>0.
\]
Suppose local $C^{m+5}$-convergence to the cylinder has been established through a time $s$.
Applying Definition~\ref{def:pseudolocality_property} at $s$ to balls centered throughout each fixed coordinate region gives survival until $\min\{s+\theta,-e^{-1}\}$ and a uniform curvature bound there.
Local derivative estimates and compactness give a complete bounded-curvature limit on this interval, which is the cylindrical flow by \cite[Theorem~1.1]{Chen_Zhu_RF_uniqueness}.
Thus the convergence extends through its right endpoint. In particular, the cylindrical bounds are recovered before the next step, so the same $r$ and $\theta$ apply.
Starting at $s=-1$, at most $\lceil(1-e^{-1})/\theta\rceil$ steps give local $C^{m+5}$-convergence on all of $[-1,-e^{-1}]$.
This proves survival of $p_i$ and convergence of the errors near it to zero, contradicting its choice.
\end{proof}

Let us now identify $\Omega'_0 ([-1,-e^{-1}])$ with $\IB^k_{1.09R} \times (\IS^{n-k}/\Gamma) \times [-1, -e^{-1}]$ via the diffeomorphism $\chi' |_{\Omega'_0}$ and express the Ricci flow on this parabolic domain by $(g_t := g_{\cyl, t} + u'_t)_{t \in [-1,-e^{-1}]}$.
We will only use the flow restricted to this product domain and our harmonic map heat flow $(\chi_t)$ will be defined on a subset of this domain as well.
So we may simplify our setting as follows:
\begin{itemize}
\item We may assume that $\MM$ is given by a conventional Ricci flow $(g_t)_{t \in [-1,-e^{-1}]}$ on $\IB^k_{1.09 R} \times (\IS^{n-k}/\Gamma)$ and $\chi' = \id$.
We write $R := R_1$ and $R_0 = 1.1 R$ and recall from \eqref{eq_chi0gtgcylt}  that $\Vert g_{t} - g_{\cyl,t} \Vert_{C^{m+5}} \leq \delta'$.
\item Our goal is to find an open subset $\Omega \subset \IB^k_{1.09 R} \times (\IS^{n-k}/\Gamma)$ and a solution $\chi_t := \chi^{(1)}_t : \Omega \to \IR^k \times (\IS^{n-k}/\Gamma)$ of the harmonic map heat flow equation \eqref{eq_HMHF_in_chi_prop} with initial condition $\chi_{-1} = \id$.
Assertion~\ref{Lem_one_step_c} is equivalent to
\begin{equation} \label{eq_double_containment}
 \IB^k_{ R} \times (\IS^{n-k}/\Gamma) \subset \chi_{t} (\Omega) \subset \IB^k_{1.1 R} \times (\IS^{n-k}/\Gamma) \qquad \text{for all} \quad t \in [-1, -e^{-1}]
\end{equation}
and Assertion~\ref{Lem_one_step_d} will follow from bounds of the form
\begin{equation} \label{eq_chigtgcyl_bound_want}
 \Vert (\chi_t)_* g_t - g_{\cyl,t} \Vert_{C^{m}} < c_0 \eta, \qquad \Vert (\chi_{-e^{-1}})_* g_{-e^{-1}} - g_{\cyl,-e^{-1}} \Vert_{C^{m+20}} < c_0 \eta  
\end{equation}
for some fixed $c_0 = c_0(m) > 0$.
The continuity statement expressed in Assertion~\ref{Lem_one_step_e} will be clear by construction.
\end{itemize}
To achieve our goal, we adapt the grafting construction from \cite{Bamler_Kleiner_uniqueness}:

\begin{Claim} \label{Cl_graft}
If $\delta' \leq \ov\delta'$, then there is a smooth family of Riemannian metrics $(g_{\graft,t})_{t \in [-1,-e^{-1}]}$ on $\IR^{k} \times (\IS^{n-k}/\Gamma)$ such that for some dimensional constant $C$:
\begin{enumerate}[label=(\alph*)]
\item \label{Cl_graft_a} $g_{\graft,t} = g_t$ on $\IB^{k}_{1.07 R} \times (\IS^{n-k}/\Gamma)$.
\item $\Vert g_{\graft,t} - g_{\cyl,t} \Vert_{C^{m+5}} \leq C\delta'$.
\item $- C\delta' g_{\graft,t} \leq \partial_t g_{\graft,t} + 2 \Ric_{g_{\graft,t}} \leq  C\delta' g_{\graft,t}$.
\item $g_{\graft,t} = g_{\cyl,t}$ on the complement of $\IB^{k}_{1.08 R} \times (\IS^{n-k}/\Gamma)$.
\item $(g_{\graft,t})$ depends continuously on $(g_t)$ and $R$ in the local smooth topology.
\end{enumerate}
\end{Claim}

\begin{proof}
Let $\omega \in C^\infty(\IR^k \times (\IS^{n-k}/\Gamma))$ be a radial cutoff function such that $\omega \equiv 1$ on $\IB^{k}_{1.07 R} \times (\IS^{n-k}/\Gamma)$ and $\omega \equiv 0$ on the complement of $\IB^{k}_{1.08 R} \times (\IS^{n-k}/\Gamma)$ and set
\[ g_{\graft,t} := \omega g_t + (1-\omega) g_{\cyl,t}. \]
The claim now follows directly from the assumptions of the lemma and it is clear that $\omega$ can be chosen in such a way that the constant $C$ is independent of $R$.
\end{proof}

In the next claim we solve the harmonic map heat flow equation with respect to the flows $(g_{\graft,t})$ and $(g_{\cyl,t})$.

\begin{Claim} \label{Cl_construct_chi}
If $\eta' > 0$ and $\delta \leq \ov\delta (m, \eta')$ and $\delta' \leq \ov\delta' (m, \eta')$, then  there is a unique smooth family of diffeomorphisms
\[ \big( \ov\chi_t : \IR^{k} \times (\IS^{n-k}/\Gamma) \lto \IR^{k} \times (\IS^{n-k}/\Gamma) \big)_{t \in [-1,-e^{-1}]} \]
which are uniformly Lipschitz and satisfy 
\begin{equation} \label{eq_chiYp_evolution}
   \partial_t \ov\chi_t = \triangle_{g_{\graft,t}, g_{\cyl,t}} \ov\chi_t, \qquad \ov\chi_{-1} = \id. 
\end{equation}
Moreover,
\begin{equation} \label{eq_chigraftgcyl}
 \Vert (\ov\chi_t)_* g_{\graft,t} - g_{\cyl,t} \Vert_{C^m} < c_0 \eta' 
\end{equation}
and we have the following displacement bound for some dimensional constant $C$
\begin{equation} \label{eq_displacement_bound}
  d_{g_{\cyl,t}}(\ov\chi_t(p), p) \leq  C \qquad \text{for all} \quad p \in \IR^k \times (\IS^{n-k}/\Gamma).
\end{equation}
The family of diffeomorphisms $(\ov\chi_t)$ depends continuously on $(g_t)$ in the local smooth topology.
\end{Claim}

\begin{proof}
The claim follows from \cite[Proposition~A.24]{Bamler_Kleiner_uniqueness}.
It is easy to see that the assumptions of this proposition are satisfied for $g_t$ and $g'_t$ replaced with $g_{\graft,t}$ and $g_{\cyl,t}$, respectively, and $\eta,\delta$ replaced with $C \delta'$.
So \cite[Proposition~A.24]{Bamler_Kleiner_uniqueness} guarantees the existence of a smooth solution $(\ov\chi_t)_{t \in [-1,T^*]}$ to \eqref{eq_chiYp_evolution} for some $T^* \in (-1,-e^{-1}]$ assuming that $\delta' \leq \ov\delta'$.
Moreover, this proposition implies that there is a dimensional constant $\delta^* > 0$ such that if the corresponding modified Ricci DeTurck flow $(\ov\chi_t)_* g_{\graft,t} = g_{\cyl,t} + u'_t$
satisfies $| u'_t| \leq \delta^*$ for all $t \in [-1,T^*]$, then we have $T^* = -e^{-1}$.

The bound \eqref{eq_chigraftgcyl} on $( u'_t)$ can be obtained from a continuity argument; compare, for example, with the proof of \cite[Lemma~A.19]{Bamler_Kleiner_uniqueness}.
We can use \cite[(A.7)]{Bamler_Kleiner_uniqueness} to conclude that there are dimensional constants $C_0 > 0$ and $0 < \delta^{**} \leq \delta^*$ such that on the maximal time-interval $[-1,T^{**}] \subset [-1,T^*]$ on which the bound $|u'_t|\leq \delta^{**}$ holds, we have
\begin{equation} \label{eq_evol_ht}
 \partial_t| u'_t|^2 \leq (g_{\cyl,t} +  u'_t)^{ij} \nabla^2_{ij} | u'_t|^2 + C_0 \delta' | u'_t| + C | u'_t|^2.
\end{equation}
So by the weak maximum principle, we obtain that if $\delta' \leq \ov\delta'$, then we even have $| u'_t|\leq C \delta' \leq \frac12 \delta^{**}$ on $[-1,T^{**}]$.
Local gradient estimates on $(u'_t)$ (for example via \cite[Lemma~A.14]{Bamler_Kleiner_uniqueness}), imply a uniform bound on $|\partial_t  u'_t|$ over the time-interval.
So if $T^{**} < T^*$, then the bound $| u'_t|\leq \delta^{**}$ holds beyond time $T^{**}$, which contradicts its maximal choice.
It hence follows that $T^{**} = T^*$ and since $\delta^{**} \leq \delta^{*}$, we have $T^* = -e^{-1}$ and $| u'_t| \leq C\delta'$ for all $t \in [-1,-e^{-1}]$.
Local derivative estimates give \eqref{eq_chigraftgcyl} if $\delta' \leq \ov\delta'(m,\eta')$.

The uniqueness statement follows from standard parabolic theory, using the uniform Lipschitz bound.
The displacement bound \eqref{eq_displacement_bound} follows as in \cite[Proposition~A.10]{Bamler_Kleiner_uniqueness}.
The continuity property follows from a limit argument via \eqref{eq_chigraftgcyl} and \eqref{eq_displacement_bound} combined with local derivative estimates (see also \cite[Lemma~A.15]{Bamler_Kleiner_uniqueness}).
\end{proof}

Set now $\Omega := \IB^k_{1.05 R} \times (\IS^{n-k}/\Gamma)$ and $\chi_t := \ov\chi_t |_{\Omega}$.
The first bound in \eqref{eq_chigtgcyl_bound_want} follows directly from \eqref{eq_chigraftgcyl} for $\eta' = c_0 \eta$, since $g_{\graft,t} = g_t$ on $\IB^k_{1.07 R} \times (\IS^{n-k}/\Gamma)$, see
Claim~\ref{Cl_graft}\ref{Cl_graft_a}; note that this bound even holds for $\ov\chi_t |_{\IB^k_{1.07 R} \times (\IS^{n-k}/\Gamma)}$.
The second bound in \eqref{eq_chigtgcyl_bound_want} follows from Lemma~\ref{Lem_RdT_estimates}\ref{Lem_RdT_estimates_a}, for a smaller $\eta'$.
Note that we can use \eqref{eq_displacement_bound} to show that for $\eta \leq \ov\eta$ and $\delta \leq \ov\delta$ we have
\[ \ov\chi_{-e^{-1}} (\IB^k_{1.06 R} \times (\IS^{n-k}/\Gamma)) \subset \ov\chi_{t} (\IB^k_{1.07 R} \times (\IS^{n-k}/\Gamma)) \]
for all $t \in [-1,-e^{-1}]$.
So the Ricci DeTurck flow is defined in a parabolic neighborhood of uniform size around every point in $\ov\chi_{-e^{-1}} (\IB^k_{1.06 R} \times (\IS^{n-k}/\Gamma))$.
The bound \eqref{eq_double_containment} also follows from \eqref{eq_displacement_bound} for $\eta \leq \ov\eta$ and $\delta \leq \ov\delta$
\end{proof}
\bigskip

\subsection{Proof of Proposition~\ref{Prop_gauging}}

\begin{proof}[Proof of Proposition~\ref{Prop_gauging}.]
By parabolic rescaling we may take $\tau_0=0$. We describe the construction for integer $\tau_1$. If the last step is shorter than one unit, we run the construction in Lemma~\ref{Lem_one_step} only up to $\tau_1$. Its domain and $C^m$ estimates and its continuous dependence hold on this shorter interval; the terminal $C^{m+20}$ estimate is needed only to start another step.
In the following we will construct a sequence of open subsets $\Omega_{\tau_0} \subset \Omega'$, $\Omega_{\tau_0+1} \subset \Omega_{\tau_0}(-e^{-(\tau_0+1)})$, $\Omega_{\tau_0+2} \subset \Omega_{\tau_0+1}(-e^{-(\tau_0+2)})$, \ldots, and a sequence of diffeomorphisms onto their images $(\chi_{j,t} : \Omega_j(t) \to \IR^k \times (\IS^{n-k}/\Gamma))_{t \in [-e^{-j},- e^{-(j+1)}]}$, for $j = \tau_0, \ldots, \tau_*-1$.
We will then set
\begin{equation}\label{eq:gauging_Omega_def}
\Omega := \bigcup_{j=\tau_0}^{\tau_*-2}  \Omega_j \big([-e^{-j},-e^{-(j+1)}) \big) \cup \Omega_{\tau_*-1}\big([-e^{-(\tau_*-1)},-e^{-\tau_*}] \big),
\end{equation}
and
\begin{equation}\label{eq:gauging_chi_def}
\chi_t := \chi_{j,t} \qquad \text{if} \quad t \in [-e^{-j}, -e^{-(j+1)}) \quad \text{or} \quad t = -e^{-\tau_*}.
\end{equation}

Let $\delta > 0$ be a constant whose value we will determine in the course of the proof.
Assume $\sigma \leq \delta$ and $\eta \leq \ov\eta$ and $\delta \leq \ov\delta(m,\eta)$, so that Lemma~\ref{Lem_one_step} is applicable
with its input threshold replaced by $1.1\delta$.
Using this lemma, we will iteratively construct radii $R_{\tau_0}, R_{\tau_0+1}, \ldots > (1.1\delta)^{-1}$ and the families of diffeomorphisms $(\chi^{(j)}_t)$.

In the first step we set $R_{\tau_0} := \frac1{1.1} R'$ and apply Lemma~\ref{Lem_one_step} to $\chi' \leftarrow \chi'$, $R_0 \leftarrow 1.1R'$ and $R_1 \leftarrow R_{\tau_0}$ to obtain $(\chi_{\tau_0,t})_{t \in [T_{\tau_0}, T_{\tau_0+1}]}$.
Set $\Omega := \Omega_{\tau_0} ([-e^{-\tau_0}, -e^{-(\tau_0+1)}])$ and $\chi_t := \chi_{\tau_0, t}$ as in \eqref{eq:gauging_Omega_def} and \eqref{eq:gauging_chi_def} for $\tau_* = \tau_0+1$ and define $g'_t$ for $t \in [-e^{-\tau_0},-e^{-(\tau_0+1)}]$ and $\td g'_\tau = g_{\cyl} + u_\tau$ on $\DD_\tau$ for $\tau \in [\tau_0, \tau_0+1]$ as explained in Assertions~\ref{Prop_gauging_c} and \ref{Prop_gauging_d}.
Note that Lemma~\ref{Lem_one_step}\ref{Lem_one_step_c} implies the following inclusion for $j' = \tau_0$
\begin{equation} \label{eq_contdoubDD}
 \IB^k_{ e^{(\tau-j')/2} R_{j'}} \times (\IS^{n-k}/\Gamma) \subset \DD_{\tau} \subset \IB^k_{1.1 e^{(\tau-j')/2} R_{j'}} \times (\IS^{n-k}/\Gamma) \qquad \text{for all} \quad \tau \in [j', j'+1].
\end{equation}

Suppose now that $R_{\tau_0}, \ldots, R_{j-1} > (1.1\delta)^{-1}$ and $(\chi_{\tau_0,t})_{t \in [-e^{-\tau_0}, - e^{-(\tau_0+1)}]}, \ldots, (\chi_{j-1,t})_{t \in [-e^{-(j-1)}, -e^{-j}]}$ have been constructed such that, if we define $\Omega$ and $\chi$ according to \eqref{eq:gauging_Omega_def} and \eqref{eq:gauging_chi_def} and $\DD$ according to Assertion~\ref{Prop_gauging_d}, then \eqref{eq_contdoubDD} holds for $j' = \tau_0, \ldots, j-1$ and assume that $j \leq \tau_1$.
This implies that
\[ \IB^k_{ e^{1/2} R_{j-1}} \times (\IS^{n-k}/\Gamma) \subset \DD_{j}. \]

\begin{Claim} \label{Cl_choice_Rjp1}
One of the following two statements is true:
We have
\begin{equation} \label{eq_choice_Rjp1_case1}
 \Vert u_{j} \Vert_{C^{m+10} (\IB^k_{ \delta^{ -1}} \times (\IS^{n-k}/\Gamma))} \geq \tfrac12\delta
\end{equation}
or there is a unique $R_{j} \in (\frac1{1.1}\delta^{ -1}, \frac1{1.1} e^{1/2} R_{j-1}]$ such that
\begin{alignat*}{2}
\Vert u_{j} \Vert_{C^{m+10} (\IB^k_{1.1 R_{j}} \times (\IS^{n-k}/\Gamma))} &= \tfrac12 \delta + \tfrac12 \delta R_{j}^{-1} \leq \delta&\qquad& \text{if} \quad R_{j} < \tfrac1{1.1} e^{1/2} R_{j-1} \\
\Vert u_{j} \Vert_{C^{m+10} (\IB^k_{1.1 R_{j}} \times (\IS^{n-k}/\Gamma))} &\leq \tfrac12 \delta + \tfrac12 \delta R_{j}^{-1} \leq \delta &\qquad&  \text{if} \quad R_{j} = \tfrac1{1.1} e^{1/2} R_{j-1}
\end{alignat*}
\end{Claim}

\begin{proof}
Suppose that \eqref{eq_choice_Rjp1_case1} is false.
Then for $R_{j} = \frac1{1.1} \delta^{-1}$ we have $$\Vert u_{j} \Vert_{C^{m+10} (\IB^k_{1.1 R_{j}} \times (\IS^{n-k}/\Gamma))} < \tfrac12 \delta + \tfrac12 \delta R_{j}^{-1}.$$
The LHS is weakly increasing in $R_{j}$, while the RHS is strictly decreasing.
So either we have equality for a unique $R_{j} \in (\frac1{1.1} \delta^{ -1}, \frac1{1.1} e^{1/2} R_{j-1}]$ or this inequality is true for all such $R_j$.
In the latter case we set $R_{j} := \frac1{1.1} e^{1/2} R_{j-1}$.
\end{proof}

If \eqref{eq_choice_Rjp1_case1} is true, then we stop and set $\tau_* := j$.
Otherwise, we can apply Lemma~\ref{Lem_one_step} for $R_0 \leftarrow e^{1/2} R_{j-1}$, $R_1 \leftarrow R_j$ and $\chi'$ equal to the restriction of $\chi_{j-1,-e^{-j}}$ to the inverse image of $\IB^k_{e^{(1-j)/2}R_{j-1}}\times(\IS^{n-k}/\Gamma)$, to obtain $(\chi_{j,t})_{t \in [-e^{-j}, -e^{-(j+1)}]}$ and repeat the construction.
Note that when updating $\Omega$ and $\chi$ via \eqref{eq:gauging_Omega_def} and \eqref{eq:gauging_chi_def} based on the extended list at each construction step, the domain $\Omega_{-e^{-j}}$ may shrink; consequently, $\chi_{-e^{-j}}$ is defined as the restriction to this updated subset.

Assertions~\ref{Prop_gauging_a}--\ref{Prop_gauging_f} follow immediately by construction and from Lemma~\ref{Lem_one_step}.
At the first restart, Assertion~\ref{Prop_gauging_f} still holds: choosing $\sigma$ small enough to give $\|u_{\tau_0+1}\|_{C^{m+20}}<\delta/2$ in the first step, Claim~\ref{Cl_choice_Rjp1} gives $R_{\tau_0+1}=e^{1/2}R_{\tau_0}/1.1>R'/2+1$.
To see Assertion~\ref{Prop_gauging_h} suppose that $\tau_* < \tau_1$.
Then \eqref{eq_choice_Rjp1_case1} must hold for $j = \tau_*$.
On the other hand, the application of Lemma~\ref{Lem_one_step}\ref{Lem_one_step_d} from the previous step implies that
\[ \Vert u_{\tau_*} \Vert_{C^{m+20} (\IB^k_{ \delta^{ -1}} \times (\IS^{n-k}/\Gamma))} \leq \eta \leq 1. \]
The lower $L^2$-bound now follows for $\sigma \leq \ov\sigma(m,\delta)$ from a basic limit argument or the Gagliardo-Nirenberg inequality.

Since the continuation condition opposite to \eqref{eq_choice_Rjp1_case1} is strict, the continuity in Lemma~\ref{Lem_one_step} implies inductively that every restart before $\tau^{*,(\infty)}$ also occurs for all sufficiently large $i$. Hence $\liminf_{i\to\infty}\tau^{*,(i)}\geq\tau^{*,(\infty)}$.
The radius in Claim~\ref{Cl_choice_Rjp1} depends continuously on $u_j$ wherever the construction continues, by the monotonicity in its proof. Lemma~\ref{Lem_one_step} thus gives convergence through the end of each common step. At a restart, \eqref{eq:gauging_chi_def} restricts the preceding terminal map to the new domain.

For Assertion~\ref{Prop_gauging_g}, suppose
$\IB_R^k\times(\IS^{n-k}/\Gamma)\subset\DD_{\tau_0}$,
$R>(\sigma')^{-1}$, and
$\|u_{\tau_0}\|_{C^{m'-1}(\IB_R^k\times(\IS^{n-k}/\Gamma))}\leq\sigma'\eta$,
where $4\leq m'\leq m$ and $\tau_2-\tau_0\in[\frac12,1]$.
By Assertion~\ref{Prop_gauging_e}, it suffices to prove
$\IB^k_{e^{\sigma'(\tau_2-\tau_0)}(R-(\sigma')^{-1})}\times(\IS^{n-k}/\Gamma)
\subset\DD_{\tau_2}$.
If $\tau_0\in[j,j+1)$, then \eqref{eq_contdoubDD} implies
$R\leq1.1e^{(\tau_0-j)/2}R_j$. If there is no integer time in
$(\tau_0,\tau_2]$ other than possibly the terminal time $\tau_*$,
then  Assertion~\ref{Prop_gauging_g} holds by taking $\sigma'<1/100$ since
\[e^{\sigma'(\tau_2-\tau_0)}\bigl(R-(\sigma')^{-1}\bigr)
<
e^{\sigma'(\tau_2-\tau_0)}R
<
0.9e^{(\tau_2-\tau_0)/2}R\le e^{(\tau_2-j)/2}R_j.\]

It remains to consider $j+1\leq\tau_2<j+2$ with $j+1<\tau_*$.
Normalize $[j,j+1]$ to logarithmic times $[0,1]$ and undo the
rescaling. Set $R':=0.9e^{(j-\tau_0)/2}R\leq0.99R_j$.
By \eqref{eq_contdoubDD}, the resulting Ricci--DeTurck flow is
defined on $\IB^k_{R'}\times(\IS^{n-k}/\Gamma)$ throughout $[-1,-e^{-1}]$.
Its perturbation has $C^2$-norm at most $C\eta$ throughout this
interval and at most $C(m)\sigma'\eta$ at the time
$-e^{j-\tau_0}$ corresponding to $\tau_0$. Choose $\ov\eta$ small
enough to apply Lemma~\ref{Lem_RdT_estimates} with the bound $C\eta$.

Let $C(m)$ also bound the changes of the $C^{m+10}$-norm under
rescaling within one normalized step. Choose $0<\delta'\leq\delta/4$.
Then choose $\sigma'$ sufficiently small that
$C(m)\sigma'\eta\leq\ov\sigma(m+10,1/100,\delta'/C(m))$
and $R'>\delta^{-1}$ exceeds the radius threshold in
Lemma~\ref{Lem_RdT_estimates}\ref{Lem_RdT_estimates_b}.
Apply the lemma with $m+10$ derivatives, $\eps=1/100$, and
$\delta'/C(m)$ in place of $\delta$. Rescaling back gives
\begin{equation} \label{eq:initial_D_2}
\lim_{\tau\nearrow j+1}
\|u_\tau\|_{C^{m+10}(\IB^k_{0.99e^{1/2}R'}\times(\IS^{n-k}/\Gamma))}
\leq\delta'.\end{equation}
Since $\delta'\leq\delta/4$ and $0.99R'\leq R_j$,
Claim~\ref{Cl_choice_Rjp1} implies
$R_{j+1}\geq\frac{0.99}{1.1}e^{1/2}R'
=0.9^2e^{(j+1-\tau_0)/2}R$.
By \eqref{eq_contdoubDD}, this proves the required inclusion, since
$0.9^2e^{(\tau_2-\tau_0)/2}>e^{\sigma'(\tau_2-\tau_0)}$
for $\tau_2-\tau_0\in[\frac12,1]$ and $\sigma'<1/100$.
\end{proof}
\bigskip

\section{Harmonic map heat flow via the PDE--ODI principle}\label{sec:pde_odi}

The main result in this section is the following theorem, which shows that if a Ricci flow is a $\delta$-neck for all times in an interval, then there exists a smooth harmonic map heat flow from this flow to the cylindrical flow.
Moreover, the geometrically insignificant modes  (conformal and affine) of the induced perturbation are dominated by geometrically significant ones.

\begin{Theorem}
[Long time existence of harmonic map heat flow]\label{thm:HMHF}
For any constants $m \geq 4, J \geq 2,  
\la <0$, there is a constant $\ov\beta(\la,J)>0$ and, for every $0<\beta\leq\ov\beta(\la,J)$, there are constants such that the following holds if
\begin{equation} \label{eq_HMHF_constants}
\eta \leq \ov\eta(\la,J,\beta), \quad 
\eps \leq \ov\eps(m,J,\la,\beta,\eta), \quad
R^* \geq \underline R^*(m,J,\la, \eta, \eps),\quad
\delta\le\ov\delta(m,J,\la,\beta, \eta, \eps,R^*).
\end{equation}

Let $(\MM,g_t)$ be a Ricci flow spacetime,
$t\in I\subset(-\infty,0)$, and let $p(t)$, $t\in I$, be a worldline
in $\MM$.  Assume $(\MM,g_t)$ is $(n,k,\delta,\Gamma)$-cylindrical at $p(t)$
at scale $\sqrt{-t}$ for all $t\in I$.
Then there exists an open subset $\Omega\subset \MM_{t\in I}$, with {$p(t)\in \Omega_t$}, and a smooth map $\chi : \Omega \lto \IR^k \times (\IS^{n-k}/\Gamma)$ whose time-slice maps $\chi_t : \Omega_t \lto \IR^k \times (\IS^{n-k}/\Gamma)$ are diffeomorphisms onto their images, and satisfy the  harmonic map heat flow equation
\[\partial_t \chi_t = \triangle_{g_{t}, g_{\cyl,t}} \chi_t ,\]
such that for all $\tau\in \td I:=\{-\log(-t),t\in I\}$, we have
\begin{equation}\label{thm:HMHF_open_2}
    \xi_{e^{\tau/2}}\circ\chi_{-e^{-\tau}}(p)\in \IB^k_{\eta}\times(\IS^{n-k}/\Gamma),
\end{equation}
Define
$u_\tau:=e^\tau(\xi_{e^{\tau/2}})_*(\chi_{-e^{-\tau}})_* g_{-e^{-\tau}}-g_{\cyl}$, and set
$\DD_\tau:=\Image(\xi_{e^{\tau/2}}\circ\chi_{-e^{-\tau}})$.
Then there is a continuous function $R : \td I\to[R^*,\infty]$, which is $C^1$ wherever finite,
such that for all $\tau \in \td I$, 
\begin{equation}\label{thm:HMHF_open_1}
\IB^{k}_{ R(\tau)} \times (\IS^{n-k}/\Gamma)  \subset \DD_\tau
\qquad\textnormal{and}\qquad
\|u_\tau\|_{C^m(\IB^k_{R(\tau)}\times(\IS^{n-k}/\Gamma))}\leq\eta,
\end{equation}
so $u_\tau \omega_{R(\tau)}$ can be extended smoothly on $\IR^k \times (\IS^{n-k}/\Gamma)$ by setting it equal to zero outside $\DD_\tau$.
Moreover, set
\begin{equation*}
 U^+_\tau :=  \PP_{\sV_{> \lambda}} (u_{\tau} \omega_{R(\tau)}) , \qquad
\mathcal U^-(\tau) := \big\|   \PP_{\sV_{\leq \lambda}}  (u_{\tau} \omega_{R(\tau)} ) \big\|_{L^2_{f}} + \eta e^{-\frac{((1-\eps)R(\tau))^2}{8}}.
\end{equation*}
If $R(\tau) = \infty$, then we set $\omega_{R(\tau)} \equiv 1$ and $e^{-\frac{((1-\eps)R(\tau))^2}{8}} = 0$.
Throughout the paper, differential estimates for locally Lipschitz functions are understood almost everywhere.
Then
the quantities $U^+_\tau$ and $\mathcal{U}^-_\tau$ are locally Lipschitz and satisfy the following evolution inequalities:
\begin{equation} \label{eq_ODI_1}
 \Big\| \partial_\tau U^+ - L U^+ - Q^+_J(U^+) \Big\| \leq C(\la, J) \| U^+ \|^{J+1} + C(\la) \sqrt{\eta} \UU^- 
\end{equation}
\begin{equation} \label{eq_ODI_2}
 \partial_\tau \UU^- \leq (\la + C(\la) \sqrt\eta) \UU^- + \| Q^-_J (U^+) \| + C(\la, J) \| U^+ \|^{J+1}. 
\end{equation}
Moreover, writing $\tau_0:=\inf\td I$ and $\tau_1:=\sup\td I$, there exist times $T_1\le T_2\le T_3$ in the extended interval $[\tau_0,\tau_1]$ and a constant $C(J,\la,\beta)>0$ such that the following holds:
\begin{enumerate}[label=(\alph*)]

\item \label{thm:HMHF_a} 
For all $\tau\in (\tau_0,T_1)$, we have $\UU^- + \| U^+_{< 0} \| \geq C \| U^+_{\geq 0} \|^2$, and 
\[ \partial_\tau ( \UU^- + \| U^+_{< 0} \|) \leq - \left(\min\left\{-\tfrac\la2,\tfrac{1}{4(n-k-1)}\right\}-C\sqrt\eta\right) ( \UU^- + \| U^+_{< 0} \|). \]

\item\label{thm:HMHF_b}  For all $\tau\in (T_1,\tau_1)$, we have
\begin{equation} \label{eq_Upp_controls_everything}
\| U^{+}_{\aff} \| + \|U^{+}_{\conf} \| \leq {C} \|U^{++} \|\quad\textnormal{and}\quad{\UU^{-}+\|U^{+}_{<0}\|}\leq {C} \|U^{++} \|^2,
\end{equation}
and moreover
\begin{equation} \label{eq:U++_evolve}
\partial_\tau \|U^{++}\|\ge -C\sqrt{\eta}\|U^{++}\|.
\end{equation}

\item\label{thm:HMHF_c} For all $\tau\in (T_1,T_3)$, we have $\|U^+\|\le C\|U^{++}_0\|$, and also
\begin{equation} \label{eq:neutral_dom_rough}
\| U^+_{> 0} \| \leq \beta \| U^{++}_0 \|. 
\end{equation}  
Moreover, we have the following behaviors depending on the dominance of $U_{0,\min}$ and $U_{0,\max}$:
\begin{itemize}
\item If $\tau \in (T_1,T_2)$, then $|U_{0,\min}|(\tau) \leq U_{0,\max}(\tau)$ and
\begin{equation} \label{eq_cmax_evol_0}
 \big| \partial_\tau U_{0,\max}(\tau) + \sqrt{2} U_{0,\max}^2(\tau) \big| \leq \beta U_{0,\max}^2(\tau). 
\end{equation}
\item If $\tau \in (T_2, T_3)$, then $-U_{0,\min}(\tau) \ge U_{0,\max}(\tau)$ and
\begin{equation} \label{eq_cmin_evol_0}
 \big| \partial_\tau U_{0,\min}(\tau) + \sqrt{2} U_{0,\min}^2(\tau) \big| \leq \beta U_{0,\min}^2(\tau). 
\end{equation}
\end{itemize}
\item\label{thm:HMHF_d} If $n-k\ge 3$, and $T_3\in(-\infty,\infty)$, we have 
$$U^+_{\conf}(T_3)=0.$$
If $n-k=2$, and $T_2\in(-\infty,\infty)$, we have 
$$U^+_{\conf}(T_2)=0.$$
At finite endpoints not contained in $\td I$, the values of $U^+_{\conf}$ are understood as one-sided limits, which exist by \eqref{eq_ODI_1} and the uniform bounds on $U^+$ and $\UU^-$.
\item\label{thm:HMHF_e} For all $\tau\in (T_3,\tau_1)$, setting $U^{++}_{>0} :=  \PP_{\sV_{1}} U^++ \PP_{\sV_{\frac12}} U^+$, we have $\| U^+ \| \leq C \| U^{++}_{>0} \|$, and
\[\partial_\tau\|U^{++}_{>0}\|\ge (\tfrac12-C\sqrt\eta)\|U^{++}_{>0}\| .\]

\item\label{thm:HMHF_f} There exists a time $T'\ge-\infty$ such that $R(\tau)=\infty$ for $\tau\le T'$ and $R(\tau)<\infty$ for $\tau>T'$.

\end{enumerate}

\end{Theorem}

The following proposition refines Theorem~\ref{thm:HMHF} when the initial perturbation is sufficiently small.

\begin{Proposition}[Appendix of Theorem~\ref{thm:HMHF}]
\label{prop:HMHF_initial}
Fix the constants in Theorem~\ref{thm:HMHF}.  
Then the following holds if we moreover choose the constants
\[\eta'\leq\eta,\quad
R_0\ge
\underline R_0(m,J,\la,\beta,\eta,\eps,R^*,\eta')\geq R^*,\quad \delta'\le \ov\delta'(m,J,\la,\beta,\eta,\eps,R^*,R_0,\eta'),
\]
Assume $t_0=\min I>-\infty$ and
$(\MM,g_{t_0})$ is $(n,k,\delta',\Gamma)$-cylindrical at $p(t_0)$ at scale
$\sqrt{-t_0}$ for some $0<\delta'\leq\ov\delta'$. Then the HMHF in
Theorem~\ref{thm:HMHF} may be chosen so that for $\tau_0=-\log(-t_0)$, $R(\tau_0)=R_0$, and
\[
\|U^+(\tau_0)\|+\UU^-(\tau_0)\leq\eta'.
\]
\end{Proposition}

We use $g_{\cyl}$ for the round cylinder in the following lemma and its proof.

\begin{Lemma}\label{claim:Lie}
For any $m\in\mathbb N$, there exists $C_m>0$ such that the following holds:
For any vector field 
$X=\tfrac12\sum_{i,j=1}^k a_{ij} \, \bx_i \partial_{\bx_j}+ \sum_{i=1}^{n-k+1}b_{i}\nabla \by_i$
 on $\IR^k\times \IS^{n-k}$, where $a_{ij}=a_{ji},b_{i}\in\IR$, we have
  \[\LL_X g_{\cyl}=\sum_{i,j=1}^k a_{ij} \, d\bx_i \otimes d\bx_j-\tfrac{1}{n-k-1}\sum_{i=1}^{n-k+1}b_{i}\,\by_i g_{S^{n-k}}\in\sV_{\aff}{\oplus}\sV_{\conf}.\]
Let $\phi$ be the self-diffeomorphism of $\IR^k\times \IS^{n-k}$ given by the flow of $-X$ at time $1$. If $\|X\|_{L^2_f}\leq1$, then
 \[\big\|\phi_* g_{\cyl}-g_{\cyl}-\LL_{X}g_{\cyl}\big\|_{C^m(\IR^k\times\IS^{n-k})}\le  C_m\|X\|^2_{L^2_f}.\]
\end{Lemma}

\begin{proof}
The first assertion is clear by computation.
Let $X_1=\tfrac12\sum_{i,j=1}^k a_{ij} \, \bx_i \partial_{\bx_j}$ be a vector field on $\IR^k$, and $X_2=\sum_{i=1}^{n-k+1}b_{i}\nabla \by_i$ a vector field on $\IS^{n-k}$.
Then the flow of $-X_1$ on $\IR^k$ is $\phi_{1,t}(\bx)=e^{-\frac12 tH}x$, where $H=(a_{ij})\in\IR^{k\times k}$. 
{So} we have $(\phi_{1,1})_*g_{\IR^k}=(\phi_{1,-1})^*g_{\IR^k}=dx^T e^{H}dx$ and $\LL_{X_1}g_{\IR^k}= dx^T Hdx$.
Since $e^{H}=I+H+O(\|H\|^2)$, it follows for $C>0$ that
\begin{align*}
\big\|(\phi_{1,1})_*g_{\IR^k}-g_{\IR^k}-\LL_{X_1}g_{\IR^k}\big\|_{C^m(\IR^k)}&
\le C\|X_1\|^2_{L^2_f(\IR^k)}.
\end{align*}

Next, let
$\phi_{2,t}$ be the flow of $-X_2$ on $\IS^{n-k}$.
Since $\frac{d}{dt}\, \phi_{2,t}^{*} g_{\IS^{n-k}} = -\phi_{2,t}^{*}(\mathcal{L}_{X_2}g_{\IS^{n-k}})$, 
we have
\[(\phi_{2,-1})^*g_{\IS^{n-k}}-g_{\IS^{n-k}}-\LL_{X_2}g_{\IS^{n-k}}=\int_0^1 \int_0^t \tfrac{d^2}{ds^2}\, \phi_{2,-s}^{*} g_{\IS^{n-k}}\,ds\,dt.\]
Seeing also $\frac{d^2}{dt^2}\, \phi_{2,t}^{*} g_{\IS^{n-k}} = \phi_{2,t}^{*}(\mathcal{L}_{X_2}\mathcal{L}_{X_2}g_{\IS^{n-k}})$, and $(\phi_{2,1})_*g_{\IS^{n-k}}=(\phi_{2,-1})^*g_{\IS^{n-k}}$, we obtain
\begin{align*}
\big\|(\phi_{2,1})_*g_{\IS^{n-k}}-g_{\IS^{n-k}}-\LL_{X_2}g_{\IS^{n-k}}\big\|_{C^m(\IS^{n-k})}
\le C_m\|X_2\|^2_{L^2(\IS^{n-k})}.
\end{align*}
Now consider $X=X_1+X_2$ as a vector field on $\IR^k\times\IS^{n-k}$, then the flow of $-X$ is $\phi=\phi_{1,1}\oplus\phi_{2,1}$. So we obtain the assertion by combining the estimates.
\end{proof}

In the following lemma, we show that given a $(n,k,\delta,\Gamma)$-neck, we can perturb the diffeomorphism so that the perturbation is dominated by any given affine and conformal modes.

We use $\mathbb B^{\aff}_r,\mathbb B^{\conf}_r$ to denote the $r$-ball {centered at the origin} in the subspaces $\sV_{\aff},\sV_{\conf}$.
	
\begin{Lemma}
\label{l:parametrize}
The following holds for any constants satisfying
\[m\in\mathbb N,\qquad \eta\le\ov\eta, \qquad r\le\eta, \qquad R_0\ge \underline{R}_0(m,r),\qquad \delta {\leq} \ov\delta(R_0,m,r).\]
Let $(M^n,g,p)$ be a {pointed} Riemannian manifold, and assume there exists a  diffeomorphism 
$\chi: \Omega\lto \IB^k_{\delta^{-1}} \times (\IS^{n-k}/\Gamma)$ for an open subset $\Omega\subset M$, such that $\chi(p)\in\bO^k\times(\IS^{n-k}/\Gamma)$ and
\[\|\chi_* g-g_{\cyl}\|_{C^{\delta^{-1}}(\IB^k_{\delta^{-1}} \times (\IS^{n-k}/\Gamma))}\le \delta.\]
Then there is a smooth family of local diffeomorphisms ({in} the local smooth topology),
\[ (\chi^{(s)} : \Omega^{(s)} \lto \IR^k \times (\IS^{n-k}/\Gamma) )_{s \in \IB^{\aff}_{r} \times \IB^{\conf}_{r}\times\IB^k_{r}}, \]
where $\Omega^{(s)}\subset M$ is a {continuous} family of open subsets,
and $\chi^{(0)}=\chi, \Omega^{(0)}=\Omega$. Let $u^{(s)}:=\chi^{(s)}_* g-g_{\cyl}$ be defined on $\IB^k_{R_0} \times (\IS^{n-k}/\Gamma)$.
Then for any $s=(v_1,v_2,v_3) { \in \IB^{\aff}_{r} \times \IB^{\conf}_{r}\times\IB^k_{r}}$,
\begin{enumerate}[label=(\alph*)]
    \item \label{i:v_3}
$\chi^{(s)} (p) \in\big\{\frac{v_3}{\eta}\big\}\times (\IS^{n-k}/\Gamma)$.
\item\label{i:v_1+v_2} $\big\|u^{(s)}\omega_{R_0}-(v_1+v_2)\big\|_{L^2_f}\le C(n,k)r^2$.
\item \label{i:Cm} $\|u^{(s)}\|_{C^{m}(\IB^k_{R_0}\times (\IS^{n-k}/\Gamma))}\le C(m,n,k)r$.
\end{enumerate}

\end{Lemma}
	
\begin{proof}
For any fixed $s=(v_1,v_2,v_3)\in \IB^{\aff}_{r} \times \IB^{\conf}_{r}\times\IB^k_{r}$, let $X^{(s)}$ be as in Lemma~\ref{claim:Lie} with $\LL_{X^{(s)}}g_{\cyl}=v_1+v_2$. For $\Gamma\ne\{1\}$, we have $v_2=0$ and $X^{(s)}$ descends to the quotient.
Let $\phi^{(s)},\psi^{(s)}$ be the diffeomorphism flows of the vector fields $-X^{(s)},\frac{v_3}{\eta}$ at time $1$ on $\IR^k\times(\IS^{n-k}/\Gamma)$.
Then Lemma \ref{claim:Lie} implies
\begin{equation} \label{eq:v_1}
\big\|\phi^{(s)}_*g_{\cyl}-g_{\cyl}-(v_1+v_2)\big\|_{C^m(\IR^k\times(\IS^{n-k}/\Gamma))}\le C_m\|v_1+v_2\|^2.
\end{equation}
Let  
$\chi^{(s)}:=\psi^{(s)}\circ\phi^{(s)}\circ\chi$. Then $\chi^{(0)}=\chi$. We will show that $\chi^{(s)}$ satisfies all assertions:
First, since ${\chi}(p)\in\bO^k\times(\IS^{n-k}/\Gamma)$ and $\phi^{(s)}(\bO^k\times(\IS^{n-k}/\Gamma))=\bO^k\times(\IS^{n-k}/\Gamma)$, Assertion \ref{i:v_3} follows from
\[\chi^{(s)}(p)\in \psi^{(s)}(\bO^k\times(\IS^{n-k}/\Gamma))=\big\{\tfrac{v_3}{\eta}\big\}\times(\IS^{n-k}/\Gamma).\]

The bound $\|v_1\|\le\ov\eta$ implies $\IB^k_{0.5\delta^{-1}} \times (\IS^{n-k}/\Gamma)\subset \phi^{(s)}(\IB^k_{\delta^{-1}} \times (\IS^{n-k}/\Gamma))$, and $\|v_3/\eta\|\le1$ implies
$\IB^k_{0.5\delta^{-1}-1} \times (\IS^{n-k}/\Gamma)\subset \psi^{(s)} \big( \phi^{(s)}(\IB^k_{\delta^{-1}} \times (\IS^{n-k}/\Gamma)) \big)$. Thus, by taking $\delta$ sufficiently small so that $0.5\delta^{-1}-1\ge R_0$, we have
$$ \IB^k_{R_0} \times (\IS^{n-k}/\Gamma)\subset \psi^{(s)} \big( \phi^{(s)}(\IB^k_{\delta^{-1}} \times (\IS^{n-k}/\Gamma)) \big) = \chi^{(s)}(\Omega) .$$
In particular, $u^{(s)}$ is defined on $\IB^k_{R_0} \times (\IS^{n-k}/\Gamma)$.
Write $g_0:=\chi_*g$.
Seeing also that {the tensor fields} $g_{\cyl},v_1+v_2$ are all invariant under the translation $\psi^{(s)}$, and taking $\delta\le r^2$ {and $\delta \leq m^{-1}$}, we obtain by \eqref{eq:v_1} (abbreviating $C^{m}(\IB^k_{R_0}\times (\IS^{n-k}/\Gamma))$ by $C^m$)
\begin{align*}
\big\|u^{(s)}-(v_1+v_2) \big\|_{C^{m}}
&=\big\|(\psi^{(s)}_*\phi^{(s)}_* g_0-g_{\cyl})-(v_1+v_2) \big\|_{C^{m}}\\
&\le\big\|(\psi^{(s)}_*\phi^{(s)}_* g_{\cyl}-g_{\cyl})-(v_1+v_2) \big\|_{C^{m}}+\big\|\psi^{(s)}_*\phi^{(s)}_* (g_0-g_{\cyl})\big\|_{C^{m}}\\
&\le \big\|\psi^{(s)}_*\big(\phi^{(s)}_* g_{\cyl}-g_{\cyl}-(v_1+v_2)\big)\big\|_{C^m} + C_m\delta\\
&\le C_m\|v_1+v_2\|^2_{{C^m}} +C_m\delta\le C_mr^2,
\end{align*}
which implies Assertion \ref{i:Cm}.
Moreover, taking $R_0$ large so that $e^{-\frac{(R_0-2)^2}{8}}{\leq} r$, we obtain
\begin{align*}
\big\|u^{(s)}\omega_{R_0}-(v_1+v_2) \big\|_{L^2_f}
&\le \big\|\big(u^{(s)}-(v_1+v_2)\big)\omega_{R_0}\big\|_{L^2_f} +\big\|(v_1+v_2)(\omega_{R_0}-1)\big\|_{L^2_f}\\
&\le Cr^2+Cre^{-\frac{(R_0-2)^2}{8}}\le Cr^2,
\end{align*}
which proves Assertion \ref{i:v_1+v_2}.
\end{proof}

For the degree argument, we approximate the lower semicontinuous stopping time from below by continuous functions.

\begin{Lemma}\label{l:lower_con}
Let $f {: X \to \IR}$ be a bounded, lower semi-continuous function on a metric space $(X,d)$. Then there exists a sequence of continuous functions $f_i$ on $X$, which satisfy $f_i\le f$ and for any sequence $x_i\to x_{\infty}\in X$, 
\[f(x_{\infty})\le \liminf_{i\to\infty}f_i(x_i).\]
\end{Lemma}

\begin{proof}
For each $i$, let
\[
f_i(x):=\inf_{y\in X}\bigl(f(y)+ i\, d(x,y)\bigr).
\]
Fix $i$ and let $x,x'\in X$. For any $y\in X$, we have
\[
{f_i(x) \leq } f(y)+ i\, d(x,y)\le f(y)+ i\, d(x',y)+ i\, d(x,x').
\]
Taking the infimum over $y$ gives
\[
f_i(x)\le f_i(x')+ i\, d(x,x').
\]
Swapping $x$ and $x'$ yields $f_i(x')\le f_i(x)+ i\, d(x,x')$, so
$f_i$ is $i$-Lipschitz and in particular continuous.

For any $x\in X$, it is easy to see $f_i(x)\le f(x)$ by definition.
To see the last statement suppose that $x_i \to x_\infty$.
By definition we can find $y_i \in X$ such that
\[ \liminf_{i \to \infty} f_i(x_i) = \liminf_{i \to \infty} \big(f(y_i) + i \, d(x_i, y_i) \big) \geq \liminf_{i \to \infty} f(y_i). \]
If the first limit is infinite, then the desired bound holds trivially; if it is finite, then $d(x_i,y_i) \to 0$ and thus $y_i \to x_\infty$, which implies that $\liminf_{i \to \infty} f(y_i) \geq f(x_\infty)$ by the lower semi-continuity.
\end{proof}

We rephrase the PDE-ODI principle for rescaled Ricci flow from \cite[Theorem~4.17]{Bamler_Lai_MCF1}.

\begin{Proposition}[PDE-ODI principle for rescaled Ricci DeTurck flow] \label{Prop_PDE_ODI_RF}
The following holds if
\[ 
m \geq 4,\quad J \geq 1, \quad 
\la \in \IR,  \quad 
\eta \leq \ov\eta(\la), \quad 
\eps \leq \ov\eps(\eta,m),  \quad
R^* \geq \underline R^*(m,J,\la, \eta,\eps). \]
Consider the setting from Proposition~\ref{Prop_gauging},
applied with regularity $\max\{m,10\}$ and restricted to
$[\tau_0,\td\tau]\subset[\tau_0,\tau_*]$ if necessary,
and assume in addition that:
\begin{enumerate}[label=(\roman*)]
\item \label{Prop_PDE_ODI_RF_ii}
For all $\tau \in [\tau_0,\td\tau]$ we have
\[ \IB_{ R^*}^k \times (\IS^{n-k}/\Gamma)  \subset \DD_{\tau} \qquad \text{and} \qquad \| u_{\tau}  \|_{L^2_f(\IB^k_{ R^*} \times (\IS^{n-k}/\Gamma))} \leq \eta e^{-(R^*)^2/8}. \]
\item \label{Prop_PDE_ODI_RF_iii} For some $R_0\in[R^*,\infty]$,
the following holds for all
$\tau\in[\tau_0,\min\{\tau_0+1,\td\tau\}]$:
\begin{equation} \label{eq_atT0bound3}
 \IB^k_{R_0+1} \times (\IS^{n-k}/\Gamma)  \subset \DD_{\tau}.
\end{equation}
\end{enumerate}
Then there is a continuous function $R:[\tau_0,\td\tau]\to[R^*,\infty]$,
which is $C^1$ where finite and satisfies $R(\tau_0)=R_0$.
For all $\tau \in [\tau_0,\td\tau]$ we have $\IB^{k}_{ R(\tau)} \times (\IS^{n-k}/\Gamma)  \subset \DD_\tau$, so $u_\tau \omega_{R(\tau)}$ can be extended smoothly on $M_{\cyl}$ by setting it equal to zero outside its domain.
So we can set
\begin{equation} \label{eq_UpUUm_def_apply}
 U^+(\tau) :=  \PP_{\sV_{> \lambda}} (u_{\tau} \omega_{R(\tau)}) , \qquad
\mathcal U^-(\tau) := \big\|   \PP_{\sV_{\leq \lambda}}  (u_{\tau} \omega_{R(\tau)} ) \big\|_{L^2_{f}} + \eta e^{-\frac{((1-\eps)R(\tau))^2}{8}}. 
\end{equation}
If $R(\tau)=\infty$, set $\omega_{R(\tau)}\equiv1$ and
$e^{-((1-\eps)R(\tau))^2/8}=0$.
The following is true:
\begin{enumerate}[label=(\alph*)]
\item \label{Prop_PDE_ODI_RF_a} We have the bound
$  \| u_{\tau} \|_{C^{m}(\IB^k_{R(\tau)} \times (\IS^{n-k}/\Gamma) )} \leq \eta$.
\item \label{Prop_PDE_ODI_RF_b}  The quantities $U^+ (\tau)$ and $\mathcal{U}^- (\tau)$ satisfy the evolution inequalities
\[ \Big\| \partial_\tau U^+  - L U^+   - Q_J^+ (U^+  ) \Big\|_{L^2_{f}}
\leq  C(\la, J) \Vert U^+  \Vert^{J+1}_{L^2_{f}} +  C(\la) \sqrt\eta\,   \mathcal{U}^-   \]
\[  \partial_\tau \mathcal U^-  
\leq  \big(\la + C(\la) \sqrt\eta  \big)\, \mathcal U^-   
+  \Vert Q^-_J(U^+ ) \Vert_{L^2_{f}}
+ C(\la,J) \Vert U^+  \Vert^{J+1}_{L^2_{f}}
  \]
\item \label{Prop_PDE_ODI_RF_c} Where $R$ is finite, we have
\begin{equation*}
-e^{R(\tau)/10} \leq \partial_\tau R(\tau) \leq 1.
\end{equation*}
For finite $R_0$, the function $R$ can be chosen in a well-defined way so that it depends continuously on the initial condition $\chi'$ and the initial value $R_0$.
Specifically, consider sequences of initial conditions $\chi^{\prime,(i)}$ and initial values $R_0^{(i)}$, for $i \leq \infty$, and let $\tau^{*,(i)}, \DD^{(i)}$ and $u^{(i)}$ be the resulting objects as explained in Assertion~\ref{Prop_gauging_i} of Proposition~\ref{Prop_gauging}.
If $\chi^{\prime,(i)} \to \chi^{\prime, (\infty)}$ locally smoothly,
$R_0^{(i)} \to R_0^{(\infty)}$, and
$\td\tau^{(i)}\to\td\tau^{(\infty)}$ with
$\td\tau^{(i)}\leq\tau^{*,(i)}$, then the resulting functions
$R^{(i)} : [\tau_0, \td\tau^{(i)}] \to [R^*,\infty)$
converge to $R^{(\infty)} : [\tau_0, \td\tau^{(\infty)}] \to [R^*,\infty)$ locally uniformly.
\end{enumerate}
\end{Proposition}

\begin{proof}
Assertions~\ref{Prop_PDE_ODI_RF_a} and
\ref{Prop_PDE_ODI_RF_b} are direct consequences of
\cite[Theorem~4.17 and Remark~4.7]{Bamler_Lai_MCF1}, with $\bY=0$,
applied to the tensor equation \eqref{eq_rescale_mod_RDTF}.
Remark~4.7 applies since Euclidean translations preserve the
cylindrical metric, its connection, and the curvature action
$u\mapsto2\Rm_{g_{\cyl}}(u)$.
Lemma~\ref{Lem_Q_precise_formula} verifies the nonlinear
bounds (4.4)--(4.5) in \cite{Bamler_Lai_MCF1}.
The nonlinear estimates are applied only on the balls in
\eqref{eq_binDueta}, where the required smallness follows from
the $C^m$-bound.
The pseudolocality hypothesis follows from
Proposition~\ref{Prop_gauging}\ref{Prop_gauging_g}; its constant
$\sigma'$ depends only on $m,\eta$ and is absorbed in the bounds
for $\eps,R^*$ above. The initial condition follows from
\eqref{eq_atT0bound3} and
Proposition~\ref{Prop_gauging}\ref{Prop_gauging_e}.

To show Assertion~\ref{Prop_PDE_ODI_RF_c}, we may restrict the flow and write $\tau_*:=\td\tau$.
We modify the radius construction in the first claim of the proof of \cite[Theorem~4.17]{Bamler_Lai_MCF1}.
First suppose that $R_0 < \infty$.

In that construction, $R$ is defined iteratively over time-intervals $[\tau'_0=\tau_0,\tau'_1], [\tau'_1, \tau'_2], \ldots$ as follows.
Set $R(\tau_0) := R_0$ and $E(\tau) :=  e^{-\frac{((1-\eps)R(\tau))^2}8}$ and define $U^+(\tau) := \PP_{\sV_{> \la}}(u_\tau \omega_{R(\tau)})$ and $\td\UU^-(\tau) := \Vert \PP_{\sV_{\leq \la}}(u_\tau \omega_{R(\tau)}) \Vert_{L^2_f}$ for all times at which $R(\tau)$ is already defined.
\begin{enumerate}[label=\arabic*.]
\item If $E(\tau'_j) > 2 ( \frac{\td\UU^-(\tau'_j)}{\eta} + (\frac{\| U^+(\tau'_j) \|}{\eta} )^{J+1} )^{1+\eps}$, then $R$ restricted to $[\tau'_j, \tau'_{j+1}]$ is chosen such that $\partial_\tau E = \la' E$, where $\la' < 0$; so we have $\partial_\tau R > 0$ on $[\tau'_j, \tau'_{j+1}]$.
The right-endpoint $\tau'_{j+1}$ is chosen maximal such that
\begin{equation} \label{eq_binDueta}
 \IB^k_{R(\tau)} \times (\IS^{n-k}/\Gamma) \subset \DD_\tau \qquad \text{and} \qquad \Vert u_\tau \Vert_{C^m (\IB^k_{R(\tau)} \times (\IS^{n-k}/\Gamma))} \leq \eta
\end{equation}
and $E \geq ( \frac{\td\UU^-}{\eta} + (\frac{\| U^+ \|}{\eta} )^{J+1} )^{1+\eps}$ hold on $[\tau'_j,\tau'_{j+1}]$.
We find a posteriori that \eqref{eq_binDueta} holds for all $\tau \in [\tau_0, \tau_1]$, so we must have $ E(\tau'_{j+1}) = ( \frac{\td\UU^-(\tau'_{j+1})}{\eta} + (\frac{\| U^+(\tau'_{j+1}) \|}{\eta} )^{J+1} )^{1+\eps}$ as long as $\tau'_{j+1}< \tau_*$.
\item If $E(\tau'_j) \leq 2 ( \frac{\td\UU^-(\tau'_j)}{\eta} + (\frac{\| U^+(\tau'_j) \|}{\eta} )^{J+1} )^{1+\eps}$, then the construction is similar, except that we solve the ODE $\partial_\tau R = -e^{R/10}$ on $[\tau'_j, \tau'_{j+1}]$ and we have $ E(\tau'_{j+1}) = 4( \frac{\td\UU^-(\tau'_{j+1})}{\eta} + (\frac{\| U^+(\tau'_{j+1}) \|}{\eta} )^{J+1} )^{1+\eps}$ as long as $\tau'_{j+1}< \tau_*$.
\end{enumerate}

We now modify this construction as follows.
Consider the quotient of the LHS and RHS in the conditions of Cases~1. and 2.:
\[ F(\tau) := E(\tau) \cdot \bigg( \frac{\td\UU^-(\tau)}{\eta} + \Big(\frac{\| U^+ (\tau)\|}{\eta} \Big)^{J+1} \bigg)^{-1-\eps}. \]
Where its denominator is positive, $F(\tau)$ depends locally
Lipschitz continuously on $\tau$ and $R(\tau)$.
Let $H : \IR \to [0,1]$ be a smooth function with $H \equiv 0$ on $(-\infty,2]$ and $H \equiv 1$ on $[4, \infty)$ and consider the following ODE:
\begin{equation} \label{eq_ODE_for_R_tau}
 \partial_\tau R(\tau) :=  H(F(\tau)) \cdot \frac{-4\la'}{(1-\eps)^2 R(\tau)} - \big( 1- H(F(\tau)) \big) e^{R(\tau)/10}. 
\end{equation}
Set $H(F)=1$ when the denominator defining $F$ vanishes.
Since $E>0$ for finite $R$ and $H\equiv1$ for $F\geq4$, this
extends $H(F)$ locally Lipschitz continuously. Thus the displayed
equation is a well-defined first-order ODE for $R$.
Let $R : [\tau_0, \tau^{**}) \to [R^*,\infty)$, for $\tau^{**} \leq \tau_*$, be its maximal solution with $R(\tau_0) = R_0$.
In addition, we restrict to the maximal $\tau^{**}$ such that for all $\tau \in [\tau_0, \tau^{**})$ we have \eqref{eq_binDueta}.
The solution is $C^1$ and unique. On compact intervals before the limiting terminal time, continuous dependence follows from \eqref{eq_ODE_for_R_tau} and Proposition~\ref{Prop_gauging}\ref{Prop_gauging_i}.

Let us now verify that the function $R(\tau)$ satisfies the assertions of the first claim in the proof of \cite[Theorem~4.17]{Bamler_Lai_MCF1}.
The bound $$e^{-\frac{(R(\tau)-2)^2}{8}} \leq \eta E(\tau)$$ follows by assuming $R^* \geq \underline R^*(\eta)$ as before.
The bound 
\[ \eta \partial_\tau E(\tau) \leq \la' \eta E(\tau) + \eta \td\UU^-(\tau) + \| U^+(\tau) \|^{J+1} \]
holds trivially if $F(\tau) \geq 4$, since this implies $H(F(\tau)) = 1$.
If $F(\tau) < 4$, then the reasoning from Case~2., as detailed in \cite{Bamler_Lai_MCF1}, applies assuming $R^* \geq \underline R^*(\la, \eta, \eps)$.
Next, consider the bound
\[ -e^{R(\tau)/10} \leq \partial_\tau R(\tau) \leq 1. \]
The first bound holds trivially by \eqref{eq_ODE_for_R_tau}.
The second bound holds whenever $F(\tau) \leq 2$; if $F(\tau) > 2$, then the reasoning from Case~1., as detailed in \cite{Bamler_Lai_MCF1}, applies assuming $R^* \geq \underline R^*(\la)$.
Assertion~(b) of the first claim in the proof of \cite[Theorem~4.17]{Bamler_Lai_MCF1} holds again by \eqref{eq_ODE_for_R_tau} and the fact that $H(F(\tau)) = 0$ if $F(\tau) \leq 1$.

It remains to show that if $\tau^{**} < \tau_*$, then $R(\tau^{**}) := \lim_{\tau \nearrow \tau^{**}} R(\tau)$ exists and
\[ \ov\IB^k_{R(\tau^{**})} \times (\IS^{n-k}/\Gamma) \subsetneq \DD_{\tau^{**}} \quad \text{and} \quad \| u_{\tau^{**}} \|_{C^m(\IB^k_{R(\tau^{**})} \times (\IS^{n-k}/\Gamma))} = \eta. \]
If this were not the case, then the maximal choice of $\tau^{**}$ would be determined by the maximality of the solution to the ODE \eqref{eq_ODE_for_R_tau} alongside the requirement that $R(\tau) \geq R^*$.
Since $\partial_\tau R(\tau) \leq 1$, this would imply that $R(\tau^{**}) := \lim_{\tau \nearrow \tau^{**}} R(\tau) = R^*$ and $\partial_\tau R(\tau^{**}) \leq 0$.
By hypothesis~\ref{Prop_PDE_ODI_RF_ii}, at $R(\tau^{**})=R^*$,
if the denominator defining $F(\tau^{**})$ is positive, then
\[
\begin{aligned}
\frac{\td\UU^-(\tau^{**})}{\eta}
+\left(\frac{\|U^+(\tau^{**})\|}{\eta}\right)^{J+1}
&\leq 2e^{-(R^*)^2/8},\\
F(\tau^{**})&\geq
2^{-1-\eps}\exp\left(\frac{(3\eps-\eps^2)(R^*)^2}{8}\right)>4
\end{aligned}
\]
for sufficiently large $R^*$. Hence $H(F(\tau^{**}))=1$, also
when the denominator vanishes, and
\eqref{eq_ODE_for_R_tau} gives
$\partial_\tau R(\tau^{**})=-4\la'/((1-\eps)^2R^*)>0$,
a contradiction.
We have $\partial_\tau R\leq1$ and
$F(\tau)<1\Longrightarrow\partial_\tau R=-e^{R/10}$, while the
domain and $C^m$-bounds in \eqref{eq_binDueta} hold before
$\tau^{**}$ and must fail there if $\tau^{**}<\tau_*$.
These properties, together with coordinate pseudolocality, allow
the argument following \cite[Claim~4.50]{Bamler_Lai_MCF1}
to apply, giving $\tau^{**}=\tau_*$.
For $R_0^{(i)}\to R_0^{(\infty)}<\infty$, the upper bound $\partial_\tau R^{(i)}\leq1$ bounds the radii uniformly on the given finite interval. The lower bound $\partial_\tau R^{(i)}\geq-e^{R^{(i)}/10}$ then gives a uniform Lipschitz bound, so convergence extends to the varying endpoints $\td\tau^{(i)}\to\td\tau^{(\infty)}$.
For $R_0=\infty$, use the same limiting
argument as in the proof of
\cite[Theorem~4.17]{Bamler_Lai_MCF1}.
Where the limiting radius is finite, passing to the limit in the
integral form of \eqref{eq_ODE_for_R_tau} gives the same ODE
and hence $C^1$-regularity.
\end{proof}

In the following lemma, we analyze the behavior of $U^+$ and $\UU^-$ from Proposition \ref{Prop_PDE_ODI_RF}. We prove that there exist unique times $T_1$ and $T_3$ which divide the time interval into three subintervals, where the stable, neutral, and unstable modes are dominant, respectively.

\begin{Lemma}\label{l:partition_1}
Consider the setting of Proposition \ref{Prop_PDE_ODI_RF}, assume $\la<0$, $\td I=[\tau_0,\tau_*]$ and $\UU^- (\tau)\neq0$ for any $\tau$.
Then for any $C_1\ge\underline{C_1}(\la,J) > 0$, any $0<\beta\leq\ov\beta(\la,J,C_1)$, and any $\eta \leq\ov\eta(\la,J, C_1,\beta)$, we can find unique times $\tau_0 \leq T_1 \leq T_3 \leq \tau_*$ such that the following is true:
\begin{enumerate}[label=(\alph*)]
\item \label{l:partition_1_a}
We have
\[ \UU^- + \| U^+_{< 0} \| \geq C_1 \| U^+_{\geq 0} \|^2 \qquad \text{on} \quad [\tau_0, T_1) \]
\[ \UU^- + \| U^+_{< 0} \| \leq C_1 \| U^+_{\geq 0} \|^2 \qquad \text{on} \quad (T_1, \tau_*] \]
Moreover, if $T_1 > \tau_0$, then
\[ \partial_\tau ( \UU^- + \| U^+_{< 0} \|) \leq - \left(\min\left\{-\tfrac\la2,\tfrac{1}{4(n-k-1)}\right\}-C\sqrt\eta\right) ( \UU^- + \| U^+_{< 0} \|)  \qquad \text{on} \quad [\tau_0, T_1] \]
\item\label{l:partition_1_b} We have
\[ \| U^+_{> 0} \| \leq \beta \| U^+_0 \|,\quad \big|\partial_\tau\|U^+_0\|\big|\le C\eta \|U^+_0\| \qquad \text{on} \quad [T_1,T_3). \]
\[ \| U^+_{> 0} \| \geq \beta \| U^+_0 \|,\quad \partial_\tau \| U^+_{> 0} \| \geq 0.1 \| U^+_{> 0} \| \qquad \text{on} \quad (T_3,\tau_*]. \]
\item\label{l:partition_1_c} 
Assume $\|U^+_{\conf}\|\le C_1\|U^{++}_{>0}\|$ holds on $(T_3,\tau_*]$ when $n-k>2$.
Then
\[ \partial_\tau \| U^{++}_{> 0} \| \geq (\tfrac12-C\sqrt\eta) \| U^{++}_{> 0} \| \qquad \text{on} \quad (T_3, \tau_*]. \]
\end{enumerate}
\end{Lemma}

\begin{proof}
As a consequence of Assertion~\ref{Prop_PDE_ODI_RF_a} of Proposition~\ref{Prop_PDE_ODI_RF}, we have
\[\|U^+_{>0}\|,\|U^+_{0}\|,\|U^+_{\ge 0}\|,\|U^+_{<0}\|,\|U^+\|,\UU^-\le C(\la,J)\eta.\]
Consider the following three conditions:
\begin{enumerate}
\item\label{i:stable_dom} $\UU^- + \| U^+_{< 0} \| \geq C_1 \| U^+_{\geq 0} \|^2$,
\item\label{i:neutral_dom} $\UU^- + \| U^+_{< 0} \| \leq C_1 \| U^+_{\geq 0} \|^2$ and $\|U^+_{>0}\|\le \beta\|U^+_0\|$.
\item\label{i:unstable_dom} $\UU^- + \| U^+_{< 0} \| \leq C_1 \| U^+_{\geq 0} \|^2$ and $\|U^+_{>0}\|\ge \beta\|U^+_0\|$.
\end{enumerate}
Set $\gamma:=\min\{-\la,1/(2(n-k-1))\}>0$.
By Lemma~\ref{Lem_cyl_spectrum}, the ODI system in
Proposition~\ref{Prop_PDE_ODI_RF}\ref{Prop_PDE_ODI_RF_b} gives
\[
\begin{aligned}
\partial_\tau(\UU^-+\|U^+_{<0}\|)
&\leq(\la+C\sqrt\eta)\UU^-
-\frac{1}{2(n-k-1)}\|U^+_{<0}\|+C\|U^+\|^2\\
&\leq-(\gamma-C\sqrt\eta)(\UU^-+\|U^+_{<0}\|)
+C\|U^+_{\geq0}\|^2.
\end{aligned}
\]
Choose $C_1\ge\underline{C_1}(\la,J)$ so large that
$C/C_1\leq\gamma/2$, fix $0<\beta\leq\ov\beta(\la,J,C_1)$,
and decrease $\ov\eta(\la,J,C_1,\beta)$ so that
$C\sqrt\eta\leq\gamma/4$. Under \eqref{i:stable_dom}, we obtain
\begin{equation} \label{eq:evo_1}
\partial_\tau ( \UU^- + \| U^+_{< 0} \|) \leq - (\tfrac12\gamma-C\sqrt\eta) ( \UU^- + \| U^+_{< 0} \|).
\end{equation}
Under \eqref{i:neutral_dom} and \eqref{i:unstable_dom}, respectively,
the same ODI system gives
\begin{equation} \label{eq:evo_2}
\big|\partial_\tau\|U^+_0\|\big|\le C\eta \|U^+_0\|,
\end{equation}
\begin{equation} \label{eq:evo_3}
\partial_\tau \| U^+_{> 0} \| \geq 0.1 \| U^+_{> 0} \|.
\end{equation}

\begin{Claim}
Consider a time $\tau \geq \tau_0$ at which at least \emph{two} of these conditions hold. We claim that for slightly smaller (resp.~larger) times, only the condition appearing earlier (resp.~later) remains valid. 
\end{Claim}

\begin{proof}
If \eqref{i:stable_dom} and \eqref{i:neutral_dom} (or \eqref{i:stable_dom} and \eqref{i:unstable_dom}) hold at time $\tau$,  then
$\UU^- + \| U^+_{< 0} \|=C_1 \| U^+_{\geq 0} \|^2>0$.
This implies 
$$\|U^+\|\le C_1\|U^+_{\ge0}\|^2+\|U^+_{\ge0}\|\le 1.1\|U^+_{\ge0}\|.$$ 
Together with Proposition \ref{Prop_PDE_ODI_RF}\ref{Prop_PDE_ODI_RF_b} it also implies $\partial_\tau \|U^+_{\ge0}\|\ge-C\|U^+\|^2$.
 So we obtain
\begin{align*}
\partial_\tau \|U^+_{\ge0}\|^2
&\ge -2C\|U^+_{\ge0}\|\|U^+\|^2
\ge -C\eta\|U^+_{\ge0}\|^2.
\end{align*}
By \eqref{eq:evo_1}, this implies $\partial_\tau (C_1\|U^+_{\ge0}\|^2)>\partial_\tau (\UU^- + \| U^+_{< 0} \|)$ in a neighborhood of $\tau$, proving the claim.

If \eqref{i:neutral_dom} and \eqref{i:unstable_dom} both hold at time $\tau$, then we have $\|U^+_{>0}\|= \beta\|U^+_0\|>0$, since $\UU^->0$ implies $\|U^+_{\ge0}\|>0$. Then by \eqref{eq:evo_2} \eqref{eq:evo_3} we have $\partial_\tau ( \beta\|U^+_0\|)<\partial_\tau \|U^+_{>0}\|$ in a neighborhood of $\tau$.
\end{proof}

It follows that the times at which conditions (1)(2)(3) hold partition the time interval $[\tau_0,\tau_*]$ into subintervals, arranged in the order in which the corresponding inequalities appear in the assertion of the lemma.
We now choose the partition points to be $T_1$ and $T_3$.
It is clear that Assertions~\ref{l:partition_1_a} and~\ref{l:partition_1_b} hold.

For Assertion \ref{l:partition_1_c}, note that $U^+_{>0}=U^{++}_{>0}$ when $n-k=2$, we may assume $C_1>1$ and $\|U^+_{\conf}\|\le C_1\|U^{++}_{>0}\|$ holds on $(T_3,\tau_*]$ in both cases.
Since when $\eta \leq\ov\eta(\la,J, C_1,\beta)$, we have $\|U^+_{>0}\|\le C\eta$ and, by the conformal domination assumption, $\|U^+_{>0}\|\le C\|U^{++}_{>0}\|$, it follows that
\begin{align*}
    \partial_\tau \| U^{++}_{> 0} \| &\geq \tfrac12 \| U^{++}_{> 0} \|-CC_1\|U^+_{>0}\|^2 \geq (\tfrac12-C\sqrt\eta) \| U^{++}_{> 0} \|.
\end{align*}
This implies Assertion \ref{l:partition_1_c}.
\end{proof}

In the next lemma, we focus on the interval $[T_1, T_3]$.
We first consider the case $n-k\ge 3$, when $\sV_0=\sV_{\aff}\oplus\sV_{\sph,0}$ and $U^+_0=U^+_{\sph,0}+U^+_{\aff}$.
We prove there exist two times dividing $[T_1, T_3]$ into three subintervals, where $U_{0,\max}$, $U^+_{\aff}$, $U_{0,\min}$ are dominant, respectively.

\begin{Lemma}\label{l:partition_3}
Assume $n-k\ge 3$ and the setting of Lemma \ref{l:partition_1}, and assume $C_1,\beta,\eta$ satisfy the conditions in Lemma \ref{l:partition_1}, after decreasing $\ov\beta(\la,J,C_1)$ if necessary.
There exists $c_2(\la,C_1,J,\beta)>0$ such that, after possibly decreasing $\ov\eta(\la,J,C_1,\beta)$, if $\eta\le\ov\eta(\la,J,C_1,\beta)$, we can find unique
times $T_2\le T_2' \in [T_1, T_3]$ such that the following is true
\[ \max \{ c_2 \| U_{\aff}^+ \|, |U_{ 0,\min}| \} \leq  U_{ 0,\max} \qquad \text{on} \quad [T_1, T_2) \]
\[ \max \{ |U_{ 0,\min}|  , U_{ 0,\max} \} \leq  c_2 \| U_{\aff}^+ \|\qquad \text{on} \quad [T_2, T_2') \]
\[ \max \{ c_2 \| U_{\aff}^+ \|, U_{ 0,\max} \} \leq  -U_{ 0,\min} \qquad \text{on} \quad [T_2', T_3). \]
\end{Lemma}

\begin{proof}
Without loss of generality, we may assume $T_1<T_3$, since otherwise we can choose $T_2=T_2'=T_1$ and the assertion holds trivially. 
We use $C(\la,J,C_1)>0$ to denote a generic constant, independent of $0<\beta\leq0.1$, and set $c_2:=\sqrt\beta$.
Note by Lemma \ref{l:partition_1} that $\UU^-\le C_1\|U^+_{\ge 0}\|^2\le C\|U^+_{0}\|^2$ on $[T_1,T_3]$.
In particular, 
$\UU^- (\tau)\neq0$ implies $U^+_{0,\tau}\neq0$ for all $\tau\in[T_1,T_3]$.

By $\PP_{\sV_{0}}Q_2(U^+_0)=\PP_{\sV_{0}}Q_2(U^+_{\sph,0})$ from \eqref{eq_PV0Q2} we have 
\begin{equation}\label{eq:star1}
    \|\PP_{\sV_{\aff}}Q_2(U^+_0)\|\le \|\PP_{\sV_{0}}Q_2(U^+_{\sph,0})\|\le C\|U^+_{\sph,0}\|^2.
\end{equation}
By Lemma \ref{l:partition_1}, we have $\|U^+_{<0}\|\le C_1\|U^+\|^2$ on $[T_1,\tau_*]$.
Using the a priori smallness of $U^+$, this implies $\|U^+_{<0}\|\le C\eta\|U^+\|$.
Using also $\|U^+_{>0}\|\le \beta\|U^+_0\|$ on $[T_1,T_3]$, it implies, after decreasing $\ov\eta(\la,J,C_1,\beta)$, that $\|U^+-U^+_0\|\le C\beta\|U^+_0\|$.
Thus it follows
\begin{equation}\label{eq:star2}
    \|Q_2(U^+)-Q_2(U^+_0)\|\le C\beta\|U^+_0\|^2.
\end{equation}

By \eqref{eq_Us2QU_4}, the positive rotational modes do not contribute to $\PP_{\sV_{\sph,0}}Q_2(U^+_{\geq0})$. The conformal cross terms have zero spherical average, while $Q_2(U^+_{\conf})$ is independent of $\bx$. Thus
\[
\begin{aligned}
\PP_{\sV_{\sph,0}}Q_2(U^+_{\geq0})
&=\PP_{\sV_{\sph,0}}Q_2(U^+_{\sph,0}),\\
\big\|\PP_{\sV_{\sph,0}}\big(Q_2(U^+)-Q_2(U^+_{\sph,0})\big)\big\|
&\leq C\|U^+\|\|U^+_{<0}\|\leq C\|U^+\|^3.
\end{aligned}
\]
Choose $\beta$ sufficiently small and then decrease $\ov\eta(\la,J,C_1,\beta)$ so that
\[
2C\sqrt\beta\leq\tfrac18,
\qquad C\big(\eta+\Psi(\eta)\big)\leq\min\{\beta^2,c_2/8\}.
\]
First, assume the first condition holds. Using \eqref{eq_Uh2QU} and the preceding projection estimate, we obtain
\begin{equation} \label{eq_cmax_evol}
\begin{split}
\big|\partial_\tau U_{0,\max}(\tau)+\sqrt{2}U_{0,\max}^2(\tau)\big|
&\le C\|U^+\|^3+\Psi(\eta)\UU^-\\
&\le C\big(\eta+\Psi(\eta)\big)c_2^{-2}U^2_{0,\max} \\
&\le \beta U_{0,\max}^2(\tau).
\end{split}
\end{equation}

Next, assume the second condition holds. Then for $c_2\le 1$, we have
\[\UU^-+\|U^+\|^2\le C\|U^+_0\|^2\le C(1+c_2^2)\|U^+_{\aff}\|^2\le C\|U^+_{\aff}\|^2.\]
By \eqref{eq:star1}\eqref{eq:star2}, we obtain
\begin{equation}\label{eq_conf}
\begin{split}
\big|\partial_\tau (c_2\|U^+_{\aff}\|)\big|&\le c_2\big(\|\PP_{\sV_{\aff}}Q_2(U^+_0)\|+C\beta\|U^+_0\|^2+C\|U^+\|^3+\Psi(\eta)\UU^-\big)\\
&\le c_2\big(C\|U^+_{\sph,0}\|^2+C(\beta+\eta+\Psi(\eta))\|U^+_{\aff}\|^2\big)\\
&\le C\left(c_2+\frac{\beta+\eta+\Psi(\eta)}{c_2}\right)(c_2\|U^+_{\aff}\|)^2
\le \tfrac14(c_2\|U^+_{\aff}\|)^2.
\end{split}
\end{equation}

Lastly, when the third condition holds, in the same way as \eqref{eq_cmax_evol}, we can show
\begin{equation} \label{eq_cmin_evol2}
\big|\partial_\tau U_{0,\min}(\tau)+\sqrt{2}U_{0,\min}^2(\tau)\big|\le \beta U_{0,\min}^2(\tau). 
\end{equation}
Equivalently, since $U_{0,\min}<0$ in this case,
\[
\partial_\tau |U_{0,\min}|(\tau)\ge (\sqrt2-\beta)U_{0,\min}^2(\tau).
\]

\begin{Claim}
Consider a time $\tau \geq T_1$ at which at least \emph{two} of the inequalities in the assertion hold.
Then, for slightly smaller (resp.~larger) times, only the inequality appearing earlier (resp.~later) in this claim can remain valid.
\end{Claim}

\begin{proof}
We consider three cases.
\medskip

\textit{Case 1: $\max \{ c_2 \| U_{\aff}^+ \|, |U_{ 0,\min}| \} \leq  U_{ 0,\max}$ and $\max \{ |U_{ 0,\min}|  , U_{ 0,\max} \} \leq  c_2 \| U_{\aff}^+ \|$} \\ 
So $|U_{ 0,\min}| \leq U_{ 0,\max} = c_2\| U^+_{\aff} \|>0$.
By \eqref{eq_conf}\eqref{eq_cmax_evol} and $\sqrt2-\beta>1/4$, we see
$\partial_\tau U_{0,\max}< \partial_\tau (c_2\|U^+_{\aff}\|)$, which implies that $U_{0,\max}(\tau')<c_2\|U^+_{\aff}(\tau')\|$ for $\tau' > \tau$ sufficiently close to $\tau$ and the reverse inequality holds for $\tau' < \tau$ sufficiently close to $\tau$.

\medskip

\textit{Case 2: $\max \{ |U_{ 0,\min}|  , U_{ 0,\max} \} \leq  c_2 \| U_{\aff}^+ \|$ and $\max \{ c_2 \| U_{\aff}^+ \|, U_{ 0,\max} \} \leq  -U_{ 0,\min}$ \quad}\\ 
This implies the claim in the same way as Case 1.

\medskip

\textit{Case 3: $\max \{ c_2 \| U_{\aff}^+ \|, |U_{ 0,\min}| \} \leq  U_{ 0,\max}$ and $\max \{ c_2 \| U_{\aff}^+ \|, U_{ 0,\max} \} \leq  -U_{ 0,\min}$} \\ 
So $U_{0,\max} = - U_{0,\min} > 0$, and by \eqref{eq_cmax_evol}\eqref{eq_cmin_evol2} we see 
$\partial_\tau U_{0,\max} < \partial_\tau (- U_{0,\min})$, which implies the claim.
\end{proof}
It follows that the times at which the inequalities from the Claim hold partition the time interval $[T_1,T_3]$ into subintervals, arranged in the order in which the corresponding inequalities appear in the assertion of the lemma.
We now choose the partition points to be $T_2$ and $T_2'$.
\end{proof}

In the next lemma, we consider the case $n-k =2$, when
$\sV_0=\sV_{\aff}\oplus\sV_{\conf}\oplus\sV_{\sph,0}$ and  $U^+_0=U^+_{\sph,0}+U^+_{\aff}+U^+_{\conf}$.
We prove that there exist two times dividing $[T_1,T_3]$ into three subintervals, where $U_{0,\max}$, $U^+_{\aff}+U^+_{\conf}$, $U^+_{0,\min}$ are dominant, respectively.

\begin{Lemma}\label{l:partition_2}
Assume $n-k=2$ and the setting of Lemma \ref{l:partition_1}, and assume $C_1,\beta,\eta$ satisfy the conditions in Lemma \ref{l:partition_1}, after decreasing $\ov\beta(\la,J,C_1)$ if necessary.
There exists $c_2(\la,C_1,J,\beta)>0$ such that, after possibly decreasing $\ov\eta(\la,J,C_1,\beta)$, if $\eta\le\ov\eta(\la,J,C_1,\beta)$, we can find unique
times $T_2\le T_2' \in [T_1, T_3]$ such that the following is true
\[ \max \{ c_2( \| U_{\aff}^+ \|+\| U_{\conf}^+ \|), |U_{ 0,\min}| \} \leq  U_{ 0,\max} \qquad \text{on} \quad [T_1, T_2) \]
\[ \max \{ |U_{ 0,\min}|  , U_{ 0,\max} \} \leq  c_2 ( \| U_{\aff}^+ \|+\| U_{\conf}^+ \|)\qquad \text{on} \quad (T_2, T_2') \]
\[ \max \{ c_2 ( \| U_{\aff}^+ \|+\| U_{\conf}^+ \|), U_{ 0,\max} \} \leq  -U_{ 0,\min} \qquad \text{on} \quad (T_2', T_3] \]
\end{Lemma}

\begin{proof}
Choose $c_2=\sqrt\beta$, decreasing the upper bounds for $\beta$ and then $\eta$ as in the proof of Lemma~\ref{l:partition_3}.
When the first or third condition holds, we can obtain \eqref{eq_cmax_evol} or \eqref{eq_cmin_evol2} in the same way as in Lemma \ref{l:partition_3}.

If the second condition holds, taking $Cc_2<1$ we have 
$$\|U^+_0\|\le (1+Cc_2)(\|U^+_{\aff}\|+\|U^+_{\conf}\|)\le 2(\|U^+_{\aff}\|+\|U^+_{\conf}\|).$$
Moreover, by \eqref{eq_PV0Q2} and Lemma \ref{l:Q_2} we have 
\begin{equation} \label{eq:Q_2_conf}
\PP_{\sV_{\conf}}Q_2(U^+_0)=\PP_{\sV_{\conf}}Q_2(U^+_{\sph,0})=0.
\end{equation}
Together with \eqref{eq:star1}\eqref{eq:star2}, this gives
\begin{equation} \label{eq:con_Q}
\begin{split}
\big|\partial_\tau (c_2\|U^+_{\aff}\|+c_2\|U^+_{\conf}\|)\big|
&\leq Cc_2\big(\|U^+_{\sph,0}\|^2+(\beta+\eta+\Psi(\eta))(\|U^+_{\aff}\|+\|U^+_{\conf}\|)^2\big)\\
&\leq C\left(c_2+\frac{\beta+\eta+\Psi(\eta)}{c_2}\right)(c_2\|U^+_{\aff}\|+c_2\|U^+_{\conf}\|)^2\\
&\leq\tfrac14(c_2\|U^+_{\aff}\|+c_2\|U^+_{\conf}\|)^2.
\end{split}
\end{equation}
So the lemma follows in the same way as Lemma \ref{l:partition_3}, when replacing \eqref{eq_conf} by \eqref{eq:con_Q}.
\end{proof}

We can now prove Theorem \ref{thm:HMHF}.

\begin{proof}[Proof of Theorem \ref{thm:HMHF} and
Proposition~\ref{prop:HMHF_initial}]
By Lemma~\ref{l:non-compact-interval}, it suffices to treat a compact interval $\td I=[\tau_0,\tau_1]$ and to establish \eqref{eq_dtR_between} below.
{Fix once and for all a number $0<\beta\leq \ov\beta(\la,J)$,
small enough for Lemmas~\ref{l:partition_1}--\ref{l:partition_3} below.  After this choice, all constants and smallness thresholds may depend on this fixed $\beta$, and we suppress this dependence from the notation.}
We use $C(m,\la,J)>0$ to denote a generic constant which may change from line to line.
It suffices to prove the assertions with $m$ replaced by
$\max\{m,10\}$, and we make this replacement in the construction.
Choose $\eta,\eps,R^*$ as in
Proposition~\ref{Prop_PDE_ODI_RF}, and then choose $\sigma_0$ for
Proposition~\ref{Prop_gauging} so that
$\sigma_0\leq\ov\sigma(m,\eta)$ and
$\sigma_0<\eta e^{-(R^*)^2/8}$.

For $\Gamma\ne\{1\}$ the conformal factor is $\{0\}$ by Lemma~\ref{l:quotient_cylinder_spectrum}, so there is no conformal boundary face.

We prove the two statements simultaneously.  Let
$C_0(\la,J)>0$ be a sufficiently large constant whose value will be
determined below.  {Let
$\ov r_1=\ov r_1(m,J,\la,\beta,\eta,\eps,R^*,\sigma_0)>0$ be smaller
than every upper bound for $r_1$ imposed below.  For
Theorem~\ref{thm:HMHF}, choose $0<r_1\leq\ov r_1$,
and set}
$\eta_1=\frac{r_1}{C_0\eta}$.  Then choose a finite number
$R_0\geq R^*$ sufficiently large for all the lower bounds involving
$R_0$ below.  Finally, choose $\ov\delta$ so
that the initial $(n,k,\delta,\Gamma)$-neck satisfies the hypotheses of
Lemma~\ref{l:parametrize} with these choices.
For Proposition~\ref{prop:HMHF_initial}, the constants from Theorem \ref{thm:HMHF} have already
been fixed independently of $\eta'$.  {After choosing $\eta'$, choose $0<r_1\leq
\min\left\{\ov r_1,\frac{3\eta'}{4C}\right\}$,
where $C$ is the constant in \eqref{eq:HMHF_initial_smallness}, and set}
$\eta_1=\frac{r_1}{C_0\eta}$, and use the prescribed value of $R_0$.
The lower bound $\underline R_0$ is chosen to imply every lower bound
for $R_0$ below, and $\ov\delta'$ is chosen so that the initial
$(n,k,\delta',\Gamma)$-neck satisfies Lemma~\ref{l:parametrize}.

We will choose $R(\tau)$ finite with
\begin{equation} \label{eq_dtR_between}
-e^{R(\tau)/10} \leq \partial_\tau R(\tau) \leq 1.
\end{equation}

Since $(M,\td g_{\tau_0})$ is an $(n,k,\delta,\Gamma)$-neck, there
exists a diffeomorphism
$\mathcal Y: \Omega\lto
\IB^k_{\delta^{-1}} \times (\IS^{n-k}/\Gamma)$ for an open subset
$\Omega\subset M$, such that
$\mathcal Y(p)\in\bO^k\times(\IS^{n-k}/\Gamma)$ and
\[
\|\mathcal Y_*\td g_{\tau_0}-g_{\cyl}\|_{C^{\delta^{-1}}
(\IB^k_{\delta^{-1}} \times (\IS^{n-k}/\Gamma))}
\leq \delta.
\]
In the proof of Proposition~\ref{prop:HMHF_initial}, we use the same
statement with $\delta$ replaced by $\delta'$.
Assume $r_1<C^{-1}_0\sigma_0$.
Then $r_1<\eta_1$ when choosing $\eta$ small so that $C_0\eta<1$.
{By the choices above,
$R_0\ge\underline R_0(m+10,r_1)$ and
$\delta\le\ov\delta(3R_0,m+10,r_1)$ in Lemma~\ref{l:parametrize}; in the
proof of Proposition~\ref{prop:HMHF_initial}, the latter inequality is
replaced by $\delta'\le\ov\delta(3R_0,m+10,r_1)$.}
By applying Lemma \ref{l:parametrize} and taking $\eta$ as $\eta_1$, and $r$ as $r_1$, with coordinate radius $3R_0$, we obtain
a smooth family of local diffeomorphisms (under the local smooth topology),
\[ (\mathcal Y^{(s)} : \Omega^{(s)} \lto \IB^k_{3R_0} \times (\IS^{n-k}/\Gamma) )_{s \in \ov \IB^{\aff}_{r_1} \times\ov \IB^{\conf}_{r_1}\times\ov\IB^k_{r_1}}, \]
where $\Omega^{(s)}\subset M$ is a continuous family of open subsets,
and $\mathcal Y^{(0)}=\mathcal Y, \Omega^{(0)}=\Omega$, such that
\[\big\|\mathcal Y^{(s)}_*\td g_{\tau_0}-g_{\cyl}\big\|_{C^{m+10}(\IB^k_{3R_0}\times (\IS^{n-k}/\Gamma))}\le C(m)r_1<\sigma_0,\]
after decreasing $r_1$, and for $u^{(s)}=\mathcal Y^{(s)}_*\td g_{\tau_0}-g_{\cyl}$ the following holds for any $s=(v_1,v_2,v_3)$,
\begin{equation} \label{eq:apply_para}
\big\|u^{(s)}\omega_{R_0}-(v_1+v_2)\big\|_{L^2_f}\le Cr_1^2,
\end{equation}
Here the estimate with cutoff $R_0$ follows from the estimate
with cutoff $3R_0$, after increasing the
lower bound for $R_0$ so that the Gaussian tail is at most $r_1^2$.
We also have $\mathcal Y^{(s)} (p) \in\big\{\tfrac{v_3}{\eta_1}\big\}\times (\IS^{n-k}/\Gamma)$.
Let $\chi^{(s)}=\xi_{e^{-\tau_0/2}}\circ\mathcal Y^{(s)}$ be a local diffeomorphism. Then $\mathcal Y^{(s)}=\xi_{e^{\tau_0/2}}\circ\chi^{(s)}$, and
\begin{equation} \label{eq:apply_para2}
\xi_{e^{\tau_0/2}}\circ\chi^{(s)} (p) \in\big\{\tfrac{v_3}{\eta_1}\big\}\times (\IS^{n-k}/\Gamma).
\end{equation}

Applying Proposition \ref{Prop_gauging} with $\sigma_0$ chosen above and $R'=2R_0$, to the restrictions of $\chi^{(s)}$ over $\IB^k_{2.2e^{-\tau_0/2}R_0}\times(\IS^{n-k}/\Gamma)$, we obtain a time $\tau_*(s)\in[\min\{\tau_0+1,\tau_1\},\tau_1]$ and a family of smooth maps (which are local harmonic map heat flows between $\td g_\tau$ and the rescaled Ricci flow of $g_{\cyl}$):
\[(\chi^{(s)}_\tau : \Omega^{(s)}_\tau \lto \DD^{(s)}_\tau )_{\tau\in [\tau_0,\tau_*(s)]},\]
with $\chi^{(s)}_{\tau_0}=\chi^{(s)}|_{\Omega^{(s)}_{\tau_0}}$, $\Omega^{(s)}_{\tau_0}\subset\Omega^{(s)}$, and 
$u^{(s)}_\tau:=(\xi_{e^{\tau/2}})_*(\chi^{(s)}_{\tau})_* \td g_{\tau}-g_{\cyl}$ with $u^{(s)}_{\tau_0}=u^{(s)}$, and by Proposition \ref{Prop_gauging}\ref{Prop_gauging_e}\ref{Prop_gauging_f}, for all $\tau\in[\tau_0,\min\{\tau_0+1,\tau_1\}]$
\begin{equation} \label{eq_atT0bound}
 \IB^k_{R_0+1} \times (\IS^{n-k}/\Gamma)  \subset \DD_{\tau}^{(s)} \qquad \text{and} \qquad  \big\| u_{\tau}^{(s)}  \big\|_{C^{m}(\IB^k_{R_0+1} \times (\IS^{n-k}/\Gamma))} < \eta.
\end{equation}
Moreover, by Proposition \ref{Prop_gauging}\ref{Prop_gauging_h}, 
if $\tau_*(s) < \tau_1$, then the following condition fails at $\tau=\tau_*(s)$:
\begin{equation} \label{eq:open_u1}
\ov{\IB^k_{\sigma_0^{-1}} \times (\IS^{n-k}/\Gamma)} \subset \DD^{(s)}_{\tau}\qquad\textnormal{and}\qquad\Vert u^{(s)}_{\tau} \Vert_{L^2 (\IB^k_{\sigma_0^{-1}} \times (\IS^{n-k}/\Gamma))} < \sigma_0.
\end{equation}
For the fixed $\sigma_0$, we can decrease $r_1$, increase the lower
bound for $R_0$, and then decrease $\delta$ for
Theorem~\ref{thm:HMHF}, or $\delta'$ for
Proposition~\ref{prop:HMHF_initial}, so that \eqref{eq:open_u1} holds for all
$\tau\in[\tau_0,\min\{\tau_0+1,\tau_1\}]$.

Let $\tau^{\#}(s)\in[\tau_0,\tau_*(s)]$ be {maximal} such that \eqref{eq:open_u1} holds for all $\tau\in[\tau_0,\tau^{\#}(s))$,
 and also
\begin{align}\label{eq:open_u2}
\xi_{e^{\tau/2}}\circ\chi^{(s)}_{\tau}(p)\in \IB^k_{\frac{r_1}{\eta_1}}\times (\IS^{n-k}/\Gamma).
\end{align}

\begin{Claim}\label{claim:tau_hash}
$\tau^{\#}(s)$ satisfies the following properties:
\begin{enumerate}[label=(\roman*)]
    \item\label{claim:tau_hash_i} If $\tau^{\#}(s)<\tau_1$, then either \eqref{eq:open_u1} or \eqref{eq:open_u2} fails at $\tau^{\#}(s)$.
    \item\label{claim:tau_hash_ii} $\tau^{\#}(s)$ is lower semi-continuous.
\end{enumerate}
\end{Claim}

\begin{proof}[Proof of Claim]
For Assertion \ref{claim:tau_hash_i}, suppose otherwise, then the maximality of $\tau^{\#}(s)$ implies $\tau^{\#}(s)=\tau_*(s)<\tau_1$, and it follows that \eqref{eq:open_u1} fails at $\tau^{\#}(s)$, a contradiction.

For Assertion~\ref{claim:tau_hash_ii}, fix $\tau_0\leq T<\tau^{\#}(s_\infty)$. Both \eqref{eq:open_u1} and \eqref{eq:open_u2} hold strictly throughout $[\tau_0,T]$. Since $T<\tau_*(s_\infty)$, Proposition~\ref{Prop_gauging}\ref{Prop_gauging_i} implies that both conditions persist on this compact interval for $s$ close to $s_\infty$. Thus $\liminf_{s\to s_\infty}\tau^{\#}(s)\geq\tau^{\#}(s_\infty)$.
\end{proof}

By Lemma \ref{l:lower_con}, we obtain a sequence of continuous functions $\{\tau^{\#}_j\}_{j=1}^{\infty}$ on $\ov \IB^{\aff}_{r_1} \times\ov \IB^{\conf}_{r_1}\times\ov\IB^k_{r_1}$, with $\tau^{\#}_j\le \tau^{\#}$ for each $j$, and for any sequence $s_j\to s_{\infty}$, 
\begin{equation} \label{c:approx_tau}
\tau^{\#}(s_{\infty})\le \liminf_{j\to\infty}\tau^{\#}_j(s_j).
\end{equation}

Assume $\sigma_0<\eta e^{-(R^*)^2/8}$, then $\sigma_0^{-1}\ge R^*$, and \eqref{eq:open_u1} implies (noting also $e^{-f}\le 1$)
\begin{equation}\label{eq:open_u}
    \IB^k_{R^*} \times (\IS^{n-k}/\Gamma) \subset \DD^{(s)}_\tau\qquad \textnormal{and}\qquad \Vert u^{(s)}_\tau \Vert_{L^2_f (\IB^k_{R^*} \times (\IS^{n-k}/\Gamma))} < \eta e^{-(R^*)^2/8}.
\end{equation}
Together with \eqref{eq_atT0bound} and assuming $R_0\ge R^*$, we can apply Proposition \ref{Prop_PDE_ODI_RF} for each $s$ and find functions, depending continuously on $s$, $R^{(s)}:[\tau_0, \tau^{\#}(s)] \longrightarrow[R^*,\infty)$ and $U^{+,(s)} : [\tau_0, \tau^{\#}(s)] \longrightarrow \sV_{> \la}$ and $\UU^{-,(s)} : [\tau_0, \tau^{\#}(s)] \longrightarrow \IR_+$, with
\begin{equation*}
 U^{+,(s)} (\tau) :=  \PP_{\sV_{> \lambda}} (u^{(s)}_{\tau} \omega_{R^{(s)}(\tau)}) , \qquad
\mathcal U^{-,(s)} (\tau) := \big\|   \PP_{\sV_{\leq \lambda}}  (u^{(s)}_{\tau} \omega_{R^{(s)}(\tau)} ) \big\|_{L^2_{f}} + \eta e^{-\frac{((1-\eps)R^{(s)}(\tau))^2}{8}}, 
\end{equation*}
which satisfy all assertions of Proposition \ref{Prop_PDE_ODI_RF}. In particular, by Proposition \ref{Prop_PDE_ODI_RF}\ref{Prop_PDE_ODI_RF_a} we have $\|U^{+,(s)} (\tau)\|,\mathcal U^{-,(s)} (\tau)\le C\eta$.
{In the setting of Proposition~\ref{prop:HMHF_initial},
\eqref{eq:apply_para} and
$\|v_1\|+\|v_2\|\leq2r_1$ imply, uniformly for every parameter $s$,
\begin{equation}\label{eq:HMHF_initial_smallness}
\begin{split}
\|U^{+,(s)}(\tau_0)\|+\mathcal U^{-,(s)}(\tau_0)
&\leq Cr_1+
\eta e^{-\frac{((1-\eps)R_0)^2}{8}}
\leq \frac34\eta'+\frac14\eta'=\eta'.
\end{split}
\end{equation}
Here we used the choice of $r_1$ above and the corresponding condition
on $R_0$ in Proposition~\ref{prop:HMHF_initial}.}
Before $\tau^{\#}(s)$, the cutoff modes depend continuously on $s$ by Proposition~\ref{Prop_gauging}\ref{Prop_gauging_i} and Proposition~\ref{Prop_PDE_ODI_RF}\ref{Prop_PDE_ODI_RF_c}. Since $R^{(s)}$ and $|\partial_\tau R^{(s)}|$ are uniformly bounded here, \eqref{eq_rescale_mod_RDTF} and the $C^m$ bound for $u^{(s)}$ give uniform Lipschitz bounds in time for the quantities in \eqref{eq_UpUUm_def_apply}. Their convergence therefore extends to the endpoint along convergent times.
Moreover, 
since $R^{(s)}(\tau_0)=R_0<\infty$, it follows by Proposition \ref{Prop_PDE_ODI_RF}\ref{Prop_PDE_ODI_RF_c} that for all $\tau\in[\tau_0,\tau^{\#}(s)]$,
\begin{equation} \label{eq:finite_R}
R^{(s)}(\tau)<\infty\qquad\textnormal{and}\qquad\UU^{-,(s)}_{\tau}>0.
\end{equation}

\begin{Claim}\label{claim:continuity_T}
Suppose that $C_1 \geq \underline C_1 (\la, J)$ and $\eta \leq \ov\eta (\la, J, C_1)$ satisfy the bounds from Lemmas~\ref{l:partition_1}, \ref{l:partition_3} (if $n-k \geq 3$) or Lemmas~\ref{l:partition_1}, \ref{l:partition_2} (if $n-k=2$), with the fixed $\beta$ chosen at the beginning of the proof, and recall the constant $c_2(\la, C_1, J) > 0$ from these lemmas.
Then for any $j\in\mathbb N$ and $s \in \ov \IB^{\aff}_{r_1} \times\ov \IB^{\conf}_{r_1}\times\ov\IB^k_{r_1}$, there are unique
$T_{1,j}^{(s)}\le T_{2,j}^{(s)}\le {T'}_{2,j}^{(s)}\le T_{3,j}^{(s)}\in[\tau_0,\tau^{\#}_j (s)]$ satisfying Lemmas~\ref{l:partition_1}, \ref{l:partition_3} (if $n-k \geq 3$) or Lemmas~\ref{l:partition_1}, \ref{l:partition_2} (if $n-k=2$). Moreover, $T_{1,j}^{(s)},T_{2,j}^{(s)},{T'}_{2,j}^{(s)},T_{3,j}^{(s)}$ depend continuously on $s$.
\end{Claim}

\begin{proof}
First, for each $s$ and $j$, by \eqref{eq:finite_R} we can apply Lemmas~\ref{l:partition_1}, \ref{l:partition_3} to $U^{+,(s)},\UU^{-,(s)}$ on the time-interval $[\tau_0, \tau^{\#}_j(s)]$, and find unique times $T_{1,j}^{(s)},T_{2,j}^{(s)},{T'}_{2,j}^{(s)},T_{3,j}^{(s)}$ satisfying these lemmas.
To show continuity, {fix} a sequence $s_i\to s_{\infty}$. Since $\tau^{\#}_j$ is continuous, it follows $\tau^{\#}_j (s_i)\to\tau^{\#}_j (s_{\infty})$ as $i\to\infty$.  
By Lemma \ref{l:partition_1}, 
$T_{1,j}^{(s_{\infty})}$ is the unique time in $[\tau_0,\tau^{\#}_j (s_{\infty})]$ such that
\begin{align*}
\UU^{-,(s_{\infty})} + \| U^{+,(s_{\infty})}_{< 0} \| &\geq C_1 \| U^{+,(s_{\infty})}_{\geq 0} \|^2 \qquad \text{on} \quad [\tau_0, T_{1,j}^{(s_{\infty})}),\\
\UU^{-,(s_{\infty})} + \| U^{+,(s_{\infty})}_{< 0} \| &\leq C_1 \| U^{+,(s_{\infty})}_{\geq 0} \|^2 \qquad \text{on} \quad (T_{1,j}^{(s_{\infty})}, \tau^{\#}_j (s_{\infty})].
\end{align*}
{Suppose} $T_1'$ is a subsequential limit of $T_{1,j}^{(s_i)}\in[\tau_0,\tau^{\#}_j(s_i)]$. Then $T_1'\in [\tau_0,\tau^{\#}_j (s_{\infty})]$. 
Due to the continuous dependence of $U^{+,(s)},\UU^{-,(s)}$ on $s$, it follows from $s_i\to s_{\infty}$ that the above two conditions also hold on $[\tau_0, T_1')$ and $(T_1', \tau^{\#}_j (s_{\infty})]$.
So by the uniqueness of $T_{1,j}^{(s_{\infty})}$ it follows $T_1'=T_{1,j}^{(s_{\infty})}$. This proves the continuity of $T_{1,j}^{(s)}$.
The continuity of $T_{2,j}^{(s)},{T'}_{2,j}^{(s)},T_{3,j}^{(s)}$ follows in the same way.
\end{proof}

We consider first the case $n-k\ge 3$.
For any fixed $j\in\mathbb N$, we define the map 
$F_j:\ov \IB^{\aff}_{r_1} \times\ov \IB^{\conf}_{r_1}\times\ov\IB^k_{r_1} \longrightarrow  \sV_{\aff}\times \sV_{\conf}\times\IR^k$ by
\[ F_j(s) = \big( U^{+,(s)}_{\aff}(T_{2,j}^{(s)}), U^{+,(s)}_{\conf}(T_{3,j}^{(s)}),\eta_1\,e^{\tau^{\#}_j(s)/2}x(\chi^{(s)}_{\tau^{\#}_j(s)}(p))\big),\]
where $x:\IR^k\times(\IS^{n-k}/\Gamma)\to \IR^k$ is the projection map.

\begin{Claim}\label{claim:F-id}
If $C_1 \geq \underline{C}_1(\la,J)$, $r_1\le\ov{r}_1(\la,J,C_1)$, $R_0\ge\underline{R}_0(\la,J,r_1,C_1)$, then for any $(v_1,v_2,v_3)\in \ov \IB^{\aff}_{r_1} \times\ov \IB^{\conf}_{r_1}\times\ov\IB^k_{r_1}$, we have
\begin{enumerate}[label=(\alph*)]
\item\label{i:v1} If $\|v_1\|=r_1$, then $T_{2,j}^{(s)}=\tau_0$, and \begin{align*}
\big\|\PP_{\sV_{\aff}} F_j(v_1,v_2,v_3)-v_1\big\|\le 0.1r_1.
\end{align*}
\item\label{i:v2} If $\|v_2\|=r_1$, then $T_{3,j}^{(s)}=T_{2,j}^{(s)}=\tau_0$ and 
\begin{align*}
\big\|\PP_{\sV_{\aff}\times\sV_{\conf}} F_j(v_1,v_2,v_3)-(v_1,v_2)\big\|\le 0.1r_1.
\end{align*}
\item\label{i:v3} If $\|v_3\|=r_1$, then $\tau^{\#}_j=T_{3,j}^{(s)}=T_{2,j}^{(s)}=\tau_0$, and
\begin{align*}
\big\|F_j(v_1,v_2,v_3)-(v_1,v_2,v_3)\big\|\le 0.1r_1.
\end{align*}
\end{enumerate}

\end{Claim}

\begin{proof}[Proof of Claim]

Fix some $s$ and $j$, and omit the dependence on them. By \eqref{eq:apply_para} we have
\begin{equation} \label{eq:F-id}
\big\|U^{+}_{\aff}(\tau_0)+ U^{+}_{\conf}(\tau_0)-(v_1+v_2)\big\| \le Cr_1^2,
\end{equation}
\begin{equation} \label{eq:other_small}
\|u_{\tau_0}\omega_{R_0}\|-\|U^+_{\conf}(\tau_0)\|-\|U^+_{\aff}(\tau_0)\|\le Cr_1^2.
\end{equation}

For Assertion \ref{i:v1}, assume $Cr_1<0.01$, then by \eqref{eq:F-id} it suffices to show $T_2=\tau_0$.
First, by $\|v_1\|=r_1$ and \eqref{eq:F-id} we have $\|U^+_{\aff}(\tau_0)\|\ge 0.9r_1$.
Assume $e^{-\frac{R_0^2}{16}}\le r_1^2$, then by \eqref{eq:other_small} we have 
$$\UU^-(\tau_0)=\|U^-_{\tau_0}\|+e^{-\frac{(1-\eps)R_0^2}{8}}\le Cr_1^2.$$
Suppose by contradiction $\tau_0<T_1$, then by Lemma \ref{l:partition_1} and \eqref{eq:other_small},
\[C_1\|U^+_{\aff}(\tau_0)\|^2\le  C_1\|U^+_{\ge 0}(\tau_0)\|^2\le \UU^-(\tau_0)+\|U^+_{<0}(\tau_0)\|\le Cr_1^2.\]
 Taking $C_1>4C$, this implies $\|U^+_{\aff}(\tau_0)\|\le 0.5r_1$, a contradiction. Thus $T_1=\tau_0$.
Next, suppose by contradiction $\tau_0=T_1<T_2$, then by Lemma \ref{l:partition_3} and \eqref{eq:other_small},
\[
c_2\|U^+_{\aff}(\tau_0)\|\le U_{0,\max}(\tau_0)\le Cr_1^2,\]
which implies $\|U^+_{\aff}(\tau_0)\|\le c_2^{-1}Cr_1^2<0.5 r_1$ by taking $r_1$ sufficiently small,
a contradiction. So $T_2=\tau_0$, which proves Assertion \ref{i:v1}.
 
 For Assertion \ref{i:v2}, by \eqref{eq:F-id} it suffices to show $T_3=\tau_0$. 
 First, taking $r_1\le\ov{r}_1(\la,J)$, then by $\|v_2\|=r_1$ and \eqref{eq:F-id} we have $\|U^+_{\conf}(\tau_0)\|\in [0.9r_1,1.1r_1]$, and also
 \[\|U^+_0\|\le \|U^+_0-U^+_{\aff}\|+\|U^+_{\aff}\|\le Cr_1^2+1.1 r_1\le 1.2 r_1.\]
Next, arguing as in Assertion~\ref{i:v1}, we obtain $T_1=\tau_0$.
Suppose by contradiction $T_1<T_3$, then by Lemma \ref{l:partition_1}\ref{l:partition_1_b} at $T_1=\tau_0$, and using $\beta<0.1$ we obtain
\[\| U^+_{\conf} (\tau_0)\| \leq\| U^+_{> 0} (\tau_0)\| \leq \beta \| U^+_0 (\tau_0)\|\le 1.2\beta r_1<0.12 r_1,\]
a contradiction. So we obtain $T_3=\tau_0$ and prove Assertion \ref{i:v2}.

For Assertion \ref{i:v3}, since  {${e^{\tau_0/2}} x(\chi^{(s)}_{\tau_0}(p))=\frac{v_3}{\eta_1}$} by \eqref{eq:apply_para2} and $\|v_3\|=r_1$, it follows from the maximality of $\tau^{\#}(s)$ that  $\tau^{\#}(s)=\tau_0=T_{2,j}^{(s)}=T_{3,j}^{(s)}$. This implies Assertion \ref{i:v3} by \eqref{eq:F-id}.
\end{proof}

In the following we fix $C_1=\underline{C}_1(\la,J)$ and omit the dependence of all constants on $C_1$ and on the fixed $\beta$.
	
\begin{Claim}\label{claim:deg=1}
For each $j$ there exists $s_j\in \IB^{\aff}_{r_1}\times \IB^{\conf}_{r_1}\times\IB^k_{r_1}$ such that $F_j(s_j)=(0,0,0)$.
\end{Claim}

\begin{proof}[Proof of Claim]
We fix $j$ and omit the dependence on $j$.
Suppose by contradiction this is not true.
Write $D:=\IB^{\aff}_1\times \IB^{\conf}_1\times\IB^k_{1}$, and let $G(z)=F(r_1z)/r_1$ for any $z\in D$. Since $\bO=(0,0,0)\notin G(D)$, we can define $g(z)=G(z)/\|G(z)\|$ for all $z\in D$.
Writing $G=(G_1,G_2,G_3)$, for each $z=(z_1,z_2,z_3)\in\partial D$ choose $i$ with $\|z_i\|=1$. By Claim~\ref{claim:F-id}, for every $t\in[0,1]$,
\[
\|(1-t)G_i(z)+tz_i\|\geq1-(1-t)\|G_i(z)-z_i\|\geq0.9.
\]
Thus $(1-t)G(z)+tz\neq\bO$ on $\partial D$.
So we can define the homotopy between the two maps $g|_{\partial D}$ and $\frac{z}{|z|}\big|_{\partial D}$ via
\[g_t:\partial D\to \big\{\sV_{\aff} \times \sV_{\conf}\times\IR^k, \| s \|=1\big\},\qquad g_t(z):=\frac{(1-t)G(z)+tz}{\|(1-t)G(z)+tz\|},\qquad t\in[0,1].\]
So the map $g|_{\partial D}$ has degree $1$.
But this contradicts the fact that $g$ can be extended continuously into the interior of $D$ by $\frac{G(z)}{\|G(z)\|}$.  
\end{proof}

For any fixed $j$, let $s_j=(v_{1,j},v_{2,j},v_{3,j})$ be from Claim \ref{claim:deg=1}, and
consider the corresponding solution $(U^{+,(s_j)}, \UU^{-,(s_j)})$, and the unique times $\tau_0\le T_{1,j}^{(s_j)} \leq T_{2,j}^{(s_j)}  \leq {T'}_{2,j}^{(s_j)}  \leq T_{3,j}^{(s_j)}  \leq \tau^{\#}_j(s_j)$ from Claim \ref{claim:continuity_T}.
Then
\begin{equation} \label{eq:T_2_vanish}
U^{+,(s_j)}_{\aff}(T_{2,j}^{(s_j)}) =0, \qquad U^{+,(s_j)}_{\conf}(T_{3,j}^{(s_j)}) = 0,\qquad x(\chi^{(s_j)}_{\tau^{\#}_j(s_j)}(p))=0.
\end{equation}

\begin{Claim}\label{eq:T_2=T_2prime}
$T_{2,j}^{(s_j)} = {T'}_{2,j}^{(s_j)}$.
\end{Claim}

\begin{proof}[Proof of Claim]
Suppose by contradiction $T_{2,j}^{(s_j)} < {T'}_{2,j}^{(s_j)}$.
Then since $U^{+,(s_j)}_{\aff}(T_{2,j}^{(s_j)}) =0$, it follows by Lemma \ref{l:partition_3} that $U^{(s_j)}_{ 0,\min}(T_{2,j}^{(s_j)}) = U^{(s_j)}_{ 0,\max}(T_{2,j}^{(s_j)}) = 0$, and {therefore} $U^{(s_j)}_0 (T_{2,j}^{(s_j)}) = 0$.
Then by Lemma \ref{l:partition_1} this implies $U^{+,(s_j)}(T_{2,j}^{(s_j)})= \UU^{-,(s_j)} (T_{2,j}^{(s_j)}) = 0$, and thus $R^{(s_j)}(T_{2,j}^{(s_j)})=\infty$, which contradicts the non-vanishing {assumption} of $\UU^{-,(s_j)}$ in \eqref{eq:finite_R}. This implies the claim.
\end{proof}

\begin{Claim}\label{claim:decay_dominance}
There is $C_2(\la,J)>0$ such that if $r_1\le\ov{r}_1(\la, J,\sigma_0)$, then
the following holds:
\begin{enumerate}[label=(\alph*)]
\item\label{claim:decay_dominance_a} On $[\tau_0,T^{(s_j)}_{1,j})$, we have
\begin{equation} \label{eq:C_2U-}
\|U^{+,(s_j)}\|+\UU^{-,(s_j)}\le  C_2 r_1,
\end{equation}
 and also
\begin{equation} \label{eq:open_u1/2}
\Vert u^{(s_j)}_{\tau} \Vert_{L^2 (\IB^k_{2\sigma_0^{-1}} \times (\IS^{n-k}/\Gamma))} < \tfrac12\sigma_0\qquad\textnormal{and}\qquad\IB^k_{2\sigma_0^{-1}} \times (\IS^{n-k}/\Gamma) \subset \DD^{(s_j)}_{\tau}.
\end{equation}

\item\label{claim:decay_dominance_b} On $[T^{(s_j)}_{1,j},\tau^{\#}_j(s_j))$, we have
\begin{equation} \label{eq:small_conf_aff}
\| U^{+,(s_j)}_{\aff} \| + \|U^{+,(s_j)}_{\conf} \| \leq {C_2} \|U^{++,(s_j)} \|\quad\textnormal{and}\quad{\UU^{-,(s_j)}+\|U^{+,(s_j)}_{<0}\|}\leq {C_2} \|U^{++,(s_j)} \|^2.
\end{equation}
Moreover, we have
\begin{equation} \label{eq:U++_evolve_1}
\partial_\tau \|U^{++,(s_j)}\|\ge -C\sqrt{\eta}\|U^{++,(s_j)}\|.
\end{equation}
\end{enumerate}
\end{Claim}

\begin{proof}[Proof of Claim]
In the following we fix some $j$ and $s_j$, and omit the dependence on them.
First, by Lemma \ref{l:partition_1} we have $$\partial_\tau (\|U^+_{< 0}\|+\UU^- )\le0\qquad\textnormal{and}\qquad C_1\|U^+_{\ge0}\|^2\le\|U^+_{< 0}\|+\UU^-,\qquad\textnormal{on}\qquad[\tau_0,T_1).$$
Since by \eqref{eq:apply_para} we have $\|U^+_{< 0}(\tau_0)\|+\UU^-(\tau_0)\le Cr_1^2$, this continues to hold on $[\tau_0,T_1)$, and it implies $\|U^+_{\ge 0}\|\le C r_1$. Together, we obtain for $r_1\le\ov{r}_1(\la, J)$,
\begin{align*}
\|U^+ (\tau)\|+\UU^- (\tau)&\le Cr_1.
\end{align*}
This proves \eqref{eq:C_2U-} for some $C_2(\la,J)>0$ after the fixed choice of $\beta$.

To verify \eqref{eq:open_u1/2} on $[\tau_0,T_1)$, by choosing $r_1\le\ov{r}_1(C_2,\sigma_0)$ and using \eqref{eq:C_2U-}, the upper bound on $\UU^-$ implies $R(\tau)\ge 2\sigma_0^{-1}+1$. So
 thus it follows that 
 $$\IB^k_{2\sigma_0^{-1}}\times(\IS^{n-k}/\Gamma)\subset\IB^k_{R(\tau)}\times(\IS^{n-k}/\Gamma)\subset \DD_\tau,$$
which verifies the second part of \eqref{eq:open_u1/2}. The first part of
\eqref{eq:open_u1/2} also follows from \eqref{eq:C_2U-} by taking $r_1\le\ov{r}_1(C_2,\sigma_0)$:
\begin{align*}
\|u\|_{L^2(\IB^k_{2\sigma_0^{-1}}\times(\IS^{n-k}/\Gamma))}
&\le \|u\|_{L^2_f(\IB^k_{2\sigma_0^{-1}}\times(\IS^{n-k}/\Gamma))}\max_{\IB^k_{2\sigma_0^{-1}}\times(\IS^{n-k}/\Gamma)} \sqrt{e^f}
\le C_2 r_1 e^{\frac{\sigma_0^{-2}}{2}}\le\tfrac12\sigma_0.
\end{align*}

Next, we prove Assertion \ref{claim:decay_dominance_b}. We first consider the interval $[T_1,T_3)$.	
Since $T_2=T_2'$ by Claim \ref{eq:T_2=T_2prime}, it follows by Lemma \ref{l:partition_3} and Lemma \ref{l:partition_1}\ref{l:partition_1_b} that on $[T_1,T_3)$,
\[ c_2\|U^+_{\aff}\|\leq\max\{|U_{0,\min}|,U_{0,\max}\},\qquad \|U^+_{\conf}\|\leq\beta\|U^+_0\|. \]
Since $\|U^+_{\sph,0}\|\leq C\max\{|U_{0,\min}|,U_{0,\max}\}$, we obtain
\begin{equation} \label{eq:dominantT_1-T_3}
\| U^+_{\aff} \| + \|U^+_{\conf} \| \leq C\max\{|U_{0,\min}|,U_{0,\max}\} \qquad \text{on} \quad [T_1, T_3).
\end{equation}

Next, we consider the interval $[T_3,\tau^{\#}_j(s_j))$.
Since the interval is empty when $T_3=\tau^{\#}_j(s_j)$, we may assume $T_3<\tau^{\#}_j(s_j)$. 
Assuming $\eta\le\ov\eta(\la,J)$, then Lemma \ref{l:partition_1} implies on $[T_3,\tau^{\#}_j(s_j))$,
\begin{align*}
\|U^+\|&\le \|U^+_{\ge 0}\|+\|U^+_{< 0}\|\le (1+C_1\|U^+_{\ge 0}\|)\|U^+_{\ge 0}\|\\
&\le 1.1(\|U^+_{0}\|+\|U^+_{> 0}\|)\le C\beta^{-1}\|U^+_{> 0}\|.
\end{align*}
Let $U^{++}_{>0}=U^{+}_{\rot,\frac12}+U^{+}_{\rot,1}$. By \eqref{eq:T_2_vanish}, $U^+_{>0}(T_3)=U^{++}_{>0}(T_3)$, so for some $C_2(\la,J)>1$ after the fixed choice of $\beta$ we have
\begin{equation} \label{eq:domi_U++}
\|U^+_{\aff}(T_3)\|+\|U^+_{\conf}(T_3)\|\le 2\|U^+(T_3)\|\le C\beta^{-1}\|U^+_{>0}(T_3)\|\le C_2\|U^{++}_{> 0}(T_3)\|
\end{equation}
The eigenvalues of $L$ on $\sV_{\conf}$ and $\sV_{\aff}$ are $c_{n-k}:=1-\frac{n-k}{2(n-k-1)}$ and $0$, respectively. With $C_2$ fixed, choose $\eta$ small enough that $C(C_2)\|U^+\|^2\le \beta(\|U^{++}_{>0}\|+\|U^+_{\conf}\|)$. Since $\beta<\frac12(\frac12-c_{n-k})$ and $C_2\ge1$, we have $\frac12C_2-\beta\ge(c_{n-k}+\beta)C_2$. The first ODI in Proposition~\ref{Prop_PDE_ODI_RF}\ref{Prop_PDE_ODI_RF_b} therefore gives
\begin{align*}
\partial_\tau(C_2\|U^{++}_{>0}\|-\|U^+_{\conf}\|-\|U^+_{\aff}\|)&\ge \tfrac12C_2\|U^{++}_{>0}\|-c_{n-k}\|U^+_{\conf}\|-C(C_2)\|U^{+}\|^2\\
&\ge (\tfrac12C_2-\beta)\|U^{++}_{>0}\|-(c_{n-k}+\beta)\|U^+_{\conf}\|\\
&\ge (c_{n-k}+\beta)(C_2\|U^{++}_{>0}\|-\|U^+_{\conf}\|-\|U^+_{\aff}\|).
\end{align*}
By \eqref{eq:domi_U++} and comparison, $\|U^+_{\aff}\|+\|U^+_{\conf}\|\le C_2\|U^{++}_{>0}\|$ on $[T_3,\tau^{\#}_j(s_j))$. Together with \eqref{eq:dominantT_1-T_3}, this implies the first inequality in \eqref{eq:small_conf_aff} on $[T_1,\tau^{\#}_j(s_j))$. Since $U^+_{\ge0}=U^{++}+U^+_{\aff}+U^+_{\conf}$, it follows by the first inequality on $[T_1,\tau^{\#}_j(s_j))$ that
$\|U^+_{\ge0}\|\le C(C_2)\|U^{++}\|$.
Thus, the second inequality in \eqref{eq:small_conf_aff} follows by Lemma \ref{l:partition_1} for a possibly larger $C_2$.
\end{proof}

After passing to a subsequence, we may assume $s_j\to s_{\infty}$.	
By \eqref{c:approx_tau} and the continuity just proved, the mode quantities converge uniformly on compact intervals before $\tau^{\#}(s_\infty)$ and at times in their domains tending to $\tau^{\#}(s_\infty)$.
In the following we will show $\tau^{\#}(s_\infty)=\tau_1$, and verify that $\chi^{(s_{\infty})}$ and $u^{(s_\infty)}_\tau,U^{+,(s_\infty)}(\tau),\UU^{-,(s_{\infty})}(\tau)$ satisfy all assertions in Theorem \ref{thm:HMHF} on $[\tau_0,\tau_1]$.

\begin{Claim}\label{claim:u=eta}
If $\tau^{\#}(s_\infty)<\tau_1$, then
$\tau^{\#}(s_{\infty})>\tau_0+1$, and \eqref{eq:open_u1} fails at $s_{\infty}$, that is,
\[\ov{\IB^k_{\sigma_0^{-1}} \times (\IS^{n-k}/\Gamma)} \not\subset \DD^{(s_\infty)}_{\tau^{\#}(s_{\infty})}\qquad \textnormal{or}\qquad\big\Vert u^{(s_\infty)}_{\tau^{\#}(s_{\infty})} \big\Vert_{L^2 (\IB^k_{\sigma_0^{-1}} \times (\IS^{n-k}/\Gamma))} = \sigma_0.\]
\end{Claim}

\begin{proof}[Proof of Claim]
Suppose by contradiction that \eqref{eq:open_u1} holds at $\tau^{\#}(s_{\infty})$. Then it follows by Claim \ref{claim:tau_hash}\ref{claim:tau_hash_i} that \eqref{eq:open_u2} must fail at $\tau^{\#}(s_{\infty})$, that is,
\begin{align}\label{eq:max_xi}
\eta_1e^{\tau^{\#}(s_{\infty})/2}\|x(\chi^{(s_\infty)}_{\tau^{\#}(s_{\infty})}(p))\|=r_1.
\end{align}
Since $\chi^{(s_\infty)}_\tau$ satisfies the harmonic map heat flow equation $\partial_\tau(\chi^{(s_\infty)}_\tau)=\Delta_{\td g_{\tau},\td g_{\cyl,\tau}}\chi^{(s_\infty)}_\tau$, we have
\begin{align*}
    \partial_\tau\big(\xi_{e^{\tau/2}}\circ\chi^{(s_\infty)}_\tau(p)\big)&=\tfrac12e^{\tau/2}x(\chi^{(s_\infty)}_\tau(p)) + (\xi_{e^{\tau/2}})_*\partial_\tau(\chi^{(s_\infty)}_\tau(p))\\
    &=\tfrac12e^{\tau/2}x(\chi^{(s_\infty)}_\tau(p)) + (\xi_{e^{\tau/2}})_*(\Delta_{\td g_{\tau},\td g_{\cyl,\tau}}\chi^{(s_\infty)}_\tau)(p)\\
    &=\tfrac12e^{\tau/2}x(\chi^{(s_\infty)}_\tau(p)) + \big(\Delta_{g_{\cyl}+u^{(s_\infty)}_\tau,g_{\cyl}}\id\big)\circ\xi_{e^{\tau/2}}\circ\chi^{(s_\infty)}_\tau(p).
\end{align*}
Projecting this equation to $\IR^k$ and multiplying with $\eta_1$,
seeing {$\eta_1= \frac{r_1}{C_0\eta}$}, $\|u^{(s_\infty)}_\tau\|_{C^1(\IB^k_1\times(\IS^{n-k}/\Gamma))}\le C\eta$, \eqref{eq:max_xi}, and choosing $C_0(\la,J)$ sufficiently large, this implies
 at $\tau=\tau^{\#}(s_{\infty})$
\begin{align*}
\partial_\tau \big\|\eta_1e^{\tau/2}x(\chi^{(s_\infty)}_\tau(p))\big\|
&\geq\tfrac12\|\eta_1e^{\tau/2}x(\chi^{(s_\infty)}_\tau(p))\| - C\eta_1{\|u^{(s_\infty)}_\tau\|_{C^1(\IB^k_1\times(\IS^{n-k}/\Gamma))}}
\ge \tfrac12r_1-C\eta\eta_1>0,
\end{align*}
and thus $\xi_{e^{\tau/2}}\circ\chi^{s_\infty}_\tau(p)\notin\overline{\IB}^k_{\frac{r_1}{\eta_1}}\times (\IS^{n-k}/\Gamma)$ for $\tau$ slightly larger than $\tau^{\#}(s_{\infty})$.
Since for each $j$, we have $\xi_{e^{\tau/2}}\circ\chi_{\tau}^{(s_j)}(p)\in {\IB}^k_{\frac{r_1}{\eta_1}}\times (\IS^{n-k}/\Gamma)$ for all $\tau\in[\tau_0,\tau^{\#}_j(s_j))$, this implies $\tau^{\#}(s_{\infty})\ge\limsup_{j\to\infty}\tau^{\#}_j(s_j)$.
Moreover, since $\tau^{\#}(s_{\infty})\le\liminf_{j\to\infty}\tau^{\#}_j(s_j)$ by \eqref{c:approx_tau}, we obtain
\[\tau^{\#}(s_{\infty})=\lim_{j\to\infty}\tau^{\#}_j(s_j).\]
Thus $x(\chi_{\tau^{\#}_j(s_j)}^{(s_j)}(p))=0$ passes to the limit, contradicting \eqref{eq:max_xi}. Hence \eqref{eq:open_u1} fails at $\tau^{\#}(s_\infty)$. Since it holds throughout the first step and $\tau^{\#}(s_\infty)<\tau_1$, we have $\tau^{\#}(s_\infty)>\tau_0+1$.
\end{proof}

Now we suppose by 
contradiction $\tau^{\#}(s_\infty)<\tau_1$.
After passing to a further subsequence, we may assume $T^{(s_j)}_{1,j}\to T_{\infty}\in[\tau_0,\tau_1]$.
We claim $T_{\infty}<\tau^{\#}(s_{\infty})$. Indeed, if not, then 
\eqref{eq:open_u1/2} passing to limit implies 
$U^{+,(s_\infty)},\UU^{-,(s_\infty)}$ satisfy \eqref{eq:open_u1} at $\tau^{\#}(s_{\infty})$, which contradicts with Claim \ref{claim:u=eta}.
So $T_{\infty}<\tau^{\#}(s_{\infty})$, and we can pass Claim \ref{claim:decay_dominance}\ref{claim:decay_dominance_b} to the limit and obtain for all $\tau\in(T_1^{(s_\infty)},\tau^{\#}(s_\infty))$:
\begin{equation} \label{eq:star}
\| U^{+,(s_{\infty})}_{\aff} \| + \|U^{+,(s_{\infty})}_{\conf} \| \leq C_2 \|U^{++,(s_{\infty})} \|,\qquad{\UU^{-,(s_{\infty})}+\|U^{+,(s_{\infty})}_{<0}\|}\leq C_2 \|U^{++,(s_{\infty})} \|^2,
\end{equation}
and also
\begin{equation} \label{eq:U++_evolve_2}
\partial_\tau \|U^{++,(s_{\infty})}\|\ge -C\sqrt{\eta}\|U^{++,(s_{\infty})}\|.
\end{equation}
These prove Theorem \ref{thm:HMHF}\ref{thm:HMHF_b}.

In the following we will only consider $u^{(s_\infty)}_{\tau^{\#}(s_\infty)},U^{+,(s_\infty)}_{\tau^{\#}(s_\infty)},\UU^{-,(s_\infty)}_{\tau^{\#}(s_\infty)}$ and abbreviate them as $u,U^+,\UU^-$.
Since $U^+_{\ge 0}$ is finite dimensional, it follows by \eqref{eq:star}
\begin{equation*} 
\|U^+_{\ge 0} \|_{C^m({\IB^k_{200} \times (\IS^{n-k}/\Gamma)})}
\leq C_m\| U^+_{\ge 0}  \|\leq C_m\| U^{++}  \|. 
\end{equation*}
Moreover, since $\tau^{\#}(s_{\infty})> \tau_0+1$ by Claim \ref{claim:u=eta}, by the {standard parabolic estimate} we have
\[\|u \|_{C^m({\IB^k_{200} \times (\IS^{n-k}/\Gamma)})} \leq C_m\|U^{++} \|. \]
Denote $U^-+U^+_{<0}$ by $\td U^-$.
Since $u = U^+_{\ge0} +\td U^-$ on $\IB^k_{200} \times (\IS^{n-k}/\Gamma)$, by triangle inequality we obtain
\[ \|\td U^- \|_{C^m({\IB^k_{200} \times (\IS^{n-k}/\Gamma)})}  \leq C_m\| U^{++}  \|. \]
Then we can use the integral bound $\|\td U^-\|\leq C \|U^{++} \|^2$ to improve this bound at the cost of losing one derivative and shrinking the domain:
\[ \|\td U^- \|_{C^{m-1}({\IB^k_{100} \times (\IS^{n-k}/\Gamma)})} \leq C_m\| U^{++}  \|^{1+\eps_0(m,n,k)}.\]
This step follows from a Gagliardo--Nirenberg inequality; see \cite[Lemma~5.107]{Bamler_Lai_MCF1}.
Since $u= U^+_{\conf}+U^+_{\aff}+U^{++}+\td U^-$, it follows 
\begin{equation} \label{eq:dom_modulo_Lie}
\|u -U^+_{\conf} -U^+_{\aff} -U^{++} \|_{C^3(\IB^k_{100}\times(\IS^{n-k}/\Gamma))}\le C\|U^{++} \|^{1+\eps_0}.
\end{equation}

Assume $X=\tfrac12\sum_{i,j=1}^k a_{ij} \, \bx_i \partial_{\bx_j}+ \sum_{i=1}^{n-k+1}b_{i}\nabla \by_i$ is a vector field  on $\IR^k\times (\IS^{n-k}/\Gamma)$ (with $b_i=0$ when $\Gamma$ is nontrivial) such that $U^+_{\conf}+U^+_{\aff}=\LL_{X}g_{\cyl}$. Then by \eqref{eq:star} we have $\|X\|\le C\|\LL_{X}g_{\cyl}\|\le C\|U^{++}\|$.
Let $\phi$ be the flow generated by the vector field $-X$ at time $1$. Then
\[\|\phi-\id\|_{C^4(\IB^k_{100}\times (\IS^{n-k}/\Gamma))}\le C\|X\|\le  C\|U^{++}\|,\]
and by Lemma \ref{claim:Lie}  we have
$$\big\|\phi^* g_{\cyl}-g_{\cyl}+U^+_{\conf}+U^+_{\aff}\big\|_{C^3(\IB^k_{100}\times (\IS^{n-k}/\Gamma))}\le  C\|U^{++}\|^2.$$
Letting $\td u=\phi^*(g_{cyl}+u)-g_{cyl}$ and writing $\mathcal B=\IB^k_{100}\times (\IS^{n-k}/\Gamma)$, this together with
 \eqref{eq:dom_modulo_Lie} implies
\begin{align*}
\big\|\td u-U^{++}\big\|_{C^3(\mathcal B)}&\le\big\|\td u-(u-U^+_{\conf}-U^+_{\aff})\big\|_{C^3(\mathcal B)}+\Psi(\eta)\|U^{++}\|\\
&=\big\|\phi^*(g_{cyl}+u)-(g_{cyl}+u)+U^+_{\conf}+U^+_{\aff}\big\|_{C^3(\mathcal B)}+\Psi(\eta)\|U^{++}\|\\
&\le\big\|\phi^*g_{\cyl}-g_{\cyl}+U^+_{\conf}+U^+_{\aff}\big\|_{C^3(\mathcal B)}+\big\|\phi^*u-u\big\|_{C^3(\mathcal B)}+\Psi(\eta)\|U^{++}\|\\
&\le C\|U^{++}\|^2+C\|U^{++}\|\|u\|_{C^4(\IB^k_{200}\times(\IS^{n-k}/\Gamma))}+\Psi(\eta)\|U^{++}\|
\le C\Psi(\eta)\|U^{++}\|.
\end{align*}
where $\Psi(\eta)\to 0$ as $\eta\to 0$.

Since $g_{\cyl}+u=(\xi_{e^{\tau^{\#}(s_\infty)/2}})_*(\chi_{\tau^{\#}(s_\infty)})_* \td g_{\tau^{\#}(s_\infty)}$, and by assumption $\td g_{\tau^{\#}(s_\infty)}$ is an $(n,k,\delta,\Gamma)$-cylinder at $p$ at scale $1$, it follows $g_{\cyl}+u$ is an $(n,k,\delta,\Gamma)$-cylinder centered at $$\xi_{e^{\tau^{\#}(s_\infty)/2}}\circ\chi_{\tau^{\#}(s_\infty)}(p)\in\ov\IB^k_{\frac{r_1}{\eta_1}}\times(\IS^{n-k}/\Gamma)$$
at scale $1$.
Since $\frac{r_1}{\eta_1}\le 1$, this implies $g_{\cyl}+u$ is an $(n,k,2\delta,\Gamma)$-cylinder at some point on $\bO^k\times(\IS^{n-k}/\Gamma)$  at scale $1$.
Since $\phi(\bO^k\times(\IS^{n-k}/\Gamma))= \bO^k\times(\IS^{n-k}/\Gamma)$, it follows that $\phi^*(g_{\cyl}+u)$ is also an $(n,k,2\delta,\Gamma)$-cylinder at some point on $\bO^k\times(\IS^{n-k}/\Gamma)$  at scale $1$. Therefore, applying Lemma \ref{l:dom_curvature} below to $\phi^*(g_{\cyl}+u)$, we obtain $$\|U^{++}\|\le {C{\delta}},$$
which by \eqref{eq:star} implies
\[\|U^+\|+\UU^-\le C\|U^{++}\|\le C\delta.\]
{By the choice of $\ov\delta$ at the beginning
of the proof}, the same argument in Claim \ref{claim:decay_dominance}\ref{claim:decay_dominance_a} implies that \eqref{eq:open_u1} holds at $\tau^{\#}(s_{\infty})$,
contradicting Claim \ref{claim:u=eta}. Thus $\tau^{\#}(s)=\tau_1$.

Moreover, \eqref{eq:open_u1}\eqref{eq:open_u2} imply \eqref{thm:HMHF_open_1}\eqref{thm:HMHF_open_2} (by replacing $\eta$ by $\eta/C_0$ in the argument above), 
\eqref{eq:star} verifies Theorem \ref{thm:HMHF}\ref{thm:HMHF_b}. Theorem \ref{thm:HMHF}\ref{thm:HMHF_a} is clear by the definition of $T_1^{(s_\infty)}$ and Lemma \ref{l:partition_1}. Claim \ref{eq:T_2=T_2prime} implies $T_2^{(s_\infty)}={T'}_2^{(s_\infty)}$, which verifies the dominance alternatives in Theorem \ref{thm:HMHF}\ref{thm:HMHF_c} by Lemma \ref{l:partition_3}. \eqref{eq:T_2_vanish} implies $U^{+,(s_\infty)}_{\conf}(T_3)=0$, verifying Theorem \ref{thm:HMHF}\ref{thm:HMHF_d}.
Finally, Theorem \ref{thm:HMHF}\ref{thm:HMHF_e} follows from Theorem \ref{thm:HMHF}\ref{thm:HMHF_b} and Lemma \ref{l:partition_1}.
So we prove Theorem \ref{thm:HMHF} in the case $n-k\ge3$.

\bigskip

Next, we consider
the case when $n-k=2$.
We only need to slightly modify the argument of the case when $n-k\ge3$.
First, we consider the map 
$F_j:\ov \IB^{\aff}_{r_1} \times\ov \IB^{\conf}_{r_1}\times\ov\IB^k_{r_1} \longrightarrow  \sV_{\aff}\times \sV_{\conf}\times\IR^k$ given by
\[ F_j : s \mapsto \big( U^{+,(s)}_{\aff}(T_{2,j}^{(s)}) , U^{+,(s)}_{\conf}(T_{2,j}^{(s)}), \eta_1\,e^{\tau^{\#}_j(s)/2}  x(\chi^{(s)}_{\tau^{\#}_j(s)}(p)) \big). \]
We can verify that Claim \ref{claim:F-id} still holds.
In fact, we can show $T_2^{(s)}=\tau_0$ by using Lemma \ref{l:partition_1} and \ref{l:partition_2} when $\|s_1\|=r_1$ or $\|s_2\|=r_1$, and show $\tau^{\#}_j=\tau_0$ when $\|s_3\|=r_1$.
So by Claim \ref{claim:deg=1} there exists $s_j=(v_{1,j}, v_{2,j},v_{3,j})$ such that $F_j(s_j) = (0,0,0)$ and the corresponding solution $(U^{+,(s_j)}, \UU^{-,(s_j)})$ satisfies
\[ U^{+,(s_j)}_{\aff}(T^{(s_j)}_{2,j}) =0, \qquad U^{+,(s_j)}_{\conf}(T^{(s_j)}_{2,j}) = 0,\qquad x(\chi^{(s_j)}_{\tau^{\#}_j(s_j)}(p))=0, \]
and we also have $T^{(s_j)}_{2,j} = {T'}_{2,j}^{(s_j)}$ by Lemma \ref{l:partition_2}.
The rest of the argument is similar, and the estimates for $\|U^+_{\conf}\|+\|U^+_{\aff}\|$ are obtained in the same way as those for $\|U^+_{\aff}\|$ when $n-k\ge3$, using Lemmas~\ref{l:partition_1} and~\ref{l:partition_2}.

The construction gives
$\|U^+_{>0}\|\leq C\sqrt\beta\,\|U^{++}_0\|$,
where $C$ is independent of $\beta$.
Replacing $\beta$ in the construction by a sufficiently small
multiple of $\beta^2$ proves \eqref{eq:neutral_dom_rough},
\eqref{eq_cmax_evol_0}, and \eqref{eq_cmin_evol_0}.

{In the setting of Proposition~\ref{prop:HMHF_initial},
\eqref{eq:HMHF_initial_smallness} passes to the limit and
$R^{(s)}(\tau_0)=R_0$ together prove
Proposition~\ref{prop:HMHF_initial}.}
\end{proof}

\bigskip

In the following lemma, we consider a metric with dominant $U^{++}\in\sV^{++}$ mode.
We show that a sufficiently good local $(n,k,\delta)$-neck forces $\|U^{++}\|\leq C\delta$.
In the proof, we first remove the affine component by a linear change of Euclidean coordinates.
The Ricci curvature on a fixed ball then controls all coefficients of $U^{++}$.
This curvature bound is intrinsic, so the neck coordinates need not agree with the coordinates used to express $U^{++}$.

\begin{Lemma}\label{l:dom_curvature}
There is a dimensional constant $C>0$ such that the following holds if $\delta\le\ov\delta$.
Suppose $g$ is a smooth metric on $\IR^k\times\IS^{n-k}$ and that there exist $U^{++}\in\sV^{++}$ and $U_{\aff}\in\sV_{\aff}$ with $\|U^{++}\|+\|U_{\aff}\|\le C^{-1}$, satisfying
\begin{align*}
\big\|g-g_{\cyl}-U^{++}-U_{\aff}\big\|_{C^{2}(\IB^k_{1}\times\IS^{n-k})}\le C^{-1}\|U^{++}\|.
\end{align*}
Suppose also that some point $\mathbf p \in \IB^k_1 \times \IS^{n-k}$ is a center of an $(n,k,\delta)$-neck at scale $1$ for $g$, or more generally, that a $2$-ball around $\mathbf p$ is $\delta$-close to an open subset of the round cylinder.
Then $\|U^{++}\|\le C{\delta}$.
\end{Lemma}

\begin{proof}
We will determine $C$ in the course of the proof and let $C'$ be a generic constant.
Write $U_{\aff}=\sum h_{ij}\,d\bx_i\otimes d\bx_j$.
Pushing forward by $(\bx,\by)\mapsto((I+(h_{ij}))^{1/2}\bx,\by)$ reduces to $U_{\aff}=0$, since this map preserves $\sV^{++}$ and changes the relevant norms by bounded factors.
After restricting to $\IB^k_{0.9}\times\IS^{n-k}$, the transformed metric still satisfies the $C^2$ estimate in the statement up to a dimensional factor, and we retain the same notation.
The change of coordinates moves $\mathbf p$ by $O(\|U_{\aff}\|)$.
Since the comparison region has radius $\delta^{-1}$ (respectively, $2$) and $g$ is close to $g_{\cyl}$ on $\IB^k_1\times\IS^{n-k}$, it contains $\IB^k_{0.9}\times\{\by_{\mathbf p}\}$ when $\delta$ is sufficiently small, where $\by_{\mathbf p}$ is the spherical coordinate of $\mathbf p$.
After an orthogonal change of Euclidean coordinates, we may write $U^{++}=(a+\sum_{i=1}^k b_i \mathfrak p^{(1)}_{i}
+\sum_{i=1}^k c_i \mathfrak p^{(2)}_{ii})g_{\IS^{n-k}}$.
Set $\Phi = 1 + a+\sum_{i=1}^k b_i \mathfrak p^{(1)}_{i}
+\sum_{i=1}^k c_i \mathfrak p^{(2)}_{ii}$.
For sufficiently large $C$ the metric $g_0 = g_{\cyl} + U^{++}=g_{\IR^k}+\Phi(\bx)g_{\IS^{n-k}}$ is a Riemannian metric on $\IB^k_1\times\IS^{n-k}$.
Recall that $\IS^{n-k}$ has radius $\sqrt{2(n-k-1)}$.
With respect to a $g_0$-orthonormal frame adapted to the product, write $i,j$ for Euclidean indices and $\alpha,\beta$ for spherical indices.
Up to curvature symmetries, the nonzero components are represented by the following formulas, with $\alpha\ne\beta$:
\begin{align*}
(\Rm_{g_0})_{i\alpha \alpha j} &= -\frac{\nabla^2_{ij} \Phi}{2\Phi} + \frac{\nabla_i \Phi \nabla_j \Phi}{4\Phi^2} = -\frac{c_i \delta_{ij}}{2\sqrt2\,\Phi}
+\frac{(b_i+c_i \bx_i)(b_j+c_j \bx_j)}{8\,\Phi^2} = -\frac{c_i \delta_{ij}}{2\sqrt2\,\Phi} + O(\| U^{++} \|^2), \\ 
(\Rm_{g_0})_{\alpha\beta\beta\alpha} &=\frac{1}{2(n-k-1)\Phi} - \frac{|\nabla \Phi|^2}{4\Phi^2} = \frac{1}{2(n-k-1)\,\Phi}
-\frac{1}{8\Phi^2}\sum_{i=1}^k (b_i+c_i \bx_i)^2 \\
&= \frac1{2(n-k-1)} - \frac{\Phi-1}{2(n-k-1)\Phi} + O(\| U^{++} \|^2),
\end{align*}
where $O(\| U^{++} \|^2)$ denotes a term bounded by $C' \| U^{++} \|^2$.
Consequently,
\begin{align}
(\Ric_{g_0})_{ij}&=-(n-k)\frac{c_i}{2\sqrt2\,\Phi}\delta_{ij}
+O(\| U^{++} \|^2), \label{eq_Ric1_approx}
\\
(\Ric_{g_0})_{\alpha\alpha}&=\frac12-\frac{\Phi-1}{2\Phi}-\frac{1}{2\sqrt2\,\Phi}\sum_{i=1}^k c_i+O(\|U^{++}\|^2). \label{eq_Ric2_approx}
\end{align}

By either comparison hypothesis, the eigenvalues of $\Ric_g$ on $\IB^k_{0.9}\times\{\by_{\mathbf p}\}$ are $C'\delta$-close to $0$ (with multiplicity $k$) and $\frac12$ (with multiplicity $n-k$).
Since
\[
\|{\Ric_{g_0}-\Ric_g}\|_{C^0(\IB^k_{0.9}\times\IS^{n-k})}
\leq C'\|g_0-g\|_{C^2(\IB^k_{0.9}\times\IS^{n-k})}
\leq C'C^{-1}\|U^{++}\|,
\]
the eigenvalues of $\Ric_{g_0}$ at these points are $C'(\delta+C^{-1}\|U^{++}\|)$-close to $0$ (with multiplicity $k$) and $\frac12$ (with multiplicity $n-k$).
The Ricci tensor of $g_0$ preserves the product splitting, and the Euclidean and spherical blocks remain separated by a fixed spectral gap.
Thus the $k$ eigenvalues near $0$ belong to the Euclidean block, and \eqref{eq_Ric1_approx} first gives $|c_i|\leq C'(\delta+C^{-1}\|U^{++}\|+\|U^{++}\|^2)$.
Substituting this bound into \eqref{eq_Ric2_approx} gives the same estimate for $|\Phi(\bx)-1|$ when $|\bx|\leq0.9$.
Since $\|U^{++}\|\leq C^{-1}$, we conclude that for $|\bx|\leq0.9$,
\[ |c_i| + |\Phi(\bx)-1| \leq C'(\delta+C^{-1}\|U^{++}\|) + C' \| U^{++} \|^2  \leq C'(\delta+2C^{-1}\|U^{++}\|) . \]
The values of $\Phi$ at finitely many points of $\IB^k_{0.9}$ determine $a,b_i,c_i$, so
\[
\|U^{++}\|\leq C'\delta+2C'C^{-1}\|U^{++}\|.
\]
So choosing $C \geq 4 C'$ implies the desired bound.
\end{proof}

We now extend the HMHF construction from compact to noncompact time intervals.

\begin{Lemma}\label{l:non-compact-interval}
{Assume Theorem \ref{thm:HMHF} and
Proposition~\ref{prop:HMHF_initial} hold for compact intervals}, using the radius choice \eqref{eq_ODE_for_R_tau}, and in addition $R(\tau)$ satisfies
\begin{equation} \label{eq_dtR_between_1}
-e^{R(\tau)/10} \leq \partial_\tau R(\tau) \leq 1.
\end{equation}
{Then Theorem \ref{thm:HMHF} and
Proposition~\ref{prop:HMHF_initial} also hold for non-compact intervals.}
\end{Lemma}

\begin{proof}
Exhaust $I$ by compact intervals $I_i=[t_{0,i},t_{1,i}]$, each containing any endpoint that belongs to $I$. Use the same $R_0$ on every $I_i$ and, in the setting of Proposition~\ref{prop:HMHF_initial}, the same initial neck chart.
By applying {the corresponding compact-interval statement} on $I_i$, we obtain 
a family of open domains $\Omega_{i,t}\subset \MM_t$, with $p\in \Omega_{i,t}$, and a family of diffeomorphisms:
\[(\chi_{i,t} : \Omega_{i,t} \lto \xi_{\sqrt{-t}}(\DD_{i,-\log(-t)})\subset \IR^k \times (\IS^{n-k}/\Gamma) )_{t\in I_i},\]
which evolves by harmonic map heat flow between $(\Omega_{i,t},g_t)$ and $(\xi_{\sqrt{-t}}(\DD_{i,-\log(-t)}),g_{\cyl,t})$,
such that setting $\tau=-\log(-t)$, then for all $\tau\in \td I_i=[-\log(-t_{0,i}),-\log(-t_{1,i})]$,
\begin{equation} \label{eq:base_pt}
\xi_{e^{\tau/2}}\circ\chi_{i,-e^{-\tau}}(p)\in \IB^k_{\eta}\times(\IS^{n-k}/\Gamma),
\end{equation}
and for $u_{i,\tau}:=e^\tau(\xi_{e^{\tau/2}})_*(\chi_{i,-e^{-\tau}})_* g_{-e^{-\tau}}-g_{\cyl}$,
and there is a $C^1$ function $R_i: \td I_i  \to [R^*, \infty)$ with $R_i(-\log(-t_{0,i}))={R_0}$, which satisfies \eqref{eq_dtR_between_1},
and such that for all $\tau \in \td I_i$, 
\begin{equation} \label{eq:u_i_small}
\IB^{k}_{ R_i(\tau)} \times (\IS^{n-k}/\Gamma)  \subset \DD_{i,\tau}\qquad\textnormal{and}\qquad\| u_{i,\tau} \|_{C^{m}(\IB^k_{R_i(\tau)} \times (\IS^{n-k}/\Gamma) )} \leq \eta,
\end{equation}

By \eqref{eq_dtR_between_1}, the functions $e^{-R_i/10}$ are uniformly Lipschitz. Hence $(R_i+1)^{-1}$ are equicontinuous and uniformly bounded, with $(\infty+1)^{-1}=0$. By Arzel\`a--Ascoli, a subsequence converges locally uniformly to $(R+1)^{-1}$ for a continuous $R:\td I\to[R^*,\infty]$, locally Lipschitz where finite.
{If $\tau_1\leq\tau_2$ and $R(\tau_1)<\infty$, then
$\partial_\tau R_i\leq1$ gives
$R_i(\tau_2)\leq R_i(\tau_1)+\tau_2-\tau_1$.  Passing to the limit
shows that $R(\tau_2)<\infty$.  Since $(R+1)^{-1}$ is continuous,
there exists $T'\in[-\infty,\sup\td I]$ such that
$R(\tau)=\infty$ for $\tau\leq T'$ and $R(\tau)<\infty$ for
$\tau>T'$.  This proves Assertion~\ref{thm:HMHF_f}.}

Set $\DD_\tau:=\IB^k_{R(\tau)}\times(\IS^{n-k}/\Gamma)$.
Continuity of $R$ makes $\DD$ open relative to
$M_{\cyl}\times\td I$. Local uniform convergence of $(R_i+1)^{-1}$
shows that every compact subset of $\DD$ lies inside
$\IB^k_{R_i(\tau)}\times(\IS^{n-k}/\Gamma)\subset\DD_{i,\tau}$
with a positive spatial margin for all sufficiently large $i$.

The rescaled inverse maps
\[
F_{i,\tau}:=
(\xi_{e^{\tau/2}}\circ\chi_{i,-e^{-\tau}})^{-1}
\]
are defined on $\DD_{i,\tau}$ and satisfy
\[
F_{i,\tau}^*(e^\tau g_{-e^{-\tau}})=g_{\cyl}+u_{i,\tau}.
\]
Equations \eqref{eq:base_pt} and \eqref{eq:u_i_small} therefore give
uniform local bounds for $F_{i,\tau}$ and its inverse, and keep their
images in a fixed neighborhood of $p(-e^{-\tau})$.
Writing $h_i=g_{\cyl}+u_{i,\tau}$ and
$\bar g_\tau=e^\tau g_{-e^{-\tau}}$, the isometry $F_{i,\tau}$ satisfies
\[
\partial_a\partial_b F_{i,\tau}^\alpha
=\Gamma(h_i)^c_{ab}\,\partial_cF_{i,\tau}^\alpha
-\Gamma(\bar g_\tau)^\alpha_{\beta\gamma}(F_{i,\tau})\,
\partial_aF_{i,\tau}^\beta\,\partial_bF_{i,\tau}^\gamma.
\]
By this identity, 
the interior estimates for $u_{i,\tau}$, and the evolution
equation for $F_{i,\tau}$, after passing to a subsequence we may assume
$F_{i}:\DD_{i}\to\MM$ converges locally smoothly
on $\DD$ to $F:\DD\to\MM$, and $F_\tau$ is a
diffeomorphism onto its image for every $\tau\in\td I$.

If $\td I$ contains its left endpoint $\tau_0$, the convergence also holds
up to $\tau_0$, since the initial maps lie in the fixed compact
family used in the compact-interval proof.
Set
$\Omega_{-e^{-\tau}}:=F_\tau(\DD_\tau)$ and
$\chi_{-e^{-\tau}}:=\xi_{e^{-\tau/2}}\circ F_\tau^{-1}$.
Then $p(t)\in\Omega_t$ by \eqref{eq:base_pt}, the maps $\chi_t$
satisfy HMHF, and
$u_\tau:=e^\tau(\xi_{e^{\tau/2}})_*(\chi_{-e^{-\tau}})_*g_{-e^{-\tau}}-g_{\cyl}$
satisfies the metric bounds in Theorem~\ref{thm:HMHF}.
The associated mode quantities converge to $U^+$ and $\UU^-$
defined with $R$, uniformly on compact time intervals.

Near any $\tau\in\td I$ with $R(\tau)<\infty$, we have $R_i\to R$ locally uniformly.
The RHS of \eqref{eq_ODE_for_R_tau} therefore converges
uniformly to the same expression with $u$ and $R$; the
convention $H(F)=1$ at a zero denominator is continuous for finite
$R$. So $R$ is
$C^1$ and satisfies \eqref{eq_ODE_for_R_tau} and
\eqref{eq_dtR_between_1} wherever finite.
Moreover,
the evolution inequalities \eqref{eq_ODI_1}--\eqref{eq_ODI_2} and
Assertions~\ref{thm:HMHF_a}--\ref{thm:HMHF_e} pass to the limit
by locally uniform convergence, after taking subsequential limits of the
partition times. The equation for $u_i$, the $C^m$-bound,
\eqref{eq_dtR_between_1}, and Gaussian decay give uniform local
Lipschitz bounds for the modes, so the limiting differential
inequalities hold.

In the setting of Proposition~\ref{prop:HMHF_initial}, each $I_i$ contains the same initial time. Thus $R(\tau_0)=R_0$ and the initial estimate follows by convergence at $\tau_0$.
\end{proof}

\section{Mode analysis, rigidity, and neck propagation}\label{sec:higher_mode_analysis}

\subsection{Higher mode analysis in the general case}

{The following proposition refines Theorem~\ref{thm:HMHF} by deriving differential inequalities for $\UU^-$ and the higher components of $U^+$. Under the stated conditions on $U^+_{\conf}$, these components become controlled by successively higher powers of $\|U^{++}\|$. In the ancient case, the estimates hold on the whole rescaled interval and yield the approximation recorded in Proposition~\ref{Prop_ancient_summary}.}

\begin{Proposition} \label{Prop_mode_analysis}
Suppose that $m \geq 4, J \geq 2$, $\la\in\IR$, and $0<\beta\leq\ov\beta(\la,J)$, and suppose that $\eta, \eps, R^*, \delta$ satisfy bounds of the form \eqref{eq_HMHF_constants} from Theorem~\ref{thm:HMHF} with this choice of $\beta$.
Suppose in addition that $\la \leq -2$, $\la \in \IZ$.
Set $J_0 := \min \{ J, -\la \} - 1$.
Set $\gamma_{\osc}:=\frac{1}{2(n-k-1)}$, a lower bound for the spectral gap of $L$ on $\sV_{\osc,<0}$ by Lemmas~\ref{Lem_cyl_spectrum} and~\ref{l:quotient_cylinder_spectrum}, and assume $\beta \in (0, \gamma_{\osc}/4)$ and a bound of the form $\eta \leq \ov\eta (J,\la,\beta)$.
In this proposition $\beta$ is fixed once and for all.  All constants and smallness thresholds below may depend on this fixed $\beta$, and this dependence is suppressed from the notation.

Let $\td I \subset \IR$ be a closed interval and consider locally Lipschitz functions $U^+ : \td I \to \sV_{>\la}$ and $\UU^- : \td I \to \IR_{\geq 0}$ that satisfy the conclusions of Theorem~\ref{thm:HMHF} for the times $T_1 \leq T_2 \leq T_3$.
Suppose also that the bound $\| U^+ \|_{L^2_f}, \UU^- \leq \eta (\| 1 \|_{L^2_f} + 1)$ holds, which follows from Theorem~\ref{thm:HMHF}. 
Then there are constants $C=C(J,\la), C_i=C_i(J,\la) >0$ with the following property.

Consider the time-dependent quantities for $i = 1, \ldots, J_0$
\[ \td\UU_i := \| U^+_{\osc, < 0} \| + \| U^+_{\rot,\sph, < -i+1} \| + \| U^+_{\rot,\aff, < -i} \|  +  \UU^-. \]
Choose times $T'_1 = T_1 \leq T'_2 \leq \ldots \leq T'_{J_0} \leq \sup \td I$ successively such that $T'_i \in [T'_{i-1}, \sup \td I]$ is minimal with the property that for $i \geq 2$
\begin{equation} \label{eq_UUim1_choice}
  \td\UU_{i-1}  \leq C_i \| U^{++} \|^{i+1}, \qquad \text{on} \quad ( T'_i, \infty) \cap \td I.  
\end{equation} 
Then the following is true:
\begin{enumerate}[label=(\alph*)]
\item\label{Prop_mode_analysis_a} On $(-\infty,T_1) \cap \td I$ we have $\| U^+_{<0} \| + \UU^-  \geq C_1 \| U^+_{\geq 0} \|^2$ and
\begin{equation} \label{eq_Upl0_UUm_decay}
 \partial_\tau ( \| U^+_{< 0} \| + \UU^- ) \leq - (\tfrac12\gamma_{\osc} - C\sqrt\eta) (\UU^- + \| U^+_{<0}  \|).
\end{equation}
\item\label{Prop_mode_analysis_b} On $[T'_i, \infty ) \cap \td I$, for $i \in \{ 1, \ldots, J_0 \}$, the following ODIs hold
\begin{equation} \label{eq_UUi_ODI}
 \partial_\tau \td\UU_i \leq -(\gamma_{\osc} - C \sqrt\eta) \td\UU_i + C  \| U^{++} \|^{i+2} + C_{i+1} \| U^{++} \| \cdot \| U^+_{\conf} \|. 
\end{equation}
\begin{equation} \label{eq_Uconf_ODI_geq3}
 \| \partial_\tau U^+_{\conf} - \la_{\conf} U^+_{\conf} \|  \leq C \| U^{++} \|  \cdot \| U^+_{\conf} \|  +  {C \| U^{++} \|^{i+1}} .
\end{equation}
Specifically, if $n-k=2$, then on $[T'_i, T_3] \cap \td I$, 
\begin{equation} \label{eq_Uconf_ODI_2}
    \| \partial_\tau U^+_{\conf}  \| \leq C \beta \| U^{++} \| \cdot \| U^+_{\conf} \| + C  \| U^{++}  \|^{i+2} + C \sqrt\eta \UU^- .
\end{equation}
\item\label{Prop_mode_analysis_c} 
If $ T'_i \leq T_3$ (for $n-k > 2$) or $T'_i \leq T_2$ (for $n-k = 2$),\footnote{If $\sup\td I=\infty$, then $T_2=T_3=\infty$, so this condition is automatically satisfied.} where $i \in \{ 1, \ldots, J_0 \}$, recall from Theorem \ref{thm:HMHF}\ref{thm:HMHF_d} that $U^+_{\conf}(T_3)=0$ (for $n-k > 2$) or $U^+_{\conf}(T_2)=0$ (for $n-k = 2$), then the following bounds hold on $[T'_i, \infty) \cap \td I$
\begin{alignat}{2}
\| U^+_{< 0} \| + \UU^- &\leq C_1 \| U^{++} \|^2 \qquad &&\text{if} \quad i = 1, \label{eq_modes_power_1} \\
\| U^+_{\conf} \|&\leq C \| U^{++} \|^{i+1} \qquad&&\text{if} \quad i \geq 1. \label{eq_modes_power_3}
\end{alignat}
 Moreover, in this case the following bounds hold on $(T'_i, T'_{i+1}) \cap \td I$, if $i < J_0$,
\begin{align}
 \td\UU_i &\geq C_{i+1} \| U^{++} \|^{i+2} \label{eq_UUi_uptoTip1_1} \\
 \partial_\tau \td\UU_i &\leq -(\gamma_{\osc} - \beta) \td\UU_i.   \label{eq_dtUUidecay}
\end{align}
\item\label{Prop_mode_analysis_d} If $\inf \td I= -\infty$ and if $U^+, \UU^- \not\equiv 0$, then $T'_1 = \ldots = T'_{J_0} = T_2 = -\infty$  and thus we have the following bounds on $\td I$:
\begin{align*}
 \UU^- + \| U^+_{\osc} \| &\leq C \| U^{++} \|^{J_0} \\
 \| U^+_{\rot,<0} \| &\leq C \| U^{++} \|^2 \\
\| U^+_{\rot,\sph, < -i+2} \| + \| U^+_{\rot,\aff, < -i+1} \|  &\leq C \| U^{++} \|^{i+1}, \qquad i = 2, \ldots, J_0
\end{align*}
\item\label{Prop_mode_analysis_e} If $\sup \td I = \infty$ and if $U^+, \UU^- \not\equiv 0$, then $T_2 = \infty$.
If $T_1 = T'_1 < \infty$, then $T'_2, \ldots, T'_{J_0} < \infty$, so the bounds from Assertion~\ref{Prop_mode_analysis_c} hold for large $\tau$.
This case corresponds to a \emph{non-degenerate neckpinch.}
\item\label{Prop_mode_analysis_f} If $\td I = \IR$, then $U^+, \UU^- \equiv 0$.
\end{enumerate}

\end{Proposition}

\begin{Remark}
We briefly contextualize these bounds. 
The two most interesting cases are those in which the flow is ancient or immortal.

In the ancient case, Assertion~\ref{Prop_mode_analysis_d} essentially captures the full result, as it yields the expected order of decay for various components of $U^+$ for $\tau \to -\infty$.

In the immortal case, the behavior is largely described by Assertions~\ref{Prop_mode_analysis_a}, \ref{Prop_mode_analysis_c}, and \ref{Prop_mode_analysis_e}. 
Here, the bounds \eqref{eq_modes_power_1} and \eqref{eq_modes_power_3} describe the order of decay for different components of $U^+$ as $\tau \to \infty$. 
However, in contrast to the ancient case, these bounds only hold on time-intervals of the form $[T'_i, \infty)$. 
Prior to these times, each quantity $\td\UU_i$ satisfies the reverse inequality and decays exponentially at a positive dimensional rate; compare with \eqref{eq_Upl0_UUm_decay} and \eqref{eq_dtUUidecay}.
If $T_1 = \infty$, then this implies exponential convergence $U^+(\tau), \UU^-(\tau) \to 0$ as $\tau \to \infty$, corresponding to a \emph{degenerate neckpinch.}
This case could be studied further, beyond what is described in this proposition. For example, one could ask whether $U^+$  contains a dominant mode among the stable modes.
Such a mode could be rotationally symmetric or oscillatory.
Alternatively, it may be that $U^+$ decays much faster than $\UU^-$, in which case the overall decay cannot be characterized by a dominant mode.
If $T_1 < \infty$, then Theorem~\ref{thm:HMHF}\ref{thm:HMHF_c} implies that $\| U^{++} \|$ decays only at a polynomial rate. 
Thus $\td\UU_i$ cannot decay exponentially indefinitely while remaining bounded below by $C_i \| U^{++} \|^{i+2}$, which enforces a quantifiable upper bound on the next time $T'_{i+1}$.

The case, where the flow is defined only on a finite time-interval, is slightly more delicate. 
In this case, the behavior of $\td\UU_i$ and $\| U^+_{\conf} \|$ depends on the location of $T'_i$ relative to the times $T_2$ or $T_3$, at which Theorem~\ref{thm:HMHF}\ref{thm:HMHF_d} guarantees $U^+_{\conf} = 0$. 
If $T'_i$ is smaller than this threshold time, then the  discussion from the immortal case applies: $\td\UU_i$ decays exponentially and must eventually fall below $C_i \| U^{++} \|^{i+1}$, which again provides a quantifiable upper bound on the subsequent time $T'_{i+1}$. 
Conversely, if $T'_i$ is larger than $T_2$ or $T_3$, then we cannot guarantee the bound \eqref{eq_modes_power_3} on $U^+_{\conf}$ over $[T'_i, \infty) \cap \td I$, and consequently only the weaker ODIs from Assertion~\ref{Prop_mode_analysis_b} are available.
In this case, the ODIs \eqref{eq_Uconf_ODI_geq3} and \eqref{eq_Uconf_ODI_2} on $U^+_{\conf}$ could be used to show that the bound \eqref{eq_modes_power_3} must hold after some quantifiable time.
However, this time may fall outside the interval $\td I$.
We have omitted a detailed analysis of this case, as it is not important for this paper.
\end{Remark}

\begin{proof}
Fix the constants $J,\la$; from now on we will omit dependencies on these constants.
In the following $C$ will denote a generic constant depending only on the constants indicated in parentheses.

Assertion~\ref{Prop_mode_analysis_a} and \eqref{eq_modes_power_1} are just restatements of Theorem~\ref{thm:HMHF}\ref{thm:HMHF_a} Theorem~\ref{thm:HMHF}\ref{thm:HMHF_b} after choosing $C_1$ large enough. 

Suppose for a moment that we have $U^{++}(\tau) = 0$ for some $\tau \in [T_1, \infty) \cap \td I$.
Then $U^{+} (\tau) = 0$ and the ODI \eqref{eq_ODI_1} together with the bound \eqref{eq_modes_power_1} on $\UU^-$ imply that we must have $U^+, \UU^- \equiv 0$ on $[T_1, \infty) \cap \td I$, which trivially implies Assertions~\ref{Prop_mode_analysis_b}--\ref{Prop_mode_analysis_e}.
So we may assume in the following that $U^{++}(\tau) \neq 0$ for all $\tau \in [T_1, \infty) \cap \td I$.

We will frequently use the following characterization of $\| U^{++} \|$ restricted to the time-interval $(T_1, T_3)$.

\begin{Claim} \label{Cl_Upp_quad_1tau}
There is a dimensional constant $C_0 > 0$ and numbers $X_1, X_2 \in \IR$ such that
\begin{alignat*}{3}
\frac{C_0^{-1}}{\tau - X_1} &\leq & \| U^{++} (\tau) \| &\leq & \frac{C_0}{\tau - X_1} \qquad &\text{for} \quad \tau \in (T_1, T_2) \\
\frac{C_0^{-1}}{X_2 - \tau} &\leq & \| U^{++} (\tau) \| &\leq & \frac{C_0}{X_2 - \tau} \qquad &\text{for} \quad \tau \in (T_2, T_3) .
\end{alignat*}
\end{Claim}

\begin{proof}
We show the first bound; the second one follows similarly.
Note first that $C^{-1} U_{0,\max} \leq \| U^{++} \| \leq C U_{0,\max}$ on $(T_1, T_2)$ for some dimensional constant $C$.
Next, \eqref{eq_cmax_evol_0} implies that
\begin{equation} \label{eq_dtUm10maxsqrt2}
 | \partial_\tau U^{-1}_{0,\max} - \sqrt 2 | \leq 0.1. 
\end{equation}
This implies that we cannot have $T_1 = - \infty$ and $T_2 > -\infty$, because otherwise $U_{0,\max}^{-1}(\tau) \to -\infty$ as $\tau \to -\infty$, in contradiction to its positivity.
So if $(T_1, T_2) \neq \emptyset$, then $T_1 > -\infty$ and we can set $X_1 := T_1 -  U_{0,\max}^{-1} (T_1) $.
Then for $\tau \in (T_1,T_2)$
\begin{multline*}
 \tau - X_1  
 =  \tau-T_1 + U^{-1}_{0,\max}(T_1)  
 \leq (\sqrt 2 - 0.1) (\tau-T_1) + U^{-1}_{0,\max}(T_1) 
 \leq  U^{-1}_{0,\max}(\tau) \\
 \leq  (\sqrt 2 + 0.1) (\tau - T_1 ) + U^{-1}_{0,\max}(T_1)
 \leq  10 (\tau - T_1) + 10 U^{-1}_{0,\max}(T_1) 
 = 10 (\tau - X_1 ),
\end{multline*}
which implies the first bound.
\end{proof}

We will first derive the ODIs from Assertion~\ref{Prop_mode_analysis_b}, while assuming that $C_1, C_2, \ldots, C_{J_0}$ are constants whose values we will determine later.

\begin{Claim} \label{Cl_UUi}
For $i\in\{1,\cdots, J_0\}$,
the ODI \eqref{eq_UUi_ODI} holds on $(T'_i, \infty) \cap \td I$ for some $C = C(C_1, \ldots, C_{i})$.
\end{Claim}

\begin{proof}
For symmetry reasons
\[ \| \PP_{\sV_{\osc}} Q_J^+(U^+) \| 
\leq C \| U^+ \| \cdot \| U^+_{\osc} \| 
\leq C \| U^{++} \| \cdot \| U^+_{\osc, <0} \| +C \| U^{++} \| \cdot \| U^+_{\conf} \| 
. \]
By definition of $\gamma_{\osc}$, all eigenvalues of $L$ on $\sV_{\osc,<0}$ are at most $-\gamma_{\osc}$. Thus, by \eqref{eq_ODI_1},
\begin{align} 
 \partial_\tau \| U^+_{\osc,<0} \| &\leq - \gamma_{\osc}  \| U^+_{\osc,<0} \| + C \| U^{++} \| \cdot \| U^+_{\osc, <0} \| +C \| U^{++} \| \cdot \| U^+_{\conf} \| + C \| U^{++} \|^{J+1} + C \sqrt\eta \UU^- \notag\\
 &\leq - (\gamma_{\osc} - C \eta) \| U^+_{\osc,<0} \| +C \| U^{++} \| \cdot \| U^+_{\conf} \| + C \| U^{++} \|^{J+1} + C \sqrt\eta \UU^-.\label{eq_dtUoscl0}
\end{align}
Next, we bound the evolution of $U^+_{\rot, \sph, < -i+1}$ and $U^+_{\rot, \aff, < -i}$.
Express
\begin{align*}
 U^+_{\rot} &= (U^+_{\rot,\sph,1} + U^+_{\rot,\aff,0}) + (U^+_{\rot,\sph,\frac12} + U^+_{\rot, \aff,-\frac12}) + \ldots + (U^+_{\rot,\sph,\la+\frac32} + U^+_{\rot, \aff,\la+\frac12}) \\
&\qquad\qquad\qquad + (U^+_{\rot,\sph,\la+1}+0) + (U^+_{\rot,\sph,\la+\frac12}+0) \\
&=: U^*_0 + U^*_1 + \ldots + U^*_{-2\la - 1} + U^*_{-2\la } + U^*_{-2\la+1}
\end{align*}
and note that the component functions of each $U^*_{i'}$ are polynomials of degree at most $i'$.
By Lemma~\ref{Lem_Q_precise_formula} we know that each component function of $Q_J(U^+_{\rot})$ is a sum of terms of the form
\[ Q_{m_1, \ldots, m_l; i_1, \ldots, i_l} := \alpha_{m_1, \ldots, m_l; i_1, \ldots, i_l} (\nabla^{m_1} U^*_{i_1} , \ldots, \nabla^{m_l} U^*_{i_l}), \]
where $\alpha_{m_1, \ldots, m_l; i_1, \ldots, i_l}$ is a spatially constant multilinear form, $m_1, \ldots, m_l \in \{ 0,1,2 \}$ with
\[ m_1 + \ldots + m_l \geq 2 \]
and $i_1, \ldots, i_l \in \{ 0, 1, \ldots, -2\la+1 \}$ and $l \geq 2$.
So the degree of $ Q_{m_1, \ldots, m_l; i_1, \ldots, i_l}$ is at most
\[ i_1 + \ldots + i_l - (m_1 + \ldots + m_l) 
\leq  i_1 + \ldots + i_l-2. \]
Due to \eqref{eq_UUim1_choice} we can bound
\begin{multline*}
 \| \nabla^{m_j} U^*_{i_j} \| \leq C\|  U^*_{i_j} \| 
\leq C \| U^+_{\rot, \sph, 1- i_j/2} \| + C \| U^+_{\rot,\aff,-i_j/2} \| \\
\leq C \| U^+_{\rot, \sph, < \lfloor 1 - i_j/2 \rfloor + 1} \| + C \| U^+_{\rot, \aff, < \lfloor  - i_j/2 \rfloor + 1} \|
\leq C \td\UU_{\lceil i_j/2 \rceil -1}
\leq C(C_1, \ldots, C_i)
\| U^{++} \|^{b_j }, 
\end{multline*}
where 
\[ b_j = \begin{cases} 1 &\text{if $i_j \in \{ 0,1,2 \}$} \\
 \lceil i_j/2 \rceil + 1 &\text{if $2 < i_j \leq 2i$} \\
i+1 &\text{if $2i < i_j$}
\end{cases} \]
and we obtain
\[ \| Q_{m_1, \ldots, m_l; i_1, \ldots, i_l} \| \leq C(C_1, \ldots, C_i) \| U^{++} \|^{b_1 + \ldots + b_l}. \]
Now $\PP_{\sV_{\rot, \sph, <-i+1}} Q^+_J(U^+_{\rot})$ and $\PP_{\sV_{\rot, \aff, <-i}} Q^+_J(U^+_{\rot})$ are each sums of the  corresponding projections of $Q_{m_1, \ldots, m_l; i_1, \ldots, i_l}$.
These vanish if the degree of all component functions of $Q_{m_1, \ldots, m_l; i_1, \ldots, i_l}$ is $\leq 2i$, so if they don't vanish, then
\begin{equation} \label{eq_2ii1ilm2}
 2i <  i_1 + \ldots + i_l-2 . 
\end{equation}
We claim that for any such non-vanishing term we have
\begin{equation} \label{eq_b1blip2}
 b_1 + \ldots + b_l \geq i+2. 
\end{equation}
Indeed, if we have $2i < i_{j}$ for some $j \in \{ 1, \ldots, l \}$, then the bound follows trivially since then $b_{j} =i+1$ and since $l \geq 2$.
So suppose that $i_j \leq 2i$ for all $j \in \{ 1, \ldots, l \}$.
In this case we have $b_j \geq \lceil i_j/2 \rceil$. Since $i_1+\ldots+i_l$ is an integer, \eqref{eq_2ii1ilm2} implies $i_1+\ldots+i_l \geq 2i+3$. Hence
\begin{equation} \label{eq_b1blip1}
 b_1 + \ldots + b_l \geq \Big\lceil \frac{i_1}2 \Big\rceil + \ldots + \Big\lceil \frac{i_l}2 \Big\rceil \geq \Big\lceil \frac{i_1 + \ldots + i_l}2 \Big\rceil \geq \Big\lceil \frac{2i+3}2 \Big\rceil = i+2. 
\end{equation}
This confirms \eqref{eq_b1blip2}.

It follows that
\[ \| \PP_{\sV_{\rot, \sph, <-i+1}} Q^+_J(U^+_{\rot}) \| + \|\PP_{\sV_{\rot, \aff, <-i}} Q^+_J(U^+_{\rot}) \| \leq C(C_1, \ldots, C_i) \| U^{++} \|^{i+2}, \] 
which implies that
\begin{multline*}
 \| \PP_{\sV_{\rot, \sph, <-i+1}} Q^+_J(U^+) \| + \|\PP_{\sV_{\rot, \aff, <-i}} Q^+_J(U^+) \| \\ 
 \leq C(C_1, \ldots, C_i) \| U^{++} \|^{i+2} + C \| U^{++} \| \cdot \| U^{+}_{\osc,<0} \|+ C \| U^{++} \| \cdot \| U^{+}_{\conf} \|. 
\end{multline*}
Plugging this into the ODI system \eqref{eq_ODI_1} \eqref{eq_ODI_2}, and projecting the first ODI onto $\sV_{\rot, \sph, <-i+1}$ and $\sV_{\rot, \aff, <-i}$ implies
\begin{align*}
 \partial_\tau  \| U^+_{\rot, \sph, <-i+1} \|
&\leq -(i-\tfrac12) \| U^+_{\rot, \sph, <-i+1} \| + C(C_1, \ldots, C_i) \| U^{++} \|^{i+2}  \\
&\qquad + C \| U^{++} \| \cdot \| U^{+}_{\osc,<0} \|+ C \| U^{++} \| \cdot \| U^{+}_{\conf} \| + C \| U^{++} \|^{J+1} + C \sqrt\eta \UU^- \\
&\leq -\tfrac12 \| U^+_{\rot, \sph, <-i+1} \| + C(C_1, \ldots, C_i) \| U^{++} \|^{i+2}  \\
&\qquad + C \eta  \| U^{+}_{\osc,<0} \|+ C \| U^{++} \| \cdot \| U^{+}_{\conf} \|  + C \sqrt\eta \UU^- .
\end{align*}
Similarly, 
\begin{align*}
 \partial_\tau  \| U^+_{\rot, \aff, <-i} \|
&\leq -\tfrac12 \| U^+_{\rot, \aff, <-i} \| + C(C_1, \ldots, C_i) \| U^{++} \|^{i+2} \\
&\qquad + C \eta \cdot \| U^{+}_{\osc,<0} \|+ C \| U^{++} \| \cdot \| U^{+}_{\conf} \|  + C \sqrt\eta \UU^-, \\
 \partial_\tau  \UU^-
&\leq -(\tfrac12-C\sqrt\eta) \UU^- + C(C_1, \ldots, C_i) \| U^{++} \|^{i+2}  \\
&\qquad + C \eta  \| U^{+}_{\osc,<0} \|+ C \| U^{++} \| \cdot \| U^{+}_{\conf} \|  .
\end{align*}
Note that $\gamma_{\osc}\le \frac12$,
the claim follows by adding these three inequalities and \eqref{eq_dtUoscl0}.
\end{proof}

\begin{Claim} \label{Cl_Uconf_geq3}
The ODI \eqref{eq_Uconf_ODI_geq3} holds on $[T'_i, \infty) \cap \td I$ for some $C = C(C_1, \ldots, C_{i})$.
\end{Claim}

\begin{proof}
Using the fact that $U^+ \mapsto \PP_{\sV_{\conf}} Q_J^+(U^+)$ has no linear term and vanishes for $U^+ \in \sV_{\rot}$, we obtain
\[ \| \PP_{\sV_{\conf}} Q_J^+(U^+) \| \leq C \| U^+ \|  ( \| U^+_{\conf} \| + \| U^+_{\osc, <0} \|). \]
For $i=1$, \eqref{eq_modes_power_1} bounds $\|U^+_{\osc,<0}\|$ by $C\|U^{++}\|^2$; for $i\geq 2$, the choice \eqref{eq_UUim1_choice} gives $\td\UU_{i-1}\leq C_i\|U^{++}\|^{i+1}$ on $[T'_i,\infty)\cap\td I$. Hence
\[ \| \PP_{\sV_{\conf}} Q_J^+(U^+) \| \leq C \| U^{++} \|  \cdot \| U^+_{\conf} \|  +  C \| U^{++} \|^{i+2}. \]
Projecting the ODI \eqref{eq_ODI_1} onto $\sV_{\conf}$ also gives the error term $C\sqrt\eta\,\UU^-$. The same bounds on $\UU^-$ absorb this error into $C\sqrt\eta\|U^{++}\|^{i+1}$, which proves \eqref{eq_Uconf_ODI_geq3}.
\end{proof}

\begin{Claim} \label{Cl_Uconf_2}
If $n-k=2$, then the ODI \eqref{eq_Uconf_ODI_2} holds on $ [T'_i, T_3] \cap \td I$ for some $C = C(C_1, \ldots, C_{i})$.
\end{Claim}

\begin{proof}
In the case $n-k=2$ recall that $\sV_{0}=\sV_{\sph,0}\oplus\sV_{\conf}\oplus\sV_{\aff}$.
View $Q_2(U', U'') := \frac12 (Q_2 (U'+U'') - Q_2(U')- Q_2(U''))$ as a symmetric bilinear form via polarization.
By \eqref{eq_PV0Q2} and Lemma~\ref{l:Q_2}
\begin{equation*}
\PP_{\sV_{\conf}}Q_2(U^+_{\rot})=\PP_{\sV_{\conf}}Q_2(U^+_{\aff})=\PP_{\sV_{\conf}}Q_2(U^+_{\conf})=\PP_{\sV_{\conf}}Q_2(U^+_{\conf},U^+_{\aff}+U^+_{\sph,0})=0. 
\end{equation*}
So by writing $U^+=U^+_{\rot}+U^+_{\conf}+U^{+}_{\osc,<0}$, using \eqref{eq_modes_power_1} and $\| U^{++}_{>0} \| \leq \beta \|U^{++} \|$ from \eqref{eq:neutral_dom_rough}, we obtain that
\begin{align*}
     \| \PP_{\sV_{\conf}} &Q^+_J(U^+) \| 
\leq  \| \PP_{\sV_{\conf}} Q_2^+(U^+) \| +  \| \PP_{\sV_{\conf}} (Q^+_J - Q_2^+)(U^+) \| \\
&\leq \| \PP_{\sV_{\conf}} Q_2^+(U^{+}_{\rot} + U^+_{\conf})  \| + C \| U^+ \| \cdot \| U^+_{\osc, <0} \| + C \| U^+ \|^2  ( \| U^+_{\conf} \| + \| U^+_{\osc, <0} \|) \\
&\leq 2 \| \PP_{\sV_{\conf}} Q_2^+(U^+_{\conf} , U^{+}_{\rot} )  \|  + C \| U^{++} \|^2 \| U^+_{\conf} \| + C  \| U^{++}  \|^{i+2} \\
&{\leq 2 \| \PP_{\sV_{\conf}} Q_2^+(U^+_{\conf} , U^{+}_{\rot,  0}  )  \|}  + C \beta \| U^+ \| \cdot \| U^+_{\conf} \|  + C  \| U^{++}  \|^{i+2} \\
&\leq 2 \| \PP_{\sV_{\conf}} Q_2^+(U^+_{\conf} , U^{+}_{\rot, \sph, 0} + U^+_{\rot,\aff,0} )  \|  + C \beta \| U^+ \| \cdot \| U^+_{\conf} \| + C  \| U^{++}  \|^{i+2} \\
&\leq  C \beta \| U^+ \| \cdot \| U^+_{\conf} \| + C  \| U^{++}  \|^{i+2} .
\end{align*}
Projecting \eqref{eq_ODI_1} onto $\sV_{\conf}$ gives \eqref{eq_Uconf_ODI_2}.
\end{proof}

Claims~\ref{Cl_UUi}, \ref{Cl_Uconf_geq3} and \ref{Cl_Uconf_2} show that Assertion~\ref{Prop_mode_analysis_b} holds for some $C = C(C_1, \ldots, C_{i})$.
The next key step is to derive the bound \eqref{eq_modes_power_3} on $\| U^+_{\conf} \|$ in the case $ T'_i \leq T_3$ (for $n-k > 2$) or $T'_i \leq T_2$ (for $n-k = 2$).

\begin{Claim}\label{Claim:small_conf_n-k>2}
Fix $i \in \{ 1, \ldots, J_0 \}$ and $\eta \leq \ov\eta (C_1, \ldots, C_i)$.
Then for some $C = C(C_1, \ldots, C_{i})$ we have
\begin{alignat}{2}
 \sup_{[T'_i, \infty)  \cap \td I} \frac{\| U^+_{\conf} \|}{\| U^{++} \|^{i +1}} &\leq  C + C\frac{\| U^+_{\conf}(T'_i) \|}{\| U^{++} (T'_i)\|^{i +1}}   &&\qquad \text{if} \quad n-k > 2, \label{eq_sup_Uconf_1} \\
 \sup_{[\max \{ T'_i, T_3 \}, \infty)   \cap \td I} \frac{\| U^+_{\conf} \|}{\| U^{++} \|^{i +1}} &\leq  C  + C\frac{\| U^+_{\conf}(\max \{ T'_i, T_3 \}) \|}{\| U^{++} (\max \{ T'_i, T_3 \})\|^{i +1}} &&\qquad \text{if} \quad n-k = 2  . \label{eq_sup_Uconf_2}
\end{alignat} 
If $T'_i \leq T_3$, then the last term in \eqref{eq_sup_Uconf_1} can be omitted.
This proves \eqref{eq_modes_power_3} in the case $n-k>2$.
Moreover, if $T'_i = T_3 = -\infty$, then the last term in \eqref{eq_sup_Uconf_2} can also be omitted.
\end{Claim}

\begin{proof}
Since $\la_{\conf} < \frac12$ (see Lemma~\ref{Lem_cyl_spectrum}), the bound \eqref{eq_Uconf_ODI_geq3} from Claim~\ref{Cl_UUi} implies that if $\eta \leq \ov\eta(C_1, \ldots, C_i)$, then for some constants $C^*(C_1, \ldots, C_i) > 0$, $c^* (C_1, \ldots, C_i) \in (0,0.1)$, which we will henceforth fix, we have on {$[T'_i, \infty) \cap \td I$}
\begin{equation} \label{eq_U_conf_Cstar}
 \partial_\tau \| U^+_{\conf} \| \leq (\tfrac12 - 3c^*) \| U^+_{\conf} \| + C^* \| U^{++} \|^{i+1}. 
\end{equation}
By Theorem~\ref{thm:HMHF}\ref{thm:HMHF_e} the following bound holds on {$[T_3, \infty) \cap \td I$} assuming $\eta \leq \ov\eta (c^* (C_1, \ldots, C_i))$
\begin{equation} \label{eq_Upp_expbound}
 \partial_\tau \| U^{++}_{> 0} \| \geq (\tfrac12 - c^*) \| U^{++}_{>0} \| 
\end{equation}
Consider the function
\[ F(\tau) := \frac{\| U^+_{\conf}(\tau) \|}{\| U^{++}_{>0} \|^{i+1}}. \]
By combining \eqref{eq_U_conf_Cstar} with \eqref{eq_Upp_expbound}, we find that on $[T'', \infty) \cap \td I$, for $T'' := \max \{T'_i, T_3 \}$,
\begin{equation} \label{eq_ODI_FcsCs}
 \partial_\tau F \leq - (\tfrac12 -c^*) F + C^*. 
\end{equation}
So if $T'' := \max \{T'_i, T_3 \} > -\infty$, then
\begin{equation} \label{eq_supFTpp}
 \sup_{[T'',\infty) \cap \td I} F(\tau) \leq \frac{C^*}{\tfrac12 -c^*}  + F(T''). 
\end{equation}
This proves \eqref{eq_sup_Uconf_2}, and also \eqref{eq_sup_Uconf_1} when $T_i'\ge T_3$.

Suppose now that $T'' = -\infty$.
The bound \eqref{eq_ODI_FcsCs} implies that the bound $F \leq \frac{C^*}{\frac12 - c^*}$ is preserved forward in time.
We claim that this bound holds on all of $\td I$.
Otherwise, we would have $F > \frac{C^*}{\frac12 - c^*}$ for all sufficiently small $\tau$ and $F(\tau) \to \infty$ as $\tau \to -\infty$.
This would imply by \eqref{eq_U_conf_Cstar} that for sufficiently small $\tau$
\[ \partial_\tau \| U^+_{\conf} \| \leq (\tfrac12 - 2c^*) \| U^+_{\conf} \| , \]
and therefore, using Theorem~\ref{thm:HMHF}\ref{thm:HMHF_b}, for some $c', c'' > 0$, 
\[ \| U^{++}_{>0} (\tau) \| \geq c' \| U^+_{\conf} (\tau) \| \geq c'' e^{(\frac12 - 2c^*)\tau} .\]
This, however, contradicts the asymptotic bound $ \| U^{++}_{>0} (\tau) \| \leq C' e^{(\frac12 - c^*)\tau}$ following from \eqref{eq_Upp_expbound}.

It remains to prove \eqref{eq_modes_power_3} in the case $n-k>2$ and $T'_i < T_3$.
First, since $T'' = T_3$, by Theorem~\ref{thm:HMHF}\ref{thm:HMHF_d} we have $F(T'') =0$, and thus it follows by \eqref{eq_supFTpp} that $F(\tau)\le \frac{C^*}{\tfrac12 -c^*}$ on $[T_3, \infty)  \cap \td I$.
So we just need to estimate $F(\tau)$ on $[T_1, T_3] \cap [T'_i,\infty) \cap \td I = [T'_i,T_3] \cap \td I$.
Recall that in this case $\la_{\conf} > 0$.
As before, we obtain from \eqref{eq_Uconf_ODI_geq3} that, after possibly adjusting $C^*$ and $c^*$,
\[ \partial_\tau \| U^+_{\conf} \| \geq 3c^* \| U^+_{\conf} \| - C^* \| U^{++} \|^{i+1}. \]
If $\eta \leq \ov\eta(c^*)$, then by Theorem~\ref{thm:HMHF}\ref{thm:HMHF_c}\eqref{eq_cmax_evol_0}\eqref{eq_cmin_evol_0} we have on $[T_1,T_3]$ 
\[ \partial_\tau \| U^{++}_0 \| \leq \frac{c^*}{i+1} \| U^{++}_0 \|, \]
so the function
\[ F(\tau) := \frac{\| U^+_{\conf}(\tau) \|}{\| U^{++}_{0} (\tau)\|^{i+1}}, \]
satisfies
\[ \partial_\tau F \geq 2c^* F - C^*. \]
If $T_3 < \infty$, then we use again the fact that $F(T_3) = 0$ to conclude
\[ \sup_{\tau \in [T'_i, T_3]} F(\tau) \leq \frac{C^*}{2c^*}  . \]
If $T_3 = \infty$ and if this bound is violated, seeing that the condition $F(\tau) > \frac{C^*}{2c^*}$ is preserved forward in time, we must have $F(\tau) > \frac{C^*}{2c^*}$ for \emph{all} large enough $\tau$, and $F(\tau) \to \infty$ as $\tau \to \infty$, implying that $\partial_\tau \| U^+_{\conf} \| \geq 2c^* \| U^+_{\conf} \|$ for large $\tau$, in contradiction to the fact that this quantity is uniformly bounded.
This finishes the proof of the claim.
\end{proof}

The next claim establishes \eqref{eq_modes_power_3} in the case $n-k=2$.

\begin{Claim}\label{Claim:small_conf_n-k=2}
Suppose that $n-k=2$ and $T'_i \leq T_2$ and $\eta \leq \ov\eta (C_1, \ldots, C_i)$.
Then \eqref{eq_modes_power_3} holds on $[T'_i, \infty) \cap \td I$ for some $C=C(C_1, \ldots, C_i)$.
\end{Claim}

\begin{proof}
In the following, $C$ denotes a generic constant that may depend on $C_1, \ldots, C_i$.
Due to \eqref{eq_sup_Uconf_2} it suffices to derive \eqref{eq_modes_power_3} over the time-interval $(T'_i, T_3]$.

If $T_2$ is finite, then integrating \eqref{eq_Uconf_ODI_2} (see Claim~\ref{Cl_Uconf_2}) from $T_2$ to any $\tau \in (T'_i, T_3]$ and using the fact that $U_{\conf}(T_2) = 0$ gives
\[ \| U^+_{\conf}(\tau) \| \leq  C\bigg| \beta  \int_{T_2}^\tau \| U^{++}(\tau') \| \cdot \|U^+_{\conf}(\tau')\| d\tau' +  \int_{T_2}^\tau \| U^{++}(\tau') \|^{i+2} d\tau' +  \sqrt\eta \int_{T_2}^\tau \td\UU_i(\tau') d\tau' \bigg|. \]
The same bound also holds if $T_2 = \pm\infty$. Indeed, if $T_2=+\infty$, then by Lemma \ref{l:partition_2}, $\|U^+_{\conf}\|\leq C U_{0,\max}\to0$ as $\tau\to+\infty$. So $\|U^+_{\conf}\|\to0$ as $\tau\to+\infty$. If $T_2=-\infty$, then $\|U^+_{\conf}(\tau)\|\leq C\|U^+(\tau)\|\to0$ as $\tau\to-\infty$. Thus, $\|U^+_{\conf}(\tau)\|\to0$ as $\tau\to T_2$ in both cases.
Using the bound from Claim~\ref{Cl_Upp_quad_1tau}, we find
\begin{equation} \label{eq_intUppip2Uppip1}
 \bigg| \int_{T_2}^\tau \| U^{++}(\tau') \|^{i+2} d\tau' \bigg| \leq C \| U^{++} (\tau) \|^{i+1}, 
\end{equation}
so we obtain
\begin{equation} \label{eq_Uconfboundedint}
 \| U^+_{\conf}(\tau) \| \leq  C \beta \bigg|  \int_{T_2}^\tau \| U^{++}(\tau') \| \cdot \|U^+_{\conf}(\tau')\| d\tau' \bigg| +  C\| U^{++} (\tau) \|^{i+1}+  C\sqrt\eta \bigg| \int_{T_2}^\tau \td\UU_i(\tau') d\tau' \bigg|. 
\end{equation}
To bound the last integral, we integrate \eqref{eq_UUi_ODI} (see Claim~\ref{Cl_UUi}), assuming $\eta \leq \ov\eta (C_1, \ldots, C_i)$.
If $T_2$ is finite, then we obtain that
\[\frac12\gamma_{\osc} \bigg| \int_{T_2}^\tau \td\UU_i (\tau') d\tau' \bigg| 
\leq \big| \td\UU_i(\tau) - \td\UU_i(T_2) \big| + C \bigg| \int_{T_2}^\tau \| U^{++} (\tau')\|^{i+2}d\tau' \bigg| + C  \bigg| \int_{T_2}^\tau \| U^{++} (\tau') \| \cdot \| U^+_{\conf}(\tau') \| d\tau' \bigg|. \]
Using \eqref{eq_UUim1_choice} (note that $\td\UU_i \leq \td\UU_{i-1}$) on the first term of the RHS and again \eqref{eq_intUppip2Uppip1} on the second term, we obtain
\[  \bigg| \int_{T_2}^\tau \td\UU_i (\tau') d\tau' \bigg| 
\leq  C  \| U^{++} (\tau)\|^{i+1} + C  \bigg| \int_{T_2}^\tau \| U^{++} (\tau') \| \cdot \| U^+_{\conf}(\tau') \| d\tau' \bigg|. \]
If $T_2 = \pm\infty$, then the same bound follows, because in this case, by the same argument as before we have $\td\UU_i(\tau) \to 0$ as $\tau \to T_2$.
Combining this with \eqref{eq_Uconfboundedint} implies
\begin{equation} \label{eq_Uconf_Cint_Uppip1}
 \| U^+_{\conf}(\tau) \| \leq C \beta \bigg|  \int_{T_2}^\tau \| U^{++}(\tau') \| \cdot \|U^+_{\conf}(\tau')\| d\tau' \bigg| +  C\| U^{++} (\tau) \|^{i+1}. 
\end{equation}

Next, we estimate
\[ H(\tau) := \bigg|  \int_{T_2}^\tau \| U^{++}(\tau') \| \cdot \|U^+_{\conf}(\tau')\| d\tau' \bigg|. \]
First, using Claim~\ref{Cl_Upp_quad_1tau} and assuming $\eta \leq \ov\eta(C_1,\ldots,C_i)$, we obtain in the case $\tau \leq T_2$ that
\[ |\partial_\tau H| \leq C \beta \| U^{++} \| \cdot H + C \| U^{++} \|^{i+2} \leq \tfrac12 (\tau-X_1)^{-1} H + C (\tau-X_1)^{-i-2}. \]
This implies 
\[
 \partial_\tau\big((\tau-X_1)^{1/2}H(\tau)\big)
 \geq -C(\tau-X_1)^{-i-3/2}.
\]
{For finite $T_2$, since $H(T_2)=0$,} integrating from $\tau$ to $T_2$ yields
\[
 (\tau-X_1)^{1/2}H(\tau)
 \leq C\int_\tau^{T_2}(\tau'-X_1)^{-i-3/2}\,d\tau'
 \leq C(\tau-X_1)^{-i-1/2},
\]
and hence again
\[
 H(\tau)\leq C(\tau-X_1)^{-i-1}\leq C\|U^{++}(\tau)\|^{i+1}.
\]
If $T_2=+\infty$, Theorem~\ref{thm:HMHF}\ref{thm:HMHF_b} and Claim~\ref{Cl_Upp_quad_1tau} give
\[
 H(\tau)\leq C\int_\tau^\infty(\tau'-X_1)^{-2}\,d\tau'
 =C(\tau-X_1)^{-1}.
\]
Thus $(\tau-X_1)^{1/2}H(\tau)\to0$ as $\tau\to+\infty$, and integrating the preceding differential inequality to infinity gives the same bound for $H(\tau)$.
If $\tau \geq T_2$, then a similar argument yields the same bound.
Plugging this bound into \eqref{eq_Uconf_Cint_Uppip1} implies \eqref{eq_modes_power_3} over the time-interval $(T'_i, T_3]$ and hence the claim.
\end{proof}

The next claim establishes the remaining bounds in Assertion~\ref{Prop_mode_analysis_c}.

\begin{Claim}
If $C_{i+1} \geq \underline{C}_{i+1}(C_1,\ldots, C_i)$ and $\eta \leq \ov\eta(C_1, \ldots, C_i)$ and $i < J_0$, then \eqref{eq_UUi_uptoTip1_1} and \eqref{eq_dtUUidecay} hold on $(T'_i,T'_{i+1})$.
\end{Claim}

\begin{proof}
Combining \eqref{eq_UUi_ODI} with \eqref{eq_modes_power_3} (see Claim \ref{Claim:small_conf_n-k>2} and \ref{Claim:small_conf_n-k=2}) implies that on $[T'_i, \infty) \cap \td I$
\[ \partial_\tau \td\UU_i \leq -(\gamma_{\osc} - C(C_1, \ldots, C_i)  \sqrt\eta) \td\UU_i + C(C_1, \ldots, C_i)  \| U^{++} \|^{i+2} . \]
So if $C_{i+1} \geq \underline C_{i+1}(C_1, \ldots, C_i)$ is chosen large enough and $\eta \leq \ov\eta(C_1, \ldots, C_i)$, then on $(T'_i, \infty) \cap \td I$
\[ \td\UU_i \geq C_{i+1} \| U^{++} \|^{i+2} \qquad \Rightarrow \qquad 
\partial_\tau \td\UU_i \leq -(\gamma_{\osc} - \beta) \td\UU_i\le  - \tfrac12 \gamma_{\osc}\,\td\UU_i.\]
On the other hand, we obtain from Theorem~\ref{thm:HMHF}\ref{thm:HMHF_b}\eqref{eq:U++_evolve}, assuming $\eta \leq \ov\eta(\gamma_{\osc})$,
\[ {\partial_\tau \| U^{++} \|^{i+2} \geq - \tfrac12 \gamma_{\osc} \| U^{++} \|^{i+2}.} \]
Therefore, the following is true on $(T'_i, \infty) \cap \td I$
\[ \td\UU_i \geq C_{i+1} \| U^{++} \|^{i+2} \qquad \Rightarrow \qquad 
\partial_\tau \td\UU_i <  \partial_\tau \big( C_{i+1} \| U^{++} \|^{i+2} \big). \] 
It thus follows that the set of times within $(T'_i, \infty) \cap \td I$ at which $\td\UU_i \geq C_{i+1} \| U^{++} \|^{i+2}$ is either empty or a single interval of the form $(T'_i, T'')$.
Recall that $T'_{i+1} \in [T'_{i}, \sup \td I]$ is minimal such that $ \td\UU_{i}  \leq C_{i+1} \| U^{++} \|^{i+2}$ on $(T'_{i+1}, \infty) \cap \td I$.
Hence, $T'' = T'_{i+1}$ and the claim follows.
\end{proof}

The previous claims show that we can successively choose the constants $C_1, \ldots, C_{J_0}$ so that Assertions~\ref{Prop_mode_analysis_a}--\ref{Prop_mode_analysis_c} hold.
Assertion~\ref{Prop_mode_analysis_d} follows from Assertions~\ref{Prop_mode_analysis_a}--\ref{Prop_mode_analysis_c}: because the quantities $\UU^- + \| U^+_{\geq 0} \|$ and $\td\UU_i$ are uniformly bounded, they cannot decay exponentially on a time-interval of the form $(-\infty,T'_i)$.
The fact that $T_2 =-\infty$ follows similarly from \eqref{eq_cmax_evol_0}.
The first statement of Assertion~\ref{Prop_mode_analysis_e} follows from \eqref{eq_cmin_evol_0}.
If $T_1 < \infty$, then Claim~\ref{Cl_Upp_quad_1tau} implies that $\| U^{++} \|$ decays only at a polynomial rate, so the quantities $\td\UU_i$ cannot decay exponentially indefinitely while staying above $C_{i+1} \| U^{++} \|^{i+2}$ (see \eqref{eq_UUi_uptoTip1_1} and \eqref{eq_dtUUidecay}).
This implies $T_2',\cdots,T_{J_0}'<\infty$.
Assertion~\ref{Prop_mode_analysis_f} follows from the fact that $T_2=-\infty$ in the ancient case and $T_2=\infty$ in the immortal case, unless $U^+, \UU^- \equiv 0$.
This finishes the proof of the proposition.
\end{proof}
\medskip

\subsection{Eternal flows and rigidity of cylinders}  \label{subsec_eternal}
As a direct corollary of Proposition~\ref{Prop_mode_analysis}, we obtain
an alternative proof of~\cite[Theorem~1.1]{Colding_Minicozzi_RF_2025}.
When $k=1$, this rigidity was proved independently by
Li--Wang~\cite[Theorem~1.1]{Li_Wang_cylinder_rigidity}.
\begin{Corollary} \label{Cor_eternal}
There is a constant $\delta=\delta(n,k,\Gamma)>0$ with the following property.
Let $(M,g,p)$ be a complete connected $n$-dimensional
normalized shrinking Ricci soliton with potential $f$, where
$\nabla f(p)=0$, and assume that $p$ is a center of an
$(n,k,\delta,\Gamma)$-neck at scale~$1$.
Then $(M,g)$ is isometric to $(M_{\cyl},g_{\cyl})$.

More generally, the following is true:
Let $\MM$ be either a connected $n$-dimensional Ricci flow with complete
time-slices or the regular part $\RR$ of a noncollapsed $\IF$-limit
$\XX$ of $n$-dimensional Ricci flows.
Assume $\MM$ is defined on $\IR_-$ and let $p:(-\infty,0)\to\MM$
be a worldline.
Suppose that $p(t)$ is a center of an $(n,k,\delta,\Gamma)$-neck at scale
$\sqrt{-t}$ for every $t<0$.
Then $\MM$ is isometric to the cylindrical flow $\MM_{\cyl}$ and,
in the limit case, $\XX_{<0}$ itself is the cylindrical metric flow.
\end{Corollary}

\begin{proof}
We first prove the more general assertion.  Fix
\[
m=J=10,
\qquad
\la=-10,
\qquad
\gamma_{\osc}:=\frac{1}{2(n-k-1)}.
\]
Let $\ov\beta(\la,J)$ be the minimal threshold from Theorem \ref{thm:HMHF} and \ref{Prop_mode_analysis}. 
Choose and fix 
\[
0<\beta\le 
\min\left\{\tfrac12\ov\beta(\la,J),
\tfrac18\gamma_{\osc}\right\}.
\]
We next choose $\eta>0$ below both the threshold
$\ov\eta(\la,J,\beta)$ in Theorem~\ref{thm:HMHF} and
Proposition~\ref{Prop_mode_analysis}.  Choose successively
\[
\eps\leq\ov\eps(m,J,\la,\beta,\eta),
\qquad
R^*\geq\underline R^*(m,J,\la,\eta,\eps),
\]
and finally choose
\[
0<\delta\leq\ov\delta(m,J,\la,\beta,\eta,\eps,R^*)
\]
sufficiently small that Theorem~\ref{thm:HMHF} and
Proposition~\ref{Prop_mode_analysis} apply.
Thus the parameters can be fixed in the order
\[
(m,J,\la),\quad \beta,\quad \eta,\quad \eps,\quad
R^*,\quad \delta.
\]
Since all parameters preceding $\delta$ have now been fixed in terms of
$n,k,\Gamma$, the resulting $\delta$ depends only on these parameters.

Apply Theorem~\ref{thm:HMHF} to the worldline $p(t)$ on
$I=(-\infty,0)$.  This produces a harmonic map heat flow $\chi$ to the cylinder and functions
\[
U^+:\IR\longrightarrow\sV_{>\la},
\qquad
\UU^-:\IR\longrightarrow\IR_{\geq0},
\]
satisfying all the conclusions of that theorem.  Our choices of
$\la,\beta$, and $\eta$ also give all the additional hypotheses of
Proposition~\ref{Prop_mode_analysis}.  Since the rescaled time interval
is the entire real line, Proposition~\ref{Prop_mode_analysis}
\ref{Prop_mode_analysis_f} yields
\[
U^+\equiv0,
\qquad
\UU^-\equiv0.
\]
By the definition of $\UU^-$ in Theorem~\ref{thm:HMHF}, we have $R(\tau)=\infty$.  Hence
$\omega_{R(\tau)}\equiv1$, and $u_\tau\equiv0$ for every
$\tau\in\IR$.  
Thus $(\xi_{(-t)^{-1/2}}\circ\chi_t)_*\big((-t)^{-1}g_t\big)\equiv g_{\cyl}$ for all $t<0$, or equivalently, $(\chi_t)_*g_t=g_{\cyl,t}$.
The harmonic map heat flow equation gives
\[
\partial_t\chi_t=\Delta_{g_t,g_{\cyl,t}}\chi_t=(\Delta_{(\chi_t)_*g_t,g_{\cyl,t}}\id)\circ\chi_t=0.
\]
Thus $\chi_t$ preserves worldlines. Since $\chi_t^{-1}$ is an isometric
embedding of the complete cylinder, its image is a connected component
of $\MM_t$. In the smooth case, connectedness gives $\Omega_t=\MM_t$.
In the limit case, \cite[Theorem~2.4(c)]{Bamler_HK_RF_partial_regularity}
implies that $\RR_t$ is connected and $(\XX_t,d_t)$ is the metric
completion of $(\RR_t,d_{g_t})$. Hence $\Omega_t=\RR_t=\XX_t$.
The heat-flow characterization in
\cite[Theorem~2.4(d)]{Bamler_HK_RF_partial_regularity} then identifies
$\XX_{<0}$ with the cylindrical metric flow.

Finally, consider the shrinker in the first assertion. Since $p$ is a
critical point of its potential $f$, the diffeomorphisms defining its
canonical self-similar flow $g_t$, with $g_{-1}=g$, fix $p$. Hence
\[
(M,(-t)^{-1}g_t,p)\cong(M,g,p)\qquad(t<0).
\]
By the two-sided pseudolocality theorem of Li--Wang~\cite[Theorem~1.6]{MR4791450}, this flow satisfies the pseudolocality property of Definition~\ref{def:pseudolocality_property}.
Thus the worldline $p(t)=(p,t)$ satisfies the neck hypothesis of the
general assertion for every $t<0$, which proves the first assertion.
\end{proof}
\medskip

\subsection{Neck propagation and uniqueness of tangent flows} \label{subsec_stability_neck}
Theorem~\ref{Thm_propagation_necks} gives a neck along the worldline through the later center at every intermediate time, assuming sufficiently accurate endpoint necks and the $H_n$-center condition.

\begin{Theorem}[Propagation of necks]\label{Thm_propagation_necks}
For every $\eps>0$ there is a $\delta=\delta(n,k,\Gamma,\eps)>0$ with the following property.

Let $\XX$ be a metric flow arising as a noncollapsed $\IF$-limit of
$n$-dimensional Ricci flows in the sense of
\cite{Bamler_RF_compactness,Bamler_HK_RF_partial_regularity},
defined on an interval containing $[-T_-,-T_+]$, where
$0<T_+<T_-$. Denote its regular part by $\RR$.
Suppose that $p_\pm\in\RR_{-T_\pm}$ are centers of
$(n,k,\delta,\Gamma)$-necks at scale $\sqrt{T_\pm}$ and that $p_-$ is
an $H_n$-center of $p_+$.
Then the worldline $p_+(t)$ through $p_+$ is defined in $\RR$
for every $t\in[-T_-,-T_+]$ and is a center of an
$(n,k,\eps,\Gamma)$-neck at scale $\sqrt{-t}$. Moreover,
\begin{equation} \label{eq_dTmpmpp}
d_{-T_-}(p_-,p_+(-T_-))\leq C(n)\sqrt{T_-}.
\end{equation}
The same statement holds for Ricci flows with complete time-slices
and uniformly bounded curvature on $[-T_-,-T_+]$.
\end{Theorem}

The following corollaries give pointed smooth convergence to the cylinder at every scale, with basepoints along a single worldline in the regular part, both at singular points and at $t=-\infty$.

\begin{Corollary}[Uniqueness of cylindrical blow-ups] \label{Cor_unique_tangent}
Let $\XX$ be a metric flow over a time-interval $I$, arising as a noncollapsed
$\IF$-limit of $n$-dimensional Ricci flows in the sense of
\cite{Bamler_RF_compactness,Bamler_HK_RF_partial_regularity}, and let
$p_0\in\XX_{t_0}$, where $t_0\in I$ and $t_0>\inf I$.
Denote its regular part by $\RR$.
Suppose that a tangent flow at $p_0$, in the sense of
\cite{Bamler_RF_compactness}, is the round shrinking cylinder
$\IR^k\times(\IS^{n-k}/\Gamma)$.  Then there are $t_-<t_0$ with $(t_-,t_0]\subset I$ and a worldline
$p:(t_-,t_0)\to\RR$ such that $p(t)\to p_0$ as $t\nearrow t_0$ and the following holds.
For $\lambda>0$, write
\[
 \RR^{(t_0,\lambda)}:=
 (\RR,\lambda^2(\mathfrak t-t_0),
 \lambda^{-2}\partial_{\mathfrak t},\lambda^2g).
\]
Then
\[
 \big(\RR^{(t_0,\lambda)},p(t_0-\lambda^{-2})\big)
 \longrightarrow (\MM_{\cyl},(o,-1))
 \qquad\text{as }\lambda\to\infty
\]
in the pointed smooth Cheeger--Gromov topology on compact subsets of
$\IR^k\times(\IS^{n-k}/\Gamma)\times(-\infty,0)$, where
$o\in\{\bO\}\times(\IS^{n-k}/\Gamma)$ is fixed (and may depend on the flow).
In particular, the limit is independent of the sequence of scales,
and the tangent flow at $p_0$ is unique up to isometry.
\end{Corollary}

\begin{Corollary}[Uniqueness of cylindrical blow-downs] \label{Cor_unique_tangent_infinity}
Let $\XX$ again be a metric flow arising as a
noncollapsed $\IF$-limit of $n$-dimensional Ricci flows, now ancient and defined over
$I=(-\infty,t_0]$, and let $p_0\in\XX_{t_0}$.
Denote its regular part by $\RR$.
Suppose that a tangent flow at infinity, based at the conjugate heat
flow $(\nu_{p_0;t})_{t<t_0}$, is the round shrinking cylinder
$\IR^k\times(\IS^{n-k}/\Gamma)$.
Then there are $t_-<t_0$ and a worldline $p:(-\infty,t_-]\to\RR$ such that, with the notation of
Corollary~\ref{Cor_unique_tangent},
\[
 \big(\RR^{(t_0,\lambda)},p(t_0-\lambda^{-2})\big)
 \longrightarrow (\MM_{\cyl},(o,-1))
 \qquad\text{as }\lambda\to0
\]
in the pointed smooth Cheeger--Gromov topology on compact subsets of
$\IR^k\times(\IS^{n-k}/\Gamma)\times(-\infty,0)$.
In particular, the limit is independent of the sequence of scales,
and the tangent flow at infinity based at $p_0$ is unique up to isometry.
\end{Corollary}

For non-collapsed $\IF$-limits, Definition~\ref{defn:asymptotically_cylindrical} agrees with the intrinsic formulation in terms of tangent flows at infinity; see~\cite[Definition~6.55]{Bamler_RF_compactness}.
The following lemma makes this correspondence precise.

\begin{Lemma}\label{lem:asymptotically_cylindrical_tangent_flows}
Let $\XX$ be an ancient non-collapsed $\IF$-limit of a sequence of $n$-dimensional Ricci flows with complete time-slices and bounded curvature on compact time-intervals in the sense of~\cite{Bamler_RF_compactness,Bamler_HK_RF_partial_regularity}, and let $k$ and $\Gamma$ be as in Definition~\ref{defn:asymptotically_cylindrical}.
Then the regular part $\RR \subset \XX$ is asymptotically $(n,k,\Gamma)$-cylindrical as in Definition~\ref{defn:asymptotically_cylindrical} if and only if $\XX$ has a non-collapsed tangent flow at infinity and every such tangent flow is isometric to the round shrinking flow on $\IR^k \times (\IS^{n-k}/\Gamma)$.

Moreover, if either condition holds and $\XX$ is induced by a complete smooth Ricci flow with bounded curvature on compact time-intervals, then this flow is $\kappa$-noncollapsed on all scales for some $\kappa>0$ depending only on $n,k,\Gamma$.
\end{Lemma}

\begin{proof}[Proof of Theorem~\ref{Thm_propagation_necks}]
For compact flows, identify $p_\pm$ with their spatial points.
We first prove \eqref{eq_dTmpmpp} under the additional assumption that
$(p_+,t)$ is a center of an $(n,k,\eps,\Gamma)$-neck at scale $\sqrt{-t}$ for
every $t\in(-T_-,-T_+]$, with $\eps$ sufficiently small depending only
on the dimension.  The neck bounds and distance distortion give the
parabolic Ricci bounds needed to apply
\cite[Proposition~9.5]{Bamler_HK_entropy_estimates} on the dyadic
intervals $[-2^{-j}T_-,-2^{-j-1}T_-]$, truncating the last one at
$-T_+$.  By the semigroup property and Wasserstein contraction, we obtain
\begin{equation} \label{eq_propagation_W1_worldline}
 \begin{aligned}
 d_{W_1}^{g_{-T_-}}
 \big(\nu_{p_+,-T_+;-T_-},\delta_{p_+}\big)
 &\leq d_{W_1}^{g_{-T_-/2}}
 \big(\nu_{p_+,-T_+;-T_-/2},\delta_{p_+}\big)+C\sqrt{T_-}\\
 &\leq d_{W_1}^{g_{-T_-/4}}
 \big(\nu_{p_+,-T_+;-T_-/4},\delta_{p_+}\big)
 +2^{-1/2}C\sqrt{T_-}+C\sqrt{T_-}\\
 &\leq\cdots\leq C\sqrt{T_-}(1+2^{-1/2}+2^{-1}+\cdots)
 \leq C(n)\sqrt{T_-}.
 \end{aligned}
\end{equation}
Here dyadic times greater than $-T_+$ are understood as $-T_+$.
The $H_n$-center assumption gives
\[
 d_{W_1}^{g_{-T_-}}
 \big(\delta_{p_-},\nu_{p_+,-T_+;-T_-}\big)
 \leq\sqrt{H_n(T_--T_+)}.
\]
The triangle inequality now proves \eqref{eq_dTmpmpp}.

It remains to prove the neck assertion.  Suppose that it was false for some
fixed $\eps>0$ and consider a sequence of counterexamples
$(M_i,(g_{i,t})_{t\in[-T_{-,i},-T_{+,i}]})$ and points $p_{\pm,i}$
that satisfy the assumptions of the theorem for a sequence
$\delta_i\to0$.  After parabolic rescaling, we may assume for now
that $T_{-,i}=1$.  By replacing $\eps$ by a smaller
constant, we may also assume that $4\eps$ is smaller than the
dimensional constants required above and in Corollary~\ref{Cor_eternal}.

By pseudolocality and uniqueness of complete bounded-curvature Ricci
flows~\cite[Theorem~1.1]{Chen_Zhu_RF_uniqueness}, the neck charts at
$(p_{-,i},-1)$ give smooth convergence on compact subsets to the round
shrinking cylinder on every compact subinterval of $[-1,0)$ on which
the flows are defined, including the final time if its limit is
negative.  Here we first obtain convergence on a short time-interval
and then repeat the argument, as in the proof of
\cite[Theorem~5.89]{Bamler_Lai_MCF1}.
If $T_{+,i}$ stayed bounded away from zero along a subsequence, then
the $H_n$-center assumption and
\cite[Corollary~9.6(b)]{Bamler_HK_entropy_estimates} would bound
$d_{g_{i,-1}}(p_{-,i},p_{+,i})$ uniformly.  The entire worldline through
$p_{+,i}$ would then lie in a bounded part of the cylindrical charts,
contradicting the choice of counterexamples.  Thus $T_{+,i}\to0$.

Set $\mu_{\cyl}:=\mu(g_{\cyl},1)$ for the chosen model.
The cylindrical minimizer on the cover descends after subtracting $\log|\Gamma|$, since entropy uses the Riemannian volume.
Using the neck chart at $(p_{+,i},-T_{+,i})$, choose a smooth potential $\ov f_{i,-T_{+,i}}$ with
$(4\pi T_{+,i})^{-n/2}\int_{M_i}e^{-\ov f_{i,-T_{+,i}}}dg_{i,-T_{+,i}}=1$ and, for some $\eps'_i\to0$,
\begin{equation} \label{eq_propagation_W_upper}
 \WW[g_{i,-T_{+,i}},\ov f_{i,-T_{+,i}},T_{+,i}]
 \leq \mu_{\cyl}+\eps'_i.
\end{equation}
Extend $\ov f_{i,-T_{+,i}}$ to a time-dependent potential
$\ov f_i\in C^\infty(M_i\times[-1,-T_{+,i}])$ such that
$v_{i,t}:=(4\pi(-t))^{-n/2}e^{-\ov f_{i,t}}$ satisfies the conjugate heat
equation $\partial_tv_i=-\triangle_{g_{i,t}}v_i+R_{g_{i,t}}v_i$,
and consider the probability measures $\nu_{i,t}:=v_{i,t}dg_{i,t}$.
Since $e^{-\ov f_{i,-T_{+,i}}}$ is close to a Gaussian on the cylinder,
we may arrange
$\Var_{-T_{+,i}}(\delta_{p_{+,i}},\nu_{i,-T_{+,i}})\leq CT_{+,i}$
for some dimensional constant $C$.

Using the heat kernel $K_i$ of the flow, we can write
\begin{equation} \label{eq_propagation_heat_kernel_representation}
 v_{i,t}(x)=\int_{M_i}K_i(q,-T_{+,i};x,t)\,d\nu_{i,-T_{+,i}}(q).
\end{equation}
By the monotonicity of the Wasserstein distance and its bound by
the square root of the variance, see
\cite[Lemmas~2.7 and~3.2]{Bamler_HK_entropy_estimates}, we obtain
\begin{equation} \label{eq_propagation_heat_kernel_close}
 d_{W_1}^{g_{i,t}}\big(\nu_{i,t},\nu_{p_{+,i},-T_{+,i};t}\big)
 \leq C\sqrt{T_{+,i}}\longrightarrow0
 \qquad\text{for all }t\in[-1,-T_{+,i}].
\end{equation}
By local parabolic regularity, after passing to a subsequence,
$v_i$ converges locally smoothly to a nonnegative smooth solution
$v_\infty$ of the conjugate heat equation on the cylinder.
By \eqref{eq_propagation_heat_kernel_close}, the fundamental
conjugate heat kernels converge locally smoothly to the same limit;
see \cite[Corollary~7.5 and Theorem~9.31(e),(f)]{Bamler_RF_compactness}.
The $H_n$-center assumption and the variance monotonicity
\cite[Corollary~3.7]{Bamler_HK_entropy_estimates} also give
$\Var_{-1}(\delta_{p_{-,i}},\nu_{i,-1})\leq C$ and
$\Var_t(\nu_{i,t})\leq C(-t)$. The first inequality
implies that $\int_{M_{\cyl}}v_{\infty,-1}\,dg_{\cyl,-1}=1$.
At any fixed time $t>-1$, the second inequality implies that
either $v_{\infty,t}\equiv0$ or
$\int_{M_{\cyl}}v_{\infty,t}\,dg_{\cyl,t}=1$: once a fixed ball
carries a positive amount of mass, that bound controls the tails.
The zero alternative is excluded by uniqueness of nonnegative
integrable conjugate heat flows on the cylinder, since it would force
$v_{\infty,-1}\equiv0$.  Thus $v_\infty$ has total mass $1$ at every
time and is positive by the strong maximum principle.
In particular, $\ov f_i\to\ov f_\infty$ locally smoothly.

Fix $t\in(-1,0)$.  The maximum principle gives
$R_{g_{i,t}}\geq-\frac{n}{2(t+1)}$, and by Perelman's monotonicity formula
and \eqref{eq_propagation_W_upper} we have
$\WW[g_{i,t},\ov f_{i,t},-t]\leq\mu_{\cyl}+\eps'_i$.
We claim that
\begin{equation} \label{eq_propagation_W_semicontinuity}
 \WW[g_{\cyl,t},\ov f_{\infty,t},-t]
 \leq\liminf_{i\to\infty}\WW[g_{i,t},\ov f_{i,t},-t].
\end{equation}
The scalar curvature lower bound, local smooth convergence and
tightness give
\[
 \int_{M_i}(R_{g_{i,t}})_-\,d\nu_{i,t}
 \longrightarrow
 \int_{M_{\cyl}}(R_{g_{\cyl,t}})_-v_{\infty,t}\,dg_{\cyl,t}=0,
\]
where $a_-:=\max\{-a,0\}$.  The nonnegative gradient and positive
scalar curvature terms are lower semicontinuous by Fatou's lemma.
To see \eqref{eq_propagation_W_semicontinuity}, it therefore remains
to show that
\[
 \int_{M_i}(\ov f_{i,t})_-\,d\nu_{i,t}
 \longrightarrow
 \int_{M_{\cyl}}(\ov f_{\infty,t})_-v_{\infty,t}\,dg_{\cyl,t}.
\]
Local cylindrical convergence and
\cite[Theorem~8.1]{Bamler_HK_entropy_estimates} give
$\mathcal N_{p_{-,i},t}(\frac{t+1}{2})\geq-C(t)$.
By \cite[Theorems~5.9 and~7.2]{Bamler_HK_entropy_estimates}, using
$\mathcal N_{x,t}(\frac{t+1}{2})$ at the earlier endpoint in the
Gaussian heat-kernel bound, we obtain
$K_i(q,-T_{+,i};x,t)\leq C(t)\exp(C(t)d_{g_{i,t}}(p_{-,i},x))$
uniformly in $q$.  Thus \eqref{eq_propagation_heat_kernel_representation}
gives the pointwise bound
\[
 (\ov f_{i,t})_-(x)\leq C(t)(1+d_{g_{i,t}}(p_{-,i},x)).
\]
The inequality $\Var_t(\nu_{i,t})\leq C(-t)$ and the positive mass in a fixed ball imply
$\Var_t(\delta_{p_{-,i}},\nu_{i,t})\leq C(t)$.
Hence the integral of $(\ov f_{i,t})_-$ outside
$B_{g_{i,t}}(p_{-,i},R)$ is at most $C(t)/R$ for $R\geq1$,
which proves the desired convergence.
Fatou's lemma now gives \eqref{eq_propagation_W_semicontinuity},
which, together with the logarithmic Sobolev inequality on the
cylinder and the upper bound, implies
$\WW[g_{i,t},\ov f_{i,t},-t]\to\mu_{\cyl}$ for every fixed
$t\in(-1,0)$.  In particular, after enlarging $\eps'_i$, we have
\begin{equation} \label{eq_propagation_W_lower}
 \WW[g_{i,-1/2},\ov f_{i,-1/2},\tfrac12]
 \geq \mu_{\cyl}-\eps'_i.
\end{equation}

Let now
\[
 t^*_i:=\sup\big\{t\in[-1,-T_{+,i}]:
 (p_{+,i},t)\text{ is not a center of an $(n,k,\eps,\Gamma)$-neck
 at scale }\sqrt{-t}\big\}.
\]
After rescaling the flows, we may assume henceforth that $t^*_i=-1$.
We keep the same notation for the rescaled quantities.
Applying \eqref{eq_propagation_W1_worldline} on $[-1,-T_{+,i}]$ gives
\begin{equation} \label{eq_propagation_W1_critical}
 d_{W_1}^{g_{i,-1}}(\delta_{p_{+,i}},\nu_{p_{+,i},-T_{+,i};-1})
 \leq C.
\end{equation}
By \eqref{eq_propagation_heat_kernel_close} and $T_{+,i}\leq1$,
\begin{equation}\label{eq_propagation_potential_first_moment}
\int_{M_i}d_{g_{i,-1}}(p_{+,i},x)\,d\nu_{i,-1}(x)
\leq C+C\sqrt{T_{+,i}}\leq C(n).
\end{equation}
We claim that $T_{-,i}\to\infty$.  Otherwise, after passing to a
subsequence, the rescaling factors $T_{-,i}$ stay bounded, so the
preceding cylindrical and heat-kernel convergence still applies.
In particular, the limit of $\nu_{p_{+,i},-T_{+,i};-1}$ is a probability
measure on the cylinder.
Since a fixed compact subset of the cylindrical charts carries at
least half of the heat-kernel mass, \eqref{eq_propagation_W1_critical}
keeps $(p_{+,i},-1)$
in a bounded part of these charts.  Smooth cylindrical convergence
then gives $(n,k,\eps,\Gamma)$-necks along this worldline at all times
sufficiently close to $-1$, contradicting the definition of $t^*_i$.

We also have $T_{+,i}\to0$. Otherwise, pass to a subsequence with
$T_{+,i}\to T_{+,\infty}>0$. The terminal neck charts converge
on every fixed ball to the complete cylinder at time
$-T_{+,\infty}$. Applying backward pseudolocality
\cite[Theorem~2.47]{Bamler_HK_RF_partial_regularity} on balls of
radius $c(n)\sqrt{T_{+,i}}$ covering these charts gives
\[
|\Rm_{g_i}|\leq C(n)T_{+,i}^{-1}
\quad\textnormal{on}\quad
B_{g_{i,-T_{+,i}}}(p_{+,i},A\sqrt{T_{+,i}})
\times[-(1+\theta)T_{+,i},-T_{+,i}]
\]
for every fixed $A<\infty$ and all sufficiently large $i$, with
$\theta=\theta(n)>0$ independent of $A$. The spatial balls are
followed along worldlines.
Metric comparison gives
$e^{-C\theta}g_{i,-T_{+,i}}\leq g_{i,t}\leq
e^{C\theta}g_{i,-T_{+,i}}$ there, so the terminal volume lower
bounds persist. Together with the curvature bound, they give
$\inj_{g_{i,t}}\geq c(n)\sqrt{T_{+,i}}$ on smaller balls.
Shi's estimates then give a smooth limit
on the latter half of this interval, including its terminal time.
The limiting metrics are complete by the same metric comparison
with the complete terminal cylinder, and have uniformly bounded curvature.
Its terminal metric is cylindrical, so by
\cite[Theorem~1]{Kotschwar_RF_backward_uniqueness} the limit is
the shrinking cylinder throughout that interval.
Repeat this argument from its earlier endpoint. At each step the
limiting metric is exactly cylindrical, so the same $\theta$ applies.
Finitely many steps reach time $-1$, since $T_{+,\infty}>0$.
This gives cylindrical convergence at $(p_{+,i},-1)$, a contradiction.

By scale invariance and Perelman's monotonicity formula~\cite[(3.4)]{Perelman1},
\eqref{eq_propagation_W_upper} and \eqref{eq_propagation_W_lower} imply
\begin{equation} \label{eq_propagation_W_defect}
 \int_{-T_{-,i}/2}^{-T_{+,i}}(-t)
 \int_{M_i}\left|
 \Ric_{g_{i,t}}+\nabla^2\ov f_{i,t}
       -\frac{1}{2(-t)}g_{i,t}
 \right|^2d\nu_{i,t}\,dt
 \longrightarrow0.
\end{equation}
The point $(p_{+,i},-1)$ is a center of an $(n,k,2\eps,\Gamma)$-neck,
so $|B_{g_{i,-1}}(p_{+,i},1)|_{g_{i,-1}}\geq c(n)>0$.
Together with \eqref{eq_propagation_W1_critical}, this gives
uniform noncollapsing by
\cite[Corollary~5.11 and Theorem~8.1]{Bamler_HK_entropy_estimates}.
By~\cite[Theorem~7.8]{Bamler_RF_compactness}, after passing to a
subsequence, we obtain an $\IF$-limit on $(-\infty,0]$, translating
the final times to zero by shifts tending to zero.
\cite[Theorem~6.40]{Bamler_RF_compactness}, applied using
\eqref{eq_propagation_W1_critical}, allows us to view
$(p_{+,i},-1)$ as a convergent sequence of points in this limit.
Backward pseudolocality and the neck charts imply that its
limit $p_\infty$ lies in the regular part. Since $T_{+,i}\to0$,
the terminal points converge to a point $x_\infty$ at time $0$,
and \eqref{eq_propagation_heat_kernel_close} identifies the
limiting measures as $\nu_{\infty,t}=\nu_{x_\infty;t}$.
It follows from
\eqref{eq_propagation_W_defect} that the
limit satisfies
\[
 \Ric+\nabla^2\ov f-\frac{1}{2(-t)}g=0
\]
on its regular part.

\begin{Claim}\label{claim:propagation_homotheties}
On the regular component containing $p_\infty$, the flow
$\Phi_{s,t}$ of $\partial_t-\nabla\ov f$ is defined between any two
interior times $s,t$ of the limiting interval and satisfies
$\Phi_{s,t}^*g_t=\frac{-t}{-s}g_s$.
\end{Claim}

\begin{proof}
The limiting density is the pointed conjugate heat kernel based at
$x_\infty$. Its Nash entropy is locally Lipschitz by
\cite[Proposition~4.28 and Theorem~15.45(a)]{Bamler_HK_RF_partial_regularity}.
For almost every $\tau=-t>0$, the conjugate heat equation and the
trace of the soliton equation give
\[
\begin{aligned}
\frac{d}{d\tau}\mathcal N_{x_\infty}(\tau)
&=\int\left(|\nabla\ov f|^2+R-\frac{n}{2\tau}\right)
\,d\nu_{\infty,-\tau}\\
&=\int\left(|\nabla\ov f|^2-\triangle\ov f\right)
\,d\nu_{\infty,-\tau}=0.
\end{aligned}
\]
Here the heat-kernel entropy estimates and the cutoffs in
\cite[Lemma~15.27]{Bamler_HK_RF_partial_regularity} justify
differentiation and integration by parts on the regular part;
the tails are controlled as in
\cite[Proposition~3.1]{Fang_Li_strong_uniqueness}.
Thus the pointed Nash entropy is constant, and
\cite[Theorem~15.69]{Bamler_HK_RF_partial_regularity} gives the
asserted global homotheties between interior times.
\end{proof}

We next find a critical point of $\ov f$ in this neck. At time $-1$,
$d\nu_{\infty,-1}=(4\pi)^{-n/2}e^{-\ov f}\,dg_{\infty,-1}$.
By \eqref{eq_propagation_potential_first_moment} and local convergence,
for a sufficiently large fixed $A=A(n)$,
\[
\nu_{\infty,-1}(B(p_\infty,A))\geq1-\frac{C(n)}{A}\geq\frac12.
\]
Choose $0\leq\zeta\leq1$ supported in $B(p_\infty,2A)$, equal to
$1$ on $B(p_\infty,A)$, with $|\nabla\zeta|\leq C(n)$.
The neck bounds and the soliton equation give
$|\nabla^2\ov f|\leq C(n)$ and $\triangle\ov f=n/2-R$.
Integration by parts yields
\[
\begin{aligned}
\int\zeta^2|\nabla\ov f|^2\,d\nu_{\infty,-1}
&=\int\zeta^2\triangle\ov f\,d\nu_{\infty,-1}
+2\int\zeta\langle\nabla\zeta,\nabla\ov f\rangle
\,d\nu_{\infty,-1}\\
&\leq C(n)+\frac12\int\zeta^2|\nabla\ov f|^2
\,d\nu_{\infty,-1}.
\end{aligned}
\]
Thus $|\nabla\ov f(z)|\leq C(n)$ at some
$z\in B(p_\infty,A)$, and $|\nabla^2\ov f|\leq C(n)$ gives
$|\nabla\ov f(p_\infty)|\leq C(n)$.
In the neck coordinates, on any fixed
$\IB_R^k\times(\IS^{n-k}/\Gamma)$ we consequently have
\[
|\nabla\ov f|\leq C(n)(1+R),\qquad
\partial^2_{\bx}\ov f=\tfrac12 I+O_R(\eps),\qquad
|\partial_{\bx}\ov f(0,\by)|\leq C(n).
\]
Fix $R>8C(n)$ and then take $\eps$ sufficiently small that
$R+(4\eps)^{-1}<(2\eps)^{-1}$ and
$\partial^2_{\bx}\ov f\geq\frac14 I$ on
$\overline{\IB}_R^k\times(\IS^{n-k}/\Gamma)$. On $|\bx|=R$,
\[
\bx\cdot\partial_{\bx}\ov f(\bx,\by)
\geq\tfrac14R^2-C(n)R>0.
\]
The minimum of $\ov f$ on
$\overline{\IB}_R^k\times(\IS^{n-k}/\Gamma)$ therefore occurs at an
interior point $p_*$, with
$\nabla\ov f(p_*)=0$ and $d_{g_{\infty,-1}}(p_*,p_\infty)\leq C(n)$.
Recentering the neck chart makes $p_*$ a center of an
$(n,k,4\eps,\Gamma)$-neck at scale $1$.

The conjugate heat equation and the soliton equation give
$\partial_t\ov f=|\nabla\ov f|^2$ and
$\partial_t\nabla\ov f=(-t)^{-1}\nabla\ov f$.
By Claim~\ref{claim:propagation_homotheties}, the trajectory through
$p_*$ of $\partial_t-\nabla\ov f$ is defined for every $t<0$;
these identities show that it is a worldline.
It is a neck center at scale $\sqrt{-t}$ for every $t<0$.
Corollary~\ref{Cor_eternal} therefore identifies the limit on
$(-\infty,0)$ with the cylindrical flow, so $p_\infty$ is a center
of an exact cylindrical neck.
Smooth
convergence near $(p_\infty,-1)$ now implies that $(p_{+,i},-1)$ is a
center of an $(n,k,\frac12\eps,\Gamma)$-neck for all sufficiently large $i$.  This is
the desired contradiction and proves the neck assertion, and hence
also \eqref{eq_dTmpmpp} by the initial argument.

For complete flows with bounded curvature on compact time intervals,
the same proof applies using conservative heat kernels and spatial
cutoffs in the entropy identities.
For a noncollapsed $\IF$-limit, decrease $\delta$ once more and
approximate the two endpoint necks smoothly by
\cite[Theorem~2.5]{Bamler_HK_RF_partial_regularity}.
Choose $H_n$-centers $z_i$ of the approximating later points
at the earlier time. The $H_n$-center condition on $p_-$ and
convergence of the conjugate heat measures
\cite[Theorems~6.40 and~9.31]{Bamler_RF_compactness} give
\[
\nu_{p_{+,i},-T_+;-T_-}
\big(B(p_{-,i},C\sqrt{T_-})\big)\geq\tfrac12,
\qquad d_{-T_-}(z_i,p_{-,i})\leq C\sqrt{T_-}.
\]
The second inequality follows from
$\Var_{-T_-}(\delta_{z_i},\nu_{p_{+,i},-T_+;-T_-})\leq H_n(T_--T_+)$.
Thus $z_i$ lies in the initial neck chart and is a sufficiently
accurate neck center. Apply the proved statement with accuracy
$\eps/2$ to $(z_i,p_{+,i})$.
Fix $t_*\in(-T_-,-T_+)$ sufficiently close to $-T_+$ that
$p_+(t_*)$ is defined in the regular part.
The neck bounds and distance distortion give a parabolic neighborhood
of $p_{+,i}(t_*)$ of fixed spatial radius, unscathed throughout
$(-T_-,-T_+)$, with curvature bounded by $C/T_+$ and a uniform
volume lower bound at $t_*$. By
\cite[Theorem~9.58]{Bamler_RF_compactness}, the limiting worldline
therefore extends backward throughout $(-T_-,-T_+]$.
Applying the same theorem on smaller neck charts gives
$(n,k,\eps,\Gamma)$-necks in the limit.
The distance bound already proved and the estimate for $z_i$ place
$p_{+,i}(-T_-)$ in a fixed compact subset of the initial neck chart.
Smooth convergence there extends the worldline to $-T_-$ and gives
\eqref{eq_dTmpmpp} for the prescribed $p_-$.
\end{proof}

\medskip

\begin{proof}[Proof of Corollaries~\ref{Cor_unique_tangent} and~\ref{Cor_unique_tangent_infinity}.]
Both statements follow from Theorem~\ref{Thm_propagation_necks}, with the following care concerning the $H_n$-centers.
At times $s_1<s_2<t_0$ from the sequence of cylindrical scales, choose $H_n$-centers $q_i\in\XX_{s_i}$ of $p_0$, and let $\td q_1\in\XX_{s_1}$ be an $H_n$-center of $q_2$.
The semigroup property, Wasserstein contraction and the $H_n$-center bounds give
\[
 \begin{aligned}
 d_{s_1}(q_1,\td q_1)
 &=d_{W_1}^{s_1}(\delta_{q_1},\delta_{\td q_1})\\
 &\leq d_{W_1}^{s_1}(\delta_{q_1},\nu_{p_0;s_1})
       +d_{W_1}^{s_1}(\nu_{p_0;s_1},\nu_{q_2;s_1})
       +d_{W_1}^{s_1}(\nu_{q_2;s_1},\delta_{\td q_1})\\
 &\leq \sqrt{H_n(t_0-s_1)}+\sqrt{H_n(t_0-s_2)}
       +\sqrt{H_n(s_2-s_1)}
 \leq 3\sqrt{H_n(t_0-s_1)}.
 \end{aligned}
\]
Thus replacing $q_1$ by $\td q_1$ preserves arbitrarily accurate
necks and ensures the $H_n$-center hypothesis.
Fix a small neck accuracy and discard finitely many
cylindrical scales. In the blow-up case, propagate the worldlines
through the later centers back to a fixed earlier cylindrical
scale. The distance bound in the theorem keeps their positions
in a compact subset of that neck chart. The uniform neck bounds
and pseudolocality then give a subsequential limiting worldline
$p(t)$ on an interval ending at $t_0$.
In the blow-down case, propagate the worldline through a fixed
later center to successively earlier cylindrical scales. Uniqueness
of worldlines makes these extensions agree.

For this worldline, \eqref{eq_propagation_W1_worldline}, the
semigroup property, and the center bounds above give
\[
d_{s_j}(p(s_j),q_j)\leq C\sqrt{t_0-s_j}
\]
at every cylindrical scale $s_j$. In the blow-up case, the same
estimates give
\[
d_{W_1}^{s}(\delta_{p(s)},\nu_{p_0;s})
\leq C\sqrt{t_0-s},
\]
which implies $p(s)\to p_0$ by Wasserstein contraction to any
fixed earlier time.
The bounded rescaled displacement shows that $p(s_j)$ are neck
centers with accuracy tending to zero. For consecutive cylindrical
times $s_1<s_2$, an $H_n$-center of $p(s_2)$ at $s_1$ is likewise
within $C\sqrt{t_0-s_1}$ of $p(s_1)$, by
\eqref{eq_propagation_W1_worldline}; hence it is an arbitrarily
accurate neck center. Applying Theorem~\ref{Thm_propagation_necks}
with any prescribed accuracy between consecutive such scales proves
cylindrical convergence along $p(t)$ at all scales.
Pseudolocality and Shi's estimates give pointed smooth
spacetime convergence on every compact subset of the cylinder.
For $\Gamma\ne\{1\}$, apply Theorem~\ref{thm:HMHF} once on a sufficiently remote time tail, with $\tau=-\log(t_0-t)$ and $m,J,-\la$ sufficiently large. Here $U^+_{\conf}=0$ by Lemma~\ref{l:quotient_cylinder_spectrum}.
Proposition~\ref{Prop_mode_analysis}\ref{Prop_mode_analysis_d} for $\tau\to-\infty$ and \ref{Prop_mode_analysis_e} for $\tau\to\infty$, together with Theorem~\ref{thm:HMHF}\ref{thm:HMHF_a}, \ref{thm:HMHF_c}, \ref{thm:HMHF_e} and Claim~\ref{Cl_Upp_quad_1tau}, give $\|U^+(\tau)\|=O(|\tau|^{-1})$ in both directions; in the forward case $T_1=\infty$, the decay is exponential instead.
The higher-order bounds in Proposition~\ref{Prop_mode_analysis}\ref{Prop_mode_analysis_d}--\ref{Prop_mode_analysis_e} and \eqref{eq_UUim1_choice}, or the preceding exponential bound, imply $R(\tau)\to\infty$ through the Gaussian term in $\UU^-$. Interpolation with the $C^m$-bound in \eqref{thm:HMHF_open_1} then gives, on each fixed compact cylindrical region $K$,
\[
\|u_\tau-U^+_{\rot}(\tau)\|_{C^2(K)}\leq C_K|\tau|^{-2}.
\]
Write $(\bx(\tau),\by(\tau))$ for the rescaled HMHF coordinates of $p(t)$.
The tension field of $g_{\cyl}+U^+_{\rot}$ has zero spherical component, since its mixed block vanishes and its spherical block is a function of $\bx$ times $g_{\IS^{n-k}/\Gamma}$. The HMHF evolution computed in the proof of Claim~\ref{claim:u=eta} thus gives
\[
|\partial_\tau\by(\tau)|\leq C|\tau|^{-2},\qquad \partial_\tau\bx(\tau)=\tfrac12\bx(\tau)+o(1).
\]
Integrating the first bound gives $\by(\tau)\to\by_\infty$ in either direction. The second equation and the boundedness of $\bx$ from \eqref{thm:HMHF_open_2} give $\bx(\tau)\to0$. Hence the limiting basepoint is the fixed point $o=(0,\by_\infty)$.
\end{proof}
\bigskip

\begin{proof}[Proof of Lemma~\ref{lem:asymptotically_cylindrical_tangent_flows}.]
Suppose first that $\RR$ is asymptotically $(n,k,\Gamma)$-cylindrical, and choose the worldline $t \mapsto p(t)$ from Definition~\ref{defn:asymptotically_cylindrical}, defined for $t \leq t_0$.
The cylindrical convergence gives a bound of the form $R(p(t)) \leq C/(1+t_0-t)$, where $C$ denotes a generic constant independent of $t$.
Thus the restriction of this worldline to $[t,t_0]$ has $\LL$-length at most
\[
\int_0^{t_0-t} \sqrt{\tau}\,\frac{C}{1+\tau}\,d\tau \leq C\sqrt{t_0-t}.
\]
Consequently, the reduced distance based at $p(t_0)$ satisfies $\ell_{p(t_0)}(p(t))\leq C$; see~\cite[Section~7]{Perelman1}.
By perturbing these paths within the cylindrical region near $p(t)$ and using geometric closeness to a round cylinder near $p(t)$, we obtain, for all sufficiently negative $t$
\[
\sup_{B_{g_t}(p(t),\sqrt{t_0-t})}\ell_{p(t_0)}\leq C.
\]
The same path estimates hold in the approximating Ricci flows.
On these flows, Perelman's differential Harnack inequality gives the conjugate heat kernel lower bound in terms of reduced distance; see~\cite[Corollary~9.5]{Perelman1}.
Passing to the limit using~\cite[Theorem~9.31(e),(f)]{Bamler_RF_compactness} yields
\[
\inf_{B_{g_t}(p(t),\sqrt{t_0-t})}K(\cdot\,;p(t_0))
\geq (4\pi(t_0-t))^{-n/2}e^{-C}.
\]
Together with the volume lower bound supplied by cylindrical convergence, this shows that $B_{g_t}(p(t),\sqrt{t_0-t})$ carries a uniformly positive amount of the conjugate heat kernel measure $\nu_{p(t_0);t}$.
The concentration estimate for $H_n$-centers therefore implies that every such center at time $t$ lies within distance $C\sqrt{t_0-t}$ of $p(t)$; see~\cite[Definition~3.10 and Proposition~3.13]{Bamler_HK_entropy_estimates}.

For each sufficiently negative $t$, the same distance estimate and smooth convergence give a lower volume bound of $c(t_0-t)^{n/2}$ for balls of radius $\sqrt{t_0-t}$ about the $H_n$-centers of the points approximating $p(t_0)$.
By~\cite[Theorem~1.8]{Chan_Ma_Zhang_2021}, this gives a uniform pointed Nash entropy lower bound, which passes to $p(t_0)$ by~\cite[Corollary~15.47]{Bamler_HK_RF_partial_regularity}; $\IF$-compactness then gives a non-collapsed tangent flow at infinity~\cite[Section~2.7]{Bamler_HK_RF_partial_regularity}.

{Since the $H_n$-centers stay within bounded rescaled distance of $p(t)$,} the time-slices of any tangent flow at infinity agree with the pointed smooth limits in Definition~\ref{defn:asymptotically_cylindrical}, which shows that every such tangent flow is the round shrinking cylindrical flow.
Finally, the tangent flows at infinity do not depend on the basepoint by~\cite[Theorem~1.5]{Chan_Ma_Zhang_closed_tangent_flows}.

Conversely, suppose that $\XX$ has a non-collapsed tangent flow at infinity and that every such tangent flow is the round shrinking cylinder $\IR^k\times(\IS^{n-k}/\Gamma)$.
Corollary~\ref{Cor_unique_tangent_infinity} gives a backward worldline $p(t)$ satisfying \eqref{eq:asymptotically_cylindrical_convergence}.
By~\cite[Theorem~2.4(c)]{Bamler_HK_RF_partial_regularity}, the regular time-slices are connected, and hence so is the regular spacetime $\RR$.

It remains to verify the last condition in Definition~\ref{defn:asymptotically_cylindrical}.
Consider two worldlines $p_1(t),p_2(t)$ satisfying~\eqref{eq:asymptotically_cylindrical_convergence}, and choose a common time $t_0$ in their domains.
For $t<t_0$, let $z_i(t)\in\XX_t$ be an $H_n$-center of $p_i(t_0)$, $i=1,2$.
The reduced-distance and heat-kernel estimates at the beginning of the proof imply $d_t(p_i(t),z_i(t))\leq C\sqrt{t_0-t}$ for sufficiently negative $t$.
The definition of $H_n$-centers and Wasserstein contraction therefore imply (see~\cite[Definition~3.35 and Proposition~3.24(c)]{Bamler_RF_compactness})
\begin{multline*}
d_t(p_1(t),p_2(t))
\leq 2C\sqrt{t_0-t}+d_t(z_1(t),z_2(t))
\leq 2(C+\sqrt{H_n})\sqrt{t_0-t}
   +d_{W_1}^{\XX_t}(\nu_{p_1(t_0);t},\nu_{p_2(t_0);t})\\
\leq 2(C+\sqrt{H_n})\sqrt{t_0-t}
   +d_{t_0}(p_1(t_0),p_2(t_0))
\leq C(p_1,p_2)\sqrt{-t}.
\end{multline*}
By~\cite[Theorem~2.4(c)]{Bamler_HK_RF_partial_regularity}, $d_t=d_{g_t}$ on $\RR_t$.

For a complete smooth flow as in the last assertion, the pointed Nash entropy at every basepoint converges as the backward time scale tends to infinity to the entropy of the normalized cylinder, by~\cite[Section~2.7 and Theorem~2.40]{Bamler_HK_RF_partial_regularity}. By monotonicity, this gives a lower bound at every scale, and the volume estimate in~\cite[Theorem~6.1]{Bamler_HK_entropy_estimates} (see also~\cite[Appendix~A]{Bamler_2021_fundamental_group}) gives the claimed $\kappa$-noncollapsing.
\end{proof}

\subsection{Ancient asymptotically cylindrical flows}
In this subsection, we summarize the results obtained so far for ancient, asymptotically cylindrical Ricci flows.
The main outcome is Proposition~\ref{Prop_ancient_summary}, which is the Ricci flow analogue of \cite[Proposition~5.101]{Bamler_Lai_MCF1}.
After choosing an appropriate harmonic map heat flow gauge, the rescaled metric perturbation can be approximated by a finite-dimensional mode $U^+(\tau) \in \sV_{\rot, \geq - J}$ of rotationally symmetric tensors on cylindrical regions of radii $R(\tau) \to \infty$ as $\tau \to -\infty$.
The evolution of $U^+(\tau)$ is described by a nonlinear ODE up to an error bounded by $C\|U^+\|^{J+1}$, and the same bound controls the error in the metric approximation, where $J$ can be chosen arbitrarily large.
Moreover, the size of $U^+$ is controlled by its leading component $U^{++}$, while the higher rotationally symmetric modes are bounded by successively higher powers of $\|U^{++}\|$.
This allows us to study the asymptotics of the flow by analyzing a finite-dimensional system, which will be the starting point for the following sections.

\begin{Proposition} \label{Prop_ancient_summary}
For $J, m \geq 10$, $\eta > 0$ and $\delta \leq \ov\delta(J,m,\eta)$, the following is true.

Let $\MM$ be an asymptotically $(n,k,\Gamma)$-cylindrical Ricci flow spacetime that satisfies the pseudolocality property.
Suppose that $q_0 \in \MM_{t_0- r_0^2}$, $r_0 > 0$, is a point that survives until time $-\infty$.
Suppose that $q_0(t)$ is a center of an $(n,k,\delta,\Gamma)$-neck at scale $(t_0 - t)^{1/2}$ for all $t \leq t_0 - r_0^2$. 
Set $\td\tau := -2\log r_0$.
Then either $\MM$ is homothetic to the round shrinking cylinder or there is
\begin{itemize}
\item an open subset $\Omega \subset \MM_{(-\infty,t_0-r_0^2]}$,
\item  a smooth map $\chi : \Omega \to \IR^k \times (\IS^{n-k}/\Gamma)$ and
\item a locally Lipschitz function
\[ U^+  : (-\infty, \td\tau ] \lto \sV_{\rot, \geq -J}  \]
\end{itemize}
such that the following is true (we will use the norm $\| \cdot \| := \| \cdot \|_{L^2_f}$ on $\sV_{\rot, \geq -J}$):
\begin{enumerate}[label=(\alph*)]
\item \label{Prop_ancient_summary_a} For all $\tau\le \td\tau$,
$$\xi_{e^{\tau/2}}\circ\chi_{t_0-e^{-\tau}}(q_0(t_0-e^{-\tau})) \in \IB_\eta^k \times (\IS^{n-k}/\Gamma).$$
\item \label{Prop_ancient_summary_b} $\chi$ evolves by harmonic map heat flow.
\item \label{Prop_ancient_summary_b2} Set $U^{++}(\tau) := \PP_{\sV^{++}} U^+(\tau)$, where $\sV^{++} = \sV_{\rot,1} \oplus \sV_{\rot, \frac12} \oplus \sV_{\sph,0}$. 
We have the following coarse bound 
\[ \| U^{++} \| \leq \|U^+\| \leq C \| U^{++} \| \leq C(\td\tau- \tau+10)^{-1}. \]
\item \label{Prop_ancient_summary_c} Set $e^{-R^2(\tau)} := \Vert U^+(\tau) \Vert^J$. Then
\begin{equation} \label{eq:R(tau)}
\IB^k_{R(\tau)} \times (\IS^{n-k}/\Gamma) \subset \xi_{e^{\tau/2}}\circ\chi_{t_0-e^{-\tau}}(\Omega_{t_0-e^{-\tau}}),
\end{equation}
and 
\[ \| e^\tau(\xi_{e^{\tau/2}})_* (\chi_{t_0-e^{-\tau}})_* g_{t_0-e^{-\tau}}-g_{\cyl} - U^+(\tau) \|_{C^m(\IB^k_{R(\tau)} \times (\IS^{n-k}/\Gamma))} \leq C(J,m) \| U^+ \|^{J+1}. \]
We also have the following bounds
\[ \| e^\tau(\xi_{e^{\tau/2}})_* (\chi_{t_0-e^{-\tau}})_* g_{t_0-e^{-\tau}}-g_{\cyl} \|_{C^m(\IB^k_{R(\tau)}\times(\IS^{n-k}/\Gamma))}
+ \| U^+(\tau) \| + \| U^+ \|_{C^m(\IB^k_{R(\tau)} \times (\IS^{n-k}/\Gamma))} \leq \eta. \]
\item \label{Prop_ancient_summary_d} The evolution of $U^+$ is controlled by the following ODI
\[ \| \partial_\tau U^+ - L U^+ - Q^+_J(U^+) \| \leq C(J) \| U^{++} \|^{J+1}. \]
\item \label{Prop_ancient_summary_e} We have
\begin{align*}
  \| U^+_{\rot, < -i} \| &\leq C(J) \| U^{++} \|^{i +2}, \qquad i = 0, \ldots, J-1.
\end{align*}
\item \label{Prop_ancient_summary_f} There is a time $T_3 \in [-\infty, \td\tau]$ such that the following is true:
\begin{enumerate}[label=(g\arabic*)]
\item If $T_3 > -\infty$, consider the minimum and maximum spectral values $U_{0,\min}(\tau)$ and $U_{0,\max}(\tau)$ of the matrix $(U_{0,ij} (\tau))$, where $\PP_{\sV_{\sph,0}} U^{++}(\tau) = \sum_{i,j} U_{0,ij}(\tau) \mathfrak p_{ij}^{(2)}g_{\IS^{n-k}/\Gamma}$.
Then on $(-\infty,T_3]$ we have $-U_{0,\min} \ge U_{0,\max}$ and
\[
\|U^+(\tau)\|\leq C|U_{0,\min}|,
\qquad
\big| \partial_\tau U_{0,\min} + \sqrt{2} U_{0,\min}^2 \big|
\leq 0.1 U_{0,\min}^2.
\]
\item \label{Prop_ancient_summary_f2} If $T_3 < \td\tau$, then on $(T_3, \td\tau]$ we have
\[ \|U^+\|\leq C\|U^{++}_{>0}\|,\qquad
\partial_\tau \Vert U^{++}_{>0} \Vert \geq 0.49 \Vert U^{++}_{>0} \Vert. \]
Consequently, for every $\tau\in (T_3,\td\tau]$,
\[ \Vert U^{++}_{>0}(\tau) \Vert \leq e^{0.49(\tau-\td\tau)}\Vert U^{++}_{>0}(\td\tau)\Vert. \]
\end{enumerate}
\item \label{Prop_ancient_summary_g} If $q_0$ survives until all times of the interval $(-\infty,t_0)$ and if $q_0(t)$ is a center of an $(n,k,\delta,\Gamma)$-neck at scale $(t_0 - t)^{1/2}$ for all $t < t_0$, then $\MM$ is homothetic to a round shrinking cylinder.
\end{enumerate}
\end{Proposition}

\begin{proof}[Proof of Proposition~\ref{Prop_ancient_summary}.]
Fix $J, m,   \eta$ as in the statement of the proposition and choose
\[ J' := 2(n+m+4) (J+1)+1, \qquad \la' := -  J' - 1   \]
Fix $0<\beta<\gamma_{\osc}/4$ below the thresholds for
these parameters in Theorem~\ref{thm:HMHF} and
Proposition~\ref{Prop_mode_analysis}.
Let $\eta', \eps', R^* > 0$ be constants whose values we will determine in the course of the proof.
Throughout the argument, we will impose  bounds on these constants and on the constant $\delta$ while respecting the following order
\[ \eta,\; m,\; J',\; \la',\; \beta,\; \eta',\;\eps',\; R^*,\; \delta. \]
For each of the other constants, we will only require bounds involving constants appearing earlier in this list.
We will frequently omit the ``$L^2_f$'' subscript when there is no chance of confusion.

Let us now apply Theorem~\ref{thm:HMHF} with the parameters $m+1, J', \la', \eta' \leq \ov\eta(\la'), \eps' \leq \ov\eps(m+1,J',\la',\beta,\eta')$  and some $R^* \geq \underline R^* (m+1, \lb J', \lb \la', \lb, \eta', \lb \eps')$.
Denote the resulting functions by
\begin{alignat*}{2}
R' &\colon (-\infty,\td\tau] \to [R^*,\infty], &\qquad
\UU^{-,\prime} &\colon (-\infty,\td\tau] \to \IR_{\geq 0},\\
U^{+,\prime} &\colon (-\infty,\td\tau] \to \sV_{>\la'}, &
U^{-,\prime} &\colon (-\infty,\td\tau] \to \sV_{\leq\la'}.
\end{alignat*}
Set $u_\tau=e^\tau(\xi_{e^{\tau/2}})_* (\chi_{t_0 - e^{-\tau}})_* g_{t_0 - e^{-\tau}} - g_{\cyl}$ and
\[ U^+ (\tau) := \PP_{\sV_{\rot, \geq - J}} U^{+,\prime}(\tau),  \qquad U^{++} (\tau) := \PP_{\sV^{++}} U^{+}(\tau). \]
By \eqref{thm:HMHF_open_1}, we have for $\tau \leq \td\tau$
\[ \IB^k_{R'(\tau)} \times (\IS^{n-k}/\Gamma) \subset \DD'_\tau, \qquad
\Vert u'_\tau \Vert_{C^{m+1} (\IB^k_{R'(\tau)} \times (\IS^{n-k}/\Gamma))} \leq \eta'. \]
Proposition~\ref{Prop_ancient_summary}~\ref{Prop_ancient_summary_a} follows directly from \eqref{thm:HMHF_open_2} in Theorem~\ref{thm:HMHF}.

Since $U^{\prime,+}$ and $\UU^{-,\prime}$ satisfy the conclusions of Theorem~\ref{thm:HMHF}, we can apply Proposition~\ref{Prop_mode_analysis} with the fixed $\beta$ above, assuming $\eta' \leq \ov\eta(\la', J', \beta)$.
Proposition \ref{Prop_mode_analysis}\ref{Prop_mode_analysis_d} (note $J_0=\min\{J',-\la'\}-1=J'-1$) implies the following bounds, which hold on the time-interval $(-\infty, \td\tau]$,
\begin{align}
 \UU^- + \| U^{\prime,+}_{\osc} \| &\leq C \| U^{++} \|^{J'-1} \label{eq:small_1}\\
 \| U^{\prime,+}_{\rot,<0} \| &\leq C \| U^{++} \|^2 \label{eq:small_2}\\
\| U^{\prime,+}_{\rot,\sph, < -i+2} \| + \| U^{\prime,+}_{\rot,\aff, < -i+1} \|  &\leq C \| U^{++} \|^{i+1}, \qquad i = 2, \ldots, J_0. \label{eq:small_3}
\end{align}
Proposition~\ref{Prop_ancient_summary}~\ref{Prop_ancient_summary_e} now follows directly from \eqref{eq:small_3}. The bound 
\begin{equation} \label{eq_Up10Uppproof}
\Vert U^{+} \Vert\le \Vert U^{\prime,+} \Vert \leq C \Vert U^{++} \Vert
\end{equation}
follows from Theorem~\ref{thm:HMHF}\ref{thm:HMHF_b} and
\eqref{eq:small_1}--\eqref{eq:small_3} for $\eta' \leq \ov\eta' (\la', J')$.
Proposition~\ref{Prop_ancient_summary}~\ref{Prop_ancient_summary_d} is a consequence of the first PDE-ODI equation \eqref{eq_ODI_1} in Theorem \ref{thm:HMHF}, via \eqref{eq:small_1} and the following bound, which is implied by \eqref{eq:small_1} and \eqref{eq:small_3}:
\begin{equation} \label{eq_UppUpJp2} 
\Vert U^{+,\prime} - U^+ \Vert \leq \| U^{+,\prime}_{\osc} \| +  \| U^{+,\prime}_{\rot, \leq - J-\frac12} \| \leq C(J') \Vert U^{++} \Vert^{J+2}. 
\end{equation}
Proposition~\ref{Prop_ancient_summary}~\ref{Prop_ancient_summary_f} follows from
Theorem~\ref{thm:HMHF}\ref{thm:HMHF_c}\ref{thm:HMHF_e} with
$T_1=-\infty$. The projection defining $U^+$ preserves the
components of $U^{+,\prime}$ in $\sV^{++}$ with positive eigenvalues of $L$, so the bound
$\|U^+\|\leq C\|U^{++}_{>0}\|$ also follows from
Assertion~\ref{thm:HMHF_e}.  On the interval where
$-U_{0,\min}\geq U_{0,\max}$, all eigenvalues of $(U_{0,ij})$ have
absolute value at most $|U_{0,\min}|$.  Hence
$\|U^{++}_0\|\leq C|U_{0,\min}|$, and the first bound in
Theorem~\ref{thm:HMHF}\ref{thm:HMHF_c} gives
$\|U^+\|\leq C|U_{0,\min}|$.
Since
\[
\|U^{++}_{>0}(\tau)\|
\leq
e^{0.49(\tau-\td\tau)}
\|U^{++}_{>0}(\td\tau)\|,
\qquad \tau\in(T_3,\td\tau],
\]
the integration argument in Claim~\ref{Cl_Upp_quad_1tau} on $(-\infty,T_3]$, combined with this and \eqref{eq_Up10Uppproof}, gives the coarse bound in Proposition~\ref{Prop_ancient_summary}~\ref{Prop_ancient_summary_b2}
\[
\|U^+(\tau)\|
\leq C(\td\tau-\tau+10)^{-1}.
\]

Next, we prove Proposition~\ref{Prop_ancient_summary}~\ref{Prop_ancient_summary_c}.

\begin{Claim} \label{Cl_upUP}
If $R^* \geq \underline R^*(J', \eta')$, then for any $\tau \leq  \td\tau$ we have $\IB^k_{2R(\tau)} \times (\IS^{n-k}/\Gamma) \subset \DD'_\tau$ and
\begin{align}
 \Vert u_\tau - U^+(\tau) \Vert_{C^{m}(\IB^k_{R(\tau)} \times (\IS^{n-k}/\Gamma))} &\leq C(J,m) \Vert U^+ (\tau) \Vert^{J+1} \label{eq_upmUp} \\
  \Vert U^+(\tau) \Vert_{C^{m}(\IB^k_{R(\tau)} \times (\IS^{n-k}/\Gamma))} &\leq C(J,m) \sqrt{\eta'} \label{eq_Upsqrt}
\end{align}
\end{Claim}

\begin{proof}
Fix $\tau \leq \td\tau$.
Recall that by the statement of the proposition we have
$e^{-R^2(\tau)} := \Vert U^+ (\tau) \Vert^J$.  If $R_{++}(\tau)$ is
defined by $e^{-R_{++}^2(\tau)}=\|U^{++}(\tau)\|^J$, then
$R(\tau)\leq R_{++}(\tau)$ because $\|U^{++}\|\leq\|U^+\|$.
The bound on $\UU^-$ from \eqref{eq:small_1} and $J'-1\ge 4J$ gives
$R'(\tau)\geq2R_{++}(\tau)$ whenever
$R(\tau) \geq R^* \geq \underline R^*(J', \eta')$.  Therefore,
\[ R'(\tau) \geq 2 R(\tau). \]
Therefore if we write $B_a := \IB^k_{a R(\tau)} \times (\IS^{n-k}/\Gamma)$, then
\[ u_{\tau} = u_{\tau} \omega_{R'(\tau)} = U^{+,\prime} + U^{-,\prime}  \qquad \text{on} \quad B_2. \]
The bound \eqref{eq:small_1} implies that
\[ \Vert U^{-,\prime} \Vert_{L^2_f (B_2)}  \leq C( J) \Vert U^{++} \Vert^{2(n+m+4) (J+1)}. \]
So since $e^{-f} \geq e^{-(4R(\tau))^2/4}$ on $B_2$, we obtain
\begin{multline*}
 \Vert U^{-,\prime} \Vert_{L^2(B_2)} = \bigg( \int_{B_2} |U^{-,\prime}|^2 dg_{\cyl}  \bigg)^{1/2} 
\leq e^{(2R(\tau))^2/8} \bigg( \int_{B_2} |U^{-,\prime}|^2 e^{-f} dg_{\cyl}  \bigg)^{1/2}  \\
\leq C(J) \Vert U^+ \Vert^{-J/2} \Vert U^{++} \Vert^{2(n+m+4) (J+1)}
\leq C(J)  \Vert U^{++} \Vert^{(n+m+4) (J+1)}.
\end{multline*}
On the other hand, using $\dim\sV_{>\la'}\le C(J)$ and \eqref{eq_Up10Uppproof}, we find that (recall $\la' = - J' - 1$)
\begin{multline} \label{eq_Upp_pol}
   \Vert U^{+,\prime} \Vert_{C^{m+1}(B_2)} 
\leq C(J,m) R^{C(J,m)}(\tau) \Vert U^{+,\prime} \Vert_{L^2_f} 
\leq C(J,m) R^{C(J,m)}(\tau) \Vert U^+ \Vert_{L^2_f} \\
= C(J,m) R^{C(J,m)}(\tau) e^{-\frac{R^2(\tau)}{J}}
\leq C( J,m) \sqrt{\eta'}.
\end{multline}
In the same way we can bound $U^+ = \PP_{\sV_{\rot,\geq - J}} U^{+,\prime}$, proving \eqref{eq_Upsqrt}, and we obtain
\[  \Vert U^{-,\prime} \Vert_{C^{m+1}(B_2)} 
\leq \Vert u_\tau \Vert_{C^{m+1}(B_2)} + \Vert U^{+,\prime} \Vert_{C^{m+1}(B_2)} 
\leq \eta' + C(J,m) \sqrt{\eta'}
\leq C( J,m). \]
So using the Gagliardo--Nirenberg inequality from \cite[Lemma~5.107]{Bamler_Lai_MCF1}, we obtain that
\begin{equation} \label{eq_Ump_Cmp1}
    \Vert U^{-,\prime} \Vert_{C^{m}(B_{2})}
\leq C( J, m) \Vert U^{-,\prime} \Vert_{L^2(B_2)}^{\frac1{n/2+m+1}}  
\leq C(J,m)  \Vert U^+ \Vert^{J+1} . 
\end{equation}

Now as in \eqref{eq_Upp_pol} we can bound, using \eqref{eq_UppUpJp2},
\begin{multline*} 
    \Vert  U^{+,\prime} - U^+ \Vert_{C^{m}(B_1)}
\leq C(J,m) R^{C(J,m)}(\tau) \Vert  U^{+,\prime} - U^+ \Vert \\
\leq C(J,m) R^{C(J,m)}(\tau) \Vert U^+ \Vert^{J+2} 
\leq C(J,m) \Vert U^+ \Vert^{J+1},
\end{multline*}
where we used $R(\tau)^2=-J\log\|U^+(\tau)\|$ to absorb the polynomial factor by $ \Vert U^+ \Vert$.
The same bound holds for $\PP_{\sV_{\rot,< -J}} U^+$.
Combining this with \eqref{eq_Ump_Cmp1} implies
\begin{equation*}
 \Vert u_\tau - U^{+} \Vert_{C^{m}(B_1)}
\leq \Vert U^{+,\prime} - U^+ \Vert_{C^{m}(B_1)} + \Vert U^{-,\prime} \Vert_{C^{m}(B_1)}  
\leq  C( J,m) \Vert U^+ \Vert^{J+1}.
\end{equation*}
This proves the claim.
\end{proof}

The first estimate in
Proposition~\ref{Prop_ancient_summary}~\ref{Prop_ancient_summary_c} follows from \eqref{eq_upmUp}.
The second estimate of the $\eta$-smallness follows from
\eqref{eq_upmUp} and \eqref{eq_Upsqrt} by taking $\eta'\le\ov\eta'(J,m,\eta)$.
This proves Assertions~\ref{Prop_ancient_summary_b}--\ref{Prop_ancient_summary_f}.
\medskip

Lastly, suppose now that $q_0$ survives on $(-\infty,t_0)$ and remains a neck center
at scale $\sqrt{t_0-t}$ for every $t<t_0$. We may then apply
Theorem~\ref{thm:HMHF} on the whole interval $(-\infty,t_0)$, whose rescaled
time interval is the entire $\IR$. So Proposition~\ref{Prop_mode_analysis}
\ref{Prop_mode_analysis_f} gives
\[
\widehat U^+,\widehat{\mathcal U}^-\equiv0.
\]
Hence $\widehat R\equiv\infty$ and $u_\tau\equiv0$ for every $\tau$.
Therefore $\MM$ is homothetic to the round shrinking cylinder, proving
Proposition~\ref{Prop_ancient_summary}~\ref{Prop_ancient_summary_g}.
This finishes the proof of the Proposition.
\end{proof}

The following corollary applies Proposition~\ref{Prop_ancient_summary} to asymptotically cylindrical Ricci flow spacetimes.
It shows that, along a worldline realizing the cylindrical limit as $t\to-\infty$, the first loss of a sufficiently good neck can be detected on a ball of fixed radius at the cylindrical scale.
We use this observation in the proof of Proposition~\ref{prop:suggested_gauge_geometric_comparison}.

\begin{Corollary}[Local detection of neck loss]\label{cor:local_neck_detection}
If $\delta>0$ and $\delta'\leq\ov\delta'(\delta)$, then the following holds.
Let $\MM$ be an asymptotically $(n,k,\Gamma)$-cylindrical Ricci flow spacetime satisfying the pseudolocality property.
Fix $t_0\in\IR$ and $r_0>0$, and let $q_0\in\MM_{t_0-r_0^2}$ have a backward worldline $q_0(t)$ along which the flow converges to the cylinder in the sense of~\eqref{eq:asymptotically_cylindrical_convergence}.
Suppose that for every $t\leq t_0-r_0^2$ there is a diffeomorphism
\[
\psi_t:B_{q_0(t)}(2\sqrt{t_0-t})\longrightarrow V_t\subset M_{\cyl,\Gamma}
\]
onto an open subset such that
\begin{equation}\label{eq:local_neck_detection_comparison}
\big\|(t_0-t)^{-1}(\psi_t)_*g_t-g_{\cyl,\Gamma}\big\|_{C^2(V_t)}<\delta'.
\end{equation}
Then $q_0(t)$ is a center of an $(n,k,\delta,\Gamma)$-neck at scale $\sqrt{t_0-t}$ for every $t\leq t_0-r_0^2$.
\end{Corollary}

\begin{proof}
Fix $m=20$ and $J=10$, and let $\eta>0$ be chosen in the course of the proof.
Without loss of generality, we may assume that $\delta$ is chosen small enough so that Proposition~\ref{Prop_ancient_summary} is applicable with these parameters.
By the assumed pointed cylindrical convergence and $\sqrt{(t_0-t)/(-t)}\to1$, the point $q_0(t)$ is a center of an $(n,k,\frac12\delta,\Gamma)$-neck at scale $\sqrt{t_0-t}$ for all sufficiently negative $t$.

If $\MM$ is homothetic to the round shrinking cylinder, the conclusion follows from homogeneity and \eqref{eq:local_neck_detection_comparison}.
Otherwise, let $t^*\leq t_0-r_0^2$ be maximal such that $q_0(t)$ is a center of an $(n,k,\frac12\delta,\Gamma)$-neck at scale $\sqrt{t_0-t}$ for all $t<t^*$.
Thus Proposition~\ref{Prop_ancient_summary} applies on $(-\infty,t^*]$ with $r^*=\sqrt{t_0-t^*}$ and gives $\chi$, $U^+$, and $R$ with basepoint $q_0(t^*)$.
Write $\tau=-\log(t_0-t)$, let $u_\tau$ be the normalized metric perturbation in Assertion~\ref{Prop_ancient_summary_c}, and set $U^{++}(\tau):=\PP_{\sV^{++}}U^+(\tau)$.
Assertions~\ref{Prop_ancient_summary_b2}, \ref{Prop_ancient_summary_c}, and \ref{Prop_ancient_summary_e} imply
\[
\big\|u_\tau-U^{++}(\tau)-U^+_{\aff}(\tau)\big\|_{C^2(\IB^k_1\times(\IS^{n-k}/\Gamma))}
\leq C\eta\|U^{++}(\tau)\|.
\]
If $\eta\leq\ov\eta$ and $\delta'\leq\ov\delta'(\delta)$, then Assertion~\ref{Prop_ancient_summary_a} and \eqref{eq:local_neck_detection_comparison} allow us to apply Lemma~\ref{l:dom_curvature}, giving $\|U^{++}(\tau)\|\leq C\delta'$.
Assertion~\ref{Prop_ancient_summary_b2} then gives $\|U^+(\tau)\|\leq C\delta'$ for all $\tau\leq-\log(t_0-t^*)$.

By Assertion~\ref{Prop_ancient_summary_c}, $R^2(\tau)=-J\log\|U^+(\tau)\|$, so a bound of the form $\delta'\leq\ov\delta'(\delta)$ ensures $R(\tau)\geq4\delta^{-1}$ for all $\tau \leq -2\log r^*$.
Assertion~\ref{Prop_ancient_summary_c} and equivalence of norms on the finite-dimensional mode space then give
\[
\|u_\tau\|_{C^{20}(\IB^k_{4\delta^{-1}}\times(\IS^{n-k}/\Gamma))}\leq C(\delta)\delta'.
\]
Using local derivative estimates, this bound can be upgraded to higher derivative bounds on slightly smaller balls.
Together with Assertion~\ref{Prop_ancient_summary_a}, this implies that $q_0(t^*)$ is a center of an $(n,k,\frac14\delta,\Gamma)$-neck.
So by the maximal choice of $t^*$ we must have $t^* = t_0 - r_0^2$.
\end{proof}
\medskip

We will frequently consider HMHF data satisfying the conclusions of
Proposition~\ref{Prop_ancient_summary}, so we introduce the following
terminology.

\begin{Definition}[Admissible HMHF data]
\label{defn:admissible_HMHF}
Let $\chi:\Omega\to\IR^k\times(\IS^{n-k}/\Gamma)$ be a harmonic
map heat flow on an open subset $\Omega$ of a Ricci flow
spacetime $\MM$, whose time-slice maps $\chi_t$ are diffeomorphisms
onto their images.
We call $\chi$ \textbf{admissible} for parameters $m,J,\eta,C$
and reference time $t_0$ on the rescaled time-interval
$(-\infty,\td\tau]$ if there exist a locally Lipschitz function
$U^+:(-\infty,\td\tau]\to\sV_{\rot,\geq-J}$ and a function
$R:(-\infty,\td\tau]\to[1000,\infty]$ satisfying the assertions of
Proposition~\ref{Prop_ancient_summary}
\ref{Prop_ancient_summary_b}--\ref{Prop_ancient_summary_f},
with $\tau=-\log(t_0-t)$ and with the following modifications:
the identity defining $R$ in
assertion~\ref{Prop_ancient_summary_c} is replaced by
$R(\tau)\to\infty$ as $\tau\to-\infty$, and every constant
denoted by $C$, $C(J)$, or $C(J,m)$ is replaced by the
prescribed constant $C$.
We call $(\chi,\Omega,R,U^+)$ the corresponding
\textbf{admissible HMHF data.}
\end{Definition}

We introduce the following definition to distinguish the two possible asymptotic behaviors as $\tau\to-\infty$ characterized in Proposition~\ref{Prop_ancient_summary}~\ref{Prop_ancient_summary_f}: either the positive modes dominate and decay exponentially, or the
neutral quadratic mode dominates and decays at order $|\tau|^{-1}$.

\begin{Definition}\label{defn:dominant_linear}
Let $\MM$ be an ancient asymptotically $(n,k,\Gamma)$-cylindrical
Ricci flow spacetime satisfying the pseudolocality property.
We say that $\MM$ has \textbf{dominant quadratic mode} if for some
admissible HMHF data $(\chi,\Omega,R,U^+)$ with parameters
$m,J\geq10$, $0<\eta<0.1$, and $C>0$, the transition time $T_3$
in Proposition~\ref{Prop_ancient_summary}\ref{Prop_ancient_summary_f}
satisfies $T_3>-\infty$,
with the convention that $T_3=-\infty$ if $U^+\equiv0$.
We say that $\MM$ has \textbf{dominant linear mode} if it does not
have dominant quadratic mode and is not homothetic to a round
shrinking cylinder.
\end{Definition}

The following lemma shows that the alternative in
Definition~\ref{defn:dominant_linear} is independent of the choice
of admissible HMHF data, parameters, and reference times.

\begin{Lemma}[Independence of the dominant mode]\label{lem:dominant_mode_independence}
Let $\MM$ be an ancient asymptotically $(n,k,\Gamma)$-cylindrical
Ricci flow spacetime satisfying the pseudolocality property.

For $i=1,2$, let $(\chi^{(i)},\Omega^{(i)},R^{(i)},U^{+,(i)})$ be
admissible HMHF data for $\MM$ with parameters $m_i,J_i\geq10$,
$0<\eta_i<0.1$, and $C_i>0$, and reference time $t_0^{(i)}$.
Let $T_3^{(i)}$ be the transition times in
Proposition~\ref{Prop_ancient_summary}\ref{Prop_ancient_summary_f},
with $T_3^{(i)}=-\infty$ when $U^{+,(i)}\equiv0$. Then
\[
T_3^{(1)}=-\infty\quad\Longleftrightarrow\quad T_3^{(2)}=-\infty.
\]
If $\MM$ is homothetic to a round shrinking cylinder, then both
transition times equal $-\infty$.

Moreover, if $\MM$ is not homothetic to a round shrinking cylinder,
the dominant linear mode condition is geometrically characterized by
the requirement that, along some worldline $p(t)$ satisfying
\eqref{eq:asymptotically_cylindrical_convergence} (or equivalently along
every such worldline), the $k$ smallest eigenvalues of the Ricci curvature
$\lambda_1(t),\ldots,\lambda_k(t)$ satisfy
$(-t)\log(-t)|\lambda_i(t)|\to0$ as $t\to-\infty$ for all $i=1,\ldots,k$.
\end{Lemma}

\begin{proof}
Suppose, by way of contradiction and after interchanging the two
data sets if necessary, that $T_3^{(1)}=-\infty$ and
$T_3^{(2)}> -\infty$. Set
$\tau^{(i)}(t):=-\log(t_0^{(i)}-t)$.

For $i=1,2$, choose $q^{(i)}$ on the inverse image of
$\{0\}\times(\IS^{n-k}/\Gamma)$ under the rescaled HMHF chart at a
sufficiently negative time, and denote its backward worldline by
$q^{(i)}(t)$. In either rescaled chart, worldlines are integral curves of the
vector field $\tfrac12\bx+o(1)$, where the error tends to zero
uniformly on fixed cylindrical regions as $\tau\to-\infty$.
Thus $\IB^k_1\times(\IS^{n-k}/\Gamma)$ is invariant backward
in time, so these worldlines
exist for all sufficiently negative $t$ and satisfy
$\xi_{e^{\tau^{(i)}(t)/2}}\circ\chi_t^{(i)}(q^{(i)}(t))
\in\IB^k_1\times(\IS^{n-k}/\Gamma)$.
Their spherical components have integrable speed
$O(|\tau|^{-J_i-1})$ by
Proposition~\ref{Prop_ancient_summary}\ref{Prop_ancient_summary_b2},
\ref{Prop_ancient_summary_c}, so these worldlines satisfy
\eqref{eq:asymptotically_cylindrical_convergence}.

Let $\lambda_1(x,t)\leq\ldots\leq\lambda_n(x,t)$ denote the
eigenvalues of the Ricci tensor of $g_t$ at $x\in\MM_t$.
By Proposition~\ref{Prop_ancient_summary}\ref{Prop_ancient_summary_f},
$\|U^{+,(1)}(\tau^{(1)})\|=O(e^{0.49\tau^{(1)}})$.
The metric approximation in
Proposition~\ref{Prop_ancient_summary}\ref{Prop_ancient_summary_c}
therefore gives the same decay in $C^2$ on every fixed rescaled
cylindrical region.
For any fixed $A>0$, the bound just established and
Proposition~\ref{Prop_ancient_summary}\ref{Prop_ancient_summary_c}
imply that the rescaled map $\xi_{e^{\tau^{(1)}/2}}\circ\chi_t^{(1)}$ sends
$B_{g_t}(q^{(1)}(t),A\sqrt{-t})$ into a fixed such region for all
sufficiently negative $t$.
Since the cylinder has $k$ zero Ricci eigenvalues, we obtain
\begin{equation}\label{eq:dominant_mode_linear_Ricci_decay}
\sup_{x\in B_{g_t}(q^{(1)}(t),A\sqrt{-t})}
\max_{1\leq j\leq k}|\lambda_j(x,t)|
=O((-t)^{-1.49})
\qquad\text{as }t\to-\infty.
\end{equation}

By Definition~\ref{defn:asymptotically_cylindrical}, for
sufficiently large $A$ and all sufficiently negative $t$, we have
$q^{(2)}(t)\in B_{g_t}(q^{(1)}(t),A\sqrt{-t})$.
Thus \eqref{eq:dominant_mode_linear_Ricci_decay} applies at $q^{(2)}(t)$.

For the second data set, write $\tau=\tau^{(2)}$ and recall
$U_{0,\min}^{(2)}$ from Proposition~\ref{Prop_ancient_summary}.
It is negative for sufficiently negative $\tau$, since
$U^{+,(2)}\not\equiv0$.
Its differential inequality in
Proposition~\ref{Prop_ancient_summary}\ref{Prop_ancient_summary_f}
implies $|\partial_\tau(U_{0,\min}^{(2)})^{-1}-\sqrt2|\leq0.1$.
Integrating gives, for sufficiently negative $\tau$,
\[
\frac{c}{|\tau|}\leq-U_{0,\min}^{(2)}(\tau)\leq\frac{C}{|\tau|}.
\]
Since this mode dominates for sufficiently negative $\tau$,
$\|U^{++,(2)}\|\leq\|U^{+,(2)}\|=O(|\tau|^{-1})$, and
Proposition~\ref{Prop_ancient_summary}\ref{Prop_ancient_summary_e}
gives $\|U_{\rot,<0}^{+,(2)}\|\leq C\|U^{++,(2)}\|^2=O(|\tau|^{-2})$.

Let $\td g_2(\tau):=e^\tau(\xi_{e^{\tau/2}})_*(\chi_t^{(2)})_*g_t$,
where $t=t_0^{(2)}-e^{-\tau}$.
The metric approximation in
Proposition~\ref{Prop_ancient_summary}\ref{Prop_ancient_summary_c}
and the computation in the proof of Lemma~\ref{l:dom_curvature}
give, on every fixed cylindrical region in a
$g_{\cyl}$-orthonormal product frame,
\begin{align*}
(\Ric_{\td g_2})_{ij}
&=-\frac{n-k}{2\sqrt2}U_{0,ij}^{(2)}(\tau)+O(|\tau|^{-2}),\\
(\Ric_{\td g_2})_{i\alpha}
&=O(|\tau|^{-2}),\\
(\Ric_{\td g_2})_{\alpha\beta}
&=\bigg(\frac12-\frac1{2\sqrt2}
\sum_{l=1}^k U_{0,ll}^{(2)}(\tau)\bigg)
\delta_{\alpha\beta}+O(|\tau|^{-2}),
\end{align*}
where $i,j$ denote Euclidean and $\alpha,\beta$ spherical directions.
Combining these formulas with the bound for $U_{0,\min}^{(2)}$
and undoing the rescaling, we obtain
\[
\max_{1\leq j\leq k}|\lambda_j(q^{(2)}(t),t)|
\geq\frac{c}{(-t)\log(-t)}
\]
for all sufficiently negative $t$, contradicting
\eqref{eq:dominant_mode_linear_Ricci_decay}.

For the final assertion, apply Proposition~\ref{Prop_ancient_summary}
along any worldline $p(t)$ satisfying
\eqref{eq:asymptotically_cylindrical_convergence}, with reference time
zero and a sufficiently negative terminal time.
If $T_3=-\infty$, the argument leading to
\eqref{eq:dominant_mode_linear_Ricci_decay} gives
$\max_{1\leq i\leq k}|\lambda_i(t)|=O((-t)^{-1.49})$.
If $T_3>-\infty$, the quadratic-mode computation above gives
$\max_{1\leq i\leq k}|\lambda_i(t)|\geq c/((-t)\log(-t))$
for all sufficiently negative $t$.
The characterization therefore follows from Definition~\ref{defn:dominant_linear}.
\end{proof}

We study the dominant linear case in Section \ref{sec:dominant_linear}, and the dominant quadratic case in Section \ref{sec:dominant_quadratic}.

\section{Flows with dominant linear mode}\label{sec:dominant_linear}
\subsection{Statement of the main results}\label{subsec:linear_main_results}
In this section, we study ancient asymptotically cylindrical Ricci flow spacetimes with dominant linear mode.
Our main goal is to describe the asymptotic geometry of their cylindrical regions by comparison with the Bryant soliton times a Euclidean factor, $\IR^{k-1}\times\MM^{n-k+1}_{\Bry}$, or a quotient of this model.
We do this by constructing a time- and $\partial_{\mathfrak t}$-preserving comparison map (see~\cite[Definitions~6.10 and~6.12]{Bamler_Kleiner_uniqueness}) under which the metric difference and all its derivatives decay exponentially in the distance from the tip of the Bryant factor.

Our main result is the following spacetime version of Theorem~\ref{thm:intro_bryant_infinity}.

\begin{Theorem}[Exponential Bryant asymptotics]\label{thm:linear_bryant_asymptotics}
Let $\MM$ be an ancient Ricci flow spacetime defined over the time-interval $(-\infty,T)$, where $T\in\IR\cup\{\infty\}$.
Suppose that $\MM$ is asymptotically $(n,k,\Gamma)$-cylindrical, satisfies the pseudolocality property, and has dominant linear mode in the sense of Definition~\ref{defn:dominant_linear}.
Let
\[
 \MM'=\IR^{k-1}\times\big(\MM^{n-k+1}_{\Bry}/\Gamma\big),
\]
where $\Gamma$ acts by rotations on the Bryant factor.
Write $g'$ for the evolving metric of $\MM'$ and let $y:\MM'\to[0,\infty)$ denote the time-slice distance to the tip in the second factor.
Then there is $Y_0=Y_0(n,k,\Gamma)>0$ such that, after parabolically rescaling $\MM$, there is a smooth diffeomorphism onto its image
\[
 \psi:\MM'_{<T,\,y>Y_0}\longrightarrow\MM
\]
with
\[
 \mathfrak t\circ\psi=\mathfrak t',
 \qquad
 d\psi(\partial_{\mathfrak t'})=\partial_{\mathfrak t}\circ\psi.
\]
For every integer $m\geq0$, there are constants $C_m,c_m>0$, depending only on $n,k,\Gamma,m$, such that
\begin{equation}\label{eq:linear_bryant_asymptotics}
 \big|\nabla_{g'}^m(\psi^*g-g')\big|_{g'}
 \leq C_m e^{-c_m y}
 \qquad\textnormal{on }\MM'_{<T,\,y>Y_0}.
\end{equation}
Moreover, for sufficiently small $\eps=\eps(n,k,\Gamma)>0$, the image of $\psi$ contains every point $p\in\MM$ for which there is $t_p>\mathfrak t(p)$ such that its backward worldline is defined for all $t\leq\mathfrak t(p)$ and $p(t)$ is the center of an $(n,k,\eps,\Gamma)$-neck at scale $\sqrt{t_p-t}$ for every such $t$.
\end{Theorem}

For $\Gamma\ne\{1\}$ the tip set is singular; the comparison is on the smooth end.

If the time-slices of $\MM$ are complete and its curvature operator is nonnegative, we even obtain a global classification.

\begin{Corollary}\label{cor:complete_linear_bryant}
Under the hypotheses of Theorem~\ref{thm:linear_bryant_asymptotics}, suppose in addition that the time-slices of $\MM$ are complete and that $\Rm\geq0$.
Then $\Gamma=\{1\}$ and $\MM$ is homothetic to a restriction of $\IR^{k-1}\times\MM^{n-k+1}_{\Bry}$.
\end{Corollary}

In particular, we obtain the following classification of $\kappa$-solutions, stated in terms of conventional Ricci flows.

\begin{Corollary}\label{thm:unique_bryant}
Let $(M,(g_t)_{t<T})$ be an ancient asymptotically $(n,k)$-cylindrical $\kappa$-solution with dominant linear mode.
Then this flow is homothetic to the restriction of $\IR^{k-1}\times\MM^{n-k+1}_{\Bry}$ to the time-interval $(-\infty,T)$.
\end{Corollary}

We now describe the argument and the structure of this section.
In Subsection~\ref{subsec:linear_asymptotics}, we first study the asymptotics of the modes $U^+$ from Proposition~\ref{Prop_ancient_summary} in the Ricci--DeTurck gauge given by an admissible harmonic map heat flow $\chi$.
Proposition~\ref{p:coefficient} identifies two asymptotic constants, a scalar $a$ and a vector $\mathbf b$, and Proposition~\ref{prop:diff_of_U} shows that they determine $U^+$ up to an error of arbitrarily high exponential order as the rescaled time $\tau$ tends to $-\infty$.
This parallels the asymptotic expansion and mode comparison in our mean curvature flow work; see~\cite[Proposition~6.2 and the proof of Lemma~6.29]{Bamler_Lai_MCF1}.
In Subsection~\ref{subsec:linear_transformations}, we study admissibility and the asymptotic constants under rotations, spatial translations, and parabolic rescalings of the HMHF data, analogously to~\cite[Propositions~6.4 and~6.6]{Bamler_Lai_MCF1}.

Subsection~\ref{subsec:linear_gauges} addresses an additional difficulty in the Ricci flow setting.
The HMHF gauge depends on time and on the choice of initial data, so it does not directly provide a comparison map preserving the time vector fields.
We therefore introduce a \emph{gauge at infinity}, which is constant along worldlines and describes the cylindrical asymptotics as $t\to-\infty$.
Such a gauge is unique up to an isometry of the cylindrical model; see Lemma~\ref{l:isometry_gauge}.
To construct it, we compensate for the drift generated by the linear mode by a time-dependent translation and take a limit along backward worldlines, as in Proposition~\ref{prop:suggested_gauge_geometric}.

In Subsection~\ref{subsec:linear_bryant_comparison}, we match the asymptotic constants with those of the Bryant model and use the gauges at infinity to construct the comparison map $\psi$.
Proposition~\ref{p:bryant_infinity} gives decay of arbitrarily high polynomial order in the spatial variable $y$.
Along backward worldlines of the Bryant model, this corresponds to decay of arbitrarily high exponential order in the rescaled time.
These polynomial estimates already suffice for the classification results in Corollaries~\ref{cor:complete_linear_bryant} and~\ref{thm:unique_bryant}, which we prove in Subsection~\ref{subsec:linear_classification} using comparison geometry and splitting.
Finally, in Subsection~\ref{subsec:linear_exponential}, we improve the polynomial estimates to exponential decay in $y$ and conclude the proof of Theorem~\ref{thm:linear_bryant_asymptotics}.

\subsection{Asymptotics at $-\infty$}\label{subsec:linear_asymptotics}
We first derive an asymptotic expansion of the rescaled metric perturbation in HMHF coordinates up to an error of order $e^{1.2\tau}$. This identifies a scalar $a$ and a vector $\mathbf b$,
which will be used later to compare different harmonic map heat flows.
In particular, Proposition~\ref{prop:diff_of_U} shows that two such data sets with the same $a,\mathbf b$ have modes agreeing to any prescribed order in $\|U^+\|$.

\begin{Proposition}\label{p:coefficient}
For every $A\geq1000$, integers $J,m_0\geq10$ and $m\geq0$, and $C_0>0$, the following holds.
Let $\MM$ be an ancient asymptotically $(n,k,\Gamma)$-cylindrical
Ricci flow spacetime satisfying the pseudolocality property and
without dominant quadratic mode.
Let $\chi:\Omega\to\IR^k\times(\IS^{n-k}/\Gamma)$ be an admissible
harmonic map heat flow on $\MM$, with parameters $m_0,J,\eta,C_0$ and
reference time $t_0$ on $(-\infty,\td\tau]$, where $\td\tau\in\IR$
and $0<\eta<0.1$.
Choose corresponding admissible data $(\chi,\Omega,R,U^+)$ with $R$ nonincreasing in $\tau$.
For $t=t_0-e^{-\tau}$, let
$u_\tau:=e^\tau(\xi_{e^{\tau/2}})_*(\chi_t)_*g_t-g_{\cyl}$
be the rescaled metric perturbation in the HMHF coordinates.
Then there are unique constants $a\in\IR$ and
$\mathbf b=(b_1,\ldots,b_k)\in\IR^k$ such that, whenever
$\tau\leq\td\tau$ and $R(\tau)\geq A+10$,
\begin{align}
u_\tau
=&
\tfrac{n-k}{4}e^\tau
\sum_{i,j=1}^k
b_ib_j\,d\bx_i\otimes d\bx_j
+
\Big(
\big(
a-\tfrac12|\mathbf b|^2\,\tau
\big)e^\tau\mathfrak p^{(0)}
+
\sum_{i=1}^k
b_i e^{\tau/2}\mathfrak p_i^{(1)}
\Big)g_{\IS^{n-k}/\Gamma}
+
\mathcal E_\tau,
\label{eq:coefficient_metric_expansion}
\end{align}
where
\[
\|\mathcal E_\tau\|_{C^m(\IB^k_A\times(\IS^{n-k}/\Gamma))}
\leq C(A,m,J,C_0)\eta e^{1.2(\tau-\td\tau)}.
\]
Moreover,
\begin{equation}\label{eq:coefficient_asymptotic_constants_bound}
e^{\td\tau/2}|\mathbf b|
+e^{\td\tau}\bigl|a-\tfrac12|\mathbf b|^2\td\tau\bigr|
\leq C(J,C_0)\eta.
\end{equation}
If $a=0$ and $\mathbf b=0$, then $\MM$ is isometric to the round
shrinking cylinder $\IR^k\times(\IS^{n-k}/\Gamma)$ with extinction
time $t_0$, restricted to the time-interval of $\MM$.
\end{Proposition}

We can therefore associate the following asymptotic data to an admissible harmonic map heat flow.

\begin{Definition}[Asymptotic constants]\label{defn:HMHF_asymptotic_constants}
We call the constants $a,\mathbf b$ in
Proposition~\ref{p:coefficient} the \textbf{asymptotic constants}
of $\chi$ with reference time $t_0$, and write
$a=a(\chi,t_0)$ and $\mathbf b=\mathbf b(\chi,t_0)$.
\end{Definition}

In fact, we may simply write $\mathbf b=\mathbf b(\chi)$,
since the following lemma shows that $\mathbf b$ is independent
of the reference time.
We will later use the fact that $\mathbf b(\chi)\ne0$ for flows
with dominant linear mode in our asymptotic characterization of such flows.

\begin{Lemma}[Change of reference time]\label{lem:HMHF_constants_reference_time}\label{l:dominant_linear_well_defined}
Suppose that $\chi$ satisfies the hypotheses of
Proposition~\ref{p:coefficient} with reference time $t_0$.
If $\chi$ also satisfies these hypotheses with reference time $t_1$,
then
\[
\mathbf b(\chi,t_1)=\mathbf b(\chi,t_0),
\qquad
a(\chi,t_1)=a(\chi,t_0)+t_0-t_1.
\]
Moreover, if $\MM$ has dominant linear mode, then $\mathbf b(\chi)\ne0$.
\end{Lemma}

The following proposition shows that these constants determine the asymptotics of $U^+$ to arbitrarily high exponential order, provided $J$ is chosen sufficiently large.

\begin{Proposition}\label{prop:diff_of_U}
For every $J\geq10$ and $C_0>0$, there are constants
$\ov\eta=\ov\eta(J,C_0)\in(0,0.1)$ and $C_1=C_1(J,C_0)>0$
such that the following holds.
For $i=0,1$, let $\MM^{(i)}$ be an ancient asymptotically
$(n,k,\Gamma)$-cylindrical Ricci flow spacetime satisfying the
pseudolocality property and without dominant quadratic mode.
Let $\chi^{(i)}:\Omega^{(i)}\to\IR^k\times(\IS^{n-k}/\Gamma)$
be an admissible harmonic map heat flow on $\MM^{(i)}$, with
parameters $10,J,\eta,C_0$ and reference time $t_{0,i}$ on the
common rescaled time-interval $(-\infty,\td\tau]$, where
$\td\tau\in\IR$ and $0<\eta\leq\ov\eta$.
Write $(\chi^{(i)},\Omega^{(i)},R^{(i)},U^{+,(i)})$ for the
corresponding admissible data, with $\tau=-\log(t_{0,i}-t)$.
If the asymptotic constants from Proposition~\ref{p:coefficient} satisfy
\[
a(\chi^{(0)},t_{0,0})=a(\chi^{(1)},t_{0,1}),
\qquad
\mathbf b(\chi^{(0)})=\mathbf b(\chi^{(1)}),
\]
then, for all $\tau\leq\td\tau$,
\[
\big\|U^{+,(1)}(\tau)-U^{+,(0)}(\tau)\big\|
\leq C_1\big\|U^{+,(0)}(\tau)\big\|^{J+1},
\]
where $\|\cdot\|$ denotes the $L^2_f$-norm.
\end{Proposition}

\medskip

We now prove these results in the order stated.

\begin{proof}[Proof of Proposition~\ref{p:coefficient}]
Given any admissible radius $R_0$, the function $R(\tau):=\inf_{s\leq\tau}R_0(s)$ is nonincreasing and tends to infinity as $\tau\to-\infty$.
Shrinking the radius preserves the admissibility assertions, so the choice of $R$ in the statement is always possible.
By Lemma~\ref{lem:dominant_mode_independence}, we have $T_3=-\infty$.
Proposition~\ref{Prop_ancient_summary}\ref{Prop_ancient_summary_b2},
\ref{Prop_ancient_summary_c}, \ref{Prop_ancient_summary_f2} therefore gives
\begin{equation}\label{eq:coefficient_coarse_decay}
\|U^{++}(\tau)\|\leq\|U^+(\tau)\|
\leq C(C_0)\eta e^{0.49(\tau-\td\tau)}
\qquad\text{for all }\tau\leq\td\tau.
\end{equation}
By Proposition~\ref{Prop_ancient_summary}\ref{Prop_ancient_summary_c},
using $J\geq10$ and $\|U^+\|\leq\eta<0.1$, we also have
\[
\|u_\tau-U^+(\tau)\|_{C^{m_0}(\IB^k_{R(\tau)}\times(\IS^{n-k}/\Gamma))}
\leq C_0\|U^+\|^{J+1}
\leq C_0\|U^+\|^3
\leq C(C_0)\eta e^{1.2(\tau-\td\tau)}.
\]
We first establish the expansion for $U^+$ and the coefficient bounds.

By \eqref{eq:coefficient_coarse_decay} and the ODI in
Proposition~\ref{Prop_ancient_summary}\ref{Prop_ancient_summary_d},
\[
\partial_\tau U^+-LU^+=O(\|U^+\|^2)
=O(\eta^2e^{0.98(\tau-\td\tau)}).
\]
Integrating the $\frac12$-mode equation gives constants $b_i\in\IR$ such that
\[
\lim_{\tau\to-\infty}e^{-\frac12\tau}
\PP_{\sV_{\rot,\frac12}}U^+(\tau)
=\sum_{i=1}^k b_i\mathfrak p_i^{(1)}g_{\IS^{n-k}},
\qquad |\mathbf b|\leq C(C_0)\eta e^{-\frac12\td\tau}.
\]
Integrating the remaining components from $-\infty$ for eigenvalues
$\leq0$ and backward from $\td\tau$ for eigenvalue $1$, we obtain
\[
U^+(\tau)=\sum_{i=1}^k b_i e^{\frac12\tau}\mathfrak p_i^{(1)}g_{\IS^{n-k}}
+O(\eta e^{0.98(\tau-\td\tau)}).
\]
In particular, $\|U^+(\tau)\|\leq C(C_0)\eta e^{\frac12(\tau-\td\tau)}$.
Since $Q_J^+-Q_2^+$ vanishes to third order, substitution into the ODI
and Lemma~\ref{Lem_Q2_computation} give
\begin{align*}
\partial_\tau U^+-LU^+
&=Q_2^+\bigg(\sum_{i=1}^k b_i e^{\frac12\tau}\mathfrak p_i^{(1)}g_{\IS^{n-k}}\bigg)
+O(\eta^2e^{1.48(\tau-\td\tau)})\\
&=\frac{n-k}{4}e^\tau\sum_{i,j=1}^k b_ib_j\,d\bx_i\otimes d\bx_j
-\frac12|\mathbf b|^2e^\tau\mathfrak p^{(0)}g_{\IS^{n-k}}
+O(\eta^2e^{1.48(\tau-\td\tau)}).
\end{align*}
The two explicit terms lie in $\sV_{\aff}$ and $\sV_{\rot,1}$, respectively.
Integrating the projected equations once more therefore gives
\eqref{eq:coefficient_metric_expansion} with $u_\tau$ replaced by $U^+(\tau)$
and an $L^2_f$ remainder $O(\eta^2e^{1.48(\tau-\td\tau)})$.
The factor $\tau e^\tau$ comes from integrating the second term on the right-hand side.
Evaluating at $\tau=\td\tau$ gives
$|a|\leq C(C_0)\eta(1+|\td\tau|)e^{-\td\tau}$.
Comparing two such expansions successively at orders $e^{\frac12\tau}$
and $e^\tau$ proves uniqueness of $\mathbf b$ and $a$.

Let $P_\tau$ denote the sum of the explicit terms preceding $\mathcal E_\tau$ in \eqref{eq:coefficient_metric_expansion}.
The difference $U^+(\tau)-P_\tau$ lies in the finite-dimensional space $\sV_{\rot,\geq-J}$, so its $L^2_f$ bound also gives a $C^{10}$ bound on $\IB^k_{A+2}\times(\IS^{n-k}/\Gamma)$, with a constant depending on $A,J,C_0$.
Now fix $\tau\leq\td\tau$ with $R(\tau)\geq A+10$.
Since $R$ is nonincreasing, $R(s)\geq R(\tau)\geq A+10$ for every $s\in[\tau-1,\tau]$.
As $m_0\geq10$, combining the preceding bound with the admissible estimate for $u_s-U^+(s)$ gives a $C^{10}$ bound for $u_s-P_s$ on $\IB^k_{A+2}\times(\IS^{n-k}/\Gamma)$ throughout this interval, with a constant independent of $m_0$.
At $s=\td\tau$, the admissible bound and the expansion of $U^+$ give $e^{\td\tau}|a-\frac12|\mathbf b|^2\td\tau|\leq C(C_0)\eta$.
Together with the bound on $\mathbf b$, the quadratic computation in Lemma~\ref{Lem_Q2_computation} shows that $P_s$ solves the rescaled Ricci--DeTurck equation up to a residual of order $O_{C^r(\IB^k_{A+2}\times(\IS^{n-k}/\Gamma))}(\eta e^{1.2(s-\td\tau)})$ for every fixed $r$.
Local parabolic derivative estimates for the Ricci--DeTurck equation first give higher derivative bounds for $u_s$ on the smaller ball $\IB^k_{A+1}\times(\IS^{n-k}/\Gamma)$ after a short waiting time; see~\cite[Lemma~A.14]{Bamler_Kleiner_uniqueness}.
Interior estimates for the equation of $u_s-P_s$, using its $C^{10}$ bound and the residual of $P_s$, then give the claimed $C^m$ bound on $\IB^k_A\times(\IS^{n-k}/\Gamma)$.
The condition $R(\tau)\to\infty$ ensures the required radius bound on $[\tau-1,\tau]$ for all sufficiently negative $\tau$.

Suppose finally that $a=0$ and $\mathbf b=0$.
By the mode expansion just proved,
$\|U^+(\tau)\|=O(e^{1.2\tau})$ as $\tau\to-\infty$.
Since the largest eigenvalue of $L$ is $1$ and $Q_J^+$ has no
constant or linear terms, the ODI in
Proposition~\ref{Prop_ancient_summary}\ref{Prop_ancient_summary_d}
implies $\partial_\tau\|U^+\|\leq1.1\|U^+\|$ for sufficiently
negative $\tau$.
Integrating from $-\infty$ yields $U^+=0$ there, and Gronwall's
inequality extends this to all $\tau\leq\td\tau$.
Thus $u_\tau=0$ on $\IB^k_{R(\tau)}\times(\IS^{n-k}/\Gamma)$.
On these regions, the HMHF equation gives $\partial_t\chi_t=0$.
Fix $y_0\in\IS^{n-k}/\Gamma$ and let
$p(t):=\chi_t^{-1}(0,y_0)$, which is therefore a worldline.
The maps $\chi_t$ identify the balls
$B_{g_t}(p(t),\sqrt{t_0-t}\,R_*(t))$ with open subsets of the
round shrinking cylinder with extinction time $t_0$, where
$R_*(t):=\frac12 R(-\log(t_0-t))\to\infty$ as $t\to-\infty$.
Apply Proposition~\ref{Prop_ancient_summary} along this worldline,
with reference time $t_0$ and a sufficiently negative terminal time.
Unless $\MM$ is already a round shrinking cylinder, this gives new
data $(\chi',\Omega',R',U^{\prime,+})$, with leading component
$U^{\prime,++}$.
By Assertions~\ref{Prop_ancient_summary_b2},
\ref{Prop_ancient_summary_c}, and~\ref{Prop_ancient_summary_e}
of Proposition~\ref{Prop_ancient_summary},
the new rescaled metric equals
$g_{\cyl}+U^{\prime,++}+U_{\aff}^{\prime,+}+O(\|U^{\prime,++}\|^2)$ in $C^3$ near
the image of $p(t)$.
Its Ricci eigenvalues there are exactly $0$ and $\frac12$, and
its scalar curvature has zero gradient.
The curvature computation in the proof of Lemma~\ref{l:dom_curvature}
therefore gives $\|U^{\prime,++}\|\leq C\|U^{\prime,++}\|^2$.
Since $\|U^{\prime,++}\|\to0$, we obtain $U^{\prime,++}=U^{\prime,+}=0$ for all
sufficiently negative times, and the identity
$e^{-(R')^2}=\|U^{\prime,+}\|^J$ in
Proposition~\ref{Prop_ancient_summary}\ref{Prop_ancient_summary_c}
gives $R'=\infty$.
Thus $\chi'$ is an isometry between an open subset of $\MM$ and a round shrinking cylinder.
By pseudolocality and forward uniqueness for complete bounded-curvature
Ricci flows~\cite[Theorem~1.1]{Chen_Zhu_RF_uniqueness}, this cylindrical
identification extends throughout the time-interval of $\MM$, proving
the assertion.
\end{proof}

\begin{proof}[Proof of Lemma~\ref{lem:HMHF_constants_reference_time}]
Set $\Delta t:=t_1-t_0$ and write
$\tau_i:=-\log(t_i-t)$ for the rescaled times corresponding to
the same physical time $t$.
Thus
\[
t=t_0-e^{-\tau_0}=t_1-e^{-\tau_1},
\qquad
\tau_0=\tau_1-\log(1-\Delta t\,e^{\tau_1}).
\]
Let $u^{(i)}$ denote the metric perturbation with reference time
$t_i$, and set
\[
\lambda:=e^{(\tau_1-\tau_0)/2}
=\sqrt{\frac{t_0-t}{t_1-t}}.
\]
For sufficiently negative $t$, we have on
$\IB^k_1\times(\IS^{n-k}/\Gamma)$
\[
u^{(1)}_{\tau_1}
=e^{\tau_1}(\xi_{e^{\tau_1/2}})_*(\chi_t)_*g_t-g_{\cyl}
=\lambda^2(\xi_\lambda)_*(g_{\cyl}+u^{(0)}_{\tau_0})-g_{\cyl}
=-\Delta t\,e^{\tau_1}g_{\IS^{n-k}/\Gamma}
+\lambda^2(\xi_\lambda)_*u^{(0)}_{\tau_0}.
\]
Substituting the asymptotic expansions
\eqref{eq:coefficient_metric_expansion} for $u^{(0)}_{\tau_0}$
and $u^{(1)}_{\tau_1}$ and comparing coefficients yields the
asserted relations between $a(\chi,t_0),\mathbf b(\chi,t_0)$
and $a(\chi,t_1),\mathbf b(\chi,t_1)$.

For the final assertion, assume that $\MM$ has dominant linear mode.
Suppose that $\mathbf b(\chi)=0$ and set $t_0':=t_0+a(\chi,t_0)$.
Choose a point $p$ at a sufficiently negative time whose rescaled
HMHF image lies in $\{0\}\times(\IS^{n-k}/\Gamma)$.
Along its worldline, the rescaled axial coordinate satisfies
$\partial_\tau\bx=\frac12\bx+X_\tau$, where $X_\tau\to0$
uniformly on $\IB^k_2\times(\IS^{n-k}/\Gamma)$ by the HMHF equation.
Since this fixed region lies in the rescaled image of $\chi$ at
every time, backward continuation keeps the worldline in $\Omega$
with $|\bx|\leq1$ and gives $\bx(\tau)\to0$ as $\tau\to-\infty$.

The expansion \eqref{eq:coefficient_metric_expansion} and the curvature
computation in the proof of Lemma~\ref{l:dom_curvature} give
\[
|\nabla R_{g_t}|(p(t))=o((-t)^{-2})
\]
as $t\to-\infty$.
Since $p(t)$ is an arbitrarily accurate neck center for sufficiently
negative $t$, we may apply Proposition~\ref{Prop_ancient_summary}
at $p(t')$ for sufficiently negative $t'$, with reference time $t_0'$.
Let $\chi'$ be the resulting harmonic map heat flow.
The new data have $T_3'=-\infty$ by
Definition~\ref{defn:dominant_linear}.
Their metric expansion gives
\[
\lim_{t\to-\infty}(-t)^2|\nabla R_{g_t}|(p(t))
=\frac{n-k}{2\sqrt2}|\mathbf b(\chi')|,
\]
so $\mathbf b(\chi')=0$.
Comparing the scalar curvatures at $p(t)$ given by
\eqref{eq:coefficient_metric_expansion} for $\chi$ and $\chi'$, and using
$t_0'=t_0+a(\chi,t_0)$, gives $a(\chi',t_0')=0$.
Proposition~\ref{p:coefficient} implies that $\MM$ is a round
shrinking cylinder, a contradiction.
\end{proof}

\begin{proof}[Proof of Proposition~\ref{prop:diff_of_U}]
By Definition~\ref{defn:dominant_linear}, both data sets have transition time $T_3=-\infty$.
Set
\[
V(\tau):=U^{+,(1)}(\tau)-U^{+,(0)}(\tau).
\]
By Proposition~\ref{p:coefficient}, the equality of the asymptotic
constants implies
\begin{equation} \label{eq:V_coarse}       
\|V(\tau)\|\leq C(J,C_0)e^{1.2(\tau-\td\tau)}.
\end{equation}
The components below $\sV_{\rot,\geq-10}$ satisfy the same bound by
Proposition~\ref{Prop_ancient_summary}\ref{Prop_ancient_summary_e}.

Subtracting the two ODIs of $U^{+,(i)}(\tau)$ in
Proposition~\ref{Prop_ancient_summary}\ref{Prop_ancient_summary_d},
and using that the largest eigenvalue of $L$ is $1$ and that $Q_J^+$
has no constant or linear term, gives, in the upper barrier sense,
\[
\partial_\tau\|V\|
\leq
(1+C(J,C_0)\eta)\|V\|
+C(J,C_0)\|U^{+,(0)}\|^{J+1}
+C(J,C_0)\|V\|^{J+1}.
\]
After decreasing $\ov\eta(J,C_0)$, the last term can be absorbed into the
first one, and hence
\[
\partial_\tau\|V\|
\leq
1.1\|V\|+C(J,C_0)\|U^{+,(0)}\|^{J+1}.
\]

Suppose first that $U^{+,(0)}\equiv0$. Integrating the inequality
from $-\infty$ to $\tau$, using the coarse bound \eqref{eq:V_coarse} on
$V$, we obtain $V\equiv0$.
Suppose now that $U^{+,(0)}\not\equiv0$. Since $T_3=-\infty$, by
Proposition~\ref{Prop_ancient_summary}\ref{Prop_ancient_summary_f} we have for $s\leq\tau$,
\[
\|U^{+,(0)}(s)\|
\leq
C e^{0.49(s-\tau)}\|U^{+,(0)}(\tau)\|.
\]
Integrating the differential inequality for $V$ from $-\infty$ to
$\tau$, and again using \eqref{eq:V_coarse}, yields
\[
\begin{split}
\|V(\tau)\|
&\leq
C(J,C_0)\int_{-\infty}^{\tau}
e^{1.1(\tau-s)}\|U^{+,(0)}(s)\|^{J+1}\,ds  \\
&\leq
C(J,C_0)\|U^{+,(0)}(\tau)\|^{J+1}
\int_{-\infty}^{\tau}
e^{(0.49(J+1)-1.1)(s-\tau)}\,ds
\leq
C(J,C_0)\|U^{+,(0)}(\tau)\|^{J+1}.
\end{split}
\]
This implies the asserted bound.
\end{proof}

\subsection{Transformation laws for the HMHF data}\label{subsec:linear_transformations}
We record how rotations, parabolic rescalings, and spatial translations affect admissibility and the asymptotic constants.
We will use the translation lemma when aligning an HMHF with a prescribed gauge at infinity.

\begin{Lemma}[Rotation]\label{l:rotation}
Let $\MM$ be an ancient asymptotically $(n,k,\Gamma)$-cylindrical Ricci flow spacetime satisfying the pseudolocality property and having dominant linear mode.
Let $\chi:\Omega\to\IR^k\times(\IS^{n-k}/\Gamma)$ be admissible with parameters $m,J\geq10$, $0<\eta<0.1$, $C_0>0$, and reference time $t_0$ on $(-\infty,\td\tau]$.
For $Q\in O(k)$ and $C\in\operatorname{Isom}(\IS^{n-k}/\Gamma)$, set $F(\bx,\by):=(Q\bx,C(\by))$.
Then $\chi':=F\circ\chi$ is admissible with the same parameters, reference time, and time-interval, and
\[
a(\chi',t_0)=a(\chi,t_0),\qquad
\mathbf b(\chi')=Q\mathbf b(\chi).
\]
\end{Lemma}

\begin{proof}
The data $(F\circ\chi,\Omega,R,F_*U^+)$ are admissible with the same parameters since $F$ is an isometry of $M_{\cyl}$.
For the asymptotic constants, write $a:=a(\chi,t_0)$ and $\mathbf b:=\mathbf b(\chi)$.
Since $u'_\tau=F_*u_\tau$, Proposition~\ref{p:coefficient} gives
\begin{align*}
u'_\tau
&=\frac{n-k}{4}e^\tau\sum_{i,j=1}^k
b_i b_j\,F_*(d\bx_i\otimes d\bx_j)\\
&\qquad+\Big(\big(a-\tfrac12|\mathbf b|^2\tau\big)e^\tau\mathfrak p^{(0)}
+e^{\tau/2}\mathbf b\cdot\mathfrak p^{(1)}(Q^{-1}\bx)\Big)
g_{\IS^{n-k}/\Gamma}+O(e^{1.2\tau})\\
&=\frac{n-k}{4}e^\tau\sum_{i,j=1}^k
(Q\mathbf b)_i(Q\mathbf b)_j\,d\bx_i\otimes d\bx_j\\
&\qquad+\Big(\big(a-\tfrac12|Q\mathbf b|^2\tau\big)e^\tau\mathfrak p^{(0)}
+e^{\tau/2}(Q\mathbf b)\cdot\mathfrak p^{(1)}(\bx)\Big)
g_{\IS^{n-k}/\Gamma}+O(e^{1.2\tau}).
\end{align*}
Here $\mathfrak p^{(1)}$ denotes the column vector of linear Hermite polynomials, and we used $Q^{-1}=Q^{\mathsf T}$ and $|Q\mathbf b|=|\mathbf b|$.
The assertion follows from the uniqueness of $a$ and $\mathbf b$ in Proposition~\ref{p:coefficient}.
\end{proof}
\bigskip

\begin{Lemma}[Parabolic rescaling]\label{l:rescaling}
Let $\MM$ be an ancient asymptotically $(n,k,\Gamma)$-cylindrical Ricci flow spacetime satisfying the pseudolocality property and having dominant linear mode.
Let $\chi:\Omega\to\IR^k\times(\IS^{n-k}/\Gamma)$ be admissible with parameters $m,J\geq10$, $0<\eta<0.1$, $C_0>0$, and reference time $t_0$ on $(-\infty,\td\tau]$.
For $\alpha>0$, consider the parabolically rescaled spacetime
$\alpha\MM:=(\MM,\alpha^2\mathfrak t,\alpha^{-2}\partial_{\mathfrak t},\alpha^2g)$ and set
$\chi'_t:=\xi_\alpha\circ\chi_{\alpha^{-2}t}$ on $\Omega'_t:=\Omega_{\alpha^{-2}t}$.
Then $\chi'$ is admissible on $\alpha\MM$ with the same parameters, reference time $t'_0:=\alpha^2t_0$, and rescaled time-interval $(-\infty,\td\tau-2\log\alpha]$, and
\[
a(\chi',t'_0)=\alpha^2\bigl(a(\chi,t_0)-|\mathbf b(\chi)|^2\log\alpha\bigr),
\qquad
\mathbf b(\chi')=\alpha\mathbf b(\chi).
\]
\end{Lemma}

\begin{proof}
After shifting the time coordinate, we may take $t_0=0$.
First, the identity $g_{\cyl,t}=\xi_{(-t)^{-1/2}}^*((-t)g_{\cyl})$ gives
\[
\alpha^2g_{\cyl,\alpha^{-2}t}
=\xi_\alpha^*g_{\cyl,t}.
\]
Thus $\xi_\alpha$ is an isometry from the rescaled cylindrical flow to the cylindrical flow.
Since $\chi$ is an HMHF, we also have
\[
\partial_t\chi_{\alpha^{-2}t}
=\alpha^{-2}\Delta_{g_{\alpha^{-2}t},g_{\cyl,\alpha^{-2}t}}\chi_{\alpha^{-2}t}
=\Delta_{\alpha^2g_{\alpha^{-2}t},\alpha^2g_{\cyl,\alpha^{-2}t}}\chi_{\alpha^{-2}t}.
\]
Hence $\chi'_t=\xi_\alpha\circ\chi_{\alpha^{-2}t}$ is an HMHF on $\alpha\MM$.
For $t=-e^{-\tau'}$, put $\tau:=\tau'+2\ln\alpha$.
The rescaled perturbations induced by $\chi$ and $\chi'$ satisfy
\[
u'_{\tau'}=u_{\tau'+2\ln\alpha}.
\]
Given admissible data $(\chi,\Omega,R,U^+)$, take
$R'(\tau'):=R(\tau'+2\ln\alpha)$ and
$U^{\prime,+}(\tau'):=U^+(\tau'+2\ln\alpha)$.
The transformed data are admissible with the same parameters on $(-\infty,\td\tau-2\ln\alpha]$.
In particular, the approximation of $u'_{\tau'}$ by $U^{\prime,+}(\tau')$ follows from that of $u_\tau$ by $U^+(\tau)$.
Write $a:=a(\chi,t_0)$ and $\mathbf b:=\mathbf b(\chi)$.
By the expansion in Proposition~\ref{p:coefficient} and $e^{1.2(\tau'+2\ln\alpha)}=\alpha^{2.4}e^{1.2\tau'}$, we have
\begin{align*}
u'_{\tau'}=u_{\tau'+2\ln\alpha}
&=\frac{n-k}{4}e^{\tau'}\sum_{i,j=1}^k
(\alpha b_i)(\alpha b_j)\,d\bx_i\otimes d\bx_j\\
&\qquad+
\Big(\big(\alpha^2(a-|\mathbf b|^2\ln\alpha)
-\tfrac12|\alpha\mathbf b|^2\tau'\big)e^{\tau'}\mathfrak p^{(0)}
+e^{\tau'/2}(\alpha\mathbf b)\cdot\mathfrak p^{(1)}\Big)
g_{\IS^{n-k}/\Gamma}
+O(e^{1.2\tau'}).
\end{align*}
The assertion follows from the uniqueness of the asymptotic constants in Proposition~\ref{p:coefficient}.
Undoing the time shift gives the stated reference time $t'_0=\alpha^2t_0$.
\end{proof}
\bigskip

\begin{Lemma}[Translation in space]\label{l:translation}
Let $\MM$ be an ancient asymptotically $(n,k,\Gamma)$-cylindrical Ricci flow spacetime satisfying the pseudolocality property and having dominant linear mode.
Let $\chi:\Omega\to M_{\cyl}$ be admissible with parameters $m,J\geq10$, $0<\eta_*<0.1$, $C_0>0$, and reference time $t_0$ on $(-\infty,\td\tau]$.
Choose corresponding data $(\chi,\Omega,R,U^+)$ with $R$ nonincreasing in $\tau$ and $R(\td\tau)\geq1010$.
Suppose that $\bx_0\in\IR^k$ with $e^{\td\tau/2}|\bx_0|\leq10$, let $T_{\bx_0}(\bx,\by):=(\bx+\bx_0,\by)$ denote translation in the Euclidean factor and set $\chi':=T_{\bx_0}\circ\chi$.
For every $0<\eta<0.1$, if $\eta_* \leq \ov\eta_* (m,J,C_0,\eta)$, then $\chi'$ is admissible with parameters $m,J,\eta,C(m,J,C_0)$, the same reference time and interval, and radius $R'(\tau):=R(\tau)-10$.
\end{Lemma}

\begin{proof}
Since $T_{\bx_0}$ is an isometry of the cylinder, $\chi'$ is a
harmonic map heat flow. Choose admissible data
$(\chi,\Omega,R,U^+)$, and write $U:=U^+$.
By Lemma~\ref{lem:dominant_mode_independence}, we have $T_3=-\infty$.
For $\bx_0(\tau):=e^{\tau/2}\bx_0$, consider the action $S_\tau:=(T_{\bx_0(\tau)})_*$ on tensor fields by translation.
This action preserves
$\sV_{\rot,\geq-J}$ and for the rescaled metric perturbations it satisfies $u'_\tau=S_\tau u_\tau$.
We set $U'(\tau):=S_\tau U(\tau)\in\sV_{\rot,\geq-J}$.
On $\sV_{\rot,\geq-J}$, the norms of $S_\tau$ and its inverse are bounded by a constant of the form $C(J)$.
In particular,
\[
C^{-1}(J)\|U\|\leq\|U'\|\leq C(J)\|U\|.
\]
Since $R$ is nonincreasing, $R(\td\tau)\geq1010$, and $|\bx_0(\tau)|\leq10$, we have $R'(\tau) = R(\tau) - 10 \geq1000$ and $R'(\tau)\to\infty$ as $\tau\to-\infty$, and the translated image contains $\IB^k_{R'(\tau)}\times(\IS^{n-k}/\Gamma)$.
The approximation bound for $U$ gives
\[
\|u'_\tau-U'(\tau)\|_{C^m}
\leq C_0\|U(\tau)\|^{J+1}
\leq C(J)C_0\|U'(\tau)\|^{J+1}
\]
on this smaller ball.
The other $\eta$-smallness bounds follow if $\eta_*\leq\ov\eta_*(m,J,C_0,\eta)$.
This establishes Assertions~\ref{Prop_ancient_summary_b}, \ref{Prop_ancient_summary_b2}, and~\ref{Prop_ancient_summary_c} of Proposition~\ref{Prop_ancient_summary} for $\chi'$, $U'$ and $R'$, with $C_0$ replaced by $C(J)C_0$.

\begin{Claim}\label{cl:translation_full_residual}
We have
\[
\|\partial_\tau U-LU-Q_J(U)\|
\leq C(J,C_0)\|U\|^{J+1}.
\]
\end{Claim}

\begin{proof}
By Proposition~\ref{Prop_ancient_summary}\ref{Prop_ancient_summary_d}, it suffices to bound $\|Q_J^-(U)\|$.
Write $U=\sum_i U_i$, where $U_i$ collects the rotational modes whose coefficients are polynomials of degree $i$ in the axial variables, and recall that $\|U\|<1$.
Proposition~\ref{Prop_ancient_summary}\ref{Prop_ancient_summary_e} and $\|U^{++}\|\leq\|U\|$ give
\[
\|U_i\|\leq C(J,C_0)\|U\|^{\max\{1,\lceil i/2\rceil\}}.
\]
As in the degree count in the proof of Proposition~\ref{Prop_mode_analysis}, Lemma~\ref{Lem_Q_precise_formula} shows that each term of $Q_J(U)$ is a polynomial in the axial variables of degree at most the sum $i_1+\cdots+i_l$ of the degrees of its input factors.
A rotationally symmetric mode with eigenvalue below $-J$ has degree greater than $2J$, so any monomial contributing to $Q_J^-(U)$ has $i_1+\cdots+i_l>2J$.
Since $\sum_{j=1}^l\max\{1,\lceil i_j/2\rceil\}\geq\lceil(i_1+\cdots+i_l)/2\rceil\geq J+1$, summing the finitely many monomials gives $\|Q_J^-(U)\|\leq C(J,C_0)\|U\|^{J+1}$.
\end{proof}

\begin{Claim}\label{cl:translation_projected_residual}
The ODI in Proposition~\ref{Prop_ancient_summary}\ref{Prop_ancient_summary_d} holds for $U'$ in the following form:
\[
\|\partial_\tau U'-LU'-Q_J^+(U')\|
\leq C(J,C_0)\|U'\|^{J+1}.
\]
\end{Claim}

\begin{proof}
Since $\partial_\tau\bx_0(\tau)=\tfrac12\bx_0(\tau)$, the drift term in $L$ gives
\[
(\partial_\tau-L)(S_\tau U)
=S_\tau\partial_\tau U-\tfrac12\nabla_{\bx_0(\tau)}(S_\tau U)
-S_\tau LU+\tfrac12\nabla_{\bx_0(\tau)}(S_\tau U)
=S_\tau(\partial_\tau-L)U.
\]
Moreover, the formula in Lemma~\ref{Lem_Q_precise_formula} shows that
$Q$ and its Taylor approximation are translation invariant.
Thus $Q_J(S_\tau U)=S_\tau Q_J(U)$, and hence
\[
\partial_\tau U'-LU'-Q_J(U')
=S_\tau\big(\partial_\tau U-LU-Q_J(U)\big).
\]
Projecting this identity onto $\sV_{\rot,\geq-J}$ gives the left-hand side of the claim.
Claim~\ref{cl:translation_full_residual}, the bound on $S_\tau$, and the norm comparison above give the stated estimate.
\end{proof}

For each $i=0,\ldots,J-1$, translations preserve $\sV_{\rot,\geq-i}$, so
\[
\PP_{\sV_{\rot,<-i}}U'
= \PP_{\sV_{\rot,<-i}} S_\tau \big( \PP_{\sV_{\rot,<-i}}  U + \PP_{\sV_{\rot,\geq-i}}  U \big)
=\PP_{\sV_{\rot,<-i}}S_\tau\PP_{\sV_{\rot,<-i}}U.
\]
The bound on $S_\tau$, Proposition~\ref{Prop_ancient_summary}\ref{Prop_ancient_summary_e} for $U$, and the comparison $\|U^{++}\|\leq\|U\|\leq C(J,C_0)\|U'^{++}\|$ therefore imply Assertion~\ref{Prop_ancient_summary_e} for $U'$.

Assertion~\ref{Prop_ancient_summary_f2} of Proposition~\ref{Prop_ancient_summary} and the norm comparison give the coarse bound $\|U'(s)\|\leq C(J,C_0)e^{0.49(s-\tau)}\|U'(\tau)\|$ for $s\leq\tau$.
Integrating the nonpositive projection of the ODI in Claim~\ref{cl:translation_projected_residual} from $-\infty$ gives $\|\PP_{\sV_{\rot,\leq0}}U'\|\leq C(J,C_0)\|U'\|^2$, since $Q_J^+$ has no linear term.
If $\eta_*\leq\ov\eta_*(m,J,C_0,\eta)$, the positive modes therefore dominate, and projecting the same ODI onto the eigenspaces with eigenvalues $1$ and $\frac12$ gives
\[
\|U'\|\leq 2\|\PP_{\sV^{++}_{>0}}U'\|,
\qquad
\partial_\tau\|\PP_{\sV^{++}_{>0}}U'\|
\geq0.49\|\PP_{\sV^{++}_{>0}}U'\|.
\]
Integrating the second inequality gives the asserted exponential in Assertion~\ref{Prop_ancient_summary_f2}.
Thus $\chi'$ is admissible on the same interval with constant $C(m,J,C_0)$.
\end{proof}
\bigskip

\subsection{Gauges at infinity}\label{subsec:linear_gauges}

In this subsection, we construct gauges at infinity for asymptotically cylindrical Ricci flow spacetimes with dominant linear mode. These gauges are constant along worldlines and unique up to an isometry of the cylindrical model. We will describe their asymptotics and compare the metrics of two flows whose asymptotic constants agree. This will be used to construct the comparison map to the Bryant model in Subsection~\ref{subsec:linear_bryant_comparison}. We first state the results and then give their proofs.

Write $M_{\cyl}:=\IR^k\times(\IS^{n-k}/\Gamma)$, with its standard product metric $g_{\cyl}$. Norms on this cylinder are taken with respect to $g_{\cyl}$; constants may depend on $n,k,\Gamma$ without further indication.

\begin{Definition}\label{defn_gauge}
Let $\MM$ be a Ricci flow spacetime over a time-interval of the form $(-\infty,T)$, where $T\in\IR\cup\{\infty\}$.
A smooth map $\phi:\Omega\to\IR^k\times(\IS^{n-k}/\Gamma)$, where $\Omega\subset\MM$ is open, is called a \textbf{gauge at infinity} if the following conditions hold:
\begin{enumerate}[label=(\alph*)]
\item\label{defn_gauge_a} The time-slices $\phi_t:\Omega_t\to\IR^k\times(\IS^{n-k}/\Gamma)$ are diffeomorphisms onto their images.
\item\label{defn_gauge_b} For $t_1<t_2<T$, the backward flow map
$\Phi^{\partial_{\mathfrak t}}_{t_2,t_1}:p\mapsto p(t_1)$ of the time vector field is defined on $\Omega_{t_2}$ and satisfies
\[
 \Phi^{\partial_{\mathfrak t}}_{t_2,t_1}(\Omega_{t_2})\subset\Omega_{t_1},
 \qquad
 \phi_{t_1}\circ\left.\Phi^{\partial_{\mathfrak t}}_{t_2,t_1}\right|_{\Omega_{t_2}}=\phi_{t_2}.
\]
Consequently, $\phi_{t_2}(\Omega_{t_2})\subset\phi_{t_1}(\Omega_{t_1})$.
\item\label{defn_gauge_c} For every compact subset $K\subset\IR^k\times(\IS^{n-k}/\Gamma)$, we have $K\subset\xi_{(-t)^{-1/2}}\big(\phi_t(\Omega_t)\big)$ for all sufficiently negative $t$.
\item\label{defn_gauge_d} We have local smooth convergence $(-t)^{-1}(\xi_{(-t)^{-1/2}})_*(\phi_t)_*g_t\longrightarrow g_{\cyl}$ as $t\to-\infty$ on $\IR^k\times(\IS^{n-k}/\Gamma)$.
\item\label{defn_gauge_mixed} The mixed components of $(\phi_t)_*g_t$ with respect to the splitting $\IR^k\times(\IS^{n-k}/\Gamma)$ converge to $0$ locally smoothly as $t\to-\infty$.
\item\label{defn_gauge_e} For every $p\in\Omega$, the maximal worldline through $p$ is contained in $\Omega$.
\end{enumerate}
\end{Definition}

{The following lemma shows that two gauges at infinity of the same flow have the same domain and only differ by a cylindrical isometry. 
This allows us to compare gauges obtained from different HMHF data.}
\begin{Lemma}\label{l:isometry_gauge}
Let $\MM$ be an ancient asymptotically $(n,k,\Gamma)$-cylindrical
Ricci flow spacetime over the time-interval $(-\infty,T)$, where
$T\in\IR\cup\{\infty\}$.
Write $M_{\cyl}:=\IR^k\times(\IS^{n-k}/\Gamma)$, equipped with
the metric $g_{\cyl}$.
For $i=0,1$, let $\phi^{(i)}:\Omega^{(i)}\to M_{\cyl}$
be gauges at infinity of $\MM$.
Then there is a unique isometry
$\psi:M_{\cyl}\to M_{\cyl}$ such that
\[
\Omega^{(0)}=\Omega^{(1)},\qquad \phi^{(1)}=\psi\circ\phi^{(0)}.
\]
Conversely, composing a gauge at infinity with an isometry of
$M_{\cyl}$ again gives a gauge at infinity.
Moreover, if a smooth map $\phi':\Omega'\to M_{\cyl}$ on an open
subset $\Omega'\subset\MM$ satisfies Properties~\ref{defn_gauge_a}--\ref{defn_gauge_mixed}
of Definition~\ref{defn_gauge}, then it extends uniquely to a gauge at infinity
$\phi:\Omega\to M_{\cyl}$ with $\Omega\supset\Omega'$ and
$\phi|_{\Omega'}=\phi'$.
\end{Lemma}

The following proposition gives an asymptotic expansion of the metric in a gauge at infinity.

\begin{Proposition}[Geometry in a gauge at infinity]\label{prop:suggested_gauge_geometric}
Let $\MM$ be an ancient asymptotically $(n,k,\Gamma)$-cylindrical Ricci flow spacetime satisfying the pseudolocality property and having dominant linear mode.
Then the following is true:
\begin{enumerate}[label=(\alph*)]
\item\label{prop:suggested_gauge_geometric_existence}
There is a gauge at infinity $\phi : \Omega^\phi \to M_{\cyl}$.
\item\label{prop:suggested_gauge_geometric_asymptotics}
For any such gauge and any $t_0\in\IR$ there are unique $a=a(\phi,t_0)\in\IR$ and $\mathbf b=\mathbf b(\phi)\in\IR^k$ such that the following holds for every $A>0$ and $m\geq0$.
For sufficiently negative $\tau$, the perturbation
\[
\td u_\tau:=e^\tau(\xi_{e^{\tau/2}}\circ\phi_t)_*g_t-g_{\cyl},
\qquad t=t_0-e^{-\tau},
\]
is defined on $\IB^k_A\times(\IS^{n-k}/\Gamma)$.
Writing $d:=n-k$, we have on this region
\begin{multline*}
\qquad\qquad \td u_\tau={}-\frac d4e^\tau\sum_{i,j=1}^k b_i b_j\,d\bx_i\otimes d\bx_j \\
+\left[\frac{e^{\tau/2}}{\sqrt2}\langle\mathbf b,\bx\rangle
+e^\tau\left(a+\frac{d-2}{4}|\mathbf b|^2\tau
+\frac d4|\mathbf b|^2\log|\mathbf b|^2\right)\right]g_{\IS^{n-k}/\Gamma}
+\mathcal E_\tau,
\end{multline*}
where
\[
\|\mathcal E_\tau\|_{C^m(\IB^k_A\times(\IS^{n-k}/\Gamma))}
\leq C(A,m,\phi,t_0)e^{1.2\tau}.
\]
\item\label{prop:suggested_gauge_geometric_control}
If $m,J\geq100$, $0<\eps<0.1$, and $\delta\leq\ov\delta(m,J,\eps)$, then the following is true.
Suppose that $\phi : \Omega^\phi \to M_{\cyl}$ is a gauge at infinity, $t_0 \in \IR$, $\td\tau \in \IR$ and $\td t := t_0 -e^{-\td\tau}$.
Let $p \in \MM_{\td t}$ and assume that $p(t)$ is a center of an $(n,k, \delta, \Gamma)$-neck at scale $\sqrt{t_0 - t}$ for all $t \leq \td t$ and that 
\[
\xi_{e^{\td\tau/2}}\circ\phi_{\td t}(p)\in\IB^k_1\times(\IS^{n-k}/\Gamma).
\]
Then for all $\tau\leq\td\tau$, with $t=t_0-e^{-\tau}$, the perturbation $\td u_\tau$ is defined on $\IB^k_{100}\times(\IS^{n-k}/\Gamma)$ and satisfies
\begin{equation}\label{eq:gauge_control_metric_bound}
\|\td u_\tau\|_{C^m(\IB^k_{100}\times(\IS^{n-k}/\Gamma))}\leq\eps.
\end{equation}
Moreover, there is an admissible HMHF $\chi : \Omega \to M_{\cyl}$ with parameters $m, J, \eps, C(m,J)$, reference time $t_0$ and rescaled time-interval $(-\infty,\td\tau]$ such that the following holds.
We have $a(\chi,t_0) = a(\phi,t_0)$ and $\mathbf b(\chi) = \mathbf b (\phi)$.
Set $\mathbf v:=\frac{n-k}{2\sqrt2}\mathbf b(\phi)$, and let $T_{\mathbf z}$ denote translation by $\mathbf z$ in the Euclidean factor.
For every compact subset $K\subset\Omega^\phi_t$, the following maps are defined for sufficiently negative $s$ and satisfy
\begin{equation}\label{eq:geometric_gauge_HMHF_limit}
T_{\mathbf v\log\left(\frac{t_0-s}{|\mathbf b(\phi)|^2}\right)}
\circ\chi_s\circ\left.\Phi^{\partial_{\mathfrak t}}_{t,s}\right|_K
\longrightarrow\left.\phi_t\right|_K
\quad\text{in }C^\infty(K)\quad\text{as }s\to-\infty.
\end{equation}
\end{enumerate}

\end{Proposition}
Along any worldline in the domain of a gauge, the neck and normalization hypotheses hold for sufficiently negative times. Thus the expansion defines $a(\phi,t_0)$ and $\mathbf b(\phi)$ independently of the point and interval used. The following proposition shows that if these constants agree for two gauges, then the difference of the rescaled metrics decays to arbitrarily high exponential order. Cylindrical control is required only on the second flow, whose controlled interval may extend beyond that of the first flow.

\begin{Proposition}[Comparison in gauges at infinity]\label{prop:suggested_gauge_geometric_comparison}
For $m, J' \geq 10$, $0 < \eps < 0.1$ and $\delta \leq \ov\delta(m,J', \eps)$ there is a constant $C_1(m,J', \eps) > 0$ such that the following is true.

Let $\MM,\MM'$ be ancient asymptotically $(n,k,\Gamma)$-cylindrical Ricci flow spacetimes satisfying the pseudolocality property and having dominant linear mode, and suppose that $\MM$ and $\MM'$ are defined over $(-\infty,T)$ and $(-\infty,T')$, respectively. 
Fix $t_0,\td\tau\in\IR$ and set $\td t:=t_0-e^{-\td\tau}$. 
We only require that $\td t < T'$, but do not impose a similar bound on $T$.
Let $\phi,\phi'$ be gauges at infinity of $\MM$ and $\MM'$, respectively and suppose that
\[
a(\phi,t_0)=a(\phi',t_0),\qquad \mathbf b(\phi)=\mathbf b(\phi').
\]
Let $p'\in\Omega^{\phi'}_{\td t}$ and suppose that $p'(t)$ is a center of an $(n,k,\delta,\Gamma)$-neck at scale $\sqrt{t_0-t}$ for all $t\leq\td t$ and that
\[
\xi_{e^{\td\tau/2}} \circ \phi'(p')\in\IB^k_1\times(\IS^{n-k}/\Gamma).
\]
Then for all $\tau \leq \td\tau - C_1$ with $t := t_0 - e^{-\tau}<T$ we have
\[ \xi_{e^{\tau/2}}  \circ \phi ( \Omega^{\phi}_{t}), \;\; \xi_{e^{\tau/2}}  \circ \phi' ( \Omega^{\phi'}_{t}) \supset \IB^k_{50} \times (\IS^{n-k} / \Gamma) \]
and setting
\[ \td u_\tau := e^\tau (\xi_{e^{\tau/2}} \circ \phi_{t})_* g_{t} - g_{\cyl}, \qquad \td u'_\tau := e^\tau (\xi_{e^{\tau/2}} \circ \phi'_{t})_* g'_{t} - g_{\cyl}, \]
we have
\begin{equation}\label{eq:gauge_comparison_decay}
\|\td u_\tau-\td u'_\tau\|_{C^m(\IB^k_{50}\times(\IS^{n-k}/\Gamma))}
\leq\eps e^{J'(\tau-\td\tau)}.
\end{equation}
\end{Proposition}

The next lemma describes how $a(\phi,t_0)$ and $\mathbf b(\phi)$ change under isometries, parabolic rescalings, and changes of reference time. These formulas will allow us to arrange that the constants of two flows agree.

\begin{Lemma}\label{l:gauge_asymp}
Let $\phi$ be a gauge at infinity of a flow as in Proposition~\ref{prop:suggested_gauge_geometric}, and write $a:=a(\phi,t_0)$ and $\mathbf b:=\mathbf b(\phi)$. Then the following holds.
\begin{enumerate}[label=(\alph*)]
\item\label{l:gauge_asymp_isometry} If $S(\bx,\by)=(A\bx+\bx_0,\rho(\by))$ is an isometry of $M_{\cyl}$, then $\phi':=S\circ\phi$ is a gauge on the same domain, and
\begin{equation}\label{eq:gauge_constants_isometry}
\mathbf b(\phi')=A\mathbf b,\qquad
a(\phi',t_0)=a-\frac{1}{\sqrt2}\langle A\mathbf b,\bx_0\rangle.
\end{equation}
\item\label{l:gauge_asymp_rescaling} On the parabolically rescaled flow $\alpha\MM$, the maps $\phi'_s:=\xi_\alpha\circ\phi_{\alpha^{-2}s}$ form a gauge, and
\begin{equation}\label{eq:gauge_constants_rescaling}
\mathbf b(\phi')=\alpha\mathbf b,\qquad
a(\phi',\alpha^2t_0)=\alpha^2\bigl(a-|\mathbf b|^2\log\alpha\bigr).
\end{equation}
\item\label{l:gauge_asymp_time} Changing only the reference time gives
\begin{equation}\label{eq:gauge_constants_reference}
a(\phi,t_1)=a(\phi,t_0)+t_0-t_1,
\qquad \mathbf b(\phi)\text{ is unchanged}.
\end{equation}
If instead the time coordinate is shifted by $\Delta t$, then $\phi'_t:=\phi_{t-\Delta t}$ is a gauge of the shifted flow, with
\begin{equation}\label{eq:gauge_constants_time}
a(\phi',t_0+\Delta t)=a,\qquad
a(\phi',t_0)=a+\Delta t,\qquad
\mathbf b(\phi')=\mathbf b.
\end{equation}
\end{enumerate}
\end{Lemma}

\bigskip

\begin{proof}[Proof of Lemma~\ref{l:isometry_gauge}]
We first show:
\begin{Claim}\label{cl:gauge_common_domain}
$\Omega^{(0)} = \Omega^{(1)}$.
\end{Claim}

\begin{proof}
Fix arbitrary points $p_i\in\Omega^{(i)}$, $i=0,1$.
By Definition~\ref{defn_gauge}\ref{defn_gauge_b}--\ref{defn_gauge_d},
their backward worldlines satisfy
\eqref{eq:asymptotically_cylindrical_convergence}.
Moreover, for every fixed $A>0$ and all sufficiently negative $t$,
\[
B_{g_t}\big(p_i(t),A\sqrt{-t}\big)\subset\Omega^{(i)}_t.
\]
The last condition in Definition~\ref{defn:asymptotically_cylindrical}
gives $d_{g_t}(p_0(t),p_1(t))\leq C(p_0,p_1)\sqrt{-t}$ for all
sufficiently negative $t$.
Choosing $A>C(p_0,p_1)$ in the preceding ball inclusion, we obtain
$p_0(t)\in\Omega^{(1)}_t$ and
$p_1(t)\in\Omega^{(0)}_t$.
Property~\ref{defn_gauge_e} of Definition~\ref{defn_gauge} implies that
the entire worldlines through $p_0$ and $p_1$ are contained in
$\Omega^{(1)}$ and $\Omega^{(0)}$, respectively, and hence
$\Omega^{(0)}=\Omega^{(1)}$.
\end{proof}

Write $\Omega=\Omega^{(0)}=\Omega^{(1)}$ and define
\[
\psi_t:=\phi_t^{(1)}\circ(\phi_t^{(0)})^{-1}:\phi_t^{(0)}(\Omega_t)\longrightarrow\phi_t^{(1)}(\Omega_t).
\]
By Definition~\ref{defn_gauge}\ref{defn_gauge_b} and \ref{defn_gauge_c}, the maps $\psi_t$ are compatible and their domains and images exhaust $M_{\cyl}$ as $t\to-\infty$.
They therefore extend to a global smooth diffeomorphism $\psi:M_{\cyl}\to M_{\cyl}$ satisfying
\[
\phi^{(1)}=\psi\circ\phi^{(0)}.
\]
It remains to prove that $\psi$ is an isometry.
To do this, set
\[
\td g^{(i)}_t:=\big(\phi_t^{(i)}\big)_*g_t,\qquad i=0,1.
\]
Then $\psi_*\td g^{(0)}_t=\td g^{(1)}_t$ on $\phi_t^{(1)}(\Omega_t)$.
Property~\ref{defn_gauge_d} of Definition~\ref{defn_gauge} gives
\begin{equation}\label{eq:isometry_gauge_scaled_conv}
(-t)^{-1}(\xi_{(-t)^{-1/2}})_*\td g^{(i)}_t\longrightarrow g_{\cyl}
\end{equation}
locally smoothly as $t\to-\infty$. In fixed cylindrical coordinates,
this implies that, for vectors $v_1,v_2$ tangent to the $\IR^k$-factor,
\begin{equation}\label{eq:isometry_gauge_eucl_conv}
\td g^{(i)}_t(v_1,v_2)\longrightarrow g_{\IR^k}(v_1,v_2),
\end{equation}
while for arbitrary vectors $v$,
\begin{equation}\label{eq:isometry_gauge_sph_conv}
(-t)^{-1}\td g^{(i)}_t(v,v)\longrightarrow g_{\IS^{n-k}/\Gamma}(v,v).
\end{equation}
Here $g_{\IS^{n-k}/\Gamma}(v,v)$ denotes the squared length of the spherical component of $v$.

Using \eqref{eq:isometry_gauge_sph_conv} and $\psi_*\td g^{(0)}_t=\td g^{(1)}_t$, we obtain, for every tangent vector $v$,
\begin{equation}\label{eq:sphere_iso}
g_{\IS^{n-k}/\Gamma}(d\psi(v),d\psi(v))=g_{\IS^{n-k}/\Gamma}(v,v).
\end{equation}
Thus $d\psi$ preserves the kernel of $g_{\IS^{n-k}/\Gamma}$, namely the Euclidean distribution, and the same is true for $d\psi^{-1}$.
Using \eqref{eq:isometry_gauge_eucl_conv}, for all Euclidean tangent vectors $v_1,v_2$ we obtain
\[
g_{\IR^k}(d\psi(v_1),d\psi(v_2))
=\lim_{t\to-\infty}\td g^{(1)}_t(d\psi(v_1),d\psi(v_2))
=\lim_{t\to-\infty}\td g^{(0)}_t(v_1,v_2)
=g_{\IR^k}(v_1,v_2).
\]
Therefore $d\psi$ restricts to an isometry on the Euclidean distribution.

To show that the spherical distribution is also preserved, let $w$ be tangent to the $\IS^{n-k}/\Gamma$-factor and write $d\psi(w)=v'+w'$ according to the product splitting.
Since $d\psi$ maps the Euclidean distribution isometrically onto itself, there is a Euclidean vector $v$ with $d\psi(v)=v'$.
By Property~\ref{defn_gauge_mixed} of Definition~\ref{defn_gauge}, the mixed terms $\td g^{(0)}_t(v,w)$ and $\td g^{(1)}_t(v',w')$ tend to zero.
Together with \eqref{eq:isometry_gauge_eucl_conv}, this gives
\[
0=\lim_{t\to-\infty}\td g^{(0)}_t(v,w)
=\lim_{t\to-\infty}\td g^{(1)}_t(d\psi(v),d\psi(w))
=g_{\IR^k}(v',v').
\]
Thus $v'=0$, so $d\psi$ preserves the spherical distribution.
Equation~\eqref{eq:sphere_iso} shows that it is an isometry on this distribution as well.
Hence $d\psi$ preserves the product metric $g_{\cyl}=g_{\IR^k}+g_{\IS^{n-k}/\Gamma}$, and $\psi$ is an isometry of $M_{\cyl}$.

The uniqueness of $\psi$ follows because $\phi_t^{(0)}(\Omega_t)$ exhausts $M_{\cyl}$ as $t\to-\infty$.
Conversely, if $\psi$ is an isometry and $\phi$ is a gauge at infinity, then the defining properties of a gauge at infinity are invariant under replacing $\phi$ by $\psi\circ\phi$.

Finally, suppose $\phi':\Omega'\to M_{\cyl}$ satisfies Properties~\ref{defn_gauge_a}--\ref{defn_gauge_mixed} of Definition~\ref{defn_gauge}.
Let $\Omega$ be the set of all points whose backward worldline meets $\Omega'$.
For $p\in\Omega$, choose $t\leq\mathfrak t(p)$ with $p(t)\in\Omega'_t$ and define $\phi(p):=\phi'_t(p(t))$.
This is independent of the choice of $t$ by Property~\ref{defn_gauge_b}, and it clearly extends $\phi'$.
Property~\ref{defn_gauge_b} gives uniqueness of such an extension.
\end{proof}

\medskip

\begin{proof}[Proof of Proposition~\ref{prop:suggested_gauge_geometric}]
Fix $A\geq200$ and integers $m,J\geq10$, and let $\delta>0$ be a constant whose value will be chosen.
The estimates for smaller balls and lower derivative orders follow by restriction.
Choose $0<\eta<0.1$ sufficiently small depending on $A,m,J$.
Choose $t_0,\td\tau\in\IR$, set $\td t:=t_0-e^{-\td\tau}$, and let $p\in\MM_{\td t}$ be such that $p(t)$ is a center of an $(n,k,\delta,\Gamma)$-neck at scale $\sqrt{t_0-t}$ for every $t\leq\td t$.
For Assertion~\ref{prop:suggested_gauge_geometric_existence}, such a point exists for sufficiently negative $\td\tau$ by asymptotic cylindricality; in Assertion~\ref{prop:suggested_gauge_geometric_control}, we use the point and endpoint in its hypotheses.
After a time shift, we may assume that $t_0=0$; we keep $\td\tau$ arbitrary.
Apply Proposition~\ref{Prop_ancient_summary} with parameters $J,m+2,\eta$, assuming $\delta\leq\ov\delta(J,m,\eta)$, and let $(\chi,\Omega,R,U^+)$ be the resulting admissible HMHF data based at $p$ on $(-\infty,\td\tau]$.
By Lemma~\ref{lem:dominant_mode_independence}, these data have $T_3=-\infty$.
The radius identity and the bound $\|U^+(\tau)\|\leq\eta$ imply $R(\tau)\geq7A$ for every $\tau\leq\td\tau$ after decreasing $\eta$.
Replacing $R$ by its nonincreasing minorant $R(\tau) := \inf_{\tau' \leq \tau} R(\tau')$ preserves this lower bound.
Set
\[
u_\tau:=e^\tau(\xi_{e^{\tau/2}}\circ\chi_{-e^{-\tau}})_*g_{-e^{-\tau}}-g_{\cyl},
\qquad
X_\tau+Y_\tau:=\Delta_{g_{\cyl}+u_\tau,g_{\cyl}}\id,
\]
where the splitting into $X_\tau$ and $Y_\tau$ is with respect to $T\IR^k\oplus T(\IS^{n-k}/\Gamma)$.
For $t=-e^{-\tau}$, along worldlines we have
\begin{align*}
\partial_\tau\chi_{-e^{-\tau}}\circ\chi_{-e^{-\tau}}^{-1}
&=e^{-\tau}\partial_t\chi_t\circ\chi_t^{-1}
=e^{-\tau}\Delta_{(\chi_t)_*g_t,g_{\cyl}}\id\\
&=e^{-\tau}\Delta_{e^{-\tau}(\xi_{e^{-\tau/2}})_*(g_{\cyl}+u_\tau),g_{\cyl}}\id\\
&=(\xi_{e^{-\tau/2}})_*\bigl(\Delta_{g_{\cyl}+u_\tau,(\xi_{e^{\tau/2}})_*g_{\cyl}}\id\bigr)\\
&=(\xi_{e^{-\tau/2}})_*\bigl(\Delta_{g_{\cyl}+u_\tau,g_{\cyl}}\id\bigr)
=(\xi_{e^{-\tau/2}})_*(X_\tau+Y_\tau).
\end{align*}
The penultimate equality holds because $(\xi_{e^{\tau/2}})_*g_{\cyl}$ and $g_{\cyl}$ have the same Levi-Civita connection.

Write $a:=a(\chi,0)$, $\mathbf b:=\mathbf b(\chi) \neq 0$ (Lemma~\ref{l:dominant_linear_well_defined}), $d:=n-k$, and set
\[
\mathbf v:=\frac{d}{2\sqrt2}\mathbf b,
\qquad
\mathsf B:=\frac d4\mathbf b\mathbf b^{\mathsf T},
\qquad
\ell_\tau:=\tau+\log|\mathbf b|^2.
\]
For estimates uniform up to $\td\tau$, we will use $\ell_\tau=(\tau-\td\tau)+\log(e^{\td\tau}|\mathbf b|^2)$.
Let $T_{\mathbf c}$ denote translation by $\mathbf c$ in the Euclidean factor, and set
\[
\widehat\chi_\tau:=T_{-\ell_\tau\mathbf v}\circ\chi_{-e^{-\tau}},
\qquad
\widehat\Omega_\tau:=\widehat\chi_\tau(\Omega_{-e^{-\tau}})\subset M_{\cyl}.
\]
Proposition~\ref{p:coefficient} applies on $\IB^k_{6A}\times(\IS^{n-k}/\Gamma)$ with regularity order $m+2$ for every $\tau\leq\td\tau$.
Its remainder in \eqref{eq:coefficient_metric_expansion} is bounded by $C(A,m,J)\eta e^{1.2(\tau-\td\tau)}$ in $C^{m+2}$ on this ball.
Consequently,
\[
X_\tau(\bx,\by)=e^{\tau/2}\mathbf v-e^\tau\mathsf B\bx+\mathcal E^X_\tau(\bx,\by),
\qquad
Y_\tau=\mathcal E^Y_\tau,
\]
where, for every $\tau\leq\td\tau$,
\[
\|\mathcal E^X_\tau\|_{C^{m+1}(\IB^k_{5A}\times(\IS^{n-k}/\Gamma))}
+\|\mathcal E^Y_\tau\|_{C^{m+1}(\IB^k_{5A}\times(\IS^{n-k}/\Gamma))}
\leq C(A,m,J)\eta e^{1.2(\tau-\td\tau)}.
\]
Along worldlines we have $\partial_\tau\widehat\chi_\tau=V_\tau\circ\widehat\chi_\tau$, where $V_\tau$ is a time-dependent vector field on $\widehat\Omega_\tau$.
Viewing $\mathsf B$ as the vector field $(\mathsf B\bx,0)$ and constant Euclidean vectors as vector fields on the cylinder, we have
\begin{equation*}
V_\tau
=-\mathbf v+(T_{-\ell_\tau\mathbf v})_*
\bigl(\partial_\tau\chi_{-e^{-\tau}}\circ\chi_{-e^{-\tau}}^{-1}\bigr) =-\mathbf v+(T_{-\ell_\tau\mathbf v}\circ\xi_{e^{-\tau/2}})_*(X_\tau+Y_\tau).
\end{equation*}
Thus, writing $w_\tau(\bx):=e^{\tau/2}(\bx+\ell_\tau\mathbf v)$, we find that the Euclidean and spherical components of $V_\tau$ are
\begin{align}
V_\tau^{\IR^k}(\bx,\by)
&=-e^\tau\mathsf B(\bx+\ell_\tau\mathbf v)
+e^{-\tau/2}\mathcal E^X_\tau(w_\tau(\bx),\by),\label{eq:gauge_adjusted_velocity}\\
V_\tau^{\sph}(\bx,\by)
&=\mathcal E^Y_\tau(w_\tau(\bx),\by).\notag
\end{align}
Equation~\eqref{eq:coefficient_asymptotic_constants_bound} gives
\begin{equation}\label{eq:gauge_shift_basepoint_bound}
\sup_{s\leq\td\tau}e^{s/2}|\ell_s\mathbf v|\leq C(A,m,J)\eta(1+|\log\eta|).
\end{equation}
For $\tau'\leq\tau\leq\td\tau$ and $|\bx|\leq4A$, we have
\[
w_{\tau'}(e^{-\tau/2}\bx)
=e^{(\tau'-\tau)/2}\bx+e^{\tau'/2}\ell_{\tau'}\mathbf v\in\IB^k_{5A}.
\]
It follows that for $\tau'\leq\tau\leq\td\tau$
\[
\big((\xi_{e^{\tau/2}})_*V_{\tau'}\big)(\bx,\by)
=-e^{\tau'}\mathsf B\bx+\mathcal E^V_{\tau,\tau'}(\bx,\by),
\]
where, on $\IB^k_{4A}\times(\IS^{n-k}/\Gamma)$,
\begin{equation}\label{eq:gauge_rescaled_velocity_error}
\|\mathcal E^V_{\tau,\tau'}\|_{C^{m+1}}
\leq C(A,m,J)\eta e^{(\tau-\td\tau)/2+0.7(\tau'-\td\tau)}.
\end{equation}
Since $R(\tau')\geq7A$ for every $\tau'\leq\td\tau$, the vector field $V_{\tau'}$ is defined on $\IB^k_{4Ae^{-\tau/2}}\times(\IS^{n-k}/\Gamma)$ throughout this interval.
The bound \eqref{eq:gauge_rescaled_velocity_error} is integrable backward in $\tau'$, so for $\eta$ sufficiently small it gives a flow starting at time $-\infty$, for every $\tau\leq\td\tau$,
\begin{equation*}
G_\tau:\IB^k_{3Ae^{-\tau/2}}\times(\IS^{n-k}/\Gamma)
\longrightarrow\IB^k_{4Ae^{-\tau/2}}\times(\IS^{n-k}/\Gamma),\qquad
\partial_\tau G_\tau=V_\tau\circ G_\tau,
\end{equation*}
such that the image of $G_\tau$ contains $\IB^k_{2Ae^{-\tau/2}}\times(\IS^{n-k}/\Gamma)$ and $G_\tau\to\id$ as $\tau\to-\infty$ pointwise.
Fix $\tau\leq\td\tau$ and notice that $G_{\tau,\tau'}:=\xi_{e^{\tau/2}}\circ G_{\tau'}\circ\xi_{e^{-\tau/2}}$ satisfies
\[
\partial_{\tau'}G_{\tau,\tau'}
=\big((\xi_{e^{\tau/2}})_*V_{\tau'}\big)\circ G_{\tau,\tau'},
\qquad G_{\tau,\tau'}\longrightarrow\id\quad\text{as }\tau'\to-\infty.
\]
Integrating from $-\infty$ to $\tau$ and using \eqref{eq:gauge_rescaled_velocity_error}, we therefore obtain by standard ODE theory
\begin{equation}\label{eq:gauge_transition_coarse_bound}
\|G_{\tau,\tau'}-\id\|_{C^{m+1}(\IB^k_{3A}\times(\IS^{n-k}/\Gamma))}
\leq C(A,m,J)\eta e^{0.7(\tau'-\td\tau)}.
\end{equation}
For $(\bx,\by)\in\IB^k_{3A}\times(\IS^{n-k}/\Gamma)$, the ODE for $G_{\tau,\tau'}$ can now be written as
\begin{equation*}
\partial_{\tau'}G_{\tau,\tau'}(\bx,\by)
=-e^{\tau'}\mathsf B\bx
+\mathcal E^V_{\tau,\tau'}\bigl(G_{\tau,\tau'}(\bx,\by)\bigr)
-e^{\tau'}\mathsf B\bigl(\pi_{\IR^k}G_{\tau,\tau'}(\bx,\by)-\bx\bigr).
\end{equation*}
By \eqref{eq:gauge_transition_coarse_bound} and \eqref{eq:gauge_rescaled_velocity_error}, the last two terms satisfy
\begin{align*}
\big\|e^{\tau'}\mathsf B(\pi_{\IR^k}G_{\tau,\tau'}-\bx)\big\|_{C^{m+1}(\IB^k_{3A}\times(\IS^{n-k}/\Gamma))}
&\leq C(A,m,J)\eta e^{1.7(\tau'-\td\tau)},\\
\big\|\mathcal E^V_{\tau,\tau'}\circ G_{\tau,\tau'}\big\|_{C^{m+1}(\IB^k_{3A}\times(\IS^{n-k}/\Gamma))}
&\leq C(A,m,J)\eta e^{(\tau-\td\tau)/2+0.7(\tau'-\td\tau)}.
\end{align*}
Integrating from $-\infty$ to $\tau$ gives
\begin{align}
\|G_{\tau,\tau}-\id+e^\tau\mathsf B\|_{C^{m+1}(\IB^k_{3A}\times(\IS^{n-k}/\Gamma))}
&\leq C(A,m,J)\eta\int_{-\infty}^\tau e^{1.7(\tau'-\td\tau)}\,d\tau'\notag\\
&\quad+C(A,m,J)\eta e^{(\tau-\td\tau)/2}
\int_{-\infty}^\tau e^{0.7(\tau'-\td\tau)}\,d\tau'\notag\\
&\leq C(A,m,J)\eta e^{1.2(\tau-\td\tau)}.\label{eq:suggested_geometric_transition}
\end{align}
Differences of maps are taken in the standard local embedding of the cylinder in $\IR^{n+1}$, and $\mathsf B$ acts only on the Euclidean component.

Since $\widehat\chi_\tau$ and $G_\tau$ satisfy the same ODE, we find that $G_{\tau}^{-1} \circ \widehat\chi_\tau (q( - e^{-\tau}))$ is locally constant in $\tau$ along worldlines, as long as it is defined.
Let $\Omega^\phi \subset \MM$ be the set of points $q \in \MM$ whose worldlines extend to time $-\infty$ and such that for all sufficiently negative $\tau$ the point $\widehat\chi_\tau(q(-e^{-\tau}))$ lies in the image of $G_\tau$.
Define $\phi(q) := G_{\tau}^{-1} \circ \widehat\chi_\tau(q(-e^{-\tau}))$, which is independent of $\tau$.
For every $\tau\leq\td\tau$, the rescaled image of $\phi_{-e^{-\tau}}$ contains $\IB^k_{2A}\times(\IS^{n-k}/\Gamma)$ and we have
\begin{equation}\label{eq:suggested_geometric_metric_transition}
(\phi_{ - e^{-\tau}})_* g_{ - e^{-\tau}} = G_{\tau}^* (T_{-\ell_\tau \mathbf v})_* (\chi_{ - e^{-\tau}})_* g_{ - e^{-\tau}}.
\end{equation}
Set $\mathsf B^\flat:=g_{\IR^k}(\mathsf B\,\cdot,\cdot)$.
Using the asymptotics from \eqref{eq:coefficient_metric_expansion} in Proposition~\ref{p:coefficient} and \eqref{eq:suggested_geometric_transition}, we obtain
\begin{align*}
g_{\cyl}+\td u_\tau
&=e^\tau(\xi_{e^{\tau/2}})_*(\phi_{-e^{-\tau}})_*g_{-e^{-\tau}}\\
&=G_{\tau,\tau}^*(T_{-e^{\tau/2}\ell_\tau\mathbf v})_*(g_{\cyl}+u_\tau)\\
&=g_{\cyl}+u_\tau+(G_{\tau,\tau}^*g_{\cyl}-g_{\cyl})
 +\big((T_{-e^{\tau/2}\ell_\tau\mathbf v})_*u_\tau-u_\tau\big)\\
&\qquad+G_{\tau,\tau}^*(T_{-e^{\tau/2}\ell_\tau\mathbf v})_*u_\tau
 -(T_{-e^{\tau/2}\ell_\tau\mathbf v})_*u_\tau\\
&=g_{\cyl}+u_\tau-2e^\tau\mathsf B^\flat
 +\frac d4e^\tau|\mathbf b|^2\ell_\tau g_{\IS^{n-k}/\Gamma}
 +\mathcal R_\tau\\
&=g_{\cyl}-\frac d4e^\tau\sum_{i,j=1}^k b_ib_j\,d\bx_i\otimes d\bx_j\\
&\qquad+\bigg[\frac{e^{\tau/2}}{\sqrt2}\langle\mathbf b,\bx\rangle
 +e^\tau\bigg(a+\frac{d-2}{4}|\mathbf b|^2\tau
 +\frac d4|\mathbf b|^2\log|\mathbf b|^2\bigg)\bigg]g_{\IS^{n-k}/\Gamma}+\mathcal E^\phi_\tau.
\end{align*}
For every $\tau\leq\td\tau$, the errors satisfy
\begin{equation}\label{eq:gauge_constructed_metric_error}
\|\mathcal R_\tau\|_{C^m(\IB^k_{2A}\times(\IS^{n-k}/\Gamma))}
+\|\mathcal E^\phi_\tau\|_{C^m(\IB^k_{2A}\times(\IS^{n-k}/\Gamma))}
\leq C(A,m,J)\eta e^{1.2(\tau-\td\tau)}.
\end{equation}
Indeed, \eqref{eq:suggested_geometric_transition} gives
\[
\|G_{\tau,\tau}^*g_{\cyl}-g_{\cyl}+2e^\tau\mathsf B^\flat\|_{C^m(\IB^k_{2A}\times(\IS^{n-k}/\Gamma))}
\leq C(A,m,J)\eta e^{1.2(\tau-\td\tau)}.
\]
Since $\mathfrak p_i^{(1)}(\bx)=\bx_i/\sqrt2$, the translation adds $e^\tau\ell_\tau\langle\mathbf b,\mathbf v\rangle/\sqrt2=\frac d4e^\tau|\mathbf b|^2\ell_\tau$ to the spherical component.
The cross term satisfies
\[
\|G_{\tau,\tau}^*(T_{-e^{\tau/2}\ell_\tau\mathbf v})_*u_\tau-(T_{-e^{\tau/2}\ell_\tau\mathbf v})_*u_\tau\|_{C^m(\IB^k_{2A}\times(\IS^{n-k}/\Gamma))}
\leq C(A,m,J)\eta e^{1.5(\tau-\td\tau)},
\]
while Proposition~\ref{p:coefficient} bounds its remainder by $C(A,m,J)\eta e^{1.2(\tau-\td\tau)}$ on this ball.
The scalar coefficient in the explicit expansion above satisfies
\begin{multline*}
e^\tau\left(a+\frac{d-2}{4}|\mathbf b|^2\tau
+\frac d4|\mathbf b|^2\log|\mathbf b|^2\right) \\
=e^{\tau-\td\tau}\left(e^{\td\tau}\left(a-\frac12|\mathbf b|^2\td\tau\right)
+\frac{d-2}{4}e^{\td\tau}|\mathbf b|^2(\tau-\td\tau)
+\frac d4e^{\td\tau}|\mathbf b|^2\log(e^{\td\tau}|\mathbf b|^2)\right).
\end{multline*}
Combining this identity with \eqref{eq:coefficient_asymptotic_constants_bound} and \eqref{eq:gauge_constructed_metric_error}, we obtain
\begin{equation}\label{eq:gauge_constructed_uniform_metric}
\|\td u_\tau\|_{C^m(\IB^k_{2A}\times(\IS^{n-k}/\Gamma))}
\leq C(A,m,J)\eta(1+|\log\eta|)
\qquad(\tau\leq\td\tau).
\end{equation}
Assertion~\ref{Prop_ancient_summary_a} of Proposition~\ref{Prop_ancient_summary}, together with Equations~\eqref{eq:gauge_shift_basepoint_bound} and~\eqref{eq:gauge_transition_coarse_bound}, shows that $p\in\Omega^\phi_{\td t}$ and
\begin{equation}\label{eq:gauge_constructed_basepoint_bound}
|e^{\td\tau/2}\pi_{\IR^k}\phi(p)|
\leq C(A,m,J)\eta(1+|\log\eta|)\leq2.
\end{equation}
This proves the expansion in Assertion~\ref{prop:suggested_gauge_geometric_asymptotics} for the constructed gauge with $t_0=0$.
As $R(\tau)\to\infty$, Proposition~\ref{p:coefficient} allows the construction for this same $\chi$ on arbitrarily large fixed balls and in every derivative order, and uniqueness of the flow of $V_\tau$ makes the resulting maps compatible.
Thus $\phi$ satisfies Definition~\ref{defn_gauge}, which proves Assertion~\ref{prop:suggested_gauge_geometric_existence}.

For any other gauge $\phi'$, Lemma~\ref{l:isometry_gauge} gives $\phi'=\psi\circ\phi$, where $\psi(\bx,\by)=(Q\bx+\mathbf c,\psi''(\by))$ is a cylindrical isometry, where $Q \in O(k)$, $\mathbf c \in \IR^k$ and $\psi'' \in O(n-k+1)$.
In rescaled coordinates this map is $\psi_\tau(\bx,\by)=(Q\bx+e^{\tau/2}\mathbf c,\psi''(\by))$.
Substituting into the expansion replaces $\mathbf b$ by $Q\mathbf b$ and $a$ by $a-\langle Q\mathbf b,\mathbf c\rangle/\sqrt2$.
These constants are uniquely determined by the coefficients of $e^{\tau/2}$ and $e^\tau$, respectively.
Changing the reference time from $0$ to $t_0$ replaces $a$ by $a-t_0$ and leaves $\mathbf b$ unchanged, by the same change of variables as in the proof of Lemma~\ref{lem:HMHF_constants_reference_time}.
Undoing the initial time shift preserves the form of the expansion.
For each fixed ball and derivative order, its remainder satisfies the bound in Assertion~\ref{prop:suggested_gauge_geometric_asymptotics}, with a constant depending on the gauge and reference time.

We finally prove Assertion~\ref{prop:suggested_gauge_geometric_control}.
Let $\phi'$ denote the prescribed gauge in this assertion.
Take the point $p$ and the endpoint $\td\tau$ in the construction above to be those given in Assertion~\ref{prop:suggested_gauge_geometric_control}.
Run that construction with $A=200$.
Assuming $\eta \leq \ov\eta(m,J)$, the bounds \eqref{eq:gauge_constructed_uniform_metric} and \eqref{eq:gauge_constructed_basepoint_bound} imply the analogue of \eqref{eq:gauge_control_metric_bound} for the constructed gauge $\phi_\chi$ on $\IB^k_{400}\times(\IS^{n-k}/\Gamma)$.
Lemma~\ref{l:isometry_gauge} gives $\phi'=\psi\circ\phi_\chi$ with $\psi(\bx,\by)=(Q\bx+\mathbf c,\psi''(\by))$, and the normalization of $\phi'(p)$ implies $e^{\td\tau/2}|\mathbf c|\leq3$.
Since $|e^{\tau/2}\mathbf c|\leq3$, composing with $\psi$ transfers the bound to $\IB^k_{100}\times(\IS^{n-k}/\Gamma)$, as asserted.

Set $\chi':=\psi\circ\chi$ and $\mathbf v':=Q\mathbf v$.
Since $T_{-\ell_{\tau'}\mathbf v'}\circ\psi=\psi\circ T_{-\ell_{\tau'}\mathbf v}$, the limit relation \eqref{eq:geometric_gauge_HMHF_limit} for $(\phi_\chi,\chi)$ also holds for $(\phi',\chi')$.
It remains to prove that $\chi'$ is admissible.
The hypotheses of Lemmas~\ref{l:rotation} and~\ref{l:translation} hold throughout $(-\infty,\td\tau]$, since $e^{\tau/2}|\mathbf c|\leq3$ and $R(\tau)\geq 7A \geq 1400$.
These lemmas show that $\chi'$ is admissible with parameters $m,J,\eps,C(m,J)$.
Applying the isometry change of coefficients computed above to the expansion of $\chi$ in Proposition~\ref{p:coefficient} gives $a(\chi',0)=a(\phi',0)$ and $\mathbf b(\chi')=\mathbf b(\phi')$.
Undoing the time shift proves Assertion~\ref{prop:suggested_gauge_geometric_control} and completes the proof.
\end{proof}

\medskip

\begin{proof}[Proof of Proposition~\ref{prop:suggested_gauge_geometric_comparison}]
We may assume without loss of generality that $t_0=0$.
We first prove the following weaker statement.

\begin{Claim}\label{cl:gauge_comparison_two_sided}
If $\delta\leq\ov\delta(m,J',\eps)$, then the following holds.
Let $\tau_*\leq\td\tau$ with $t_*:=-e^{-\tau_*}<T$.
Suppose that the corresponding worldline $p(t):=\phi_t^{-1}(\phi'(p'))$ exists for all $t\leq t_*$ and that \emph{both} $p(t)$ and $p'(t)$ are centers of $(n,k,\delta,\Gamma)$-necks at scale $\sqrt{-t}$ for every $t\leq t_*$.
Then \eqref{eq:gauge_comparison_decay} holds for all $\tau\leq\tau_*$.
\end{Claim}

\begin{proof}
It suffices to establish \eqref{eq:gauge_comparison_decay} for $\tau = \tau_*$, since we may always reduce $\tau_*$.
After parabolic rescaling, we may assume that $\tau_*=0$ and $t_*=-1$; by Lemma~\ref{l:gauge_asymp}\ref{l:gauge_asymp_rescaling}, the asymptotic constants of $\phi$ and $\phi'$ remain equal.
We retain the notation for the rescaled quantities.
Fix $m,J',\eps$ as in the statement, and suppose that $\delta\leq\ov\delta(m,J',\eps)$ is small enough that Proposition~\ref{prop:suggested_gauge_geometric} is applicable; we will impose further bounds on $\delta$ in the course of this proof.
Set $M:=\max\{m+10,100\}$, choose $J\geq100$ with $J>10J'$, and let $0<\eta_0<0.1$ be a constant that we will determine later.
Assume that $\delta\leq\ov\delta(M,J,\eta_0)$ so that Assertion~\ref{prop:suggested_gauge_geometric_control} of Proposition~\ref{prop:suggested_gauge_geometric} supplies admissible HMHF data $(\chi,\Omega,R,U^+)$ and $(\chi',\Omega',R',U^{\prime,+})$ with reference time $0$ and parameters $M,J,\eta_0,C_0$, where $C_0=C_0(M,J)$.
These data are also admissible at lower regularity orders by restriction.
We choose the unprimed data on $(-\infty,0]$ and the primed data on the full model interval $(-\infty,\td\tau]$, where we restrict the latter to $(-\infty,0]$ for the comparison.
Notice that we have $a(\chi,0)=a(\chi',0)$ and $\mathbf b(\chi)=\mathbf b(\chi')=:\mathbf b$.

For $\tau\leq0$ and $t=-e^{-\tau}$, recall the metric perturbations
\[
u_\tau:=e^\tau(\xi_{e^{\tau/2}}\circ\chi_t)_*g_t-g_{\cyl},
\qquad
\td u_\tau:=e^\tau(\xi_{e^{\tau/2}}\circ\phi_t)_*g_t-g_{\cyl},
\]
and define $u'_\tau,\td u'_\tau$ analogously.
By admissibility, $u_\tau,u'_\tau$ are defined on $\IB^k_{1000}\times(\IS^{n-k}/\Gamma)$, while Proposition~\ref{prop:suggested_gauge_geometric} implies that $\td u_\tau,\td u'_\tau$ are defined on $\IB^k_{100}\times(\IS^{n-k}/\Gamma)$.
Set $q(\tau):=\|U^{\prime,+}(\tau)\|$ and $q_0:=q(0)$.
We absorb the dependence on $C_0=C_0(M,J)$ into the constants $C(m,J)$ below.
Proposition~\ref{prop:diff_of_U}, applied with the primed data indexed by $0$, together with admissibility and equivalence of norms on the fixed mode space, implies
\begin{equation}\label{eq:gauge_comparison_HMHF_bound}
\|u_\tau-u'_\tau\|_{C^{m+3}(\IB^k_{800}\times(\IS^{n-k}/\Gamma))}
\leq C(m,J)q^{J+1}(\tau).
\end{equation}
Moreover, Assertion~\ref{Prop_ancient_summary_f2} of Proposition~\ref{Prop_ancient_summary} implies $q(\tau)\leq C(m,J)e^{0.49\tau}q_0$ for $\tau\leq0$.
Since the primed data are defined up to $\td\tau\geq0$, the same assertion gives $q_0\leq C(m,J)\eta_0e^{-0.49\td\tau}$.

Set $\ell_\tau:=\tau+\log|\mathbf b|^2$ and
\[
\widehat\chi_t:=T_{-\ell_\tau\mathbf v}\circ\chi_t,
\qquad
\widehat\chi'_t:=T_{-\ell_\tau\mathbf v}\circ\chi'_t,
\qquad \mathbf v:=\frac{n-k}{2\sqrt2}\mathbf b.
\]
Applying \eqref{eq:geometric_gauge_HMHF_limit} to the preimages of $\overline{\IB}^k_{50}\times(\IS^{n-k}/\Gamma)$ under $\phi_{-1}$ and $\phi'_{-1}$ gives the maps
\begin{alignat*}{2}
P_\tau:=\widehat\chi_{-e^{-\tau}}\circ\phi_{-e^{-\tau}}^{-1}
&:\quad &\IB^k_{50}\times(\IS^{n-k}/\Gamma)
&\longrightarrow\IR^k\times(\IS^{n-k}/\Gamma),\\
P'_\tau:=\widehat\chi'_{-e^{-\tau}}\circ(\phi'_{-e^{-\tau}})^{-1}
&:\quad &\IB^k_{50}\times(\IS^{n-k}/\Gamma)
&\longrightarrow\IR^k\times(\IS^{n-k}/\Gamma)
\end{alignat*}
for all sufficiently negative $\tau$.
By \eqref{eq:geometric_gauge_HMHF_limit}, both maps converge smoothly to the identity map as $\tau\to-\infty$.
We will show that $P_\tau,P'_\tau$ remain defined for every $\tau\leq0$, with images contained in $\IB^k_{400} \times (\IS^{n-k}/\Gamma)$, and that
\begin{equation}\label{eq:gauge_comparison_transition_bound}
\|P_0-P'_0\|_{C^{m+1}(\IB^k_{50}\times(\IS^{n-k}/\Gamma))}\leq C(m,J)q_0^{J+1}.
\end{equation}
Proposition~\ref{p:coefficient} gives $|\mathbf b|\leq C(m,J)\eta_0$, so the translation satisfies $|\ell_0\mathbf v|\leq C(m,J)\eta_0(1+|\log\eta_0|)$ and is small for sufficiently small $\eta_0$.
We first establish \eqref{eq:gauge_comparison_transition_bound}.
We have
\begin{equation}\label{eq:gauge_comparison_transition_ode}
\partial_\tau P_\tau=V_\tau\circ P_\tau,
\end{equation}
where, for  $\Delta_{g_{\cyl}+u_\tau,g_{\cyl}}\id=X_\tau+Y_\tau\in T\IR^k\oplus T(\IS^{n-k}/\Gamma)$, 
\[
V_\tau^{\IR^k}(\bx,\by)=e^{-\tau/2}X_\tau(e^{\tau/2}(\bx+\ell_\tau\mathbf v),\by)-\mathbf v,
\qquad
V_\tau^{\sph}(\bx,\by)=Y_\tau(e^{\tau/2}(\bx+\ell_\tau\mathbf v),\by).
\]
Analogous identities hold for $P'_\tau$.
Since $|\mathbf b|\leq C(m,J)\eta_0$ by Proposition~\ref{p:coefficient}, we have
\[
e^{\tau/2}|\ell_\tau\mathbf v|
\leq C(m,J)\eta_0(1+|\log\eta_0|).
\]
Thus $e^{\tau/2}(\bx+\ell_\tau\mathbf v)\in\IB^k_{800}$ for $\bx\in\IB^k_{500}$, provided $\eta_0 \leq \ov\eta_0(m,J)$.
We obtain the bounds
\begin{align}
\|V_\tau-V'_\tau\|_{C^{m+1}(\IB^k_{500}\times(\IS^{n-k}/\Gamma))}
&\leq C(m,J)e^{-\tau/2}q^{J+1}(\tau),\label{eq:gauge_comparison_velocity_difference}\\
\|V_\tau\|_{C^{m+2}(\IB^k_{500}\times(\IS^{n-k}/\Gamma))}
+\|V'_\tau\|_{C^{m+2}(\IB^k_{500}\times(\IS^{n-k}/\Gamma))}
&\leq C(m,J)\eta_0 e^{-\tau/2}e^{1.2\tau}.\label{eq:gauge_comparison_velocity_bound}
\end{align}
The first bound follows from \eqref{eq:gauge_comparison_HMHF_bound}; the factor $e^{-\tau/2}$ in $V_\tau^{\IR^k}$ gives its first inequality.
The second follows via \eqref{eq:coefficient_metric_expansion}; note that this gives a leading term $e^{\tau/2}\mathbf v$, which cancels the subtraction of $\mathbf v$ in the definition of $V^{\IR^k}_\tau$.
The metric remainder contributes at most $C(m,J)\eta_0e^{-\tau/2}e^{1.2\tau}$, and the remaining explicit terms decay faster by \eqref{eq:coefficient_asymptotic_constants_bound}.

Starting sufficiently far back, solve \eqref{eq:gauge_comparison_transition_ode} for $P_\tau$ and its analogue for $P'_\tau$ forward to $0$.
The bound \eqref{eq:gauge_comparison_velocity_bound} implies a displacement bound of at most $C(m,J)\eta_0$.
Since the initial maps are close to the identity, their images remain in $\IB^k_{400}\times(\IS^{n-k}/\Gamma)$, provided $\eta_0\leq\ov\eta_0(m,J)$.
By admissibility and the translation bound above, this region lies in the image of $\widehat\chi_{-e^{-\tau}}$, where $V_\tau$ is defined, so the trajectories correspond to worldlines and $P_\tau$ extends to all $\tau\leq0$.
The same argument applies to $P'_\tau$.

Equations~\eqref{eq:gauge_comparison_velocity_difference} and~\eqref{eq:gauge_comparison_velocity_bound}, together with Gr\"onwall's inequality, now imply
\begin{align*}
\|P_0-P'_0\|_{C^{m+1}(\IB^k_{50}\times(\IS^{n-k}/\Gamma))}
&\leq C(m,J)\int_{-\infty}^0 e^{-\tau/2}q(\tau)^{J+1}\,d\tau\\
&\leq C(m,J)q_0^{J+1}\int_{-\infty}^0e^{(0.49(J+1)-1/2)\tau}\,d\tau
\leq C(m,J)q_0^{J+1}.
\end{align*}
This proves \eqref{eq:gauge_comparison_transition_bound}.
Integrating the $C^{m+2}$ bound in \eqref{eq:gauge_comparison_velocity_bound} also gives uniform $C^{m+1}$ bounds for $P_0$ and $P'_0$, since both maps converge to $\id$ as $\tau\to-\infty$.
We now use the pullback identities
\[
g_{\cyl}+\td u_0=P_0^*(T_{-\ell_0\mathbf v})_*(g_{\cyl}+u_0),
\qquad
g_{\cyl}+\td u'_0=(P'_0)^*(T_{-\ell_0\mathbf v})_*(g_{\cyl}+u'_0).
\]
Subtracting them gives
\begin{equation*}
\td u_0-\td u'_0
=P_0^*(T_{-\ell_0\mathbf v})_*(u_0-u'_0)
+\big(P_0^*-(P'_0)^*\big)(T_{-\ell_0\mathbf v})_*(g_{\cyl}+u'_0).
\end{equation*}
Equation~\eqref{eq:gauge_comparison_HMHF_bound} and the uniform bounds for $P_0$ control the first term in $C^m$ by $C(m,J)q_0^{J+1}$.
Admissibility and \eqref{eq:gauge_comparison_transition_bound} give the same bound for the second term.
Thus
\begin{equation}\label{eq:suggested_gauge_comparison_mode_bound}
\|\td u_0-\td u'_0\|_{C^m(\IB^k_{50}\times(\IS^{n-k}/\Gamma))}
\leq C(m,J)q_0^{J+1}
\leq C(m,J)\eta_0^{J+1}e^{-0.49(J+1)\td\tau}
\leq\eps e^{-J'\td\tau},
\end{equation}
provided $\eta_0\leq\ov\eta_0(m,J,\eps)$, since $J>10J'$.
This proves the claim.
\end{proof}

We now remove the extra neck hypothesis from the claim.
Fix $m,J',\eps$ as in the proposition and fix $\delta_0 \leq \ov\delta (m,J',\eps)$ such that the claim applies with $\delta$ replaced by $\delta_0$.
Assume moreover that $\delta\leq\delta_0$; we will impose further bounds on $\delta$ later.
Set $z:=\phi'(p')$ and define $p(t):=\phi_t^{-1}(z)$ for sufficiently negative $t$; by the gauge asymptotics, $p(t)$ is a center of an $(n,k,\frac12\delta_0,\Gamma)$-neck for $t\ll0$.
Let $C_1>0$ be fixed below, and let $t^*\leq\min\{-e^{C_1-\td\tau},T\}$ be maximal such that $p(t)$ exists and is a center of an $(n,k,\frac12\delta_0,\Gamma)$-neck for every $t<t^*$.
If $t^*<\min\{-e^{C_1-\td\tau},T\}$, then the claim applies at $t^*$ and gives \eqref{eq:gauge_comparison_decay} for all $\tau\leq\tau_*$, where $\tau_*:=-\log(-t^*)$.

Since $p'(t)$ is a center of an $(n,k,\delta,\Gamma)$-neck, \eqref{eq:gauge_control_metric_bound} and \eqref{eq:gauge_comparison_decay} allow us to verify the assumption of Corollary~\ref{cor:local_neck_detection} for $p(t)$ so that the left-hand side of \eqref{eq:local_neck_detection_comparison} is bounded by $C(m,J',\eps)(\delta+\eps e^{-J'C_1})$ for all $t\leq t^*$ with $t<T$.
Assuming $C_1 \geq \underline C_1(m,J',\eps)$ and $\delta\leq\ov\delta(m,J',\eps)$, this bound is less than $\ov\delta'(\frac14\delta_0)$ from Corollary~\ref{cor:local_neck_detection}.
The corollary therefore implies that $p(t)$ is a center of an $(n,k,\frac14\delta_0,\Gamma)$-neck at these times, so by the maximal choice of $t^*$ we must have $t^*=\min\{-e^{C_1-\td\tau},T\}$, and the claim proves \eqref{eq:gauge_comparison_decay} for every $\tau\leq\td\tau-C_1$ with $t<T$.
\end{proof}

\medskip

\begin{proof}[Proof of Lemma~\ref{l:gauge_asymp}]
The defining gauge properties are preserved by these transformations. For a target isometry, the normalized perturbation transforms by $(S_\tau)_*$, where $S_\tau(\bx,\by)=(A\bx+e^{\tau/2}\bx_0,\rho(\by))$. In the expansion from Proposition~\ref{prop:suggested_gauge_geometric}, the rotation replaces $\mathbf b$ by $A\mathbf b$, and the translation contributes $-e^\tau\langle A\mathbf b,\bx_0\rangle g_{\IS^{n-k}/\Gamma}/\sqrt2$. This proves Assertion~\ref{l:gauge_asymp_isometry}.

For parabolic rescaling, the normalized perturbations satisfy $\td u'_{\tau'}=\td u_{\tau'+2\log\alpha}$. Substituting this identity into the expansion proves Assertion~\ref{l:gauge_asymp_rescaling}. For a change of reference time, set $\tau_i=-\log(t_i-t)$ and $\lambda^2=e^{\tau_1-\tau_0}=1-(t_1-t_0)e^{\tau_1}$. At the same physical time,
\[
\td u^{(1)}_{\tau_1}=(\lambda^2-1)g_{\IS^{n-k}/\Gamma}
+\lambda^2(\xi_\lambda)_*\td u^{(0)}_{\tau_0}.
\]
The linear spherical term is unchanged, and the background term adds $(t_0-t_1)e^{\tau_1}g_{\IS^{n-k}/\Gamma}$. Comparing coefficients proves \eqref{eq:gauge_constants_reference}. Shifting the time coordinate and reference time together leaves the normalized metrics unchanged, which proves the remaining formulas.
\end{proof}
\bigskip

\subsection{Comparison with the Bryant solution times a Euclidean factor}\label{subsec:linear_bryant_comparison}
In this subsection, we prove the polynomial version of Theorem~\ref{thm:linear_bryant_asymptotics}. After matching the asymptotic constants, Proposition~\ref{prop:suggested_gauge_geometric_comparison} gives the required comparison using cylindrical control on the model. We will show that this includes the entire Bryant end, outside a region of uniformly bounded radius.

\begin{Proposition}\label{p:bryant_infinity}
Suppose we are in the setting of Theorem~\ref{thm:linear_bryant_asymptotics}, with model flow $\MM'$ and radial function $y$ as defined there. There is $Y_0=Y_0(n,k,\Gamma)>0$ such that the following holds after parabolically rescaling $\MM$.
There is a time- and $\partial_{\mathfrak t}$-preserving diffeomorphism onto its image
\[
\psi:\MM'_{<T,\,y\geq Y_0}\longrightarrow\MM,
\]
defined on an open neighborhood of the displayed domain, with the following properties:
\begin{enumerate}[label=(\alph*)]
\item\label{p:bryant_infinity_b} For integers $m,J\geq10$ and $0\leq j\leq m$,
\begin{equation}\label{eq:polynomial_comparison_local}
\big|\nabla_{g'}^j(\psi^*g-g')\big|_{g'}\leq C(m,J)R^{J+1}
\qquad\textnormal{on }\MM'_{<T,\,y\geq Y_0},
\end{equation}
where $R$ is the scalar curvature of $g'$.
\item\label{p:bryant_infinity_c} For sufficiently small $\eps=\eps(n,k,\Gamma)>0$, the image of $\psi$ contains every point $p\in\MM$ for which there is $t_p>\mathfrak t(p)$ such that its backward worldline is defined for all $t\leq\mathfrak t(p)$ and $p(t)$ is the center of an $(n,k,\eps,\Gamma)$-neck at scale $\sqrt{t_p-t}$ for every such $t$.
\end{enumerate}
\end{Proposition}

\begin{proof}
We first fix the comparison map $\psi$. By Lemma~\ref{lem:dominant_mode_independence}, the model $\MM'$ has dominant linear mode, since the $k$ smallest eigenvalues of the Ricci curvature decay like $O((-t)^{-2})$ along every non-tip backward worldline. Choose gauges at infinity $\phi:\Omega\to M_{\cyl}$ and $\phi':\Omega'\to M_{\cyl}$ on $\MM,\MM'$ by Proposition~\ref{prop:suggested_gauge_geometric}. Their vectors $\mathbf b(\phi)$ and $\mathbf b(\phi')$ are nonzero by the same proposition. Lemma~\ref{l:gauge_asymp} therefore allows us to make them equal by parabolically rescaling $\MM$ and rotating its gauge. A translation in their common direction then makes $a(\phi,0)$ and $a(\phi',0)$ equal; see again Lemma~\ref{l:gauge_asymp}. Fix these choices for the rest of the proof. For $t<T$, write $V_t:=\phi_t(\Omega_t)\cap\phi'_t(\Omega'_t)$ and set
\[
W_t:=(\phi_t)^{-1}(V_t),\qquad
W'_t:=(\phi'_t)^{-1}(V_t),\qquad
\psi_t:=(\phi_t)^{-1}\circ\phi'_t:W'_t\longrightarrow W_t.
\]
The sets $W\subset\Omega$ and $W'\subset\Omega'$ are open and contain the entire backward worldline through each of their points. Corresponding worldlines remain in these sets for as long as both are defined, and $\psi:W'\to W$ preserves time and the time vector fields. By Lemma~\ref{l:gauge_asymp}, the constants agree at every common reference time.

We must show that $\MM'_{<T,\,y\geq Y_0}\subset W'$, that the estimate in Assertion~\ref{p:bryant_infinity_b} of the proposition holds there, and that those points characterized in Assertion~\ref{p:bryant_infinity_c} lie in the image of this end. For the metric estimate, fix $m,J\geq10$. We will first allow $Y_0$ to depend on $m,J$ and then show that it can be chosen uniformly.

The following claim gives the comparison on a sufficiently far-out part of the model.
\begin{Claim}\label{cl:polynomial_comparison_local}
For every $m,J\geq10$, there is $Y_{m,J}>0$ such that $\MM'_{<T,\,y\geq Y_{m,J}}\subset W'$ and \eqref{eq:polynomial_comparison_local} holds there for an appropriate choice of $C(m,J)$.
\end{Claim}

\begin{proof}
Choose $\delta>0$ sufficiently small that Proposition~\ref{prop:suggested_gauge_geometric_comparison} applies with parameters $m,J+1$ and $\eps=0.01$, and let $C_1=C_1(m,J+1,0.01)>0$ be the constant in that proposition.
Choose $\delta$ small enough also for Assertion~\ref{prop:suggested_gauge_geometric_control} of Proposition~\ref{prop:suggested_gauge_geometric} with regularity $\max\{m,100\}$, approximation order $100$, and $\eps=0.01$.
Consider constants $Y_{m,J}>Y'_{m,J} > 0$ to be chosen in the course of the proof.
Let $p'\in\MM'_{t_1}$ with $t_1<T$ and $y(p')\geq Y_{m,J}$.
Follow the worldline of $p'$ forward to the time $\td t\geq t_1$ at which $y(p'(\td t))=Y'_{m,J}$ (note that $\MM'$ is defined for all times), and set
\[
t_0:=\td t+\frac{n-k}{2R(p'(\td t))},\qquad
\tau_1:=-\log(t_0-t_1),\qquad \td\tau:=-\log(t_0-\td t).
\]
For $Y'_{m,J}$ sufficiently large, the Bryant asymptotics imply that $p'(t)$ is the center of an $(n,k,\delta,\Gamma)$-neck at scale $\sqrt{t_0-t}$ for every $t\leq\td t$, and that $c\leq R(p'(t))(t_0-t)\leq C$ on this interval.
Choosing $Y_{m,J}$ sufficiently large depending on $Y'_{m,J}$ and $C_1$, we also have $\td\tau-\tau_1\geq C_1$ and $e^{\tau_1-\td\tau}\leq C(m,J)R(p')$.

Proposition~\ref{prop:suggested_gauge_geometric} and Lemma~\ref{l:isometry_gauge} give $p'(\td t)\in\Omega'$. 
Using Lemma~\ref{l:gauge_asymp}, we can translate both gauges by the same Euclidean vector so that $\phi'(p')\in\{0\}\times(\IS^{n-k}/\Gamma)$.
This preserves the identity $\psi = \phi^{-1} \circ \phi'$ and equality of the asymptotic constants $a$ and $\mathbf b$. Proposition~\ref{prop:suggested_gauge_geometric_comparison} therefore implies that $p'\in W'$ and
\begin{multline*}
\big\| e^{\tau_1} (\xi_{e^{\tau_1/2}} \circ \phi_{t_1})_* g_{t_1} - e^{\tau_1} (\xi_{e^{\tau_1/2}} \circ \phi'_{t_1})_* g'_{t_1} \big\|_{C^m(\IB^k_{50}\times(\IS^{n-k}/\Gamma))}
\\
\leq 0.01e^{(J+1)(\tau_1-\td\tau)}
\leq C(m,J)R(p')^{J+1}.
\end{multline*}

Since $\phi_{t_1}\circ\psi_{t_1}=\phi'_{t_1}$, pulling back the preceding estimate and undoing the rescaling gives \eqref{eq:polynomial_comparison_local}.
The bound in Assertion~\ref{prop:suggested_gauge_geometric_control} of Proposition~\ref{prop:suggested_gauge_geometric}, applied to $\phi'$ with regularity $\max\{m,100\}$ and approximation order $100$, gives the uniform comparison between the normalized model metric and $g_{\cyl}$ needed to pass from the coordinate $C^m$ estimate to the intrinsic estimate.
\end{proof}

Choose $Y_0\geq2Y_{10,10}$. 
The claim gives $\MM'_{<T,\,y\geq Y_0/2}\subset W'$, and $\psi^*g$ and $g'$ are uniformly equivalent on this region. For each $m,J$, it also gives the required estimate on $\{y\geq Y_{m,J}\}$, after increasing $Y_{m,J}$ to be at least $Y_0$. We will remove this dependence on $m,J$ at the end of the proof.

We next prove Assertion~\ref{p:bryant_infinity_c} of Proposition~\ref{p:bryant_infinity}. Let $p$ be as described there, where $\eps>0$ will be chosen below. Proposition~\ref{prop:suggested_gauge_geometric} and Lemma~\ref{l:isometry_gauge} imply $p\in\Omega$. Since the gauges $\phi, \phi'$ are constant along worldlines and the subsets $V_t$ form an exhaustion, we have $p \in W$.
Set $p':=\psi^{-1}(p)\in W'$.
Suppose, for a contradiction, that $y(p')\leq Y_0$, and let $t^*\leq\mathfrak t(p)$ be the first time at which $y(p'(t^*))=Y_0$.
Since the model $\MM'$ is well understood, the $k$-th smallest eigenvalue of the Ricci curvature at $p'(t^*)$ is bounded below by a constant of the form $c(Y_0)R(p'(t^*))$.
The estimate \eqref{eq:polynomial_comparison_local} with $m=J=10$ has $C^2$ error at most $CR^{11}$.
Thus, after possibly increasing $Y_0$, a similar statement is true at $p(t^*)$, which contradicts the fact that $p(t^*)$ is a center of an $(n,k,\eps,\Gamma)$-neck if $\eps \leq \ov\eps(Y_0)$.

Finally, we remove the dependence of $Y_0$ on $m,J$. 
For this, note that $\widehat g:=\psi^*g$ is a Ricci flow on $W' \subset \MM'$. 
The low-order estimates on $\{y\geq Y_0/2\}$ and local Shi estimates give uniform bounds for all curvature derivatives of $\widehat g$ and $g'$ on $\{y\geq Y_0\}$. 
For fixed $m,J$, every point in the remaining strip $Y_0\leq y\leq Y_{m,J}$ can be followed backward to $\{y\geq Y_{m,J}\}$ within some uniformly bounded time. 
At that earlier time, the estimates from Claim~\ref{cl:polynomial_comparison_local} bound $\widehat g - g'$ and its derivatives through order $m$. 
Integrating the evolution equations for the metrics and the difference of their Levi-Civita connections gives, by induction, $|\nabla_{g'}^j(\hat g-g')|_{g'}\leq C(m,J)$ for $j\leq m$ throughout this time-interval. Since $R\geq c(m,J)>0$ on the remaining strip, increasing $C(m,J)$ proves Assertion~\ref{p:bryant_infinity_b} of the proposition on the fixed region $\{y\geq Y_0\}$.
\end{proof}
\bigskip

\subsection{Classification of $\kappa$-solutions with dominant linear mode}\label{subsec:linear_classification}
In the case of non-negative curvature operator, the polynomial decay towards a $\IR^{k-1} \times \MM_{\Bry}^{n-k+1}$ established in Proposition~\ref{p:bryant_infinity}, modulo rescaling, can be upgraded to isometry.

\begin{proof}[Proof of Corollaries~\ref{cor:complete_linear_bryant} and \ref{thm:unique_bryant}]
Fix $Y_0>0$ and the diffeomorphism $\psi:\IR^{k-1} \times \MM_{\textnormal{Bry}, y \geq Y_0}^{n-k+1}\to W \subset \MM$ from Proposition~\ref{p:bryant_infinity}.

\begin{Claim} \label{Cl_C0_bound}
There is a constant $C_0 > 1$ such that any point in $M$ has distance smaller than $C_0$ from $W$.
\end{Claim}

\begin{proof}[Proof of Claim \ref{Cl_C0_bound}]
Fix $p \in M \setminus W$ and let $q \in W$ be its closest point.
Let $\gamma : [0,d] \to M$ be an arclength minimizing geodesic between $p$ and $q$ and fix $\eps > 0$.
Then $q$ is also the closest point to $p_\eps := \gamma(\eps)$ and $q$ does not lie in the cut locus of $p_\eps$.
So the distance sphere $S(p_\eps, d-\eps)$ is smooth near $q$ and tangent to $\partial W$, and the open ball $B(p_\eps, d-\eps)$ is disjoint from $W$.
So the mean curvature of $S(p_\eps, d-\eps)$ at $q$ cannot be smaller than that of $\partial W$ at $q$.
By comparison geometry, the former is $\leq \frac{n-1}{d-\eps}$ and by the characterization from Proposition~\ref{p:bryant_infinity}, the latter is $\geq c(Y_0)$ for some uniform $c(Y_0)> 0$.
It follows that $\frac{n-1}{d-\eps} \geq c$ for any $\eps > 0$, which implies the desired upper bound on $d$.
\end{proof}

\begin{Claim}
If $k=1$, then $\MM$ must be isometric to $\MM_{\textnormal{Bry}}^{n}$.
\end{Claim}

\begin{proof}
In the setting of Corollary~\ref{cor:complete_linear_bryant}, Claim~\ref{Cl_C0_bound} and the asymptotic characterization in Proposition~\ref{p:bryant_infinity} show that every time-slice is the union of a compact subset and an end diffeomorphic to $\IR_+ \times (\IS^{n-k}/\Gamma)$.
The flow of the time vector field $\partial_{\mathfrak t}$ therefore gives diffeomorphisms between the time-slices, so $\MM$ is a conventional Ricci flow.
Claim~\ref{Cl_C0_bound}, Proposition~\ref{p:bryant_infinity}, and standard volume comparison also show that this flow is $\kappa$-noncollapsed at all scales for some uniform $\kappa>0$.
Hamilton's Harnack inequality~\cite[Corollary~1.2]{Hamilton_Harnack} shows that the scalar curvature is pointwise nondecreasing, yielding bounded curvature on compact time-intervals; hence $\MM$ is a $\kappa$-solution.

If $n=3$, the claim follows from Brendle's classification~\cite[Theorem~1.3]{Brendle_3dBryant_unique} of noncompact ancient $\kappa$-solutions (see also \cite{MR4319063} for an alternative approach).
If $n\geq4$, then we can apply Brendle--Naff's classification~\cite[Corollary~1.6]{Brendle_Naff_ancient} once we have shown that $\MM$ is uniformly $\textnormal{PIC}$.
Suppose for a contradiction that this condition fails.
Then there are points $x_i\in\MM_{t_i}$ for which the ratio of the smallest isotropic curvature to $|\Rm|(x_i)$ tends to $0$ as $i\to\infty$.
Since the cylinder $\IR\times(\IS^{n-1}/\Gamma)$ is uniformly $\textnormal{PIC}$, Claim~\ref{Cl_C0_bound} and the asymptotic estimate from Proposition~\ref{p:bryant_infinity} imply that the points $x_i$ remain within uniformly bounded distance of $\partial W$.
After passing to a subsequence, the flows $\MM$, time-shifted by $-t_i$ and based at $x_i$, therefore converge smoothly to an ancient $\kappa$-solution with the same Bryant asymptotics and with vanishing isotropic curvature at a point at time $0$.
The strong maximum principle~\cite[Theorem~12.53]{Chow_book_series_Part_II} then implies that the curvature operator of the limit has a zero eigenvalue everywhere, contradicting its Bryant asymptotics.
(Alternatively, the recent preprint of Hebbar~\cite{Hebbar_cylindrical_kappa} gives the uniformly $\textnormal{PIC}$ condition, after which we can apply~\cite[Corollary~1.6]{Brendle_Naff_ancient}.)
\end{proof}

Assume now that $k>1$.
We will show that $\MM$ is isometric to a product $\IR^{k-1}\times\MM'$, where $\MM'$ is an $(n-k+1)$-dimensional Ricci flow spacetime.
The factor $\MM'$ is then asymptotic to $\IR\times(\IS^{n-k}/\Gamma)$ as $t\to-\infty$ and satisfies the pseudolocality property, possibly with a different constant.
The decay of the small Ricci eigenvalues of $\MM$ passes to $\MM'$, so Lemma~\ref{lem:dominant_mode_independence} shows that $\MM'$ also has dominant linear mode.
The case $k=1$ therefore identifies $\MM'$ with $\MM^{n-k+1}_{\Bry}$ and proves the desired classification.
\medskip

To show that $\MM$ splits off an $\IR^{k-1}$-factor, it suffices to establish the splitting for each time-slice.
The strong maximum principle~\cite[Theorem~12.53]{Chow_book_series_Part_II} then ensures that the splitting holds globally.
We therefore fix a time $t<T$, set $(M,g):=(\MM_t,g_t)$, and write $W_t$ as $W$.

\begin{Claim} \label{Cl_C1_bound}
There is a constant $C_1 > 0$ such that for any $A > 0$ there is an $A' (A) > 0$ with the following property.
Suppose that $p \in W$ such that $B^M(p,C_1) \subset W$.
Then
\[ B^M(p,A) \cap W \subset B^W(p, A'). \]
Here $B^M(p,A)$ denotes the radius-$A$ distance ball within $(M,g)$ while $B^W(p,A')$ denotes the radius-$A'$ distance ball with respect to the induced length metric in $(W, g|_W)$.
\end{Claim}

\begin{proof}
Set $S:=B^M(p,A)\cap W$.
We first show that $S$ is connected if $C_1$ is sufficiently large.
Otherwise, there is a component $U$ of $S$ not containing $p$, and $d(p,\cdot)$ attains its minimum on $U$ at a point $q\in U\cap\partial W$.
Writing $r:=d(p,q)$, the ball $B^M(p,r)$ lies locally outside $W$ and is tangent to $\partial W$ at $q$.
The mean curvature comparison from the proof of Claim~\ref{Cl_C0_bound} gives $r\leq C'_0$ for a uniform constant $C'_0$.
For $C_1>C'_0$, this contradicts $B^M(p,C_1)\subset W$, since $q\in\partial W$.

Choose a maximal set $q_1,\ldots,q_N\in S$ with pairwise intrinsic distances $d^W(q_i,q_j)\geq2$.
The balls $B^W(q_i,2)$ cover $S$, while the balls $B^W(q_i,1)$ are pairwise disjoint and contained in $B^M(p,A+1)$.
Proposition~\ref{p:bryant_infinity} gives a uniform positive lower bound for the volumes of the latter balls, and volume comparison bounds the volume of $B^M(p,A+1)$ from above, so $N\leq\ov N(A)$.
Since $S$ is connected, any two covering balls can be joined by a chain of at most $N$ overlapping balls.
Their centers have intrinsic distance less than $4$ whenever the balls overlap, so $d^W(p,x)<4N$ for every $x\in S$.
The claim follows with $A'(A):=4\ov N(A)$.
\end{proof}
\medskip

\begin{Claim} \label{Cl_C2_bound}
Let $C_1>0$ be from Claim \ref{Cl_C1_bound}.
There is a constant $C_2 > 0$ with the following property.
Let $p_1, p_2 \in W$ such that $B^M(p_i,C_1) \subset W$.
Then
\[  d (p_1, p_2) \leq d^W(p_1,p_2) \leq C_2 d (p_1, p_2) \]
Here $d$ and $d^W$ denote the induced distance functions on $(M,g)$ and $(W, g|_W)$, respectively.
\end{Claim}

\begin{proof}
The first bound is obvious.
We may take $C_1>1$. By Proposition~\ref{p:bryant_infinity}, moving a point of $W$ in the direction of increasing $y$ by a uniformly bounded distance gives a point $q\in W$ with $B^M(q,C_1)\subset W$. Indeed, the uniform bound for $|\nabla(y\circ\psi^{-1})|$ on $W$ implies that, once $y(\psi^{-1}(q))-(Y_0+10)$ is sufficiently large, every curve from $q$ to $\partial W$ has length greater than $C_1$. Thus, after increasing $C_0$, Claim~\ref{Cl_C0_bound} gives, for every $x\in M$, a point $q\in W$ with $d(x,q)<C_0$ and $B^M(q,C_1)\subset W$.
For the second bound we may assume $d(p_1, p_2) > 1$ and choose an arclength minimizing geodesic $\gamma :[0,d] \to M$ between $p_1$ and $p_2$, where $d = d(p_1, p_2)$.
For every integer $s \in (0,d)$ choose $q_s \in W$ with $d(\gamma(s),q_s)<C_0$ and $B^M(q_s,C_1)\subset W$, as above.
We also set $q_0 := p_1$ and $q_d := p_2$.
Then for any two $s_1, s_2 \in [0,d]$ with $|s_1-s_2| \leq 1$ for which $q_{s_1}$ and $q_{s_2}$ are defined we have
\[ d(q_{s_1}, q_{s_2}) \leq d(q_{s_1}, \gamma(s_1)) + d(\gamma(s_1), \gamma(s_2)) + d(\gamma(s_2), q_{s_2}) < 2C_0 + 1.\]
It follows, using Claim~\ref{Cl_C1_bound}, that $$q_{s_2} \in B^M(q_{s_1},2C_0+1) \cap W \subset B^W(q_{s_1},A'(2C_0+1)),$$ thus $d^W(q_{s_1},q_{s_2}) \leq A'(2C_0+1)$ and we obtain
\[ d^W(p_1,p_2) \leq d^W(q_0,q_1) + d^W(q_1, q_2) + \ldots
\leq \lceil d \rceil \cdot A'(2C_0+1) 
\leq 2  A'(2C_0+1) \cdot d(p_1,p_2). \]
This proves the claim for $C_2 := 2  A'(2C_0+1)$.
\end{proof}
\medskip

Consider now the tangent cone at infinity $(X,d)$ of $(M,g)$.
Recall that this is an Alexandrov space and the Gromov-Hausdorff limit of rescalings $(M,\la_i d_g,p)$ for any fixed $p \in M$ and any sequence $\la_i \to 0$.
The asymptotic characterization from Proposition~\ref{p:bryant_infinity} combined with Claim~\ref{Cl_C2_bound} implies that we can take $X = \IR^{k-1} \times \IR_{\geq 0}$ such that $d$ is bilipschitz to the Euclidean metric and locally isometric to it on $\IR^{k-1} \times \IR_{+}$.
So the regular part of $(X,d)$ is equal to $\IR^{k-1} \times \IR_{+}$.
Since the regular part is convex by \cite{Petrunin_1998}, it follows that $(X,d)$ is, in fact, isometric to the upper half space.
It now follows from standard Alexandrov geometry arguments (see e.g.  \cite{Burago_Burago_Ivanov_2001}) that $(M,g)$ splits off an $\IR^{k-1}$-factor.
This finishes the proof.
\end{proof}
\bigskip

\subsection{Upgrade to exponential convergence}\label{subsec:linear_exponential}

In this subsection we strengthen the arbitrary-order polynomial estimate in
Proposition~\ref{p:bryant_infinity} to exponential convergence.  The main
ingredient is Lemma~\ref{l:relative_grafting_one_step}, which provides one
step of the induction.
It is modeled on the grafting constructions of Lai
\cite[Lemma~6.5]{Lai_O2_steady} and Bamler--Kleiner
\cite[Claim~5 in the proof of Proposition~12.22]
{Bamler_Kleiner_uniqueness}, together with Lai's harmonic-map-heat-flow
argument in \cite[Theorem~6.7]{Lai_O2_steady}.

The lemma compares two Ricci flows close to the cylindrical flow on a
product ball of radius proportional to $L$, for
$-1\leq\tau\leq-e^{-1}$.  Given an initial diffeomorphism under which
the two metrics are close, it produces a final diffeomorphism on a
smaller ball.  The metric error and the difference of the map from the
identity, measured in the indicated $C^m$-norms, increase by at most a
fixed factor $C_m$, with an additional error $C_me^{-\beta L^2}$.

We use the grafting construction in
Subsection~\ref{subsec_modified_hmhf}, in particular
Claim~\ref{Cl_graft}, to extend both metrics to the complete cylinder,
using the same cutoff.  The extended metrics agree with the original
metrics on a smaller ball and with the exact cylindrical flow outside
a larger ball.  The errors introduced by this construction and by
cutting off the initial map to the identity are supported a distance
of order $L$ from the ball on which we estimate.  Their contribution
is therefore bounded by $Ce^{-\beta L^2}$ by a Gaussian estimate.

In Proposition~\ref{p:bryant_exponential}, we verify that the lemma applies
to the Bryant flow at scale $\rho$ with $L=a\sqrt\rho$, where
$a=a(m)>0$ is a small constant.  We then begin the iteration
$\ell$ steps earlier at a scale at least $2^\ell\rho$.  The order-$N$
polynomial estimate there contributes at most
$C_m^\ell(2^\ell\rho)^{-N}=\rho^{-N}(C_m2^{-N})^\ell$, because
each step multiplies the preceding error by at most $C_m$.
Choosing $N$ large makes this tend to zero as $\ell\to\infty$,
whereas the sum of the new Gaussian errors remains exponentially small.
The diffeomorphism is removed only after this iteration.

A related polynomial-to-exponential argument appears in
Deruelle--Schulze--Simon
\cite[Sections~5--6]{Deruelle_Schulze_Simon_2022b}.  They use
arbitrary-order polynomial decay to eliminate the initial contribution
and an off-diagonal Gaussian estimate to control the cutoff errors
exponentially.  By comparison,
Lemma~\ref{l:relative_grafting_one_step} permits a nonzero input error
and propagates it with the fixed factor $C_m$.

\smallskip

For a map $\Psi$ close to the identity and a metric $h$ on a
domain $U$, the norm $\|\Psi-\id\|_{C^q_h(U)}$ denotes the
$C^q$-norm of the vector field
$x\mapsto\exp_{x,h}^{-1}\Psi(x)$, where the derivatives and norms
are taken with respect to $h$.

\begin{Lemma}
\label{l:relative_grafting_one_step}
There is a dimensional constant $\beta>0$ with the
following property.  For every integer $m\geq10$, there are
$\varepsilon_m>0$ and $C_m,L_m<\infty$ such that the following
holds for all $L\geq L_m$.
Let $g_\tau,g'_\tau$, $-1\leq\tau\leq-e^{-1}$, be Ricci flows
on $\IB^k_{\frac52L}\times(\IS^{n-k}/\Gamma)\subset\IR^k\times(\IS^{n-k}/\Gamma)$.  Suppose that
\begin{equation}\label{eq:cap_cylinder_smallness}
 \sup_{-1\leq\tau\leq-e^{-1}}
 \left(\|g_\tau-g_{\cyl,\tau}\|_{C^{m+6}_{g_{\cyl,\tau}}
 (\IB^k_{\frac52L}\times(\IS^{n-k}/\Gamma))}
 +\|g'_\tau-g_{\cyl,\tau}\|_{C^{m+6}_{g_{\cyl,\tau}}
 (\IB^k_{\frac52L}\times(\IS^{n-k}/\Gamma))}\right)<\varepsilon_m.
\end{equation}
Assume that there is a smooth map $\Phi:\mathcal B_1:=\IB^k_{\frac45L}\times(\IS^{n-k}/\Gamma)\to \IB^k_{\frac52L}\times(\IS^{n-k}/\Gamma)$ which is a diffeomorphism onto the image
satisfying
\begin{equation}\label{eq:cap_step_input}
 \|\Phi^*g'_{-1}-g_{-1}\|_{C^{m+2}_{g_{\cyl,-1}}(\mathcal B_1)}
 +\|\Phi-\id\|_{C^{m+1}_{g_{\cyl,-1}}(\mathcal B_1)}
 <\varepsilon_m.
\end{equation}
Then there is a smooth map $\Phi^1:\mathcal B_0:=\IB^k_{\frac35L}\times(\IS^{n-k}/\Gamma)\to \IB^k_{\frac52L}\times(\IS^{n-k}/\Gamma)$ which is a diffeomorphism onto the image
such that
\begin{equation}\label{eq:cap_one_step}
\begin{split}
 &\|(\Phi^1)^*g'_{-e^{-1}}-g_{-e^{-1}}\|_
 {C^{m+2}_{g_{\cyl,-e^{-1}}}(\mathcal B_0)}
 +\|\Phi^1-\id\|_{C^{m+1}_{g_{\cyl,-e^{-1}}}(\mathcal B_0)}\\
 &\quad\leq C_m\left(
 \|\Phi^*g'_{-1}-g_{-1}\|_{C^{m+2}_{g_{\cyl,-1}}(\mathcal B_1)}
 +\|\Phi-\id\|_{C^{m+1}_{g_{\cyl,-1}}(\mathcal B_1)}
 +e^{-\beta L^2}\right).
\end{split}
\end{equation}
\end{Lemma}

\begin{proof}
Let $\omega$ be a smooth cutoff function as in
Claim~\ref{Cl_graft}, equal to
one on $\IB^k_{\frac52L-1}$ and compactly supported in
$\IB^k_{\frac52L}$, with
$|\nabla^j\omega|\leq C_j$.  Set
\[
 h_\tau=\omega g_\tau+(1-\omega)g_{\cyl,\tau},
 \qquad
 h'_\tau=\omega g'_\tau+(1-\omega)g_{\cyl,\tau}.
\]
These are complete metrics on $\IR^k\times(\IS^{n-k}/\Gamma)$, with
uniformly bounded geometry.  Their bounds through order $m+6$
depend only on $m$, not on $L$.  The tensors
\[
 \mathcal E_\tau=\partial_\tau h_\tau+2\Ric(h_\tau),
 \qquad
 \mathcal E'_\tau=\partial_\tau h'_\tau+2\Ric(h'_\tau)
\]
are supported where $\frac52L-1\leq|x|\leq\frac52L$.  By the assumption \eqref{eq:cap_cylinder_smallness} we have
\begin{equation}\label{eq:cap_completed_bounds}
 \sup_{-1\leq\tau\leq-e^{-1}}
 \left(\|h'_\tau-h_\tau\|_{C^{m+4}_{h_\tau}(\IR^k\times(\IS^{n-k}/\Gamma))}
 +\|\mathcal E_\tau\|_{C^{m+2}_{h_\tau}(\IR^k\times(\IS^{n-k}/\Gamma))}
 +\|\mathcal E'_\tau\|_{C^{m+2}_{h_\tau}(\IR^k\times(\IS^{n-k}/\Gamma))}\right)
 \leq C_m\varepsilon_m.
\end{equation}

Choose a smooth cutoff $\zeta$, equal to one on
$\IB^k_{\frac34L}\times(\IS^{n-k}/\Gamma)$ and compactly supported in
$\mathcal B_1$, with $|\nabla^j\zeta|\leq C_j$.  Define a map $\overline\Phi:\mathcal B_1\to \IR^k\times(\IS^{n-k}/\Gamma)$ by
\[
 \overline\Phi(x)
 =\exp_{x,g_{\cyl,-1}}
 \left(\zeta(x)\exp_{x,g_{\cyl,-1}}^{-1}\Phi(x)\right),
\]
and set $\overline\Phi=\id$ outside $\mathcal B_1$.
By \eqref{eq:cap_step_input}, this is a diffeomorphism of the
complete cylinder, and its $C^{m+1}$-difference from the identity
is bounded by $C_m\|\Phi-\id\|_{C^{m+1}_{g_{\cyl,-1}}(\mathcal B_1)}$.
\begin{Claim}\label{Cl_cap_hmhf_existence}
If $\varepsilon_m$ is sufficiently small, there is a smooth family of HMHFs by diffeomorphisms $\Theta_\tau$ on $\IR^k\times(\IS^{n-k}/\Gamma)$,
$-1\leq\tau\leq-e^{-1}$, satisfying
\begin{equation}\label{eq:cap_relative_hmhf}
 \partial_\tau\Theta_\tau=\triangle_{h'_\tau,h_\tau}\Theta_\tau,
 \qquad \Theta_{-1}=\overline\Phi^{-1}.
\end{equation}
Set $v_\tau=(\Theta_\tau^{-1})^*h'_\tau-h_\tau$ and
$F_\tau=(\Theta_\tau)_*\mathcal E'_\tau-\mathcal E_\tau$.
For a constant $C_m$ independent of $L$, we have, with all norms
taken on $\IR^k\times(\IS^{n-k}/\Gamma)$,
\begin{equation}\label{eq:cap_higher_map_bounds}
 \sup_{-1\leq\tau\leq-e^{-1}}
 \left(\|\Theta_\tau-\id\|_{C^0_{h_\tau}}
 +\|v_\tau\|_{C^{m+2}_{h_\tau}}
 +\|F_\tau\|_{C^{m+2}_{h_\tau}}\right)
 \leq C_m\varepsilon_m.
\end{equation}
\end{Claim}

\begin{proof}[Proof of Claim~\ref{Cl_cap_hmhf_existence}]
Since $\Phi$ is an isometry from $\Phi^*g'_{-1}$ to $g'_{-1}$, we have
\[
d\Phi(\nabla^{\Phi^*g'_{-1}}_{\partial_i}\partial_j)
=\nabla^{g'_{-1}}_{d\Phi(\partial_i)}d\Phi(\partial_j),
\]
which implies in local coordinates that
\[
 \partial_i\partial_j\Phi^\alpha
 =\Gamma(\Phi^*g'_{-1})^k_{ij}\partial_k\Phi^\alpha
 -\Gamma(g'_{-1})^\alpha_{\beta\gamma}(\Phi)
       \partial_i\Phi^\beta\partial_j\Phi^\gamma.
\]
Subtracting the corresponding identity with $\Phi=\id$, we can express
the RHS in terms of $\Phi-\id$, its first derivatives, and
the difference $\Phi^*g'_{-1}-g'_{-1}$.
Differentiating successively and using \eqref{eq:cap_step_input} and
\eqref{eq:cap_cylinder_smallness} imply
$\|\Phi-\id\|_{C^{m+3}_{g_{\cyl,-1}}(\mathcal B_1)}
\leq C_m\varepsilon_m$.
The cutoff construction gives the same bound for
$\overline\Phi-\id$ on the complete cylinder.  Therefore
$v_{-1}=\overline\Phi^*h'_{-1}-h_{-1}$ satisfies
$\|v_{-1}\|_{C^{m+2}_{h_{-1}}}\leq C_m\varepsilon_m$.

Fix $0<\eta_0<\eta_1<1$ as in the short-time existence result in \cite[Appendix~A.3]{Bamler_Kleiner_uniqueness}, where $\eta_0$ is the initial bound on the perturbation and $\eta_1$ is the continuation threshold. We may assume $\|v_{-1}\|_{C^0_{h_{-1}}}<\eta_0$ by taking $\varepsilon_m$ small.  By the estimates of $\partial_\tau h+2\Ric(h)$ and
$\partial_\tau h'+2\Ric(h')$ in \eqref{eq:cap_completed_bounds}, we can apply this result and obtain a short-time solution of
\eqref{eq:cap_relative_hmhf}.

We
first work on the interval on which
$\|v_\tau\|_{C^{m+2}_{h_\tau}}<\delta$, where
$0<\delta<\eta_1$ is a small constant to be chosen later.
To estimate $F_\tau$, we observe that
$\Theta_\tau^{-1}$ is an isometry from $h_\tau+v_\tau$ to
$h'_\tau$.  Thus
\[
 (\nabla^{h_\tau+v_\tau})^j
       \big((\Theta_\tau)_*\mathcal E'_\tau\big)
 =(\Theta_\tau)_*
       \big((\nabla^{h'_\tau})^j\mathcal E'_\tau\big).
\]
Since
$\nabla^{h_\tau+v_\tau}-\nabla^{h_\tau}
=(h_\tau+v_\tau)^{-1}*\nabla v_\tau$, we can rewrite the LHS
using $h_\tau$-covariant derivatives of
$v_\tau$ only through order $j$.
Here and below, all derivatives and tensor norms are with respect to
$h_\tau$, unless otherwise indicated.
The assumed $C^{m+2}$-estimate for $v_\tau$ and \eqref{eq:cap_completed_bounds} consequently give
\begin{equation*} \label{eq:eq_of_F}
 \|F_\tau\|_{C^{m+2}_{h_\tau}}\leq C_m\varepsilon_m.
 \end{equation*}

Put $a^{ij}=(h_\tau+v_\tau)^{ij}$.
The Ricci--DeTurck equation in
\cite[Appendix~A.1]{Bamler_Kleiner_uniqueness}, with the terms
from $\mathcal E,\mathcal E'$ retained, has the form
\begin{equation}\label{eq:cap_tensor_evolution}
 \partial_\tau v
 =a^{ab}\nabla_a\nabla_bv
 +\mathcal A_\tau(v)v
 +\mathcal Q_\tau(v)(\nabla v,\nabla v)+F.
\end{equation}
Explicitly, with curvature and covariant derivatives taken with respect to
$h_\tau$, we have
\[
\begin{aligned}
 (\mathcal A_\tau(v)v)_{ij}
 &=a^{pq}\bigl(
 R^\ell{}_{pij}v_{\ell q}+R^\ell{}_{pji}v_{\ell q}
 -R^\ell{}_{ipq}v_{\ell j}-R^\ell{}_{jpq}v_{\ell i}\bigr),\\
 \bigl(\mathcal Q_\tau(v)(\nabla v,\nabla v)\bigr)_{ij}
 &=-\frac12a^{pq}a^{rs}\bigl(
 -\nabla_i v_{pr}\nabla_jv_{qs}
 -2\nabla_rv_{ip}\nabla_qv_{js}
 +2\nabla_rv_{ip}\nabla_sv_{jq}\\
 &\hspace{32mm}
 +2\nabla_qv_{is}\nabla_jv_{pr}
 +2\nabla_iv_{pr}\nabla_qv_{js}\bigr).
\end{aligned}
\]
For $0\leq j\leq m+2$, differentiating
\eqref{eq:cap_tensor_evolution} we obtain $\partial_\tau\nabla^jv
 =a^{ab}\nabla_a\nabla_b\nabla^jv+\nabla^jF+\mathcal R_j$,
where
\[
\begin{split}
 \mathcal R_j={}&[\partial_\tau,\nabla^j]v
 +[\nabla^j,a^{ab}\nabla_a\nabla_b]v
 +\nabla^j\big(\mathcal A_\tau(v)v
       +\mathcal Q_\tau(v)(\nabla v,\nabla v)\big).
\end{split}
\]
For $j\geq1$, we have
$[\partial_\tau,\nabla]v=\nabla(\partial_\tau h_\tau)*v$ and
$[\partial_\tau,\nabla^j]v
=\sum_{r+s=j-1}\nabla^{r+1}(\partial_\tau h_\tau)*\nabla^sv$.
Similarly, by $[\nabla_a,\nabla_b]v=\Rm(h_\tau)*v$, and
repeated commutation we have
$[\nabla^j,\nabla_a\nabla_b]v
=\sum_{r+s=j}\nabla^r\Rm(h_\tau)*\nabla^sv$.
Together with the product rule, this shows that
$[\nabla^j,a^{ab}\nabla_a\nabla_b]v$ consists of terms
$\nabla^rv*\nabla^{j-r+2}v$, with $1\leq r\leq j$,
and $a*\nabla^r\Rm(h_\tau)*\nabla^sv$, with $r+s=j$.
Both commutators vanish for $j=0$.
The last term can be expanded into the forms
$\nabla^r(\mathcal A_\tau(v))*\nabla^sv$, with $r+s=j$,
and $\nabla^r(\mathcal Q_\tau(v))*\nabla^{s+1}v*\nabla^{t+1}v$,
with $r+s+t=j$.
Note that $\mathcal R_j$ contains no derivatives of $v$ above
order $j+1$.
Pairing the equation of $\nabla^jv$ with $2\nabla^jv$, we obtain
\[
\begin{split}
 (\partial_\tau-a^{ab}\nabla_a\nabla_b)|\nabla^jv|^2
 &=-2a^{ab}\langle\nabla_a\nabla^jv,\nabla_b\nabla^jv\rangle
 +2\langle\nabla^jv,\mathcal R_j+\nabla^jF\rangle
 +O_m(|\nabla^jv|^2)\\
 &\leq-2c|\nabla^{j+1}v|^2
 +C_m\sum_{r=0}^j|\nabla^rv|^2+C_m|\nabla^jF|^2.
\end{split}
\]
Seeing that
$\Delta_{h_\tau+v_\tau}f
=a^{ab}\nabla_a\nabla_bf
-\nabla v* \nabla f$ and
taking $f=|\nabla^jv|^2$, the first-order term is bounded by
$C|\nabla v||\nabla^jv||\nabla^{j+1}v|$.
So the smallness of $|\nabla v|$ and Young's inequality bound this by
$c|\nabla^{j+1}v|^2+C_m|\nabla^jv|^2$. So we obtain
\begin{equation}\label{eq:cap_scalar_squared}
\begin{split}
 (\partial_\tau-\Delta_{h_\tau+v_\tau})|\nabla^jv|^2
 &\leq-c|\nabla^{j+1}v|^2
 +C_m\sum_{r=0}^j|\nabla^rv|^2+C_m|\nabla^jF|^2.
\end{split}
\end{equation}
Setting $W=\sum_{j=0}^{m+2}|\nabla^jv|^2$, we obtain
\begin{equation}\label{eq:cap_derivative_scalar}
 (\partial_\tau-\Delta_{h_\tau+v_\tau})W
 \leq C_m W+C_m\sum_{j=0}^{m+2}|\nabla^jF|^2.
\end{equation}
The bound for $F_\tau$, the bound for $v_{-1}$ and the maximum principle imply $\sup_\tau\|v_\tau\|_{C^{m+2}_{h_\tau}}\leq C_m\varepsilon_m$.
For sufficiently small $\varepsilon_m$, this is less than
$\delta/2$, so the bootstrap closes.  In particular,
$|v_\tau|<\eta_1/2$, and the continuation statement in
\cite[Appendix~A.3]{Bamler_Kleiner_uniqueness} extends the solution
to $[-1,-e^{-1}]$.

Finally,
the HMHF equation \eqref{eq:cap_relative_hmhf} gives
$(\partial_\tau\Theta_\tau)\circ\Theta_\tau^{-1}
=\Delta_{h_\tau+v_\tau,\,h_\tau}\id$, namely
\begin{equation}\label{eq:cap_gauge_velocity}
 V^b_\tau:=\big((\partial_\tau\Theta_\tau)\circ\Theta_\tau^{-1}\big)^b
 =-\frac12(h_\tau+v_\tau)^{ij}(h_\tau+v_\tau)^{bq}
 (\nabla_iv_{jq}+\nabla_jv_{iq}-\nabla_qv_{ij}).
\end{equation}
Thus $\|V_\tau\|_{C^0_{h_\tau}}\leq C_m\varepsilon_m$.
Integrating $\partial_\tau\Theta_\tau=V_\tau\circ\Theta_\tau$
from $\Theta_{-1}=\overline\Phi^{-1}$ gives the map bound in
\eqref{eq:cap_higher_map_bounds} and proves the claim.
\end{proof}

We now use \eqref{eq:cap_derivative_scalar} to improve the bound for
$v$ on $\IB^k_{\frac23L}\times(\IS^{n-k}/\Gamma)$ in terms of the initial
error and an exponentially small term.  Continue to take derivatives
and tensor norms with respect to $h_\tau$, and set
$W=\sum_{j=0}^{m+2}|\nabla^jv|^2$.
Let $G(x,\tau;z,s)$ be the positive fundamental solution of the parabolic equation
$(\partial_\tau-\Delta_{h_\tau+v_\tau}-C_m)u=0$.
The evolving metrics $h_\tau+v_\tau$ are complete and uniformly
equivalent to the cylindrical metric.  The background bounds,
\eqref{eq:cap_higher_map_bounds}, and the Ricci--DeTurck equation
give uniform bounds for their curvature and time derivative.
Here we use only the uniform derivative bounds supplied by
\eqref{eq:cap_higher_map_bounds}, not their smallness.
Applying \cite[Corollary~26.26]{Chow_book_series_Part_III} gives
\begin{equation}\label{eq:cap_gaussian_tail}
 G(x,\tau;z,s)
 \leq \frac{C_m}{(\tau-s)^{n/2}}
 \exp\left(-\frac{d_{h_s+v_s}(x,z)^2}{C(\tau-s)}\right),
 \qquad -1\leq s<\tau\leq-e^{-1}.
\end{equation}
The constant potential contributes only the bounded factor
$e^{C_m(\tau-s)}$, which is absorbed into $C_m$.
Integrating \eqref{eq:cap_gaussian_tail} gives, for every $r\geq0$,
\[
\begin{aligned}
\int_{\IR^k\times(\IS^{n-k}/\Gamma)}G(x,\tau;z,s)\,d\mu_{h_s+v_s}(z)
&\leq C_m,\\
\int_{(\IR^k\times(\IS^{n-k}/\Gamma))\setminus B_{h_s+v_s}(x,r)}
G(x,\tau;z,s)\,d\mu_{h_s+v_s}(z)
&\leq C_m\exp\left(-\frac{r^2}{C(\tau-s)}\right).
\end{aligned}
\]
The constant in the Gaussian exponent
is dimensional; the prefactor may depend on $m$, but not on $L$.
For $-1\leq\tau\leq-e^{-1}$ and
$x\in\IB^k_{\frac23L}\times(\IS^{n-k}/\Gamma)$, the reproduction formula gives
\[
\begin{split}
 W_\tau(x)
 &\leq\int_{\IR^k\times(\IS^{n-k}/\Gamma)} G(x,\tau;z,-1)W_{-1}(z)\,
                d\mu_{h_{-1}+v_{-1}}(z)\\
 &\qquad\qquad+C_m\int_{-1}^\tau\int_{\IR^k\times(\IS^{n-k}/\Gamma)} G(x,\tau;z,s)
       \sum_{j=0}^{m+2}|\nabla^jF_s(z)|^2\,
                d\mu_{h_s+v_s}(z)\,ds.
\end{split}
\]
Denote the two terms on the right by $\mathcal I,\mathcal J$.
For $\mathcal I$, split the integral into
$\IB^k_{\frac34L}\times(\IS^{n-k}/\Gamma)$ and its complement.
On the former, $\Theta_{-1}^{-1}=\overline\Phi=\Phi$, and the
grafted metrics agree with the original metrics, so
$v_{-1}=\Phi^*g'_{-1}-g_{-1}$.
The complement is at $(h_{-1}+v_{-1})$-distance at least $cL$
from $x$, since the product radii differ by $L/12$.
\eqref{eq:cap_higher_map_bounds} gives
$W_{-1}\leq C_m$ on the complete cylinder.
Together with the preceding inequalities, this gives
\[
 \mathcal I\leq
 C_m\|\Phi^*g'_{-1}-g_{-1}\|_{C^{m+2}_{g_{\cyl,-1}}(\mathcal B_1)}^2
 +C_m\exp\left(-\frac{cL^2}{\tau+1}\right)\le C_m\|\Phi^*g'_{-1}-g_{-1}\|_{C^{m+2}_{g_{\cyl,-1}}(\mathcal B_1)}^2
 +C_me^{-cL^2}.
\]
For $\mathcal J$, $F_s$ is supported in the grafting annulus
or its image under $\Theta_s$.  By
\eqref{eq:cap_higher_map_bounds}, it and all its derivatives
vanish on $\IB^k_{\frac23L+cL}\times(\IS^{n-k}/\Gamma)$, for a uniform
constant $c>0$.
\eqref{eq:cap_higher_map_bounds} gives
$\sum_{j=0}^{m+2}|\nabla^jF_s|^2\leq C_m$.
Together with \eqref{eq:cap_gaussian_tail}, this implies
\[
 \mathcal J\leq C_m\int_{-1}^\tau
 \exp\left(-\frac{cL^2}{\tau-s}\right)\,ds
 \leq C_me^{-cL^2}.
\]
Combining these bounds and taking square roots, we obtain for some dimensional $\beta>0$ that
\begin{equation}\label{eq:cap_integer_smoothing}
 \sup_{-1\leq\tau\leq-e^{-1}}
 \|v_\tau\|_{C^{m+2}_{h_\tau}(\IB^k_{\frac23L}\times(\IS^{n-k}/\Gamma))}
 \leq C_m\left(
 \|\Phi^*g'_{-1}-g_{-1}\|_{C^{m+2}_{g_{\cyl,-1}}(\mathcal B_1)}
 +e^{-\beta L^2}\right).
\end{equation}

We finally estimate the map using \eqref{eq:cap_gauge_velocity}.
For $x\in\mathcal B_0$, the curve
$\tau\mapsto\Theta_\tau\circ\Theta_{-e^{-1}}^{-1}(x)$
ends at $x$ and satisfies
\[
 \partial_\tau(\Theta_\tau\circ\Theta_{-e^{-1}}^{-1})
 =V_\tau
       \circ(\Theta_\tau\circ\Theta_{-e^{-1}}^{-1}).
\]
These curves stay in $\IB^k_{\frac23L}\times(\IS^{n-k}/\Gamma)$ by
\eqref{eq:cap_higher_map_bounds}.  Integrating backward from
$-e^{-1}$, differentiating through order $m+1$, and applying
Gr\"onwall's inequality bounds
$\Theta_{-1}\circ\Theta_{-e^{-1}}^{-1}-\id$
by the time integral of the $C^{m+1}$-norm of
$V_\tau$.
Since $\Theta_{-1}^{-1}=\overline\Phi=\Phi$, this implies the bounds on $\Theta_{-e^{-1}}^{-1}-\Phi$.
Thus \eqref{eq:cap_gauge_velocity}, the $C^{m+2}$-bound for $v$ in
\eqref{eq:cap_integer_smoothing}, and the initial $C^{m+1}$-bound for
$\Phi-\id$ bound
$\Theta_{-e^{-1}}^{-1}-\id$ in $C^{m+1}$ on
$\mathcal B_0$ by the RHS of
\eqref{eq:cap_one_step}.  Let
$\Phi^1=\Theta_{-e^{-1}}^{-1}|_{\mathcal B_0}$.
Then
$(\Phi^1)^*g'_{-e^{-1}}-g_{-e^{-1}}=v_{-e^{-1}}$ on
$\mathcal B_0$, and the assertion follows from \eqref{eq:cap_integer_smoothing}.
\end{proof}

We now establish the exponential estimate in Theorem~\ref{thm:linear_bryant_asymptotics}.

\begin{Proposition}[Exponential convergence]
\label{p:bryant_exponential}
Let $\MM$, $Y_0$, and $\psi:\MM'_{<T,\,y\ge Y_0}\to\MM$ be as in Proposition~\ref{p:bryant_infinity}.
After increasing $Y_0$ by an amount depending only on $n,k,\Gamma$, for every integer $m\geq0$ there are constants
$c_m,C(m)>0$, depending only on $m,n,k,\Gamma$, such that for all $j=0,\cdots,m$,
\begin{equation}\label{eq:bryant_exponential}
 \big|\nabla^j_{g'}(\psi^*g-g')\big|_{g'}
 \leq C(m)e^{-c_m y}
 \qquad\text{on }\MM'_{<T,\,y\geq Y_0},
\end{equation}
where $y$ is the distance function to the tips on $\MM'$.
\end{Proposition}

\begin{proof}
We may assume without loss of generality that $\psi=\id$ and write $\psi^*g$ as $g$.
Then $g$ is defined on $\MM'_{<T,\,y\ge Y_0}$, and Assertion~\ref{p:bryant_infinity_b} of Proposition~\ref{p:bryant_infinity} gives, for every $j,P\geq0$,
\begin{equation}\label{eq:bryant_polynomial_decay}
 \big|\nabla_{g'}^j(g-g')\big|_{g'}
 \leq C(j,P)R^P
\end{equation}
on $\{y\geq Y_0\}$.
From now on we assume $m\geq10$.
Let $\mathbf z:\MM'\to\IR^{k-1}$ denote the Euclidean
coordinate and $\theta:\MM'_{y>0}\to(\IS^{n-k}/\Gamma)$
the spherical coordinate.

For any point $p\in\MM'_{<T,y\ge Y_0}$, we write $\rho(p):=e(n-k)/(2R(p))$ and define a map
$\kappa_p:\{q\in\MM'_{[\mathfrak t(p)-(1-e^{-1})\rho(p),\mathfrak t(p)]}:
y(q(\mathfrak t(p)))\ge Y_0\}\to\IR^k\times(\IS^{n-k}/\Gamma)$ by
\[
 \kappa_p(q)
 =\Big(\rho(p)^{-1/2}\big(\mathbf z(q(\mathfrak t(p)))-\mathbf z(p)\big),
 \rho(p)^{-1/2}\big(y(q(\mathfrak t(p)))-y(p)\big),\theta(q(\mathfrak t(p)))\Big).
\]
Then $\kappa_p$ is constant along worldlines.
For $\tau\in[-1,-e^{-1}]$, denote
$\kappa_{p,\tau}:=\kappa_p|_{\MM'_{\mathfrak t(p)+\rho(p)(\tau+e^{-1})}}$.

\begin{Claim}\label{Cl_cap_bryant_cylindrical_estimate}
There is a small constant $a=a(m)>0$ and a constant $\underline{Y}(m)>0$ such that the following holds for every $p\in\MM'_{<T,\,y> \underline{Y}(m)}$ with $t_0=\mathfrak t(p)$.
The two Ricci flows
\[
 h_\tau:=(\kappa_{p,\tau})_*(\rho(p)^{-1}g_{t_0+\rho(p)(\tau+e^{-1})}),
 \qquad h'_\tau:=(\kappa_{p,\tau})_*(\rho(p)^{-1}g'_{t_0+\rho(p)(\tau+e^{-1})})
\] 
satisfy \eqref{eq:cap_cylinder_smallness} in Lemma \ref{l:relative_grafting_one_step} with $L=a\sqrt{\rho(p)}\ge L_m$, for $L_m>0$ from Lemma \ref{l:relative_grafting_one_step}.
\end{Claim}

\begin{proof}[Proof of Claim~\ref{Cl_cap_bryant_cylindrical_estimate}]
Translating the time coordinates of both flows, we may assume
$t_0=-\rho/e$, so that $t_0+\rho(\tau+e^{-1})=\rho\tau$.
We abbreviate $\kappa=\kappa_p$,
$\kappa_\tau=\kappa_{p,\tau}$, and $\rho=\rho(p)$.
Fix $(\bx_0,\theta_0)\in\IB^k_{\frac52a\sqrt\rho}\times(\IS^{n-k}/\Gamma)$
and choose $q_0\in\MM'_{-\rho/e,\,y\ge Y_0}$ with
$\kappa_{-e^{-1}}(q_0)=(\bx_0,\theta_0)$.
Then
$\kappa_\tau(q_0(\rho\tau))=(\bx_0,\theta_0)$ for every
$\tau\in[-1,-e^{-1}]$, and $\bx(q_0)=\bx_0$ and
$y(q_0)=y(p)+\sqrt\rho\,x_{0,k}$.
The Bryant asymptotics $2yR/(n-k)\to1$ and
$|x_{0,k}|\leq\frac52a\sqrt\rho$ give
\begin{equation} \label{eq:y_q_0}
 |y(q_0)-y(p)|\leq\tfrac52a\rho,\qquad
 \frac{n-k}{2\rho R(q_0)}
 =\frac{R(p)}{eR(q_0)}=e^{-1}+O(a)+o(1).
\end{equation}
In particular, $y(q_0)\geq c\rho>Y_0$ for small $a$ and large
$\rho$, which also justifies the choice of $q_0$.
On the same domain of $\kappa$, we define a map to the cylinder by
\[
 \pi_{q_0}(q')
 =\Big(
 \sqrt{\tfrac{2R(q_0)}{n-k}}\big(\mathbf z(q'(-\rho/e))-\mathbf z(q_0)\big),
 \sqrt{\tfrac{2R(q_0)}{n-k}}\big(y(q'(-\rho/e))-y(q_0)\big),
 \theta(q'(-\rho/e))
 \Big).
\]
By the asymptotics of the Bryant soliton
\cite[Chapter~1, Section~4]{Chow_book_series_Part_I}, we see that the Ricci flow $(\pi_{q_0})_*
(\tfrac{2R(q_0)}{n-k}g'_{-\rho/e+\frac{n-k}{2R(q_0)}(t+1)})$ smoothly converges to the cylindrical flow $g_{\cyl,t}$ locally uniformly as $\rho\to\infty$ for all $t\in(-\infty,0)$.
Taking $t=-1+\frac{2\rho R(q_0)}{n-k}(\tau+e^{-1})$,
by \eqref{eq:y_q_0} we have $t\in[-e-1,-1]$ for sufficiently
small $a$ and all $\tau\in[-1,-e^{-1}]$, so
the convergence of the flow implies
\begin{equation} \label{eq:close_cyl_pi}
 (\pi_{q_0})_*\big(\tfrac{2R(q_0)}{n-k}g'_{\rho \tau}\big)
 =g_{\IR^k}
 +\left(1-\frac{2\rho R(q_0)}{n-k}(\tau+e^{-1})\right)g_{\IS^{n-k}/\Gamma}
 +o_r(1),\end{equation}
where $o_r(1)\to0$ locally uniformly in the $C^r$-norm on $\IR^k\times(\IS^{n-k}/\Gamma)$.
Since the two coordinate maps satisfy $\kappa_\tau\circ(\pi_{q_0}|_{\MM'_{\rho \tau}})^{-1}
 ({\bx},\theta)
 =\left(\bx(q_0)+\sqrt{\frac{n-k}{2\rho R(q_0)}}{\bx},\theta\right)$,
it follows from \eqref{eq:close_cyl_pi} that, on the 
$100$-neighborhood of $\{\bx(q_0)\}\times(\IS^{n-k}/\Gamma)$, uniformly for
$\tau\in[-1,-e^{-1}]$,
\[
\begin{split}
 h'_\tau=(\kappa_\tau)_*(\rho^{-1}g'_{\rho \tau})
 &=\frac{n-k}{2\rho R(q_0)}
 \big(\kappa_\tau\circ(\pi_{q_0}|_{\MM'_{\rho \tau}})^{-1}\big)_*
 \big[(\pi_{q_0})_*(\tfrac{2R(q_0)}{n-k}g'_{\rho \tau})\big]\\
 &=g_{\IR^k}+
 \left(\frac{n-k}{2\rho R(q_0)}-\tau-e^{-1}\right)g_{\IS^{n-k}/\Gamma}
 +o_r(1)=g_{\cyl,\tau}+
 O(a)
 +o_r(1),
\end{split}
\]
where we used
$\frac{n-k}{2\rho R(q_0)}(\kappa_\tau\circ(\pi_{q_0}|_{\MM'_{\rho \tau}})^{-1})_*o_r(1)=o_r(1)$
because $C^{-1}\le\rho R(q_0)\le C$, and also
\eqref{eq:y_q_0}.
The estimate holds uniformly for $\tau\in[-1,-e^{-1}]$.
Since $(\bx_0,\theta_0)$ was arbitrary, taking the supremum gives
\begin{equation}\label{eq:cap_bryant_cylindrical_estimate}
 \sup_{-1\leq\tau\leq-e^{-1}}
 \|h'_\tau-g_{\cyl,\tau}\|_
 {C^r_{g_{\cyl,\tau}}(\IB^k_{\frac52a\sqrt\rho}\times(\IS^{n-k}/\Gamma))}
 \leq Ca+o_r(1).
\end{equation}

Along the backward worldlines, the steady soliton equation $R+|\nabla f|^2=1$ gives
\begin{equation}\label{eq:cap_bryant_y_drift}
 \frac{d}{dt}y(p(t))
=-f'(y(p(t)))
=-\sqrt{1-R(p(t))}=-1+O(y(p(t))^{-1}).
\end{equation}
Together with \eqref{eq:y_q_0}, this implies $y\geq c\rho$
on $\kappa^{-1}_\tau(\IB^k_{\frac52a\sqrt\rho}\times(\IS^{n-k}/\Gamma))$ for
$-1\leq\tau\leq-e^{-1}$.  For $0\leq j\leq m+6$, \eqref{eq:bryant_polynomial_decay}
and $R\leq C\rho^{-1}$ imply
\[
 \left|\nabla^j\big(\rho^{-1}(g-g')\big)\right|_{\rho^{-1}g'}
 =\rho^{j/2}|\nabla^j(g-g')|_{g'}
 \leq C(m,P)\rho^{(m+6)/2-P}.
\]
Choose $P>(m+6)/2$, and choose $a=a(m)$ small so that $Ca<\varepsilon_m/4$ in
\eqref{eq:cap_bryant_cylindrical_estimate} with $r=m+6$, where
$\varepsilon_m$ is from Lemma~\ref{l:relative_grafting_one_step},
and then increase $\rho$.  The coordinate norms and the rescaled
Bryant norms are uniformly equivalent through this order.
The triangle inequality together with \eqref{eq:cap_bryant_cylindrical_estimate} therefore shows that $h_\tau,h'_\tau$ satisfy
\eqref{eq:cap_cylinder_smallness} with $L=a\sqrt\rho$.
\end{proof}

Choose and fix any point $p\in\MM'_{<T}$ with $y(p)\ge\underline{Y}(m)$, where $\underline{Y}(m)$ is from Claim~\ref{Cl_cap_bryant_cylindrical_estimate}.
Starting with $t_0=\mathfrak t(p)$, define
\[
 \rho_j=\frac{e(n-k)}{2R(p(t_j))},\qquad
 t_{j+1}=t_j-(1-e^{-1})\rho_j,\qquad j\geq0.
\]
By \eqref{eq:cap_bryant_y_drift} and $2yR/(n-k)\to1$, we can see
$y(p(t_{j+1}))=y(p(t_j))+(1-e^{-1})\rho_j+O(1)$.
Consequently, after increasing $\underline{Y}(m)$, we have
\begin{equation}\label{eq:cap_center_drift}
 \frac{\rho_{j+1}}{\rho_j}=e+o(1),\qquad
 2\rho_j\leq\rho_{j+1}\leq3\rho_j
\end{equation}
for all $j$.  In particular, $t_j\to-\infty$.
For each $j$, let
$\kappa_{j,\tau}:=\kappa_{p(t_j),\tau}$ and $L_j=a\sqrt{\rho_j}$.  Thus $\kappa_{j,\tau},\tau\in[-1,-e^{-1}]$
is defined on the slice at time $t_j+\rho_j(\tau+e^{-1})\in[t_{j+1},t_j]$.
Set 
\[h_{j,\tau}:=(\kappa_{j,\tau})_*
   \big(\rho_j^{-1}g_{t_j+\rho_j(\tau+e^{-1})}\big),\qquad h'_{j,\tau}:=(\kappa_{j,\tau})_*
   \big(\rho_j^{-1}g'_{t_j+\rho_j(\tau+e^{-1})}\big).\]
Then Claim \ref{Cl_cap_bryant_cylindrical_estimate} implies that these two Ricci flows
satisfy \eqref{eq:cap_cylinder_smallness} with $L=L_j$.

Fix some $\ell\geq1$.
We will construct a sequence of diffeomorphisms onto their images
\[
 \Phi_{j}:
 \mathcal B_{0,j}:=\IB^k_{\frac35L_j}\times(\IS^{n-k}/\Gamma)
 \longrightarrow\IB^k_{\frac52L_j}\times(\IS^{n-k}/\Gamma),
 \qquad j=0,\ldots,\ell,
\]
by downward induction on $j$, starting with
$\Phi_{\ell}=\id$, such that for each map,
the associated quantity
\[
\begin{split}
 S_j:=
 \|(\Phi_{j})^*h'_{j,-e^{-1}}-h_{j,-e^{-1}}\|_
 {C^{m+2}_{g_{\cyl}}(\mathcal B_{0,j})}
 +\|\Phi_{j}-\id\|_
 {C^{m+1}_{g_{\cyl}}(\mathcal B_{0,j})}.
\end{split}
\]
satisfies the inductive assumptions $S_j\le\eps'_m$ and
\begin{equation}\label{eq:cap_inductive_bound}
 S_j\leq C(m,N)C_0^{\ell-j}\rho_\ell^{-N}
 +\sum_{i=j}^{\ell-1}C_0^{i-j+1}e^{-\beta'\rho_i},
\end{equation}
where $0<\eps'_m<\eps_m$, $C_0=C_0(m)>0$, $\beta'(m)>0$, $N(m)\in\mathbb N$ will be chosen below, and $\eps_m$ is from Lemma~\ref{l:relative_grafting_one_step}.
We first verify this for the first step when $j=\ell$: Since $\Phi_\ell=\id$, the arbitrary-order polynomial estimate
\eqref{eq:bryant_polynomial_decay} gives
$S_\ell\leq C(m,N)\rho_\ell^{-N}$, which is \eqref{eq:cap_inductive_bound} and $S_\ell\le \eps'_m$ for $j=\ell$, by taking $\underline{Y}(m)$ sufficiently large depending on $\eps'_m$.

Now assume by induction that we have constructed $\Phi_{\ell},\cdots,\Phi_{j+1}$ and that $S_{\ell},\cdots,S_{j+1}$ satisfy both inductive assumptions $S_i\le\eps'_m$ and \eqref{eq:cap_inductive_bound}.
Let $D_j:=\kappa_{j,-1}\circ\kappa_{j+1,-e^{-1}}^{-1}$.
Consider the map
$D_j\circ\Phi_{j+1}\circ D_j^{-1}$, restricted to
$\IB^k_{\frac45L_j}\times(\IS^{n-k}/\Gamma)$.
Note that the initial map is well-defined since
\[
 \kappa_{j,-1}^{-1}(\IB^k_{\frac45L_j}\times(\IS^{n-k}/\Gamma))
 \subset B_{g'_{t_{j+1}}}(p(t_{j+1}),a\rho_j)
 \subset
 \kappa_{j+1,-e^{-1}}^{-1}(\mathcal B_{0,j+1}),
\]
where the second inclusion uses
$\rho_j\leq\rho_{j+1}/2$.
By the equivalence of the rescaled Bryant soliton metrics and the cylindrical metrics, the initial norms in \eqref{eq:cap_step_input} for $D_j\circ\Phi_{j+1}\circ D_j^{-1}$ are bounded by $C_0S_{j+1}$.
Thus, choosing $\eps'_m=\eps_m/(2C_0)$, the inductive assumption $S_{j+1}\le\eps'_m$ implies \eqref{eq:cap_step_input}  for $D_j\circ\Phi_{j+1}\circ D_j^{-1}$.
Here we also used the observation
$h_{j,-1}=(\rho_{j+1}/\rho_j)(D_j)_*h_{j+1,-e^{-1}}$
and the same identity for the primed metrics.
So applying
Lemma~\ref{l:relative_grafting_one_step} to $h_{j,\tau},h'_{j,\tau}$, with $L=L_j$ and initial map $D_j\circ\Phi_{j+1}\circ D_j^{-1}$, and choosing $C_0$ sufficiently large compared with $C_m$ from that lemma to account for the coordinate changes and rescaling, we obtain the map
$\Phi_{j}:\mathcal B_{0,j}\to\IB^k_{\frac52L_j}\times(\IS^{n-k}/\Gamma)$ such that
\begin{equation}\label{eq:cap_multiscale_recurrence}
 S_j\leq C_0\left(S_{j+1}+e^{-\beta'\rho_j}\right),
\end{equation}
using $L_j=a\sqrt{\rho_j}$ and setting
$\beta'=\beta a^2$.
Substituting \eqref{eq:cap_inductive_bound} for $S_{j+1}$ proves
\eqref{eq:cap_inductive_bound} for $S_j$.
Moreover, by \eqref{eq:cap_center_drift} and choosing $N$ large so that $C_0 2^{-N}<1/2$, \eqref{eq:cap_inductive_bound} implies
\begin{equation}\label{eq:cap_partial_iteration}
\begin{split}
 S_j
 &\leq C(m,N)C_0^{\ell-j}\rho_{\ell}^{-N}
 +\sum_{q=0}^{\ell-j-1}C_0^{q+1}e^{-\beta'2^q\rho_0}\le C(m,N)2^{-\ell}\rho_{0}^{-N}+ 2C_0e^{-\beta'\rho_0}.
\end{split}
\end{equation}
So by choosing $\underline{Y}_m$ sufficiently large, this verifies the inductive assumption $S_j<\eps_m'$ for $S_j$. 
Therefore, by induction we construct all the maps with the desired properties.  

Finally, the $C^{m+1}$-bound for $\Phi_0-\id$ and the $C^{m+2}$-bound for
$\Phi_0^*h'_{0,-e^{-1}}-h_{0,-e^{-1}}$ in $S_0$
allow us to remove $\Phi_0$ and bound
$h'_{0,-e^{-1}}-h_{0,-e^{-1}}$ in $C^m$.
More precisely, returning to the original coordinates and scale, we obtain
\[
\begin{split}
 |\nabla^j_{g'}(g-g')|_{g'}(p)
 \leq C(m)\rho_0^{-j/2}S_0\leq C(m)S_0,
 \qquad 0\leq j\leq m.
\end{split}
\]
Letting
$\ell\to\infty$ in \eqref{eq:cap_partial_iteration} for $S_0$ we prove \eqref{eq:bryant_exponential} {for $y\geq Y_m$, where $Y_m$ depends only on $m$ and the dimension,} since
$\rho_0=e(n-k)/(2R(p))\geq c\,y(p)$.

\medskip
Finally, Assertion~\ref{p:bryant_infinity_b} of Proposition~\ref{p:bryant_infinity} gives $|\nabla_{g'}^j(g-g')|_{g'}\leq C_m$ on the fixed region $\{y\geq Y_0\}$ for every $0\leq j\leq m$, since $R$ is uniformly bounded there. On the remaining strip $Y_0\leq y\leq Y_m$, increasing $C(m)$ to at least $C_m e^{c_mY_m}$ therefore gives \eqref{eq:bryant_exponential}. Thus the exponential estimate holds on $\{y\geq Y_0\}$, with $Y_0$ independent of $m$.
If $Y_0$ was enlarged, the image conclusion of Theorem~\ref{thm:linear_bryant_asymptotics} follows after decreasing $\eps$, by the same Ricci-eigenvalue argument used to prove Assertion~\ref{p:bryant_infinity_c} of Proposition~\ref{p:bryant_infinity}.
\end{proof}

\section{Flows with dominant quadratic mode}\label{sec:dominant_quadratic}

{In this section we study the neutral dominant case.}

\subsection{Asymptotics and the quadratic invariant}

{We first associate an exact matrix ODE solution to the
quadratic coefficients of $U^+$.}

\begin{Proposition}[Asymptotics of the quadratic mode] \label{Prop_RF_dom_qu_asymp}
For any integers $J,m \geq 10$ and any $0<\eta\leq \frac1{10}$ there is a constant $\delta(J,m,\eta)>0$ such that the following is true.

Let $\MM$ be an ancient asymptotically $(n,k,\Gamma)$-cylindrical Ricci flow.
Suppose that $q_0 \in \MM_{t_0-r_0^2}$ and that $q_0(t)$ is a center of an $(n,k,\delta,\Gamma)$-neck at scale $(t_0-t)^{1/2}$ for all $t \leq t_0-r_0^2$.
Set $\td\tau := -2\log r_0$.
Choose admissible HMHF data $(\chi,\Omega,R,U^+,q_0)$ with parameters
$m,J+2,\eta,C(J,m)$ and reference time $t_0$ on $(-\infty,\td\tau]$,
using Proposition~\ref{Prop_ancient_summary} for noncylindrical flows and fixed product coordinates for round shrinking cylinders, with the approximation
\[
 U^+ : (-\infty,\td\tau] \lto \sV_{\rot,\geq -J-2}.
\]
Write $\PP_{\sV_{\sph,0}} U^+(\tau)
 = \sum_{i,j=1}^k c_{ij}(\tau) \mathfrak p^{(2)}_{ij}(\bx) g_{\IS^{n-k}/\Gamma}$
and write $\mathsf C(\tau):=(c_{ij}(\tau))_{i,j=1}^k$ for the symmetric coefficient matrix.
Then there exists a unique maximal solution $\ov U : (-\infty,\ov T) \lto \IR^{k\times k}_{\leq 0}$
for some $\ov T\in\IR\cup\{\infty\}$ to the ODE
\begin{equation} \label{eq_ovU_ODE}
 \partial_\tau \ov U
 = - \sqrt 2\, \ov U^2 + \tr(\ov U^2)\ov U
 +(4-n+k)\ov U^3
\end{equation}
such that, for $\tau \leq \td\tau-C(J,m,\eta)$,
\begin{equation} \label{eq_RF_U0_close_ovU}
 \| \mathsf C(\tau)-\ov U(\tau) \| \leq C(J,m,\eta)(\td\tau-\tau+10)^{-3}.
\end{equation}
Moreover, we have $\ov U\equiv 0$ if and only if $T_3=-\infty$ in Proposition~\ref{Prop_ancient_summary}.
\end{Proposition}

The proof follows the finite-dimensional ODE argument
of~\cite[Section~7]{Bamler_Lai_MCF1}.

\begin{proof}
If $\MM$ is a round shrinking cylinder with extinction time $T_{\cyl}$, take fixed product coordinates centered at $q_0$ and set
\[
u_\tau=U^+(\tau)=(T_{\cyl}-t_0)e^\tau g_{\IS^{n-k}/\Gamma},\qquad
e^{-R(\tau)^2}=\|U^+(\tau)\|^{J+2}.
\]
The neck assumption gives $|T_{\cyl}-t_0|\leq C\delta r_0^2$, so these data are admissible for sufficiently small $\delta$, with $R=\infty$ when $U^+=0$, $\mathsf C=0$, and $T_3=-\infty$.

By parabolic rescaling, assume $\td\tau=-10$.
First suppose $T_3=-\infty$. Proposition~\ref{Prop_ancient_summary}
\ref{Prop_ancient_summary_f} gives exponential decay of $U^+$.
Taking $\ov U\equiv0$ gives \eqref{eq_RF_U0_close_ovU}.

Henceforth assume $T_3>-\infty$.
We derive an equation for $\mathsf C$ and construct the corresponding
exact solution. In this proof, the $O(\cdot)$-constants depend only
on $n,k,J,m,\eta$. Write
$\PP_{\sV_{\aff}}U^+=\sum\mathsf D_{ij}\,d\bx_i\otimes d\bx_j$.
By Proposition~\ref{Prop_ancient_summary}\ref{Prop_ancient_summary_d},
\begin{equation}\label{eq:RF_Uplus_evolution}
\|\partial_\tau U^+-LU^+-Q_{J+2}^+(U^+)\|
\leq C\|U^+\|^{J+3}.
\end{equation}
By Proposition~\ref{Prop_ancient_summary}\ref{Prop_ancient_summary_b2}
and \ref{Prop_ancient_summary_e} with $i=0$, we have
\begin{equation}\label{eq_RF_Upp_standard_bound}
\|U^+\|=O(|\tau|^{-1}),\qquad
\|\PP_{\sV_{\rot,<0}}U^+\|=O(|\tau|^{-2}).
\end{equation}
Since $L=0$ on $\sV_{\sph,0}$ and
$\|(Q_{J+2}-Q_3)(U^+)\|\leq C\|U^+\|^4$, projecting
\eqref{eq:RF_Uplus_evolution} onto $\sV_{\sph,0}$ gives
\begin{equation}\label{eq:RF_spherical_neutral_evolution}
\partial_\tau\bigg(\sum_{i,j}\mathsf C_{ij}
\mathfrak p_{ij}^{(2)}g_{\IS^{n-k}/\Gamma}\bigg)
=\PP_{\sV_{\sph,0}}Q_3(U^+)+O(|\tau|^{-4}).
\end{equation}
Projecting \eqref{eq:RF_Uplus_evolution} onto $\sV_{\rot,>0}$
and integrating gives $\|\PP_{\sV_{\rot,>0}}U^+\|=O(|\tau|^{-2})$.
Together with \eqref{eq_RF_Upp_standard_bound}, this gives
\begin{equation}\label{eq:RF_nonneutral_bound}
\|\PP_{\sV_{\rot,\ne0}}U^+\|=O(|\tau|^{-2}).
\end{equation}
Hence
\[
\|Q_{J+2}^+(U^+)-Q_2^+(\PP_{\sV_{\rot,0}}U^+)\|
\leq C\|U^+\|\|\PP_{\sV_{\rot,\ne0}}U^+\|+C\|U^+\|^3
=O(|\tau|^{-3}).
\]
Projecting \eqref{eq:RF_Uplus_evolution} onto $\sV_{\sph,0}$ and
$\sV_{\aff}$ and using Lemma~\ref{Lem_Q2_computation}, this implies
\begin{equation}\label{eq:RF_neutral_quadratic_system}
\partial_\tau\mathsf C=-\sqrt2\mathsf C^2+O(|\tau|^{-3}),\qquad
\partial_\tau\mathsf D=\frac{n-k}{2}\mathsf C^2+O(|\tau|^{-3}).
\end{equation}
Thus $\partial_\tau(\mathsf D+\frac{n-k}{2\sqrt2}\mathsf C)=O(|\tau|^{-3})$.
Since $\mathsf C,\mathsf D\to0$ at $-\infty$, integration gives
\begin{equation}\label{eq:RF_affine_slaving}
\mathsf D=-\frac{n-k}{2\sqrt2}\mathsf C+O(|\tau|^{-2}).
\end{equation}

By \eqref{eq_Us2QU_1}, we have
\[
E(\tau):=\PP_{\sV_{\rot,<0}}Q_2(\PP_{\sV_{\rot,0}}U^+)=\frac{n-k}{\sqrt2}
\sum_{a,b,i,j}c_{ai}c_{bj}\mathfrak p_{ij}^{(2)}
\,d\bx_a\otimes d\bx_b,\qquad LE=-E.
\]
By \eqref{eq:RF_neutral_quadratic_system},
$\partial_\tau E=O(|\tau|^{-3})$.
Set $Z:=\PP_{\sV_{\rot,<0}}U^+-E$. Using $LE=-E$ in
\eqref{eq:RF_Uplus_evolution}, we obtain
\[
\partial_\tau Z-LZ
=\PP_{\sV_{\rot,<0}}\big[Q_{J+2}^+(U^+)-Q_2(\PP_{\sV_{\rot,0}}U^+)\big]
-\partial_\tau E+O(|\tau|^{-J-3})=O(|\tau|^{-3}).
\]
Since $\Spec(L|_{\sV_{\rot,<0}})\subset(-\infty,-\frac12]$,
integrating from $-\infty$ gives
\begin{equation}\label{eq:RF_negative_slaving}
\|Z(\tau)\|
\leq C\int_{-\infty}^\tau e^{-(\tau-s)/2}|s|^{-3}\,ds
\leq C|\tau|^{-3}.
\end{equation}

\begin{Claim}\label{claim:RF_spherical_neutral_Q3}
We have
\[
\PP_{\sV_{\sph,0}}Q_3(U^+)
=\sum_{i,j}\big(-\sqrt2\mathsf C^2
+\tr(\mathsf C^2)\mathsf C+(4-n+k)\mathsf C^3\big)_{ij}
\mathfrak p_{ij}^{(2)}g_{\IS^{n-k}/\Gamma}+O(|\tau|^{-4}).
\]
\end{Claim}

\begin{proof}
By \eqref{eq:RF_negative_slaving},
$U^+=\PP_{\sV_{\rot,0}}U^++E
+\PP_{\sV_{\rot,>0}}U^++Z$,
where $\|E\|+\|\PP_{\sV_{\rot,>0}}U^+\|=O(|\tau|^{-2})$
and $\|Z\|=O(|\tau|^{-3})$.
By \eqref{eq_Us2QU_4},
$\PP_{\sV_{\sph,0}}Q_2(\PP_{\sV_{\rot,\geq0}}U^+)
=\PP_{\sV_{\sph,0}}Q_2(\PP_{\sV_{\rot,0}}U^+)$.
Using \eqref{eq_Us2QU_4} and the quadraticity of $Q_2$,
the triangle inequality gives
\[
\begin{aligned}
&\big\|\PP_{\sV_{\sph,0}}\big[
Q_2(U^+)-Q_2(\PP_{\sV_{\rot,0}}U^++E)\big]\big\|\\
&\quad\leq\big\|Q_2(U^+)
-Q_2(\PP_{\sV_{\rot,\geq0}}U^++E)\big\|
+\big\|\PP_{\sV_{\sph,0}}\big[
Q_2(\PP_{\sV_{\rot,\geq0}}U^++E)
-Q_2(\PP_{\sV_{\rot,0}}U^++E)\big]\big\|\\
&\quad\leq\big\|\PP_{\sV_{\sph,0}}\big[
Q_2(\PP_{\sV_{\rot,\geq0}}U^+)
-Q_2(\PP_{\sV_{\rot,0}}U^+)\big]\big\|
+C(\|U^+\|+\|Z\|)\|Z\|
+C\|E\|\|\PP_{\sV_{\rot,>0}}U^+\|\\
&\quad\leq C(\|U^+\|+\|Z\|)\|Z\|
+C\|E\|\|\PP_{\sV_{\rot,>0}}U^+\|
=O(|\tau|^{-4}).
\end{aligned}
\]
In the second norm, \eqref{eq_Us2QU_4} cancels the contribution
independent of $E$, leaving a bilinear term in $E$ and
$\PP_{\sV_{\rot,>0}}U^+$.
Since $Q_3-Q_2$ is cubic,
\[
\|(Q_3-Q_2)(U^+)
-(Q_3-Q_2)(\PP_{\sV_{\rot,0}}U^+)\|
\leq C\|U^+\|^2\|\PP_{\sV_{\rot,\ne0}}U^+\|
=O(|\tau|^{-4}).
\]
Combining these estimates gives
\[
\PP_{\sV_{\sph,0}}Q_3(U^+)
=\PP_{\sV_{\sph,0}}Q_2(\PP_{\sV_{\rot,0}}U^++E)
+\PP_{\sV_{\sph,0}}(Q_3-Q_2)(\PP_{\sV_{\rot,0}}U^+)
+O(|\tau|^{-4}).
\]
By \eqref{eq_Uh2QU}, we have
\[
\begin{aligned}
\PP_{\sV_{\sph,0}}Q_2(\PP_{\sV_{\rot,0}}U^++E)
&=\sum_{i,j}\bigg(-\sqrt2(\mathsf C^2)_{ij}
-\frac{n-k}{2}\sum_{a,b}\mathsf C_{ai}\mathsf C_{ab}\mathsf C_{bj}\bigg)
\mathfrak p_{ij}^{(2)}g_{\IS^{n-k}/\Gamma}\\
&=\sum_{i,j}\bigg(-\sqrt2\mathsf C^2
-\frac{n-k}{2}\mathsf C^3\bigg)_{ij}
\mathfrak p_{ij}^{(2)}g_{\IS^{n-k}/\Gamma}.
\end{aligned}
\]
Diagonalizing $\mathsf C$ to apply Lemma~\ref{Lem_Q3_computation}
and returning to the original coordinates, we obtain
\[
\PP_{\sV_{\sph,0}}(Q_3-Q_2)(\PP_{\sV_{\rot,0}}U^+)
=\sum_{i,j}\big(\tr(\mathsf C^2)\mathsf C+4\mathsf C^3
+\sqrt2\mathsf C\mathsf D\mathsf C\big)_{ij}
\mathfrak p_{ij}^{(2)}g_{\IS^{n-k}/\Gamma}.
\]
Combining these identities and substituting
\eqref{eq:RF_affine_slaving} proves the claim.
\end{proof}

Combining the claim with \eqref{eq:RF_spherical_neutral_evolution}
proves
\begin{equation}\label{eq_RF_U0_approx_ODE}
\big\|\partial_\tau\mathsf C+\sqrt2\mathsf C^2
-\tr(\mathsf C^2)\mathsf C-(4-n+k)\mathsf C^3\big\|
\leq C|\tau|^{-4}.
\end{equation}
Therefore, by
the same argument in the proof of
\cite[Proposition~7.1(d)]{Bamler_Lai_MCF1}
we find a maximal solution
$\ov U:(-\infty,\ov T)\to\IR^{k\times k}_{\leq0}$
to \eqref{eq_ovU_ODE} satisfying \eqref{eq_RF_U0_close_ovU}.
Moreover, for two exact solutions $\ov U_1,\ov U_2$ satisfying
\eqref{eq_RF_U0_close_ovU}, the argument in the proof of
\cite[Lemma~7.32]{Bamler_Lai_MCF1} gives $\ov U_1\equiv\ov U_2$,
proving uniqueness.

Lastly, Proposition~\ref{Prop_ancient_summary}
\ref{Prop_ancient_summary_f} gives
$|\lambda_{\min}(\mathsf C)|\geq c|\tau|^{-1}$ for
sufficiently negative $\tau$, so
$\ov U\ne0$.
\end{proof}

The following lemma parametrizes the space of nonpositive definite solutions of \eqref{eq_ovU_ODE} that converge to zero as $\tau\to-\infty$.
This parametrization will be used in Definition~\ref{Def_RF_Qu} to define the quadratic invariant of asymptotically cylindrical Ricci flows.

\begin{Lemma} \label{Lem_RF_ODE_UQ}
Let $\mathfrak U^k$ denote the set of maximal solutions $\ov U:(-\infty,\ov T)\to\IR^{k\times k}_{\leq0}$ to \eqref{eq_ovU_ODE} with $\lim_{\tau\to-\infty}\ov U(\tau)=0$.
There is a topology on $\mathfrak U^k$ for which convergence is equivalent to convergence on each fixed time-slice.
Moreover, there are homeomorphisms
\[
 \Qu_k:\mathfrak U^k\lto \IR^{k\times k}_{\geq0},
\]
for all $k\geq0$, such that the following is true.
\begin{enumerate}[label=(\alph*)]
\item \label{Lem_RF_ODE_UQ_a} For every orthogonal matrix $\mathsf S\in O(k)$, we have $\Qu_k(\mathsf S^T\ov U\mathsf S)=\mathsf S^T\Qu_k(\ov U)\mathsf S$.
\item \label{Lem_RF_ODE_UQ_b} The nullspaces of $\ov U(\tau)$, for any $\tau$, and $\Qu_k(\ov U)$ agree.
\item \label{Lem_RF_ODE_UQ_c} If $\ov U$ is regarded as a solution in dimension $k'>k$ by adding zero rows and columns, then $\Qu_{k'}(\ov U)$ is obtained from $\Qu_k(\ov U)$ in the same way.
\item \label{Lem_RF_ODE_UQ_d} If $\ov U'(\tau)=\ov U(\tau+2\log\alpha)$, then $\Qu_k(\ov U')=\alpha\Qu_k(\ov U)$.
\end{enumerate}
\end{Lemma}

\begin{proof}
The RHS of \eqref{eq_ovU_ODE} is a polynomial in $\ov U$ and $\tr(\ov U^2)$.
On the range of a nonzero solution $\ov U$ of rank $r$, we have
\[
\partial_\tau\ov U^{-1}
=\sqrt2\,\mathsf I_r-\tr(\ov U^2)\ov U^{-1}-(4-n+k)\ov U.
\]
Following the argument in the proof of
\cite[Proposition~7.1(d), following (7.29)]{Bamler_Lai_MCF1},
successive integration gives
\begin{equation} \label{eq_RF_Ubar_asymp}
 \ov U(\tau)=
 \big(2\sqrt2\log\mathsf A+(\sqrt2\,\tau-\tfrac{r+4-n+k}{\sqrt2}\log(-\tau))\mathsf I_r\big)^{-1}
 +O(|\tau|^{-3}\log|\tau|)
\end{equation}
for a unique positive definite matrix $\mathsf A$ on the range of $\ov U$.

As in \cite[proof of Lemma~7.6]{Bamler_Lai_MCF1}, there exists an $\varepsilon_0(n-k)>0$ such that every nonzero solution has a unique time $\tau_{\ov U}$ such that $\lambda_{\min}(\ov U(\tau_{\ov U}))=-\varepsilon_0$.
Define $\Qu_k(0):=0$ and
\[
\Qu_k(\ov U):=-\varepsilon_0^{-1}e^{-\tau_{\ov U}/2}\ov U(\tau_{\ov U})\quad(\ov U\ne0).
\]
The same proof shows that $\Qu_k$ is a homeomorphism.
Assertions~\ref{Lem_RF_ODE_UQ_a}--\ref{Lem_RF_ODE_UQ_c} follow directly from this formula and the ODE.
Assertion~\ref{Lem_RF_ODE_UQ_d} follows from $\tau_{\ov U'}=\tau_{\ov U}-2\log\alpha$.
\end{proof}

{We define the quadratic invariant $\Qu$ of an asymptotically cylindrical Ricci flow as follows.}

\begin{Definition} \label{Def_RF_Qu}
Let $\MM$ be an ancient Ricci flow spacetime satisfying the assumptions of
Proposition~\ref{Prop_RF_dom_qu_asymp} or homothetic to a round shrinking cylinder.
In the cylindrical case set $\ov U\equiv0$; otherwise let $\ov U$ be
the solution from Proposition~\ref{Prop_RF_dom_qu_asymp}.
We define
\[
 \Qu(\MM):=\big(\lambda_1(\Qu_k(\ov U)),\ldots,\lambda_k(\Qu_k(\ov U))\big)
 \in\IR^k_{\uparrow,\geq0},
\]
where $\lambda_1(A)\leq\cdots\leq\lambda_k(A)$ are the eigenvalues of a symmetric matrix $A$, and $\IR^k_{\uparrow,\geq0}$ denotes the set of vectors with nonnegative, weakly increasing entries.
We call $\Qu(\MM)$ the \textbf{quadratic mode at $-\infty$}.
\end{Definition}

{The following proposition shows that $\Qu(\MM)$ is independent of the HMHF data and the parameters used in its construction. It also records its behavior under translations and parabolic rescaling, and characterizes the flows for which $\Qu(\MM)=0$. We establish independence of the choices of HMHFs by
recovering the eigenvalues of $\ov U$ from the intrinsic Ricci
eigenvalues at a common point.}

\begin{Proposition} \label{Prop_RF_Q_basic}
The vector $\Qu(\MM)$ and the orthogonal conjugacy class of
the solution $\ov U$ depend only on $\MM$ and not on the choices of
$J,m,\eta,\delta,r_0,t_0,q_0$, or on the particular harmonic map heat
flow supplied by Proposition~\ref{Prop_ancient_summary}.
Moreover, $\Qu(\MM)$ is invariant under time-translation and spatial translation of the flow, and under parabolic rescaling it satisfies
\[
 \Qu(\alpha\MM)=\alpha\Qu(\MM), \qquad \alpha>0.
\]
Lastly, $\Qu(\MM)=0$ if and only if $\MM$ is homothetic to a round
shrinking cylinder or has dominant linear mode.
\end{Proposition}

\begin{proof}
For a round shrinking cylinder the assertions follow from $\ov U\equiv0$,
so assume that $\MM$ is not homothetic to a round shrinking cylinder.
To prove independence, let $p_1(t),p_2(t)$ be the distinguished worldlines of the two HMHFs and $\ov U_1,\ov U_2$ their corresponding solutions of the ODE \eqref{eq_ovU_ODE} constructed above.
Write $t_{0,\ell}$ for their reference times and $\tau^{(\ell)}=-\log(t_{0,\ell}-t)$, $\ell=1,2$.
In this proof, we allow
$O(\cdot)$ to depend on these two HMHF data sets.
If one of $\ov U_1,\ov U_2$ vanishes, then
Proposition~\ref{Prop_RF_dom_qu_asymp}
and Lemma~\ref{lem:dominant_mode_independence} imply that both vanish.
Henceforth assume that both are nonzero.

\begin{Claim}\label{claim:RF_intrinsic_quadratic_Ricci}
For $\ell=1,2$, the $k$ smallest Ricci eigenvalues $\lambda_j(t)$, $1\leq j\leq k$, of $g_t$ at $p_1(t)$ satisfy
\[
(t_{0,\ell}-t)\lambda_j(t)
=\lambda_j\!\left(
-\frac{n-k}{2\sqrt2}\ov U_\ell(\tau^{(\ell)})
-\frac{(n-k)(n-k+2)}8\ov U_\ell(\tau^{(\ell)})^2\right)
+O(|\tau^{(\ell)}|^{-3}).
\]
The remaining $n-k$ rescaled Ricci eigenvalues are $\frac12+O(|\tau^{(\ell)}|^{-1})$.
\end{Claim}

\begin{proof}
Fix $\ell\in\{1,2\}$ and work in the $\ell$-th HMHF coordinates, writing $\tau=\tau^{(\ell)}$.
Proposition~\ref{Prop_ancient_summary}\ref{Prop_ancient_summary_a},
\ref{Prop_ancient_summary_b2}, \ref{Prop_ancient_summary_c} and local
parabolic estimates give
\eqref{eq:asymptotically_cylindrical_convergence} along both
distinguished worldlines. Definition~\ref{defn:asymptotically_cylindrical}
then yields $(t_{0,\ell}-t)^{-1/2}d_{g_t}(p_1(t),p_\ell(t))\leq C$
for all sufficiently negative $t$.
By Proposition~\ref{Prop_ancient_summary}\ref{Prop_ancient_summary_a},
$\xi_{e^{\tau/2}}\circ\chi_{\ell,t}(p_\ell(t))\in\IB_\eta^k\times(\IS^{n-k}/\Gamma)$.
This implies
$p_1(t)$ lies in a fixed rescaled coordinate ball of the $\ell$-th HMHF
for sufficiently negative $t$. Write $(\bx(\tau),\by(\tau))$ for its coordinates
in the $\ell$-th chart. 
Using the HMHF equation, we have
\[
\partial_\tau\bx=\tfrac12\bx+X_\tau(\bx,\by),\qquad
X_\tau:=\operatorname{proj}_{\IR^k}
\big(\Delta_{g_{\cyl}+u_\tau,g_{\cyl}}\id\big).
\]
By \eqref{eq_RF_Upp_standard_bound}, the bound
$U^+-\PP_{\sV_{\rot,0}}U^+=O(|\tau|^{-2})$ proved above, and
Proposition~\ref{Prop_ancient_summary}\ref{Prop_ancient_summary_c},
\[
X_\tau(\bx,\by)=\frac{n-k}{2\sqrt2}\mathsf C(\tau)\bx+O(|\tau|^{-2}).
\]
Since $\|\mathsf C(\tau)\|=O(|\tau|^{-1})$, we obtain $|X_\tau(\bx,\by)|\leq C|\tau|^{-1}|\bx|+C|\tau|^{-2}$.
For a fixed sufficiently negative $T$, integration gives
\[
|\bx(\tau)|\leq e^{(\tau-T)/2}|\bx(T)|
+C\int_\tau^T e^{(\tau-s)/2}
\big(|s|^{-1}|\bx(s)|+|s|^{-2}\big)\,ds.
\]
Starting with bounded $\bx$, this gives successively
$|\bx(\tau)|=O(|\tau|^{-1})$ and $O(|\tau|^{-2})$.

Write
$g_{\cyl}+U^+=H_{ab}(\bx)\,d\bx_a\otimes d\bx_b
+(1+v(\bx))g_{\IS^{n-k}/\Gamma}$.
Here $H:=\sum_{a,b=1}^k H_{ab}(\bx)\,d\bx_a\otimes d\bx_b$ denotes the metric on the base $\IR^k$.
On a fixed coordinate ball, \eqref{eq:RF_negative_slaving} gives
\[
\begin{aligned}
H_{ab}
&=\delta_{ab}+\mathsf D_{ab}
+\frac{n-k}{4}\big((\mathsf C\bx)_a(\mathsf C\bx)_b
-2(\mathsf C^2)_{ab}\big)+O_{C^2}(|\tau|^{-3}),\\
v(0)&=-\frac{\tr\mathsf C}{\sqrt2}+O(|\tau|^{-2}),
\qquad \partial^2_{ab}v(0)=\frac{\mathsf C_{ab}}{\sqrt2}
+O(|\tau|^{-3}),\qquad
\nabla v(0)=O(|\tau|^{-2}).
\end{aligned}
\]
By the first formula, the Ricci curvature of $H$ at $\bx=\bO$ is $\Ric_H
=\frac{n-k}{4}
\big((\tr\mathsf C)\mathsf C-\mathsf C^2\big)
+O(|\tau|^{-3})$.
By the warped-product formula, we can compute the Ricci tensor of $g_{\cyl}+U^+$ at $\bx=\bO$, and obtain the following components under the $\bx$-coordinates
\[
\begin{aligned}
\Ric|_{T\IR^k}
&=\Ric_H-\frac{n-k}{2(1+v)}\nabla_H^2v
+\frac{n-k}{4(1+v)^2}dv\otimes dv
=-\frac{n-k}{2\sqrt2}\mathsf C
-\frac{n-k}{4}\mathsf C^2+O(|\tau|^{-3}).
\end{aligned}
\]
By using \eqref{eq:RF_affine_slaving} in
$H^{-1/2}=I-\frac12\mathsf D+O(|\tau|^{-2})$, this implies
\[H^{-1/2}(\Ric|_{T\IR^k})H^{-1/2}
=-\frac{n-k}{2\sqrt2}\mathsf C
-\frac{(n-k)(n-k+2)}8\mathsf C^2+O(|\tau|^{-3}),\]
and the same expansion holds at $\bx(\tau)=O(|\tau|^{-2})$.
By Proposition~\ref{Prop_ancient_summary}\ref{Prop_ancient_summary_c},
$\|u_\tau-U^+\|_{C^2}=O(|\tau|^{-J-3})$ on this ball.
Together with \eqref{eq_RF_U0_close_ovU}, this gives
the asserted expansion. The spherical block gives
the remaining rescaled Ricci eigenvalues $\frac12+O(|\tau|^{-1})$.
\end{proof}

Since $\frac{t_{0,2}-t}{t_{0,1}-t}=1+O(e^{\tau^{(1)}})$, $\tau^{(2)}-\tau^{(1)}=O(e^{\tau^{(1)}})$,
and the derivative of
$c\mapsto-(n-k)c/(2\sqrt2)-(n-k)(n-k+2)c^2/8$
is bounded away from zero near $c=0$, Claim~\ref{claim:RF_intrinsic_quadratic_Ricci} implies
\[
\big|\lambda_j(\ov U_1(\tau))-\lambda_j(\ov U_2(\tau))\big|
=O(|\tau|^{-3}),\qquad 1\leq j\leq k.
\]
So modulo an orthogonal transformation, we may assume $
\|\ov U_1(\tau)-\ov U_2(\tau)\|
=O(|\tau|^{-3})
$.
The uniqueness argument in the proof of Proposition~\ref{Prop_RF_dom_qu_asymp} gives $\ov U_1\equiv\ov U_2$.
By Lemma~\ref{Lem_RF_ODE_UQ}\ref{Lem_RF_ODE_UQ_a},
the ordered vector $\Qu(\MM)$ is therefore independent of the choices.

The transformation laws are obtained by applying Proposition~\ref{Prop_RF_dom_qu_asymp} to the translated or rescaled flow and using Lemma~\ref{Lem_RF_ODE_UQ}\ref{Lem_RF_ODE_UQ_d}.
The last assertion follows from Proposition~\ref{Prop_RF_dom_qu_asymp}
and Definition~\ref{defn:dominant_linear}.
\end{proof}

{The following proposition turns equality of $\Qu$ into agreement of the modes and rescaled metric perturbations to any prescribed polynomial order in $|\tau|^{-1}$. For each such order, the comparison uses the corresponding HMHF data and a fixed rotation of the Euclidean coordinates.}

\begin{Proposition} \label{Prop_RF_same_Q_close}
For any integers $J,m\geq10$ and $0<\eta\leq\frac1{10}$ there is a constant $\delta(J,m,\eta)>0$ with the following property.
Let $\MM^{(0)}$ and $\MM^{(1)}$ be two ancient Ricci flow spacetimes satisfying the assumptions of Proposition~\ref{Prop_RF_dom_qu_asymp}, with the same values of $t_0,r_0,\td\tau=-2\log r_0$ and with $\Qu(\MM^{(0)})=\Qu(\MM^{(1)})$.
Let $U^{+,(i)},\chi^{(i)},R_i(\tau)$ be the corresponding objects from Proposition~\ref{Prop_RF_dom_qu_asymp}, constructed with approximation order $J+2$.
After a fixed rotation of the Euclidean cylindrical coordinates for $\MM^{(1)}$,
we have, for $\tau\leq \td\tau-C(J)$,
\begin{equation}\label{eq_RF_same_Uplus}
 \|U^{+,(1)}(\tau)-U^{+,(0)}(\tau)\|_{L^2_f}
 \leq C(J)(\td\tau-\tau+10)^{-J}.
\end{equation}
Consequently, the corresponding rescaled perturbations satisfy
\begin{equation}\label{eq_RF_same_metric}
 \|u_\tau^{(1)}-u_\tau^{(0)}\|_{C^m(\IB^k_{\min\{R_0(\tau),R_1(\tau)\}}\times(\IS^{n-k}/\Gamma))}
 \leq C(J,m)(\td\tau-\tau+10)^{-J}.
\end{equation}
\end{Proposition}

\begin{proof}
By parabolic rescaling, assume $\td\tau=-10$.
We omit the superscript $(0)$ and use primes for the quantities
associated with $\MM^{(1)}$.
By Definition~\ref{Def_RF_Qu} and Lemma~\ref{Lem_RF_ODE_UQ}\ref{Lem_RF_ODE_UQ_a},
equality of the ordered invariants allows us to rotate the cylindrical
coordinates for $\MM^{(1)}$ so that $\Qu_k(\ov U)=\Qu_k(\ov U')$.
The injectivity of $\Qu_k$ then gives $\ov U=\ov U'$.
If $\ov U=\ov U'=0$, then
Proposition~\ref{Prop_RF_dom_qu_asymp} gives $T_3=-\infty$ for both HMHFs.
Finite-dimensional norm equivalence, the radius relation, and
Proposition~\ref{Prop_ancient_summary}\ref{Prop_ancient_summary_c},
\ref{Prop_ancient_summary_f} give
\[
\|U^+(\tau)\|\leq C\eta e^{0.49(\tau+10)},\qquad
\|u_\tau\|_{C^m(\IB^k_{R(\tau)}\times(\IS^{n-k}/\Gamma))}
\leq C(J,m)\|U^+(\tau)\|^{1/2},
\]
and similarly for the primed quantities. These bounds imply
\eqref{eq_RF_same_Uplus} and \eqref{eq_RF_same_metric}.
Henceforth assume $\ov U=\ov U'\ne0$.
Set $N=J+2$, and write
\[
V=U^+-U^{+,\prime},\qquad
V_C=\mathsf C-\mathsf C',\qquad
V_D=\mathsf D-\mathsf D'
+\frac{n-k}{2\sqrt2}V_C,\qquad
V_\perp=\PP_{\sV_{\rot,\ne0}}V.
\]
The mode estimates in the proof of
Proposition~\ref{Prop_RF_dom_qu_asymp} use only Proposition~\ref{Prop_ancient_summary}~\ref{Prop_ancient_summary_b2}, \ref{Prop_ancient_summary_d}--\ref{Prop_ancient_summary_f}, so their constants and $C_0$ in the construction of $\ov U$ depend only on $J$.
Fix a common cutoff $T=T(J)<-10$ for both HMHFs; all $O(\cdot)$-constants in the following mode comparison depend only on $J$.
Equality of the ODE solutions gives $V_C=O(|\tau|^{-3})$,
while \eqref{eq:RF_affine_slaving} and
\eqref{eq:RF_nonneutral_bound} give
$V_D,V_\perp=O(|\tau|^{-2})$.

We now subtract the evolution equations for $U^+$ and $U^{+,\prime}$,
whose errors satisfy \eqref{eq:RF_Uplus_evolution}.
By \eqref{eq_Us2QU_4} and \eqref{eq_Us2QU_1}, the
$\sV_{\sph,0}$- and $\sV_{\aff}$-projections of
$Q_2(\PP_{\sV_{\rot,0}}U^+)$ have coefficient matrices
$-\sqrt2\mathsf C^2$ and $\frac{n-k}{2}\mathsf C^2$, respectively,
and similarly for $U^{+,\prime}$.
For the $V_C$ equation, use
\[
\mathsf C^2-(\mathsf C')^2=\mathsf C V_C+V_C\mathsf C'.
\]
For the $V_D$ equation, the quadratic contributions cancel by
\eqref{eq:RF_neutral_quadratic_system} and the definition of $V_D$.
Using \eqref{eq_RF_Upp_standard_bound} and \eqref{eq:RF_nonneutral_bound}
for both HMHFs, together with
$\|\PP_{\sV_{\rot,0}}V\|\leq C(\|V_C\|+\|V_D\|)$, we obtain
\[
\begin{aligned}
\|(Q_{J+2}-Q_2)(U^+)-(Q_{J+2}-Q_2)(U^{+,\prime})\|
&\leq C|\tau|^{-2}\|V\|,\\
\big\|Q_2(U^+)-Q_2(U^{+,\prime})
-Q_2(\PP_{\sV_{\rot,0}}U^+)
+Q_2(\PP_{\sV_{\rot,0}}U^{+,\prime})\big\|
&\leq C|\tau|^{-2}(\|V_C\|+\|V_D\|)
+C|\tau|^{-1}\|V_\perp\|
\end{aligned}
\]
Together with the two errors of size $O(|\tau|^{-N-1})$ in
\eqref{eq:RF_Uplus_evolution}, these estimates give
\begin{equation}\label{eq:RF_quadratic_difference_system}
\begin{aligned}
\|\partial_\tau V_C+\sqrt2(\mathsf C V_C+V_C\mathsf C')\|
+\|\partial_\tau V_D\|
&\leq C|\tau|^{-2}(\|V_C\|+\|V_D\|)
+C|\tau|^{-1}\|V_\perp\|+C|\tau|^{-N-1},\\
\|\partial_\tau V_\perp-LV_\perp\|
&\leq C|\tau|^{-1}(\|V_C\|+\|V_D\|+\|V_\perp\|)
+C|\tau|^{-N-1}.
\end{aligned}
\end{equation}
The initial bounds make the RHS of the second inequality
$O(|\tau|^{-3})$. Integrating from $-\infty$ on
$\sV_{\rot,<0}$ and backward from $T$ on $\sV_{\rot,>0}$
gives $V_\perp=O(|\tau|^{-3})$; the boundary terms at $T$
decay exponentially backward. Substituting this improved bound into
the first inequality gives $\partial_\tau V_D=O(|\tau|^{-4})$.
Since $V_D\to0$ at $-\infty$, integration gives $V_D=O(|\tau|^{-3})$.
Thus $\|V_C\|+\|V_D\|+\|V_\perp\|\leq C(J)|\tau|^{-3}$ for $\tau\leq T$.

We now apply the comparison argument in the proof of
\cite[Proposition~7.11]{Bamler_Lai_MCF1} to
\eqref{eq:RF_quadratic_difference_system}, treating $(V_C,V_D)$
as the neutral component. Its leading equations are
$\partial_\tau V_C=-\sqrt2(\mathsf C V_C+V_C\mathsf C')$
and $\partial_\tau V_D=0$.
By Proposition~\ref{Prop_ancient_summary}\ref{Prop_ancient_summary_f},
this evolution has the same growth bound used there; the constant
$V_D$ component does not increase that bound.
The initial $O(|\tau|^{-3})$ estimate makes the contribution at
$-\infty$ vanish, so the argument gives
\begin{equation}\label{eq:RF_same_modes_higher_order}
\|V(\tau)\|\leq C(J)|\tau|^{-J-2}.
\end{equation}

To prove \eqref{eq_RF_same_metric}, set $r=\min\{R(\tau),R'(\tau)\}$.
If $r=\infty$, then $u_\tau=u'_\tau=0$ and the assertion is immediate.
The radius identities give
$\|V\|\leq\|U^+\|+\|U^{+,\prime}\|\leq2e^{-r^2/(J+2)}$.
Combining this with \eqref{eq:RF_same_modes_higher_order} and
finite-dimensional norm equivalence, we obtain, for some $d=d(J)$,
\[
\begin{aligned}
\|V\|_{C^m(\IB_r^k\times(\IS^{n-k}/\Gamma))}
&\leq C(J,m)(1+r)^d\|V\|^{J/(J+2)}\|V\|^{2/(J+2)}\\
&\leq C(J,m)|\tau|^{-J}(1+r)^d e^{-2r^2/(J+2)^2}
\leq C(J,m)|\tau|^{-J}.
\end{aligned}
\]
By Proposition~\ref{Prop_ancient_summary}\ref{Prop_ancient_summary_c},
the two metric approximation errors on this region are
$O(|\tau|^{-J-3})$. Together with
\eqref{eq:RF_same_modes_higher_order}, this proves
\eqref{eq_RF_same_Uplus} and \eqref{eq_RF_same_metric}.
\end{proof}

\subsection{Existence of prescribed quadratic modes}

In this subsection we take $\Gamma=\{1\}$.

{To construct flows with prescribed $\Qu$, we first record its continuity under smooth convergence when the limit is again asymptotically cylindrical. This property will allow us to retain the prescribed invariant when passing to a limit in the proof of Theorem~\ref{Thm_RF_Q_realization}.}

\begin{Proposition}[Continuity of the quadratic mode] \label{Prop_RF_Q_continuity}
Let $\MM^{(i)}$, $i\in\IN\cup\{\infty\}$, be ancient asymptotically
$(n,k)$-cylindrical Ricci flow spacetimes, each the regular part of a
noncollapsed $\IF$-limit as in
Lemma~\ref{lem:asymptotically_cylindrical_tangent_flows}.
Assume that $\MM^{(i)}\to\MM^{(\infty)}$ in the pointed smooth
spacetime topology on compact time-intervals.
Then
\[
 \Qu(\MM^{(i)}) \lto \Qu(\MM^{(\infty)}).
\]
\end{Proposition}

\begin{proof}
Fix $J,m\geq10$, $0<\eta<\frac1{10}$, the corresponding
$\delta$ in Proposition~\ref{Prop_RF_dom_qu_asymp} with regularity
$m+2$ and smallness $\eta/2$, and choose $\eps>0$ sufficiently small
to apply Theorem~\ref{Thm_propagation_necks} with accuracy $\delta$.
Start with an arbitrary subsequence.
Choose the worldline $p_\infty(t)$ in Definition~\ref{defn:asymptotically_cylindrical}
and take $t_*<0$ so negative that $p_\infty(t)$ is an $(n,k,\eps/2)$-neck
center at scale $\sqrt{-t}$ for every $t\leq t_*$.
Smooth convergence gives $(n,k,\eps)$-neck centers
$p_i\in\MM^{(i)}_{t_*}$ with $p_i\to p_\infty(t_*)$ for large $i$.

\begin{Claim}\label{lem:backward_cylindrical_worldline}
For every $i\leq\infty$, the backward worldline $p_i(t)$ is defined
for all $t\leq t_*$ and is a center of an $(n,k,\delta)$-neck
at scale $\sqrt{-t}$ for every such $t$.
\end{Claim}

\begin{proof}
For each $i\leq\infty$, Lemma~\ref{lem:asymptotically_cylindrical_tangent_flows}
shows that every tangent flow at infinity based at $p_i$ is cylindrical.
By
\cite[Theorem~2.5]{Bamler_HK_RF_partial_regularity}, for
every sufficiently negative $s\le T_i$ (which may depend on $\MM_i$), an $H_n$-center of $p_i$ at time $s$
lies in $\MM^{(i)}_s$ and is a center of an $(n,k,\eps)$-neck
at scale $\sqrt{-s}$.
Applying Theorem~\ref{Thm_propagation_necks} to this center and $p_i$, the claim holds on $[s,t_*].$
Since $s$ can be arbitrarily negative and the extensions agree by
uniqueness, the claim follows.
\end{proof}

For $i<\infty$, choose the HMHF data in
Proposition~\ref{Prop_RF_dom_qu_asymp} with these parameters,
reference time $t_0=0$, and $\td\tau=-\log(-t_*)$.
By the same
compactness argument in the proof of
Lemma~\ref{l:non-compact-interval}, after passing to a further subsequence,
the restricted HMHFs and mode functions converge locally to admissible data for $\MM^{(\infty)}$
with regularity $m$, smallness $\eta$, and the same approximation
order and reference time. In particular,
$\mathsf C^{(i)}\to\mathsf C^{(\infty)}$ locally uniformly.
Let $\ov U^{(i)}$, $i\leq\infty$, be the corresponding solutions
from Proposition~\ref{Prop_RF_dom_qu_asymp}.
Then by Proposition~\ref{Prop_ancient_summary}
\ref{Prop_ancient_summary_b2} and \eqref{eq_RF_U0_close_ovU}, after passing to a subsequence, $\ov U^{(i)}$
converges locally smoothly to a solution $\ov V$
of \eqref{eq_ovU_ODE}.
Passing to the limit in \eqref{eq_RF_U0_close_ovU}, whose constant
and time cutoff are uniform in $i$, gives
\[
\|\mathsf C^{(\infty)}(\tau)-\ov V(\tau)\|
\leq C(\td\tau-\tau+10)^{-3}
\qquad(\tau\leq\td\tau-C).
\]
By uniqueness of the solution satisfying \eqref{eq_RF_U0_close_ovU},
we have $\ov V\equiv\ov U^{(\infty)}$.
Lemma~\ref{Lem_RF_ODE_UQ} and independence of the HMHF choices in
Proposition~\ref{Prop_RF_Q_basic} then give
$\Qu(\MM^{(i)})\to\Qu(\MM^{(\infty)})$ along this further
subsequence. Since the initial subsequence was arbitrary, the full
sequence converges.
\end{proof}

{The following lemma establishes smooth pointed compactness for the symmetric Ricci flows below, with limits that are asymptotically cylindrical $\kappa$-solutions with extinction time $0$.
This compactness will be used to establish properness of $\Qu$ in the proof of Theorem~\ref{Thm_RF_Q_realization}.}

\begin{Lemma}\label{Lem_RF_NPC_convergence}
There is a constant $\delta_0>0$ with the following property.
Let $(M_i,(g_{i,t})_{t\in(-T_i,0)})$ be a sequence of Ricci flows with $T_i\to\infty$ and
\[
 M_i=\IR^{l_i}\times \IS^{n-l_i}, \qquad l_i \in \{ 0, \ldots, k \}.
\]
Assume that each flow has non-negative curvature operator, splits off the $\IR^{l_i}$-factor isometrically, is invariant under the standard $\IZ_2^{l_i} \times \IZ_2^{k-l_i} \times O(n-k+1)$-action on $M_i$ (here we view $M_i$ as a subset of $\IR^{l_i} \times \IR^{k-l_i} \times \IR^{n-k+1}$) and becomes singular at time $0$.
Consider the central spherical fiber
\[
 S_i:=\bO^k\times \IS^{n-k}\subset M_i
\]
which is fixed by the $\IZ_2^k = \IZ_2^{l_i} \times \IZ_2^{k-l_i}$-action and let $p :=\mathbf e_{n+1}\in S_i$.
Suppose the following hold:
\begin{enumerate}[label=(\arabic*)]
\item\label{Lem_RF_NPC_convergence_a} There is some $t_0<0$ such that for every $i$ and every $t\in(-T_i,t_0]$, $(p,t)$ is a center of an $(n,k,\delta_0)$-neck at scale $\sqrt{-t}$.
\item\label{Lem_RF_NPC_convergence_b} For each $i$, there is a time $s_i\in(-T_i,t_0]$ at which $T_pS_i$ is the sum of the eigenspaces of the Ricci endomorphism of $(-s_i)^{-1}g_{i,s_i}$ with eigenvalue greater than $\frac14$.
\end{enumerate}
Then, after passing to a subsequence, we have smooth pointed convergence
\begin{equation} \label{eq_RF_H_conv_Mi}
 (M_i,(g_{i,t})_{t\in(-T_i,0)},(p,-1))
 \lto
 (M_\infty,(g_{\infty,t})_{t<0},(p_\infty,-1))
\end{equation}
in the sense of Hamilton~\cite{Hamilton_RF_compactness}.
The limit is a $\kappa$-solution whose time-slices, rescaled by $(-t)^{-1}$ and based at $p_\infty$, converge smoothly to the normalized cylinder $\IR^k\times\IS^{n-k}$ as $t\to-\infty$.
We can take $(M_\infty, p_\infty) = (\IR^{l_\infty} \times \IS^{n-l_\infty}, \mathbf e_{n+1})$ for some $l_\infty \in \{ 0, \ldots, k \}$ and the flow $(g_{\infty,t})_{t <0}$ is invariant under the standard $\IZ^k_2 \times O(n-k+1)$-action, isometrically splits off the Euclidean $\IR^{l_\infty}$-factor and becomes extinct in a round cylindrical singularity at time $0$.
\end{Lemma}

\begin{proof}
By $O(n-k)$ symmetry at $p$, the Ricci endomorphism of $g_{i,t}$ acts on $T_pS_i$ by a scalar $\lambda_i(t)$. For small $\delta_0$, Assumption~\ref{Lem_RF_NPC_convergence_a} keeps the rescaled Ricci spectrum away from $\frac14$. By continuity, the inequality $(-t)\lambda_i(t)>\frac14$ from Assumption~\ref{Lem_RF_NPC_convergence_b} holds throughout $(-T_i,t_0]$.

Thus $T_pS_i$ tends to the spherical tangent space as $\delta_0\to0$. Sectional curvature on the round, totally geodesic $S_i$ then gives convergence of its rescaled radius to $\sqrt{2(n-k-1)}$. A half great circle in $S_i$ from $p$ converges in the expanding neck charts to one on the model sphere, so its rescaled ambient endpoint distance and the rescaled intrinsic diameter both tend to $\pi\sqrt{2(n-k-1)}$. Hence for small $\delta_0$,
\begin{equation}\label{eq:RF_central_fiber_diameter}
\frac1{10}\pi\sqrt{-2(n-k-1)t}\leq\diam_t S_i\leq10\pi\sqrt{-2(n-k-1)t}\qquad(t\in(-T_i,t_0]).
\end{equation}
Away from collapsed fibers, write $g_{i,t}=h_{i,t}+\varphi_{i,t}^2g_{\IS^{n-k}(1)}$, with $g_{\IS^{n-k}(1)}$ the unit round metric. Nonnegative fiber curvature gives $|\nabla\varphi_{i,t}|\leq1$, so $F_{i,t}(x,\theta):=\varphi_{i,t}(x)\theta$, extended by zero on collapsed fibers, is $1$-Lipschitz.

After passing to a subsequence, we may assume that $l_i$ is constant.
Since the Euclidean factors split isometrically, we may divide them out and assume that $l_i=0$ and $M_i=\IS^n$.
The general case then follows by taking the product with the limiting Euclidean factor.

Hamilton's Harnack inequality~\cite[Corollary~1.2]{Hamilton_Harnack} gives
\begin{equation} \label{eq_RF_Harnack}
  \partial_t R \geq - \frac{R}{t+T_i} .
\end{equation}
Using \cite[11.6]{Perelman1} (cf \cite[Corollary~45.1]{Kleiner_Lott_notes}) and the neck assumption, we can deduce uniform curvature bounds depending only on the  distance to $p_i$ and on time, at all times $t < t_0$.
By integrating the bound \eqref{eq_RF_Harnack} backwards we moreover obtain uniform scalar curvature bounds of the form
\[ \sup_{t \in [-T_i/2, t_0]} \sup_{B(p_i, t_0, r)} R (\cdot, t) \leq C(r). \]
Due to the non-negative curvature operator condition these bounds imply uniform bounds on the Riemannian curvature, allowing us to extract a subsequence such that we have convergence in the sense of Hamilton
\begin{equation} \label{eq_MitoMinftyt0}
 (M_i,(g_{i,t})_{t\in(-\infty,t_0]},(p,t_0))
 \lto
 (M_\infty,(g_{\infty,t})_{t\leq t_0},(p_\infty,t_0)).
\end{equation}
The time-slices of the limit are complete and have non-negative curvature operator, but the curvature need not be globally bounded.
This fact will complicate our arguments, as we cannot directly apply standard results, such as Perelman's Pseudolocality Theorem, to the limiting flow.
Since $M_i$ is compact, we can, however, apply the Pseudolocality Theorem \cite[10.3]{Perelman1} and again \cite[11.6]{Perelman1} near $(p_i,t_0)$, which allows us to extend the convergence \eqref{eq_MitoMinftyt0} to a time-interval of the form $(-\infty,t_1)$, for some $t_1 > t_0$, after passing to another subsequence.
Suppose now that $t_1$ is chosen maximal with this property.
Since the flows $(M_i, (g_{i,t}))$ go extinct in a spherical singularity \cite{Bohm2008} at time $0$, we must have $t_1 \leq 0$.
Moreover, repeating the same argument as above implies that $R(\cdot,t)$ cannot remain bounded near $p_\infty$ as $t \nearrow t_1$.
In the following, we will study the limiting flow $(M_\infty, (g_{\infty,t})_{t < t_1})$ and show that, in fact, $t_1 = 0$.
This will imply that we can choose the basepoints of the convergence at time $-1$, as asserted in \eqref{eq_RF_H_conv_Mi}.

By \eqref{eq:RF_central_fiber_diameter}, the isometric $\IZ_2^k \times O(n-k+1)$-action limits to an isometric action on $(M_\infty, g_{\infty,t})$ for all $t < t_1$.
In addition, the fibers $S_i \subset M_i$, which are totally geodesic and round, limit to a submanifold $S_\infty \subset M_\infty$, which is fixed by the limiting $\IZ_2^k$-action, invariant under the limiting $O(n-k+1)$-action and totally geodesic with respect to any $g_{\infty,t}$.

\begin{Claim}
If $\delta_0 \leq \ov\delta_0$, then $M_\infty$ is simply-connected.
\end{Claim}

\begin{proof}
Consider the universal covering flow $( \td M_\infty, ( \td g_{\infty,t})_{t < t_1})$ and choose two lifts $\td p_0, \td p_1 \in \td M$ of $p_\infty$.
By the limiting neck condition, we have curvature bounds of the form $|{\Rm}| \leq C(-t)^{-1}$ on $B(\td p_j, t, C^{-1} \sqrt{-t})$ for some uniform $C > 0$.
A well-known distance distortion bound therefore implies that $\partial_t d_t (\td p_0, \td p_1) \geq - C (-t)^{-1/2}$, after adjusting $C$.
Integrating this backwards in time implies, again after adjusting $C$,
\begin{equation} \label{eq_dtdt0}
 d_t (\td p_0, \td p_1)  \leq d_{t_0} (\td p_0, \td p_1) + C \sqrt{-t}, \qquad \text{for all} \quad t \leq t_0. 
\end{equation}
On the other hand, by the neck assumption and the fact that $(n,k,\delta_0)$-neck regions are simply-connected, we must have for some universal $c > 0$
\[ d_t (\td p_0, \td p_1) \geq c \delta_0^{-1} \sqrt{-t}, \]
which contradicts \eqref{eq_dtdt0} for $\delta_0 \leq \ov\delta_0$ and $t \ll t_0$.
\end{proof}
\medskip

\begin{Claim} \label{Cl_limitsymmetry}
Assuming $\delta_0 \leq \ov\delta_0$, the flow $(M_\infty, (g_{\infty,t})_{t < t_1})$ is asymptotic to $\IR^k \times \IS^{n-k}$ at $-\infty$.
Moreover, the flow is isometric to the Cartesian product of $\IR^{l_\infty}$, for some $l_\infty \in \{ 0, \ldots, k \}$, with a flow on $\IS^{n-l_\infty}$ that is invariant under the standard $\IZ_2^{k-l_\infty} \times O(n-k+1)$-action and which develops a spherical singularity at time $t_1$.
In addition, there is a uniform constant $C > 0$ such that for any $t < t_1$ we have
\begin{equation} \label{eq_limsupSi}
 \limsup_{i \to \infty} \vol_{t} S_i \leq C (t_1 - t)^{(n-k)/2}. 
\end{equation}
\end{Claim}

\begin{proof}
The radius bound above controls $\varphi_{i,t}(p)$ for $t\leq t_0$, and $\partial_t\varphi_{i,t}(p)\leq0$ gives a bound for later times. Thus $F_{i,t}(p)$ is uniformly bounded on compact time-intervals.
After passing to a further subsequence, these maps converge locally uniformly to equivariant maps $F_t$; the convergence is smooth where $F_t\ne0$, since the orbit radii are determined by infinitesimal rotation lengths and the angular coordinates by the smooth fixed-point sections.
If $F_t(q)\ne0$, small rotations align the spherical directions of approximating points with $F_t(q)$, so its $O(n-k)$-stabilizer fixes $q$; equivariance gives the reverse inclusion.
Where $F_t(q)=0$, the orbit diameters tend to zero, so the whole group fixes $q$.
For such a shrinking orbit in a convex ball, minimizing the averaged squared distance gives a nearby fixed point in $M_i$; its standard normal representation passes to the limit.
Thus the time-independent spherical coordinates pass to the limit without identifications, and the standard slices at the fixed points give
\begin{equation}\label{eq_prodmap}
M_\infty\cong(\Sigma'\times\IS^{n-k}(1))/\!\sim,\qquad g_{\infty,t}=h_t+\varphi_t^2g_{\IS^{n-k}(1)}\quad\textnormal{on }\Int\Sigma'\times\IS^{n-k}(1),
\end{equation}
where each boundary fiber is collapsed to a point.
Here $\Sigma'$ is the complete orbit space, with smooth totally geodesic boundary, and the collapsed locus in $M_\infty$ has codimension $n-k+1$.
The central fiber projects to a point $\bar p\in\Int\Sigma'$, where $\nabla\varphi_t=0$.
The warped-product curvature formulas give
\[
\nabla_{h_t}^2\varphi_t\leq0,\qquad K(X,V)=-\varphi_t^{-1}\nabla_{h_t}^2\varphi_t(X,X)
\]
for unit horizontal $X$ and vertical $V$.

Hamilton's strong maximum principle makes the curvature-operator kernel parallel, and its rank can decrease only finitely many times; see~\cite[Appendix~A]{Kleiner_Lott_notes}.
Work first on a terminal subinterval where this rank is constant, and set $E:=\ker\nabla_{h_t}^2\varphi_t$.
For horizontal $X$, the preceding formula and nonnegative curvature operator imply
\[
X\in E\quad\Longleftrightarrow\quad X\wedge V\in\ker\Rm\ \textnormal{for every vertical }V.
\]
Parallel transport along horizontal curves preserves the horizontal and vertical spaces, so $E$ is parallel on $\Int\Sigma'$.
The $E$-component of $\nabla\varphi_t$ is therefore parallel and vanishes at $\bar p$, giving $d\varphi_t(E)=0$.
For a horizontal lift $X\in E$ we then have $\nabla_VX=X(\log\varphi_t)V=0$, so $E$ is also parallel in the full metric.
It extends across the collapsed locus: loops avoiding this codimension-$(n-k+1)$ submanifold approximate arbitrary loops, and parallel transport is continuous.
The de Rham theorem gives a global product, and the time-invariance and backward splitting assertions of the same strong maximum principle extend it to the whole flow:
\begin{equation}\label{eq_MpMpp_splitting}
(M_\infty,(g_{\infty,t})_{t<t_1})\cong(M'',(g''_t)_{t<t_1})\times(M',(g'_t)_{t<t_1}),
\end{equation}
where $E$ is tangent to $M''$ on the terminal subinterval, the spherical action is trivial on $M''$, and $\varphi_t$ depends only on $M'$.

The orbit space of $M'$ has no nonzero null direction of $\nabla^2\varphi_t$ on this subinterval.
If it were noncompact, a minimizing ray from the image of $\bar p$ would remain in its interior, since the boundary is totally geodesic.
Along this ray, $\varphi_t$ is positive and concave, with derivative zero at the origin, so it must be constant.
This contradicts the absence of a null direction.
Hence this orbit space, and consequently $M'$, is compact.
 
Write $p_\infty=(p'',p')$ and fix any sequence $\lambda_j\to\infty$. By Assumption~\ref{Lem_RF_NPC_convergence_a}, $(p_\infty,t)$ is a center of an $(n,k,2\delta_0)$-neck at scale $\sqrt{-t}$ for $t<t_0$. On compact time intervals in $s<0$, these necks give uniform lower bounds for the injectivity radius at $p_\infty$ and the volume of balls of any fixed radius about $p_\infty$ in the metrics $\lambda_j^{-1}g_{\infty,\lambda_j s}$. Applying \cite[Corollary~45.1(b)]{Kleiner_Lott_notes} to balls of four times the desired radius gives uniform curvature bounds on such balls. Hamilton compactness gives a complete pointed smooth spacetime limit on $s<0$ after passing to a subsequence. The same bounds give subsequential limits of both factors in \eqref{eq_MpMpp_splitting}.

For each $s<0$, we have $\lambda_j s<t_0$ for large $j$, so the blowdown limit has an $(n,k,3\delta_0)$-neck at scale $\sqrt{-s}$ along its basepoint worldline.
It inherits pseudolocality under smooth convergence.
For sufficiently small $\delta_0$, Corollary~\ref{Cor_eternal} therefore identifies the limit as the round shrinking cylinder.
Thus
\begin{equation} \label{eq_cyl_splitting}
 \MM_{\cyl}\cong (M''_{\infty},(g''_{\infty,s})_{s<0})\times(M'_{\infty},(g'_{\infty,s})_{s<0}).
\end{equation}
The central sphere $S_\infty$ lies in $\{p''\}\times M'$ and is totally geodesic. Its radius bound $\varphi_t(\bar p)\leq C\sqrt{-t}$ gives positive sectional curvature in the spherical directions at $p'$ after rescaling, bounded away from zero at each fixed $s<0$. Hence the second factor in \eqref{eq_cyl_splitting} contains the sphere factor at every time, and the first is the static Euclidean flow by the de Rham decomposition and the Ricci flow equation. Since the sequence $\lambda_j$ was arbitrary, we have smooth pointed spacetime convergence
\begin{equation} \label{eq_Mpp_to_Eucl}
 (M'',(\lambda^{-1}g''_{\lambda s})_{s<0},p'')\xrightarrow[\lambda\to\infty]{}(\IR^{l''},(g_{\IR^{l''}})_{s<0},\bO),
\end{equation}
where $l''=\dim M''$.

To prove that $M''$ is Euclidean, we may assume it is nontrivial. Fix $t^*<t_1$ and consider the $\LL$-distance based at $(p'',t^*)$.
Since $\Ric \geq 0$, we have for any $C^1$-curve $\gamma : [0,  \tau] \to M''$:
\begin{multline*}
 \LL(\gamma) = \int_0^{\tau} \sqrt{\tau'} \Big(R(\gamma(\tau'), t^*-\tau') + |\gamma'(\tau')|^2_{g''_{t^*-\tau'}} \Big) d\tau '\geq \int_0^{\tau} \sqrt{\tau'} \, |\gamma'(\tau')|^2_{g''_{t^*}}  d\tau' \\
\geq \bigg( \int_0^{\tau} \frac{d\tau'}{\sqrt{\tau'}} \bigg)^{-1} \bigg(\int_0^{\tau} |\gamma'(\tau')|_{g''_{t^*}} d\tau'\bigg)^2 = \frac{L_{t^*}(\gamma)^2}{2\sqrt{\tau}}.
\end{multline*}
So curves starting at $p''$ with bounded $\LL$-length remain in a uniform compact subset, yielding minimizing $\LL$-geodesics.
Since the flow is complete and $\partial_\tau g''_{t^*-\tau}=2\Ric\geq0$, \cite[Theorem~2.1]{Yokota_reduced_volume_gap} gives monotonicity and the bound $\td V_{(p'',t^*)}(\tau)\leq1$.

We will now show that $\td V_{(p'',t^*)}(\tau) \equiv 1$.
To do this, we first claim that for any $D, \eps > 0$ and $\tau  \geq \underline \tau(D,\eps)$ the following bound holds for Perelman's reduced distance function, for all $q \in M''$ with $d_{t^*-\tau}(p'',q) \leq D \sqrt{\tau}$:
\begin{equation} \label{eq_ellpppts}
 \ell_{(p'',t^*)} (q,\tau) \leq \frac{d^2_{t^*-\tau}(p'',q)}{4 \tau} + \eps. 
\end{equation}
Fix a small $\delta>0$. By \eqref{eq_Mpp_to_Eucl}, the metrics $\tau^{-1}g''_{t^*+\tau s}$ converge smoothly to the Euclidean metric on bounded regions around $p''$ over $s\in[-1,-\delta]$, since $t^*/\tau\to0$.
We can then consider a path $\gamma : [0,\tau] \to M''$ with $\gamma(0) = p''$ and $\gamma(\tau) = q$ such that $\gamma(\tau') = p''$ for $\tau' \leq \delta\tau$ and such that $\gamma|_{[\delta\tau,\tau]}$ is close to an $\LL$-geodesic in Euclidean space.
By the neck condition and the product splitting, $R''(p'',t^*-s)\leq C/(1+s)$ for $s\geq0$, with $C$ depending on $t^*$. Thus
\[
\begin{aligned}
\ell_{(p'',t^*)}(q,\tau)
&\leq \frac{C}{2\sqrt{\tau}}\int_0^{\delta\tau}\frac{\sqrt{s}}{1+s}\,ds+\frac1{2\sqrt{\tau}}\LL(\gamma|_{[\delta\tau,\tau]})
\leq C\sqrt{\delta}+\frac{d_{t^*-\tau}(p'',q)^2}{4\tau(1-\sqrt{\delta})}+o_\tau(1).
\end{aligned}
\]
For fixed $\delta$, the last estimate uses Euclidean convergence on $[-1,-\delta]$. Choosing $\delta$ sufficiently small and then $\tau$ sufficiently large gives \eqref{eq_ellpppts}.
The bound \eqref{eq_ellpppts} implies for the reduced volume based at $(p'', t^*)$
\[ \limsup_{\tau \to \infty} \td V_{(p'', t^*)} (\tau) \geq \int_{\IR^{l''}} (4\pi\tau)^{-l''/2} e^{-|\bx|^2/4\tau} d\bx = 1. \]
By monotonicity and the upper bound, $\td V_{(p'',t^*)}(\tau)\equiv1$.
The equality case of \cite[Theorem~2.1]{Yokota_reduced_volume_gap} implies that $(M'',(g''_t)_{t\leq t^*})$ is Euclidean.
Since $t^* < t_1$ was chosen arbitrarily, we obtain that the entire flow $(M'', (g''_{t < t_1} ))$ is Euclidean.
\medskip

In summary, we have shown that
\[ (  M_\infty, (  g_{\infty,t})_{t < t_1}) \cong \IR^{l_\infty} \times ( M', ( g'_{t})_{t < t_1}) \]
for some compact flow $( M', ( g'_{t})_{t < t_1})$, where the $O(n-k+1)$-action is trivial on the first factor and the $\IZ_2^k$-action preserves each factor.
{The Euclidean tangent space at $p_\infty$ lies in the span of the $(-1)$-eigenlines of the coordinate reflections and is invariant under those reflections, so it is spanned by a subset of these lines.
Thus, after reordering the original coordinate reflections, we may assume that $\IZ_2^{l_\infty} \times \{ + 1 \}^{k-l_\infty}$ acts on the $\IR^{l_\infty}$-factor and fixes all points in $M'$, while $\{ + 1 \}^{l_\infty} \times \IZ_2^{k-l_\infty}$ fixes the normal directions of $M' \subset M_\infty$ and acts on $M'$.}
Now observe that $\bO^{l_\infty} \times \IS^{n-l_\infty}$ and $M'$ are the fixed point sets of $\IZ_2^{l_\infty} \times \{ + 1 \}^{k-l_\infty}$ for $M_i$ and $M_\infty$, respectively, and these fixed point sets contain $S_i$ and $S_\infty$, respectively.
Since we have smooth convergence of $(M_i,g_{i,t},p_i) \to (M_\infty, g_{\infty,t},p_\infty)$, we also have smooth convergence of these fixed point sets.
Thus $M' \approx \IS^{n-l_\infty}$ and by Palais's rigidity theorem we can identify $M'$ with $\IS^{n-l_\infty}$ in such a way that the $\IZ_2^{k-l_\infty} \times O(n-k+1)$-action on $M'$ agrees with the standard action on $\IS^{n-l_\infty} \subset \IR^{k-l_\infty} \times \IR^{n-k+1}$.
This implies the first part of the claim.
The bound \eqref{eq_limsupSi} follows from the fact that $( M', ( g'_t)_{t < t_1})$ has a round spherical singularity at time $t_1$ by \cite{Bohm2008}.
\end{proof}

To show that $t_1 = 0$, we use the following claim.

\begin{Claim} \label{Cl_dtSi}
There is a uniform constant $c > 0$ such that
\[ \partial_t (\vol_t S_i)^{2/(n-k)} \leq - c. \]
\end{Claim}

\begin{proof}
Since $S_i$ is totally geodesic, the ambient sectional curvatures tangent to $S_i$ are proportional to $(\vol_t S_i)^{-2/(n-k)}$.
Since the flows $(g_{i,t})$ have non-negative sectional curvature, this implies a similar lower bound on the Ricci curvature tangent to $S_i$.
Therefore,
\[ \partial_t \vol_t S_i \leq - c (\vol_t S_i)^{-2/(n-k)} (\vol_t S_i), \]
which implies the claim.
\end{proof}

Claim~\ref{Cl_dtSi} implies that $\vol_t S_i \geq c(-t)^{(n-k)/2}$, after adjusting $c$.
So due to \eqref{eq_limsupSi} we must have $t_1 = 0$.
In summary, $(M_\infty, (g_{\infty,t})_{t < 0})$ is an ancient flow, which has non-negative curvature operator and is the product of a Euclidean factor with a compact flow.
Claim~\ref{Cl_limitsymmetry} establishes the symmetry properties described in the statement of the lemma.
The compact factor gives bounded curvature on compact time-intervals, so Lemma~\ref{lem:asymptotically_cylindrical_tangent_flows} gives uniform noncollapsing at all scales.
So $(M_\infty, (g_{\infty,t})_{t < 0})$ is a $\kappa$-solution with the asserted properties.
\end{proof}

{The following theorem realizes every $\Qu'\in\IR^k_{\uparrow,\geq0}$ by a symmetric $\kappa$-solution. The number of zero entries of $\Qu'$ determines the dimension of its Euclidean factor; the remaining compact factor becomes extinct in a round spherical singularity at time zero. These solutions provide comparison models for the asymptotics described in Proposition~\ref{Prop_RF_same_Q_close}.}

\begin{Theorem}\label{Thm_RF_Q_realization}
Let $\Qu' \in \IR^k_{\uparrow, \geq 0}$.
Then there is a $\kappa$-solution $(M,(g_t)_{t <0})$ that is asymptotic to $\IR^k \times \IS^{n-k}$ at $-\infty$ such that the corresponding Ricci flow spacetime $\MM$ satisfies $\Qu(\MM) =  \Qu'$.
Moreover, the following is true:
\begin{enumerate}[label=(\alph*)]
\item $M = \IR^l \times \IS^{n-l}$, where $l \geq 0$ is the number of zero entries of $\Qu'$ and the $\IR^l$-factor splits off isometrically at all times.
\item The flow is invariant under the standard $\IZ_2^k\times O(n-k+1)$-action: $\IZ_2^k$ reflects the coordinates of $\IR^l\times\IR^{k-l}$, and $O(n-k+1)$ acts on the remaining $\IR^{n-k+1}$ coordinates of $\IR^l\times\IR^{n-l+1}$.
\item The flow becomes singular at time $0$ and converges everywhere, modulo rescaling, to a round cylinder of the form $\IR^l \times \IS^{n-l}$
\end{enumerate}
\end{Theorem}

\begin{proof}
We adapt the proof of \cite[Theorem~7.13]{Bamler_Lai_MCF1}. The round shrinking cylinder realizes $\Qu'=0$. It suffices to treat targets $\Qu'=(Q'_1,\ldots,Q'_k)$ with all entries positive: taking a product with $\IR^l$ appends $l$ zero entries to $\Qu$, by the Ricci eigenvalue expansion in Claim~\ref{claim:RF_intrinsic_quadratic_Ricci} and Lemma~\ref{Lem_RF_ODE_UQ}\ref{Lem_RF_ODE_UQ_c}.

On each symmetric flow, fix $p=\mathbf e_{n+1}$ on the central spherical fiber and mark the $(-1)$-eigenlines of the coordinate reflections at $p$ in order. Use a distance $d_{CGH}$ for which convergence to an ancient flow means that the starting times tend to $-\infty$ and the flows converge smoothly in the pointed sense, with the marked lines converging in order at time $-1$. After passing to a subsequence, the reflections themselves converge: they preserve balls about $p$, and $r\circ\exp_x=\exp_{r(x)}\circ d r|_x$ gives smooth convergence from Arzel\`a--Ascoli.
Let $\mathfrak M$ be the space of $\kappa$-solutions with this product structure, symmetry and extinction at time $0$, allowing any $0\leq l\leq k$ and any ordering of the reflection lines. We require cylindrical convergence at $-\infty$ based at $p$.

The map $\Qu:\mathfrak M\to\IR^k_{\uparrow,\geq0}$ is continuous by Proposition~\ref{Prop_RF_Q_continuity}. To prove properness, let $\MM_i\in\mathfrak M$ have bounded $\Qu(\MM_i)$. Let $t_{0,i}$ be the supremum of the terminal times for which Assumption~\ref{Lem_RF_NPC_convergence_a} in Lemma~\ref{Lem_RF_NPC_convergence} holds for $\MM_i$; this set is nonempty by cylindrical convergence at $p$.
Write $\lambda_i(t)$ for the Ricci eigenvalue on $T_pS_i$ and $\varphi_{i,t}(p)$ for the radius of $S_i$. If $(-t)\lambda_i(t)<\frac14$ at any $t<t_{0,i}$, the opening observation in the proof of Lemma~\ref{Lem_RF_NPC_convergence} gives the same inequality throughout $(-\infty,t_{0,i})$. By the Ricci flow equation, total geodesicity and nonnegative sectional curvature,
\[
\frac{d}{dt}\varphi_{i,t}(p)^2=-2\lambda_i(t)\varphi_{i,t}(p)^2,\qquad \lambda_i(t)\varphi_{i,t}(p)^2\geq n-k-1.
\]
The first equation then gives $\varphi_{i,t}(p)^2=O((-t)^{1/2})$ as $t\to-\infty$, while the second gives $\varphi_{i,t}(p)^2\geq(n-k-1)(-t)$ for sufficiently negative $t$, a contradiction. Thus $(-t)\lambda_i(t)>\frac14$ and Assumption~\ref{Lem_RF_NPC_convergence_b} holds at any terminal time below $t_{0,i}$.

If the $t_{0,i}$ are bounded below, Lemma~\ref{Lem_RF_NPC_convergence} gives subsequential convergence in $\mathfrak M$ at a common terminal time below all $t_{0,i}$. Otherwise, pass to a subsequence with $t_{0,i}\to-\infty$ and apply the lemma to $(-t_{0,i})^{-1/2}\MM_i$ at terminal time $-2$. The limit $\MM_\infty\in\mathfrak M$ satisfies
\[ \Qu(\MM_\infty) = \lim_{i \to \infty} \Qu \big((-t_{0,i})^{-1/2}  \MM_i \big) 
= \lim_{i \to \infty} (-t_{0,i})^{-1/2} \Qu(  \MM_i)  = 0. \]
By Proposition~\ref{Prop_RF_Q_basic} and Corollary~\ref{thm:unique_bryant}, the limit is a round shrinking cylinder or a rescaled product with the Bryant soliton. The latter is excluded since the limit becomes extinct at time $0$.
For large $i$, smooth convergence to the cylinder extends Assumption~\ref{Lem_RF_NPC_convergence_a} from $t\leq-2$ through $t=-1/2$ in the rescaled flows, contradicting the rescaled supremum $-1$.

Set $K:=\{\MM\in\mathfrak M:|\Qu(\MM)|\leq|\Qu'|\}$, which is compact by properness. For $\delta\in(0,1)$, let $\Delta_\delta^k:=\{x\in\IR^k_{\geq0}:x_1+\cdots+x_k\leq\delta\}$ and, for $\bb\in\Delta_\delta^k$, equip the ellipsoid
\[ \big\{ (\bb_1 x_1)^2 + \ldots + (\bb_k x_k)^2 + x_{k+1}^2 + \ldots + x_{n+1}^2 = 1 \big\} \subset \IR^{n+1}, \]
with its induced metric. Weak convexity gives nonnegative curvature operator. Let $\MM^{\bb}$ be the Ricci flow starting from this metric at time $-T_{\bb}$, with time translated so that the round extinction of the compact factor occurs at $0$ by \cite{Bohm2008}. The flow preserves the coordinate symmetries and splits off the directions with $\bb_j=0$.
The lifetimes $T_{\bb}$ and flows $\MM^{\bb}$ depend continuously on $\bb$ in the marked topology, using the common ellipsoid coordinates near $p$ to identify the marks. In particular, pseudolocality and a limiting argument give $T_{\bb}\to 1/(2(n-k-1))$ as $\bb\to\bO$.

\begin{Claim}\label{claim:RF_ellipsoid_compactness}
There are constants $\ov\eps, \delta > 0$ with the following property.
Suppose $a_i\to\infty$ and $\bb_i\in\Delta_\delta^k$ satisfy
\begin{equation}\label{eq_dcghstar}
d_{CGH}(a_i\MM^{\bb_i},K)\leq\ov\eps.
\end{equation}
Then, after passing to a subsequence, $a_i\MM^{\bb_i}\to\MM_\infty\in\mathfrak M$ in $d_{CGH}$ and $\bb_i\to\bO$.
\end{Claim}

\begin{proof}
Choose $\delta_1\in(0,\frac12\delta_0)$ so that Theorem~\ref{Thm_propagation_necks} gives $(n,k,\delta_0)$-necks from endpoint $(n,k,\delta_1)$-necks. Each flow in $K$ has a $(n,k,\delta_1/2)$-neck at $(p,t)$ at scale $\sqrt{-t}$ for some $t<0$. By smooth convergence, nearby flows have a $(n,k,\delta_1)$-neck there. Compactness of $K$ therefore lets us choose the times $t_{0,i}$ for the flows in \eqref{eq_dcghstar} from a fixed finite set if $\ov\eps$ is small enough.

For small $\delta$, the initial ellipsoid chart gives a $(n,k,\delta_1)$-neck at $(p,-a_i^2T_{\bb_i})$ at scale $a_i\sqrt{T_{\bb_i}}$ with central slice $S_i$.
Choose an $H_n$-center $z_i$ of $(p,t_{0,i})$ on the initial slice. By reflection invariance, its image under $\bx\mapsto-\bx$ is also an $H_n$-center, so their distance in the unrescaled ellipsoid is at most $2\sqrt{H_nT_{\bb_i}}$. Thus $|\bx(z_i)|\leq\sqrt{H_nT_{\bb_i}}$, and for small $\delta$ the initial charts give a $(n,k,\delta_1)$-neck centered at $z_i$ at the same scale.
Theorem~\ref{Thm_propagation_necks} gives $(n,k,\delta_0)$-necks at $(p,t)$ at scale $\sqrt{-t}$ for $t\in[-a_i^2T_{\bb_i},t_{0,i}]$.
The initial ellipsoid chart identifies $T_pS_i$ with the spherical Ricci eigenspace, so Assumption~\ref{Lem_RF_NPC_convergence_b} holds just after the initial time.
Since $T_{\bb_i}$ is bounded away from zero, $a_i^2T_{\bb_i}\to\infty$.
We can therefore apply Lemma~\ref{Lem_RF_NPC_convergence} at a fixed terminal time below all $t_{0,i}$ to obtain subsequential convergence $a_i\MM^{\bb_i}\to\MM_\infty\in\mathfrak M$.

If $\bb_i\to\bb_\infty\ne\bO$ along a subsequence, the unrescaled flows $\MM^{\bb_i}$ have $(n,k,\delta_0)$-necks at $(p,t)$ at scale $\sqrt{-t}$ up to time $a_i^{-2}t_{0,i}\to0$. These necks pass to $\MM^{\bb_\infty}$ with error $2\delta_0$ for every $t\in(-T_{\bb_\infty},0)$. But near extinction $\MM^{\bb_\infty}$ approaches a round cylinder with fewer than $k$ Euclidean directions, a contradiction for small $\delta_0$.
\end{proof}

For $1\leq j\leq k$, define
\[ X_j := \IR^k_{\geq 0} \cap \Big\{ \frac{x_j}{Q'_j} \leq \min_{1\leq i\leq k} \frac{x_i}{Q'_i} \Big\} \cap \ov\IB^k_{|\Qu'|}, \]
such that
\[ \{ \Qu' \} = X_1 \cap \ldots \cap X_k \cap \big\{ x\in\IR^k : |x|=|\Qu'| \big\}. \]
We now label the entries of $\Qu$ by the reflection lines at $p$. For sufficiently negative times these lines span the lower Ricci eigenspace by the verification of Assumption~\ref{Lem_RF_NPC_convergence_b} above. With reference time $0$ and HMHF based at $p$, apply the inverse of the quadratic polynomial in Claim~\ref{claim:RF_intrinsic_quadratic_Ricci} to the rescaled Ricci eigenvalues on these lines. The resulting values match the eigenvalues of $\ov U$, counted with multiplicity, up to $O(|\tau|^{-3})$. The error is uniform for a convergent sequence, by choosing a common terminal time as in Proposition~\ref{Prop_RF_Q_continuity} and using the bounds of Proposition~\ref{Prop_ancient_summary} at the HMHF basepoint. The case $\ov U=0$ follows from the exponential bound in Proposition~\ref{Prop_ancient_summary}\ref{Prop_ancient_summary_f}, or from the exact cylinder formula.
For distinct branches $u,v$, subtracting their scalar equations in \eqref{eq_ovU_ODE} gives $|\partial_\tau\log|u-v||\leq 5/(2|\tau|)$ for sufficiently negative $\tau$. Indeed, integrating $\partial_\tau(-u)\geq(\sqrt2-O(|\tau|^{-1}))u^2$ backwards gives $|u|\leq (1/\sqrt2+o(1))/|\tau|$ uniformly on the same tail. Thus branches separated by more than twice the matching error bound at a fixed time stay separated at all earlier times. By continuity each reflection line matches a single branch on a sufficiently negative tail. Write $q_j(\MM)$ for the corresponding entry of $\Qu_k(\ov U)$ and $q=(q_1,\ldots,q_k)$. If the Ricci eigenvalue on the $j$-th line vanishes identically, then $q_j=0$, since each nonzero branch is of order $|\tau|^{-1}$.
For marked convergence, Proposition~\ref{Prop_RF_Q_continuity} and Lemma~\ref{Lem_RF_ODE_UQ} give convergence of the branches in the ordering induced by the sorted $\Qu$. At a fixed sufficiently negative time, branches with distinct limiting entries of $\Qu$ are separated by more than twice the uniform error. The Ricci eigenvalues on the marked lines converge at that time, so each line matches the same group of branches as in the limit. The separation persists backwards by the same estimate. Hence $q$ is continuous. Under parabolic rescaling the matching shifts rescaled time, so $q(\alpha\MM)=\alpha q(\MM)$ by Lemma~\ref{Lem_RF_ODE_UQ}\ref{Lem_RF_ODE_UQ_d}.

Set $K_j:=\{\MM\in\mathfrak M:q(\MM)\in X_j\}$. The $K_j$ are closed in $K$, and $K=\bigcup_{j=1}^k K_j$. Fix $0<\eps\leq\ov\eps$. By compactness of $K$ and continuity of $\Qu$ and parabolic rescaling, choose $0<\rho<\eps/2$ such that any $\MM\in\mathfrak M$ with $d_{CGH}(\MM,K)\leq\rho$ can be rescaled into $K$ at distance less than $\eps/2$ from $\MM$: the factor $|\Qu'|/|\Qu(\MM)|$ tends to $1$ whenever rescaling is needed and $d_{CGH}(\MM,K)\to0$. For $1\leq j\leq k$, define the open sets
\begin{equation*}
 U_0^{a,\eps}:=\{\bb\in\Delta_\delta^k:d_{CGH}(a\MM^{\bb},K)>\rho\},\qquad
 U_j^{a,\eps}:=\{\bb\in\Delta_\delta^k:d_{CGH}(a\MM^{\bb},K_j)<\eps\}.
\end{equation*}
They cover $\Delta_\delta^k$ since $\rho<\eps$ and $K=\bigcup_{j=1}^kK_j$. The rescaled flows $a\MM^{\bO}$ converge to the round cylinder in $K$, so $\bO\notin\overline{U_0^{a,\eps}}$ for sufficiently large $a$.
Let $F_j := \Delta_\delta^k \cap \{ x_j = 0 \}$ and $F_0 := \Delta_\delta^k \cap \{ x_1 + \ldots + x_k = \delta \}$ be the faces of $\Delta_\delta^k$.

\begin{Claim}
If $a \geq \underline a(\eps)$, then
\[ F_0 \subset U^{a,\eps}_0 \qquad \text{and} \qquad F_j \subset U_0^{a,\eps} \cup U_j^{a,\eps} \quad \text{for all} \quad j = 1, \ldots, k. \]
\end{Claim}

\begin{proof}
If an inclusion fails along $a_i\to\infty$ and $\bb_i\in F_j$, then $d_{CGH}(a_i\MM^{\bb_i},K)\leq\rho$. Claim~\ref{claim:RF_ellipsoid_compactness} gives $\bb_i\to\bO$ and, after passing to a subsequence, $a_i\MM^{\bb_i}\to\MM_\infty\in\mathfrak M$ with $d_{CGH}(\MM_\infty,K)\leq\rho$. This rules out $j=0$. If $j>0$, the $j$-th line in $a_i\MM^{\bb_i}$ is a split direction. Its Ricci eigenvalue remains zero in the limit, so $q_j(\MM_\infty)=0$. Rescale $\MM_\infty$ into $K$ as above. The rescaled flow still has $q_j=0$, hence lies in $K_j$ and is at distance less than $\eps/2$ from $\MM_\infty$. It follows that $\bb_i\in U_j^{a_i,\eps}$ for large $i$, a contradiction.
\end{proof}

By \cite[Claim~7.46]{Bamler_Lai_MCF1}, for sufficiently large $a$ choose $\bb_{a,\eps}\in\bigcap_{j=0}^k U_j^{a,\eps}$. Letting $a\to\infty$ along a subsequence, Claim~\ref{claim:RF_ellipsoid_compactness} gives a limit $\MM_{\infty,\eps}\in\mathfrak M$ with
\[
d_{CGH}(\MM_{\infty,\eps},K)\geq\rho,\qquad d_{CGH}(\MM_{\infty,\eps},K_j)\leq\eps\quad(1\leq j\leq k).
\]
As $\eps\to0$, compactness of $K$ gives a subsequential limit $\MM'\in\bigcap_{j=1}^k K_j$. The first bound gives $|\Qu(\MM_{\infty,\eps})|>|\Qu'|$, so continuity and $\MM'\in K$ give $|\Qu(\MM')|=|\Qu'|$. By the definition of the $X_j$, this gives $q(\MM')=\Qu'$, hence $\Qu(\MM')=\Qu'$. Positivity excludes a Euclidean factor.

\end{proof}

\bibliography{references}	
\bibliographystyle{amsalpha}

\end{document}